\documentclass[12pt, openany]{report}
\usepackage[utf8]{inputenc}
\usepackage[T1]{fontenc}
\usepackage[a4paper,left=2cm,right=2cm,top=2cm,bottom=2cm]{geometry}
\usepackage[pdftex]{graphicx}
\usepackage{amsmath,amssymb,amsthm}
\usepackage{hyperref}
\usepackage{pgf,tikz,pgfplots}
\usetikzlibrary{arrows}
\usepackage{graphicx}
\usepackage{array}
\usepackage{color}

\usepackage[russian,english]{babel}

\newcommand{\HRule}{\rule{\linewidth}{0.5mm}}

\newtheorem{theorem}{Theorem}[chapter]
\theoremstyle{definition}
\newtheorem{definition}[theorem]{Definition}
\theoremstyle{remark}
\newtheorem{example}[theorem]{Example}
\theoremstyle{remark}
\newtheorem{remark}[theorem]{Remark}
\theoremstyle{theorem}
\newtheorem{lemma}[theorem]{Lemma}
\theoremstyle{theorem}
\newtheorem{conjecture}[theorem]{Conjecture}
\theoremstyle{theorem}
\newtheorem{proposition}[theorem]{Proposition}
\theoremstyle{theorem}
\newtheorem{corollary}[theorem]{Corollary}

\newcommand{\enum}{\noindent\hspace*{1cm}$-$ }
\newcommand{\linep}[2]{{#1}_{\rangle #2}}
\newcommand{\inverse}[1]{{#1}^{-1}}
\newcommand{\orth}[1]{{#1}^{\bot}}
\newcommand{\polar}[1]{{#1}^{\ast}}
\newcommand{\ensemble}[2]{\left\{#1\,\left|\,#2\right.\right\}}
\newcommand{\grass}[2]{\text{Gr}_{#1}\left(#2\right)}
\newcommand{\matspace}[2]{\mathcal{M}_{#1}\left(#2\right)}
\newcommand{\PP}[2]{\mathbb{P}^{#1}(#2)}
\newcommand{\restreint}[2]{{#1}_{|#2}}
\newcommand{\scalar}[2]{\langle#1|#2\rangle}
\newcommand{\rscalar}[2]{\left(#1\cdot#2\right)}
\newcommand{\enonce}[2]{\textbf{#1. }\emph{#2}}
\newcommand{\transp}[1]{{#1}^{\text{T}}}
\newcommand{\class}{\mathcal{C}}
\newcommand{\id}{\text{Id}\,}
\newcommand{\az}{\text{az}}
\newcommand{\im}{\text{Im}\,}
\newcommand{\rk}{\text{rk}\;}
\newcommand{\proj}{\text{proj}}

\newcommand{\NN}{\mathbb{N}}
\newcommand{\ZZ}{\mathbb{Z}}

\newcommand{\RR}{\mathbb{R}}
\newcommand{\CC}{\mathbb{C}}
\newcommand{\RP}{\mathbb{R}\mathbb{P}}
\newcommand{\CP}{\mathbb{C}\mathbb{P}}
\renewcommand{\SS}{\mathbb{S}}
\newcommand{\HH}{\mathbb{H}}
\newcommand{\dd}{\text{d}}
\newcommand{\go}{«}
\newcommand{\gf}{»}

\renewcommand{\mod}{\,\text{mod}\,}

\begin{document}

\begin{titlepage}
  \begin{sffamily}
  \begin{center}

    \includegraphics[height=2.5cm]{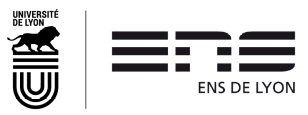}\\[0.4cm]
    \end{center}
    {Numéro National de Thèse : 2021LYSEN014}\\
    \begin{center}
    
    {\bf\large TH\`ESE de DOCTORAT DE L’UNIVERSITÉ DE LYON}\\
	{\large opérée par}\\
	{\bf\large l’École Normale Supérieure de Lyon}\\[1cm]
	
	{\bf\large École Doctorale $\text{N}^{\text{o}}\text{512}$}\\
	{\bf\large Informatique et Mathématiques de Lyon}\\[1cm]
	
	{\bf\large Discipline : Mathématiques}\\[1cm]
	
	{\large Soutenue publiquement le 25/05/2021, par}\\[0.1cm]
	{\bf\large Corentin FIEROBE}\\[0.3cm]

    \HRule \\[0.4cm]
    { \LARGE \bfseries Billards projectifs et complexes, orbites périodiques et systèmes Pfaffiens\\[0.4cm] }
    \HRule \\[0.7cm]
  \end{center}
    
    {\large Devant le jury composé de :}\\[0cm]

	\begin{tabular}{m{.25\textwidth} >{\centering\arraybackslash}m{.23\textwidth} >{\centering\arraybackslash}m{.27\textwidth} >{\centering\arraybackslash}m{.18\textwidth}}
	Stolovitch, Laurent & Directeur de recherche & Université Côte d'Azur & Rapporteur\\
	Tabachnikov, Sergei & Professor & Pennsylvannia State University & Rapporteur\\
	 
	Mazzucchelli, Marco & Chargé de recherche & École Normale Supérieure de Lyon & Examinateur\\
	 Paris-Romaskevich, Olga & Chargée de recherche & Institut de Mathématiques de Marseille & Examinatrice\\
	 Radnović, Milena & Associate professor & The University of Sydney & Examinatrice\\
	 Sorrentino, Alfonso & Professor & Universita degli Studi di Roma Tor Vergata & Examinateur\\
	 
	Zeghib, Abdelghani & Directeur de recherche & École Normale Supérieure de Lyon & Examinateur\\
	 
	 Glutsyuk, Alexey & Chargé de recherche & École Normale Supérieure de Lyon & Directeur de thèse\\
	\end{tabular}

  \end{sffamily}
\end{titlepage}

\begin{titlepage}
  \begin{sffamily}
  \begin{center}


    \textsc{\large École Normale Supérieure de Lyon}\\
    \textsc{\large Unité de Mathématiques Pures et Appliquées}\\
    \HRule \\[2cm]

    {\large\bf TH\`ESE DE DOCTORAT}\\[1.3cm]
    
    { \huge \bfseries Billards projectifs et complexes, orbites périodiques et systèmes Pfaffiens\\[0.4cm] }
    
\vspace{2cm}

\definecolor{xdxdff}{rgb}{0.49,0.49,1}
\definecolor{ffffqq}{rgb}{1,1,0}
\definecolor{tttttt}{rgb}{0.2,0.2,0.2}
\begin{tikzpicture}[line cap=round,line join=round,>=triangle 45,x=0.7cm,y=0.7cm]
\clip(-8,-6) rectangle (13,7);
\fill[dash pattern=on 2pt off 2pt,color=tttttt,fill=tttttt,fill opacity=0.1] (-7.3,-4.46) -- (3.46,5.94) -- (12.08,-2.26) -- cycle;
\draw [line width=1.2pt,color=tttttt] (-7.3,-4.46)-- (3.46,5.94);
\draw [line width=1.2pt,color=tttttt] (3.46,5.94)-- (12.08,-2.26);
\draw [line width=1.2pt,color=tttttt] (12.08,-2.26)-- (-7.3,-4.46);
\draw [line width=5.2pt,color=ffffqq] (6.54,3.01)-- (-3.44,-0.73);
\draw [line width=5.2pt,color=ffffqq] (6.54,3.01)-- (2.88,-3.3);
\draw [line width=5.2pt,color=ffffqq] (-3.44,-0.73)-- (6.48,-2.9);
\draw [line width=5.2pt,color=ffffqq] (0.82,3.39)-- (9.84,-0.13);
\draw [line width=5.2pt,color=ffffqq] (9.84,-0.13)-- (6.48,-2.9);
\draw [line width=5.2pt,color=ffffqq] (0.82,3.39)-- (2.88,-3.3);
\draw (2.88,-3.3)-- (6.54,3.01);
\draw (6.54,3.01)-- (-3.44,-0.73);
\draw (-3.44,-0.73)-- (6.48,-2.9);
\draw (6.48,-2.9)-- (9.84,-0.13);
\draw (9.84,-0.13)-- (0.82,3.39);
\draw (0.83,3.39)-- (2.88,-3.3);
\draw [dash pattern=on 5pt off 5pt] (-3.44,-0.73)-- (3.48,0.5);
\draw [dash pattern=on 5pt off 5pt] (3.48,0.5)-- (0.83,3.39);
\draw [dash pattern=on 5pt off 5pt] (3.48,0.5)-- (6.54,3.01);
\draw [dash pattern=on 5pt off 5pt] (3.48,0.5)-- (9.84,-0.13);
\draw [dash pattern=on 5pt off 5pt] (3.48,0.5)-- (6.48,-2.9);
\draw [dash pattern=on 5pt off 5pt] (3.48,0.5)-- (2.88,-3.3);
\begin{scriptsize}
\end{scriptsize}
\end{tikzpicture}
\vspace{2cm}

    {\large réalisée par}\\[0.2cm]
    {\LARGE\bf Corentin Fierobe}\\[1cm]
    {\large dirigée par}\\[0.2cm]
    {\LARGE\bf Alexey Glutsyuk}\\[2cm]
    

    \vfill
    {\large }
  \end{center}
  \end{sffamily}
\end{titlepage}

\chapter*{Remerciements}
	Je tiens à remercier très sincèrement mon directeur de thèse, Alexey Glutsyuk, qui s'est beaucoup investi dans mon travail et qui m'a beaucoup soutenu. J'ai appris beaucoup grâce à toi, en mathématiques mais aussi sur la recherche. Tes compétences, ton ouverture d'esprit et ta gentillesse m'ont permis de faire mûrir mes réflexions en toute liberté. Un grand merci de m'avoir permis de participer à de nombreux séminaires et conférences ; visiter Moscou, Nizhnyi Novgorod et Novosibirsk furent des expériences inoubliables.

J'aimerais aussi remercier tous les membres du jury de cette thèse pour leur expertise précieuse, Marco Mazzucchelli, Olga Paris-Romaskevich, Milena Radnović, Alfonso Sorrentino, Laurent Stolovitch, Sergei Tabachnikov, Abdelghani Zeghib. Un merci appuyé à Laurent Stolovitch et Sergei Tabachnikov d'avoir accepté d'être rapporteurs. Merci aussi pour vos remarques et vos conseils. 

Je remercie également Alfonso Sorrentino avec qui j'ai pu faire un peu de maths, mais aussi observer les débuts de la crise sanitaire depuis les régions glacées de Sibérie. Merci aussi à Lior Shalom et Sergyi Maksymenko d'avoir rendu le voyage plus sympathique encore.

Merci aux doctorants de l'UMPA et de l'ICJ, merci Simon Allais pour nos discussions stimulantes, merci Valentine Roos pour ton aide dans l'élaboration du cours de calcul différentiel pour économiste, merci Matthieu Joseph pour ton aide dans le déroulement du TD d'analyse complexe. Merci Mélanie Théillière, Anatole Ertul, Gauthier Clerc, Mete Demircigil, Gabriele Sbaiz, pour votre aide dans l'organisation du séminaire des doctorants deux années de suite. Merci Mendes Oulamara pour nos discussions au sujet de la thèse.

Je voudrais aussi remercier tous les membres de l'UMPA pour leur soutien inconditionnel. Un grand merci à Magalie Le Borgne, Virginia Goncalves et Laure Savetier pour leur aide précieuse et leur écoute. Merci notamment de m'avoir laissé libre de partir en mission sans prendre l'avion. Merci aussi à Micaël Calvas pour ses conseils en informatique.

Un grand merci à ma future femme, Anastasia. Mes premières pensées vont naturellement vers toi, pour ton soutien, ton écoute et ta présence. Merci à toute ma famille pour leur aide et leur soutien, à mes frères, Stéphane, Victor, Hippolyte et Achille, mes parents Florence et Thierry, mes grand-parents Jacqueline, Pierre, Claudie et Edgar, ainsi qu'à \'Elise, Eugénie, Balthazar, Lisiane, Aurore, Julie et au reste de ma famille. \foreignlanguage{russian}{Спасибо Алине и всем членам ее семьи, Мите, Ольге, Кириллу, Ксении, Артуру}. Merci Romy, Pesto, Opia, Affli, Nida, Mao, Muse, Blue.

Merci !

\tableofcontents

\chapter*{Introduction en français}
\addcontentsline{toc}{chapter}{\protect\numberline{}Introduction en français}%

 Un billard peut être décrit comme un système dynamique modélisant le comportement d'un objet sans volume ni masse, par exemple une particule infiniment petite ou un grain de lumière, qui évolue sans frottements dans un milieu homogène délimité par une paroi réfléchissante. Comme l'ont très bien résumé Valerii V. Kozlov et Dmitrii V. Treshchëv \cite{treshchev}, l'étude des billards qui « \emph{[a commencé] avec les travaux de D. Birkhoff, a été un sujet de recherche populaire combinant différents éléments de théorie ergodique, théorie de Morse, théorie KAM, etc. Les billards sont d'autant plus remarquables qu'ils apparaissent naturellement dans un grand nombre de problèmes de mécanique et de physique (systèmes vibrant à impacts, diffraction des ondes courtes, etc.).} »\footnote{\label{footnote:treshchev}
\go \foreignlanguage{russian}{Начиная с работ Дж. Биркгофа, биллиарды являются популярной темой исследования, где естественным образом переплетаются различные сюжеты из эргодический теории, теории Морса, КАМ-теории и т.д. С другой стороны, биллиардные системы замечательны еще и тем, что естественно возникают в ряде важных задач механики и физики (виброударные системы, дифракция коротких волн и др.).}\gf} La thèse ci-présente s'inscrit dans ce champ de recherche et tente d'apporter des réponses partielles à de grandes questions qui la traversent.

Le mouvement d'une particule dans un billard est régi par deux contraintes: 1) elle se déplace \textit{en ligne droite} à l'intérieur du milieu et 2) se réfléchit sur la paroi selon la loi d'optique géométrique \textit{angle d'incidence = angle de réflexion.} Le modèle mathématique le plus courant pour décrire les assertions 1) et 2) est celui d'une variété Riemannienne complète : le déplacement en lignes droites est celui qui suit les géodésiques, et la mesure des angles est donnée par la métrique. On peut donc par exemple étudier des billards dans le plan, dans l'espace, sur un hyperboloïde ou sur une sphère, ce dernier cas pouvant s'avérer utile par exemple dans une simulation où la courbure de la terre n'est plus négligeable. Il existe cependant d'autres modèles de billards que ces billards dits \textit{classiques}~: évoquons les billards extérieurs, les billards filaires, les billards dans les pavages ou les billards pseudo-Euclidiens. Dans cette thèse, une attention particulière sera portée aux billards dits \textit{projectifs} ainsi qu'aux billards \textit{complexes}. Ces deux derniers modèles généralisent les billards classiques et peuvent permettre de démontrer certains résultats liés à la théorie classique du billard, ce dont une partie de cette thèse va s'attacher à montrer.

\begin{figure}[!h]
\centering
\definecolor{qqwuqq}{rgb}{0,0.39,0}
\begin{tikzpicture}[line cap=round,line join=round,>=triangle 45,x=2.0cm,y=2.0cm]
\clip(-1.4,-0.8) rectangle (1.4,0.8);
\draw [shift={(-0.85,0.52)},color=qqwuqq,fill=qqwuqq,fill opacity=0.1] (0,0) -- (-150.66:0.15) arc (-150.66:-117.76:0.15) -- cycle;
\draw [shift={(-0.85,0.52)},color=qqwuqq,fill=qqwuqq,fill opacity=0.1] (0,0) -- (-3.56:0.15) arc (-3.56:29.34:0.15) -- cycle;
\draw [rotate around={0:(0,0)}] (0,0) ellipse (2.46cm and 1.44cm);
\draw (-0.27,-0.54)-- (-1.21,-0.15);
\draw (-1.21,-0.15)-- (-0.85,0.52);
\draw [-latex] (-0.27,-0.54) -- (-0.81,-0.31);
\draw [-latex] (-1.21,-0.15) -- (-1.01,0.23);
\draw [-latex] (-0.85,0.52) -- (-0.11,0.47);
\begin{scriptsize}
\end{scriptsize}
\end{tikzpicture}
\hspace{1.5cm}
\begin{tikzpicture}[line cap=round,line join=round,>=triangle 45,x=1.5cm,y=1.5cm]
\clip(-1.6,-1.2) rectangle (1.6,1.2);
\draw [rotate around={0:(0,0)}] (0,0) ellipse (2.12cm and 1.5cm);
\draw (-1.56,-0.1)-- (-1.27,0.11);
\draw (-1.51,0.3)-- (-1.12,0.39);
\draw (-1.24,0.73)-- (-0.96,0.53);
\draw (-0.88,0.99)-- (-0.64,0.69);
\draw (-0.3,1.18)-- (-0.2,0.79);
\draw (0.4,1.16)-- (0.31,0.78);
\draw (1.02,0.91)-- (0.77,0.64);
\draw (1.37,0.58)-- (1.09,0.42);
\draw (1.55,0.17)-- (1.28,-0.06);
\draw (1.07,-0.87)-- (0.76,-0.65);
\draw (0.66,-1.09)-- (0.37,-0.77);
\draw (0.14,-1.2)-- (0,-0.8);
\draw (-0.46,-1.15)-- (-0.39,-0.76);
\draw (-0.81,-1.03)-- (-0.79,-0.63);
\draw (-1.11,-0.85)-- (-1.03,-0.47);
\draw (-1.4,-0.54)-- (-1.24,-0.19);
\draw (1.11,-0.39)-- (1.52,-0.26);
\begin{scriptsize}
\draw[color=black] (-0.59,1.1) node {$\partial\Omega$};
\draw[color=black] (-0.8,0) node {$\Omega$};
\end{scriptsize}
\end{tikzpicture}
\caption{\`A gauche, un rayon lumineux se réfléchissant sur le bord d'un domaine selon la loi d'optique géométrique. \`A droite, un billard projectif et son champ de droites transverses.}
\label{figure:billiard_intro}
\end{figure}
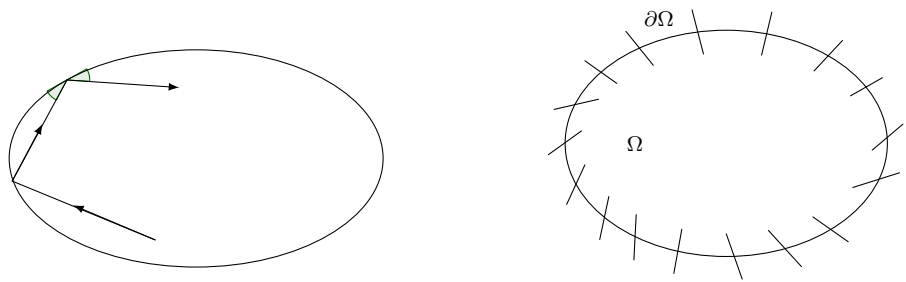

Les \textit{billards complexes} sont une extension naturelle des billards classiques au plan Euclidien complexifié, c'est-à-dire à $\CC^2$. Ils ont été introduits et étudiés par Glutsyuk \cite{glut,glut1,glut2} pour résoudre la conjecture de Ivrii à quatre réflexions, la conjecture des billards commutants en dimension $2$, ou encore la conjecture d'invisibilité de Plakhov (cas planaire à $4$ réflexions). Souvent combinés à la théorie des sytèmes Pfaffiens, ils permettent notamment d'appliquer des méthodes d'analyse complexe à la résolution de problèmes réels. Nous reviendrons plus en détails sur ces questions.

Introduits par Tabachnikov qui les a étudiés en détails \cite{taba_projectif_ball, taba_projectif}, les \textit{billards projectifs} généralisent les billards classiques. Un billard projectif est un domaine borné d'un espace euclidien dont le bord est traversé par un champ de droites transverses, dites \textit{droites projectives}. Une particule à l'intérieur du domaine se déplace le long de droites. Elle est réfléchie sur le bord de sorte que la droite incidente, la droite réfléchie, la droite projective en le point d'impact, et la droite obtenue par intersection de l'hyperplan contenant ces trois premières droites avec l'hyperplan tangent à la surface forment une famille harmonique. Lorsque la droite projective est perpendiculaire au bord, cette condition impose à la réflexion de suivre la loi d'optique géométrique. Ceci reste vrai quand la droite projective est perpendiculaire à la surface pour une métrique pseudo-Euclidienne ou encore une métrique projectivement équivalente à la métrique Euclidienne (c'est à dire dont les géodésiques sont supportées par des droites). Ainsi les billards projectifs englobent différents types de billards.

Dans le modèle du billard classique à l'intérieur d'un domaine $\Omega$ borné de frontière $\partial\Omega$ lisse, la dynamique d'une particule évoluant à l'intérieur de $\Omega$ se décrit à l'aide de deux objets. Le premier est l'\textit{espace des phases}, c'est à dire l'ensemble des morceaux de trajectoires entre deux rebonds. Il peut notamment être codé par un couple $(p,v)$, où $p$ est un point du bord $\partial\Omega$ et $v$ et un vecteur unitaire dirigé vers l'intérieur de $\Omega$ et représentant la direction de la trajectoire. Dans le plan, on peut aussi remplacer $v$ par une mesure $\theta\in[0,\pi]$ de l'angle qu'il forme avec la tangente $T_p\partial\Omega$. Dépendant de ces deux paramètres, l'espace des phases est ainsi de dimension $2$ pour les billards du plan, et de façon générale de dimension $2(d-1)$ pour les billards dans un espace de dimension $d$. Le deuxième objet modélisant la dynamique du billard est l'\textit{application de billard}, une application qui, étant donné un couple $(p,v)$ de l'espace des phases codant la trajectoire d'une particule émise du point $p$ avec une direction $v$, renvoie le couple $(q,w)$ de l'espace des phases où $q\in\partial\Omega$ est le prochain point d'impact de la particule et $w$ est le vecteur unitaire dirigeant la trajectoire après réflexion. Ces deux objets, espace des phases et application de billard, peuvent aussi être définis pour d'autres types de billards.

\section*{Conjecture de Ivrii}

L'un des enjeux de la théorie des billards est l'étude des \textit{trajectoires périodiques}, c'est-à-dire des trajectoires qui se répètent après un nombre fini de réflexions. 
Ivrii \cite{ivrii} a montré en 1980 que l'étude des orbites périodiques de billards a une application dans un problème célèbre, qui a été résumé par Kac \cite{kac} en une question : \textit{peut-on entendre la forme d'un tambour ?}\footnote{\go \textit{Can one hear the shape of a drum} \gf, titre de l'article cité, \cite{kac}.} Il s'agit de comprendre si la donnée des valeurs propres du problème de Dirichlet dans un domaine borné $\Omega\subset\RR^d$ permet de retrouver $\Omega$. Les valeurs propres du problème de Dirichlet sont les réels $\lambda$ pour lesquels le système
\begin{equation}
\label{equation:probleme_dirichlet}
\left\{\begin{array}{l}
\Delta u +\lambda u=0\\
\restreint{u}{\partial\Omega}=0
\end{array}\right.
\end{equation}
possède des solutions non-triviales. Elles peuvent être interprétées physiquement comme les différents modes de vibration d'une forme $\Omega$ donnée, ce qui explique la question de Kac. La réponse à cette question s'est avérée être négative et des exemples de domaines de formes distinctes ont été donnés pour lesquels les problèmes de Dirichlet \eqref{equation:probleme_dirichlet} correspondants ont les mêmes valeurs propres. Néanmoins se pose toujours la question de pouvoir retrouver des informations sur $\Omega$ à partir des valeurs propres du problème de Dirichlet. Weyl \cite{weyl} a montré que l'on peut \textit{entendre le volume\footnote{\go The first pertinent result is that one can hear the area of $\Omega$ \gf, \cite{kac}} de $\Omega$}, au sens ou la connaissance du spectre de Dirichlet permet de retrouver ce volume. En effet, les valeurs propres du problème de Dirichlet peuvent être énumérées par une famille $(\lambda_n)_n$ de sorte que $0\leq\lambda_1\leq\lambda_2\leq\ldots\leq\lambda_n\leq\ldots$ avec $\lambda_n\to+\infty$. On note $N(\lambda)$ le nombre de valeurs propres inférieures ou égales à $\lambda$. Alors Weyl a prouvé que
$N(\lambda)\sim(2\pi)^{-d}v_d\text{vol}({\Omega})\lambda^{d/2}$, où $v_d$ est le volume de la boule unité de $\RR^d$. Il a aussi conjecturé le second terme de ce développement asymptotique~: 
\begin{equation}
\label{conjecture:weyl}
N(\lambda)=(2\pi)^{-d}v_d\text{vol}({\Omega})\lambda^{d/2}-\frac{1}{4(2\pi)^{d-1}}\text{area}(\partial\Omega)\lambda^{(d-1)/2}+o(\lambda^{(d-1)/2}).
\end{equation}
Cette formule reste une conjecture dans sa généralité malgré de nombreuses avancées dont une notable est due à Ivrii \cite{ivrii}, qui a prouvé que \eqref{conjecture:weyl} est vérifiée sous réserve 
que le billard constitué par $\Omega$ a \textit{peu} d'orbites périodiques. Plus précisément, la condition imposée est que l'ensemble des paramètres correspondant aux orbites périodiques dans l'espace des phases du billard soit de mesure nulle. Cela a donné lieu à une célèbre conjecture portant son nom~:

\enonce{Conjecture de Ivrii}{\'Etant donné un domaine d'un espace Euclidien dont le bord est suffisamment lisse, l'ensemble de ses orbites périodiques est de mesure nulle.}

Cette conjecture, qui tient toujours, relève d'une grande complexité malgré sa simplicité apparente. Si elle est vérifiée, elle impliquerait notamment qu'un billard ne possède pas d'\textit{ouvert d'orbites périodiques}, c'est à dire que son espace des phases ne contient pas d'ouvert contenant uniquement des paramètres $(p,v)$ associés à des orbites périodiques d'une période donnée $k$. On ne sait pas encore si un tel billard, dit \textit{$k$-réfléchissant}, existe ou non. Son existence aurait la conséquence amusante suivante: elle permettrait de construire une salle dont les murs sont recouverts de miroirs et de sorte qu'il existe un endroit de la salle où un observateur regardant devant lui peut toujours voir son image de dos, même s'il se déplace un peu et/ou tourne légèrement sur lui-même.

La conjecture de Ivrii a été abordée dans de nombreux articles. Elle a d'abord été prouvée de façon générique par Petkov et Stojanov \cite{petkovstojanov} : l'ensemble des domaines de $\RR^d$ de bord $\class^{\infty}$ ayant pour tout $k\geq 2$ un nombre fini d'orbites périodiques de période $k$ contient un ensemble résiduel, c'est-à-dire une intersection dénombrable d'ouverts denses. Une autre réponse partielle à la conjecture a été donnée par Vasiliev \cite{vasiliev} qui l'a prouvée pour un domaine convexe de bord analytique. Notons aussi qu'il est possible de restreindre la conjecture à l'ensemble des orbites périodiques d'une période donnée arbitraire, et que l'ensemble de ces conjectures restreintes est équivalent à la conjecture globale. Dans cet idée, Rychlik \cite{rychlik}, puis Stojanov \cite{stojanov} ont démontré que l'ensemble des orbites de période $3$, ou \textit{triangulaires}, est de mesure nulle dans un billard du plan de frontière de classe $\class^3$, et Vorobets \cite{vorobets} a étendu ce résultat aux billards en toute dimension. Un peu plus tard, Wojtkowski \cite{wojtkowski}, puis Baryshnikov et Zharnitsky \cite{bary} ont donné de nouvelles preuves de ce résultat. Plus récemment, Glutsyuk et Kudryashov \cite{glutkud2} ont démontré la conjecture pour les orbites périodiques de période $4$ dans des billards planaires de classe $\class^4$. En toute généralité dans le cas Euclidien, la conjecture de Ivrii tient toujours pour un nombre quelconque de réflexions, même pour des classes de billards de frontière très lisse (par exemple analytique par morceaux).

\begin{figure}[!h]
\centering
\includegraphics[scale=0.4]{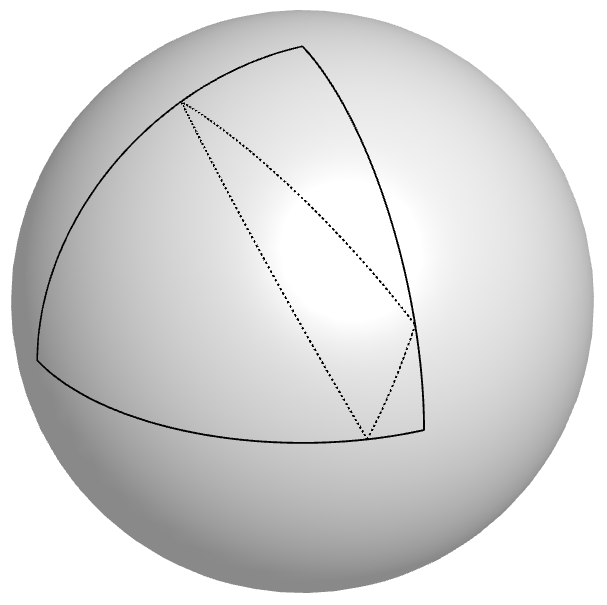}
\caption{Un exemple de billard $3$-réfléchissant sur la sphère proposé par Barychnikov. Le triangle extérieur représente le bord du billard, le triangle intérieur en pointillé est une orbite. On peut bouger arbitrairement deux points de l'orbite sans changer son caractère périodique.}
\label{figure:sphere_bary}
\end{figure}

La conjecture de Ivrii s'énonce de façon analogue pour des billards non-Euclidiens, par exemple pour les billards en courbure constante, sur une sphère ou un hyperboloïde. Des exemples remarquables \cite{bary_introuvable,VKNZ} de billards $2$- ou $3$-réfléchissants existent sur la sphère de dimension $2$, liés d'une certaine façon à l'existence de points joints par une infinité de géodésiques distinctes et contredisant la conjecture de Ivrii sur la sphère, voir la Figure \ref{figure:sphere_bary}. Les articles cités \cite{bary_introuvable,VKNZ} donnent une classification des billards sur la sphère unité $\SS^2$ ayant un ouvert d'orbites de période $3$ ainsi que la non-existence de tels billards sur l'hyperboloïde $\HH^2$.

Malgré tous ces résultats, la conjecture de Ivrii reste encore ouverte. Il semble d'ailleurs que les spécialistes sont partagé.e.s entre celleux qui pensent qu'elle est vraie, et celleux qui pensent qu'elle est fausse et qui recherchent des contre-exemples.

\section*{Billards intégrables}

Un autre enjeu de la théorie des billards est l'étude des billards dits \textit{intégrables}. Un billard $\Omega$ du plan est dit globalement intégrable si son espace des phases est feuilleté de façon lisse par une famille de courbes fermées invariantes par l'application de billard. On dit aussi que $\Omega$ est localement intégrable si seul un voisinage du bord, correspondant à la courbe $\{\theta=0\}$ dans l'espace des phases, admet un tel feuilletage. Cette propriété se manifeste par l'existence de \textit{caustiques} correspondant à ces courbes invariantes et qui se définissent de façon indépendante en toute dimension : une caustique d'un billard $\Omega$ est une hypersurface $\Gamma\subset\Omega$ telle que toute droite tangente à $\Gamma$ et intersectant la frontière $\partial\Omega$ en un point $p$ est réfléchie en une droite tangente à $\Gamma$ après réflexion en $p$ sur le bord de $\Omega$. 

Un exemple de billard globalement intégrable est le disque, puisque tout cercle concentrique inclus dans le disque est une caustique du disque. L'ellipse est un exemple de billard localement intégrable, puisque toute trajectoire de billard qui ne passe pas entre les foyers reste tangente à une même ellipse homofocale, qui dès lors est une caustique de l'ellipse initiale. La question a été posée par Birkhoff et Poritsky de savoir si ce sont les seuls exemples de billards intégrables et cela a donné lieu à la célèbre conjecture de Birkhoff, ou Birkhoff-Poritsky comme cela a été rappelé dans  \cite{kaloshinsorrentino}.

\enonce{Conjecture de Birkhoff-Poritsky}{Les seuls billards localement intégrables sont les ellipses.}

Certaines avancées majeures ont été réalisées sur cette conjecture. Citons le théorème de Bialy \cite{bialy} énonçant que si l'espace des phases d'un billard est feuilleté par des courbes fermées continues invariantes et non-homotopes à un point, alors $\partial\Omega$ est un cercle. Cela implique que le seul billard globalement intégrable est le disque. Ce résultat nécessite néanmoins l'hypothèse que le feuilletage est global et ne permet pas de conclure que la conjecture est vraie en toute généralité. Une version algébrique de la conjecture de Birkhoff-Poritsky a été démontrée conjointement par Bialy, Glutsyuk et Mironov \cite{bialymironov_poly, bialymironov_poly_bis, glut_integrability, glut_integrability_bis} pour les billards sur le plan et sur les autres hypersurfaces de courbure constante. Kaloshin et Sorrentino \cite{kaloshinsorrentino} ont prouvé la version locale de la conjecture, démontrant que toute déformation intégrable d'une ellipse est une ellipse. En dimension supérieure, l'étude des billards ayant des caustiques a été conclue par Berger \cite{berger_caustics} qui a montré que si un billard de $\RR^d$, avec $d\geq 3$, admet une caustique, alors ce dernier est une quadrique et sa caustique est une quadrique homofocale. Ainsi en dimension au moins $3$, il suffit juste d'une seule caustique, et non plus un feuilletage, pour que la conjecture de Birkhoff-Poritsky soit vérifiée. 

\section*{Résultats obtenus dans cette thèse}
\addcontentsline{toc}{section}{\protect\numberline{}Résultats obtenus dans cette thèse}%

Cette thèse présente différents résultats sur les billards complexes et projectifs, applicables pour certains à la théorie des billards classiques. Elle se divise en trois chapitres : le \textbf{Chapitre \ref{chapter_billiards}} présente en détails les modèles des billards projectifs et complexes. Le \textbf{Chapitre \ref{chapter_caustics}} étudie la notion de caustique dans ces deux modèles de billard. Le \textbf{Chapitre \ref{chapter_ivrii}} porte son attention sur l'analogue de la conjecture de Ivrii appliquée aux billards projectifs.

\subsection*{Détails du Chapitre \ref{chapter_billiards}}

Ce chapitre présente les deux classes de billards étudiées tout au long de cette thèse, les billards complexes et les billards projectifs. Nous exposons brièvement quelques aspects de ces billards pour rendre compréhensible les résumés des chapitres suivants. Plus de détails seront donnés dans le Chapitre \ref{chapter_billiards} lui-même.

Un billard \textit{projectif} est un domaine borné $\Omega$ de $\RR^d$ dont le bord est lisse et muni d'un champ de droites transverses. Ce champ de droites induit en chaque point $p\in\partial\Omega$ du bord une transformation de l'ensemble des droites orientées passant par $p$, qui permet de considérer les orbites du billard: une droite orientée $\ell_0$ intersectant $\Omega$ en $p$ est réfléchie en une droite orientée $\ell_1$ par la transformation décrite précédemment. Si $\ell_1$ intersecte le bord en un autre point, cette construction peut être répétée, et ainsi de suite.

Un billard \textit{complexe} est une courbe complexe $\gamma$ de $\CP^2$ sur laquelle on définit une loi de réflexion de droites complexes qui l'intersecte. Cette construction est réalisée en considérant la complexification de la métrique Euclidienne $\dd x^2+\dd y^2$ à $\CC^2$. \'Etant donnée une droite complexe $L\subset\CC^2$ dite \textit{non-isotrope}, on peut définir une symétrie de droites complexes par rapport à $L$ : cette symétrie est l'unique involution affine non triviale qui fixe les points de $L$ et préserve la forme quadratique complexifiée définie précédemment. Deux droites complexes $\ell,\ell'$ intersectant $\gamma$ en un point $p$ sont dites symétriques (pour cette loi de réflexion complexe) si la symétrie de droites complexes par rapport à la tangente $T_p\gamma$ envoie l'une sur l'autre. Pour les autres droites $L$, dites \textit{isotropes}, on utilise un passage à la limite.

\subsection*{Détails du Chapitre \ref{chapter_caustics}}

Ce chapitre propose l'étude de propriétés relatives aux caustiques des billards projectifs et complexes. La Section \ref{section:general_properties_on_quadrics} présente un premier résultat publié \cite{fierobe_caustics} sur les caustiques dites \textit{complexes} d'une ellipse ou d'une hyperbole. On dira qu'une conique $C'\subset\CP^2$ est une \textit{caustique complexe} d'une autre conique $C\subset\CP^2$ si toute droite $\ell$ tangente à $C'$ est réfléchie en une droite tangente à $C$ par réflexion complexe en l'un des deux points d'intersection de $\ell$ avec $C$. \'Etant donnés $a,b\in\polar{\RR}$, on introduit la famille $(\mathcal{C}_{\lambda})_{\lambda\in\CC}$ de coniques de $\CC^2$ définies par l'équation
$$\mathcal{C}_{\lambda}:\frac{x^2}{a-\lambda}+\frac{y^2}{b-\lambda}=1$$
et on étudie le billard complexe sur $\mathcal{C}_0$. Il est connu que dans le cas du billard réel, les coniques réelles $\mathcal{C}_{\lambda}$ avec $\lambda\in\RR$ sont des caustiques du billard formé par $\mathcal{C}_0$. On s'interroge sur le fait de savoir si cela reste vrai dans le cas du billard complexe sur $\mathcal{C}_0$ et quels sont les caustiques complexes inscrites dans des orbites périodiques. On prouve les deux résultats suivants:

\enonce{Proposition}{Toute conique $\mathcal{C}_{\lambda}$ est une caustique complexe de $\mathcal{C}_0$.}

\enonce{Proposition}{Pour tout entier $n\geq 3$, il existe un polynôme en $(a,b,\lambda)$, noté $\mathcal{B}^n_{a,b}(\lambda)$, dont les racines complexes en $\lambda$ correspondent aux caustiques $\mathcal{C}_{\lambda}$ qui sont inscrites dans des orbites de période $n$. Pour $(a,b)$ en dehors d'un nombre fini de valeurs $a/b$, le degré en $\lambda$ du polynôme $\mathcal{B}^n_{a,b}(\lambda)$ est $(n^2-1)/4$ si $n$ est impair, et $n^2/4-1$ si $n$ est pair.}

Ainsi les racines distinctes en $\lambda$ de $\mathcal{B}^n_{a,b}(\lambda)$ différentes de $a$ et $b$ correspondent aux caustiques complexes de $\mathcal{C}_{0}$ inscrites dans les orbites périodiques de période $n$. Nous avons pu montrer que pour un nombre générique de $(a,b)$ (au sens du résultat précédent), ni $a$ ni $b$ ne sont racines (en $\lambda$) de $\mathcal{B}^n_{a,b}(\lambda)$. Il reste à déterminer si $\mathcal{B}^n_{a,b}(\lambda)$ est génériquement à racines simples en $\lambda$ ou non. Pour l'instant le résultat n'est pas connu, mais est vérifiée pour de petites périodes. Et en effet, un phénomène surprenant se produit dans le cas des orbites de période $3$ lorsque $\mathcal{C}_0$ est une ellipse (avec des résultats similaires pour une hyperbole ou pour les orbites de période $4$)~:

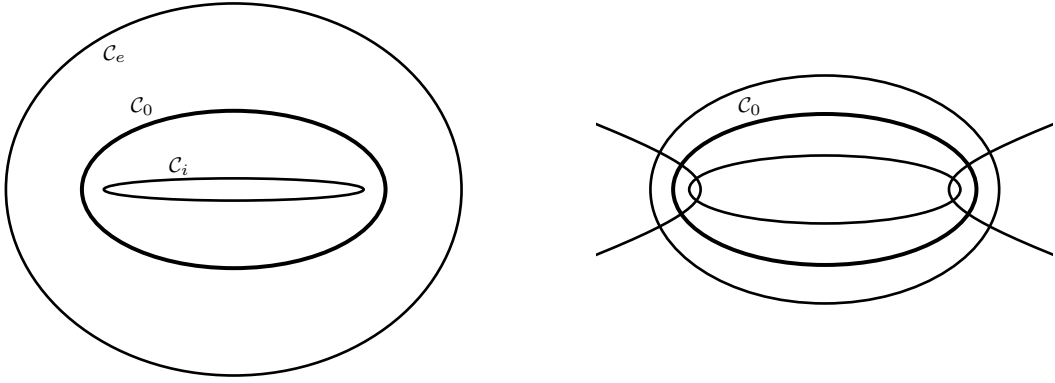
\begin{figure}
\centering
\begin{tikzpicture}[line cap=round,line join=round,>=triangle 45,x=1cm,y=1cm]
\clip(-3.5,-2.7) rectangle (3.5,2.7);
\draw [rotate around={0:(0.0025086145409116416,0)},line width=1.5pt] (0.0025086145409116416,0) ellipse (1.9974913854590879cm and 1.0432433360133084cm);
\draw [rotate around={0:(0.0025086145409115726,0)},line width=1pt] (0.0025086145409115726,0) ellipse (2.9974913854590883cm and 2.4664426668897765cm);
\draw [rotate around={0:(0.002508614540915727,0)},line width=1pt] (0.002508614540915727,0) ellipse (1.709768500390961cm and 0.1472859398657943cm);
\begin{scriptsize}
\draw[color=black] (-1.2,1.1) node {$\mathcal{C}_0$};
\draw[color=black] (-1.5597148670526204,1.8) node {$\mathcal{C}_{e}$};
\draw[color=black] (-0.7111058200015711,0.31037283751140887) node {$\mathcal{C}_{i}$};
\end{scriptsize}
\end{tikzpicture}
\hspace{1cm}
\begin{tikzpicture}[line cap=round,line join=round,>=triangle 45,x=1cm,y=1cm]
\clip(-3,-2.7) rectangle (3,2.7);
\draw [rotate around={0:(0.004465809016780549,0)},line width=1.5pt] (0.004465809016780549,0) ellipse (1.9955341909832196cm and 1.0008587655128343cm);
\draw [rotate around={0:(0.004465809016780692,0)},line width=1pt] (0.004465809016780692,0) ellipse (1.784232681737871cm and 0.45060828188388125cm);
\draw [rotate around={0:(0.004465809016780525,0)},line width=1pt] (0.004465809016780525,0) ellipse (2.2943056805364748cm and 1.5110923588129084cm);
\draw [samples=50,domain=-0.99:0.99,rotate around={0:(0.0044658090167803855,0)},xshift=0.0044658090167803855cm,yshift=0cm,line width=1pt] plot ({1.6329479850417183*(1+(\x)^2)/(1-(\x)^2)},{0.5602850319501439*2*(\x)/(1-(\x)^2)});
\draw [samples=50,domain=-0.99:0.99,rotate around={0:(0.0044658090167803855,0)},xshift=0.0044658090167803855cm,yshift=0cm,line width=1pt] plot ({1.6329479850417183*(-1-(\x)^2)/(1-(\x)^2)},{0.5602850319501439*(-2)*(\x)/(1-(\x)^2)});
\begin{scriptsize}
\draw[color=black] (-0.9777216831311774,1.1) node {$\mathcal{C}_0$};
\end{scriptsize}
\end{tikzpicture}
\caption{A gauche, une ellipse $\mathcal{C}_0$ avec ses deux caustiques complexes $\mathcal{C}_i$ et $\mathcal{C}_e$ inscrites dans des orbites triangulaires. Ce sont des ellipses complexifiées, l'une incluse dans $\mathcal{C}_0$ et l'autre la contenant. Le graphique représente leur partie réelle. A droite, les trois caustiques complexes de $\mathcal{C}_0$ pour les orbites de période 4.}
\label{figure:caustics_complex}
\end{figure}

\enonce{Proposition}{Si $a,b>0$, il existe exactement deux coniques complexes homofocales à $\mathcal{C}_0$ dont les orbites complexes qui leur sont circonscrites sont périodiques de période $3$. Ce sont des ellipses complexifiées : l'une $\mathcal{C}_i$ est incluse dans $\mathcal{C}_0$, l'autre $\mathcal{C}_e$ \textit{contient} $\mathcal{C}_0$ (voir Figure \ref{figure:caustics_complex}).}

Nous avons cherché des propriétés curieuses de ces deux ellipses qui pourraient apparaître, comme la question de savoir si $\mathcal{C}_0$ ou $\mathcal{C}_i$ sont des caustiques de la plus grande ellipse $\mathcal{C}_e$ inscrites dans des orbites périodiques du billard réel. Mais les simulations ont échoué à mettre en évidence un tel phénomène. Nous montrons alors qu'un invariant du billard elliptique réel connu sous le nom d'\textit{invariant de Joachimsthal} se généralise au cas complexe, et qu'il entretient des liens étroits avec les caustiques complexes de l'ellipse. 

Cette thèse propose ensuite une étude sur l'existence de caustiques dans les billards projectifs. Notons d'abord que de nombreux résultats ont été obtenus par Tabachnikov \cite{taba_projectif_ball,taba_projectif} sur l'existence de formes d'aire dans l'espace des phases qui sont invariantes par l'application de billard projectif, et sur les propriétés d'intégrabilités qui en découlent. Citons par exemple \cite{taba_projectif} Corollaire F : \textit{si l'application de billard dans un cercle muni d'une structure de billard projectif a une forme d'aire invariante lisse au voisinage du bord, alors le billard est intégrable}. Notons aussi qu'une nouvelle preuve de l'intégrabilité du billard elliptique dans le plan Euclidien, sur l'hyperboloïde ou sur la sphère a été donnée par des considérations sur les billards projectifs (voir Corollaire G de \cite{taba_projectif}).

Dans la Section \ref{section_projective_caustics}, nous considérons le cas des caustiques pour des quadriques munies d'une structure de billard projectif. Précisons que dans le terme \textit{quadriques} sont aussi comprises les \textit{coniques}. Nous montrons le résultat suivant qui découle d'une construction proposée dans \cite{CKS} pour généraliser le théorème de Poncelet, mais qui ne mentionne pas les billards projectifs~:

\enonce{Proposition}{Soit $Q_1$ et $Q_2$ deux coniques ou quadriques distinctes. On peut munir un ouvert dense de $Q_1$ d'une structure de billard projectif de sorte que $Q_2$ est une caustique pour le billard projectif induit sur $Q_1$.}

\'Etant données deux quadriques $Q_1$ et $Q_2$ distinctes, on peut alors considérer le faisceau $\polar{\mathcal{F}}(Q_1,Q_2)$ de quadriques qui contient $Q_1$ et $Q_2$ et est défini ainsi par dualité~: l'ensemble des quadriques duales des quadriques de $\polar{\mathcal{F}}(Q_1,Q_2)$ est une droite qui contient les quadriques duales de $Q_1$ et $Q_2$ (dans l'espace des quadriques). On peut le voir comme une généralisation des faisceaux de quadriques homofocales. On prouve alors:

\enonce{Proposition}{Les quadriques de $\polar{\mathcal{F}}(Q_1,Q_2)$ sont des caustiques de $Q_1$ pour la structure de billard projectif induite par $Q_2$ sur $Q_1$. Toute quadrique de $\polar{\mathcal{F}}(Q_1,Q_2)$ induit la même structure projective sur $Q_1$ que celle induite par $Q_2$.}

En dimension au moins $3$, l'étude des billards classiques possédant des caustiques a été conclue par Berger \cite{berger_caustics} qui a énoncé un résultat dont les hypothèses sont beaucoup plus faibles que dans la conjecture de Birkhoff-Poritsky: Berger a montré que s'il existe des hypersurfaces $S$, $U$, $V$ de $\RR^d$, avec $d\geq 3$, ayant des secondes formes fondamentales non-dégénérées et telles qu'il existe un ouvert de droites tangentes à $U$ et intersectant $S$ qui sont réfléchies sur $S$ en des droites tangentes à $V$, alors $S$ est un morceau de quadrique, et $U,V$ sont des morceaux d'une seule et même quadrique homofocale. Ainsi la conjecture de Birkhoff-Poritsky est vérifiée dès l'existence d'au moins une caustique.

Dans la Section \ref{section_berger_property}, nous prouvons qu'un argument clé de la preuve de Berger peut se généraliser au cas des billards projectifs de $\RR^d$, avec toujours $d\geq3$, et nous l'avons appliqué pour généraliser le résultat de Berger aux billards pseudo-Euclidiens convexes:

\enonce{Proposition}{Soit $\Omega\subset\RR^d$, $d\geq 3$, un billard pseudo-Euclidien strictement convexe qui admet une caustique $\Gamma$. Alors $\partial\Omega$ est un ellipsoïde et $\Gamma$ est un morceau de quadrique homofocale pour la métrique pseudo-Euclidienne.}

L'argument de Berger que nous généralisons repose sur l'idée suivante. Soit $S\subset\RR^d$ une hypersurface, et $U,V$ comme dans l'énoncé de Berger cité plus haut. Toute droite $\ell$ de l'ouvert de droites tangentes à $U$, intersectant $S$ en $p$ et réfléchie en une droite $\ell'$ tangente à $V$, est telle que l'hyperplan tangent à $U$ contenant $\ell$ et l'hyperplan tangent à $V$ contenant $\ell'$ intersectent $T_pS$ en un même hyperplan $H$ de $T_pS$. Un tel hyperplan $H\subset T_pS$ est dit \textit{autorisé}, et l'argument de Berger est que pour $p$ fixé il y a au plus $d-1$ hyperplans autorisés. Nous montrons que dans le cas projectif, l'argument est encore valable génériquement (un sens plus précis sera donné à ce mot)~:

\enonce{Proposition}{Génériquement en un point de réflexion d'un billard projectif en dimension $\geq 3$, le nombre d'hyperplans autorisés est au plus $d-1$.}

Nous pensons que ce résultat, valable pour tout billard projectif, n'est pas applicable uniquement pour caractériser les billards pseudo-Euclidiens ayant des caustiques, mais peut-être encore pour d'autres billards. Peut-être permettrait-il au moins d'affirmer que si un billard projectif admet une caustique, alors cette caustique est une quadrique. Comme ce résultat semble délicat à démontrer, une première avancée pourrait consister à le prouver pour une classe assez générale de billards projectifs, ceux ayant un champ dit \textit{exact} de droites projectives et qui contient la classe des billards pseudo-Euclidiens, voir \cite{taba_projectif_ball}.

\subsection*{Détails du Chapitre \ref{chapter_ivrii}}

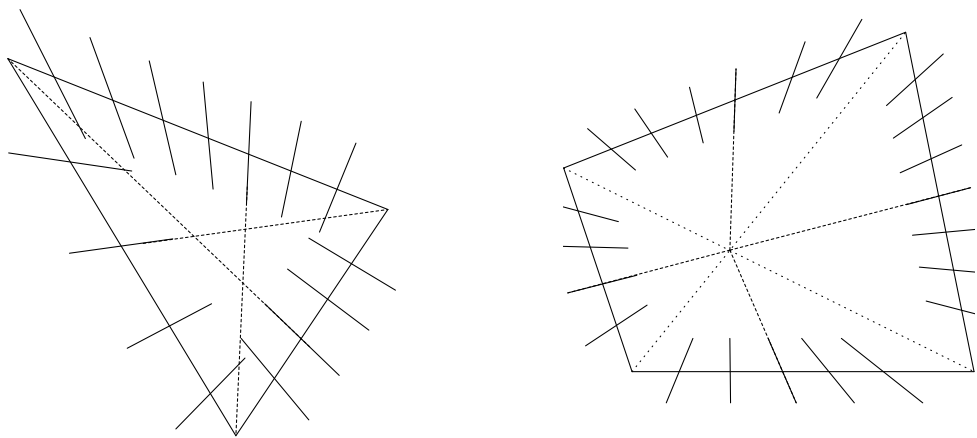
\begin{figure}[!h]
\centering
\begin{tikzpicture}[line cap=round,line join=round,>=triangle 45,x=1cm,y=1cm]
\clip(-3.4,-2.5) rectangle (2.4,4);
\draw (2,1)-- (0,-2);
\draw (0,-2)-- (-3,3);
\draw (-3,3)-- (2,1);
\draw (-2.84,3.65)-- (-1.98,1.94);
\draw (-1.92,3.28)-- (-1.34,1.68);
\draw (-1.14,2.97)-- (-0.79,1.46);
\draw (-0.43,2.69)-- (-0.3,1.27);
\draw (0.2,2.43)-- (0.14,1.09);
\draw (0.86,2.17)-- (0.6,0.9);
\draw (1.58,1.88)-- (1.1,0.7);
\draw (-2.99,1.75)-- (-1.37,1.51);
\draw (-2.19,0.42)-- (-0.83,0.61);
\draw (-1.43,-0.84)-- (-0.32,-0.25);
\draw (-0.79,-1.91)-- (0.12,-0.97);
\draw (0.07,-0.71)-- (0.96,-1.79);
\draw (1.36,-1.2)-- (0.38,-0.25);
\draw (0.68,0.21)-- (1.75,-0.6);
\draw (0.96,0.62)-- (2.1,-0.07);
\draw [line width=0.4pt,dash pattern=on 1pt off 1pt] (0.97,-0.82)-- (-3,3);
\draw [line width=0.4pt,dash pattern=on 1pt off 1pt] (-1.22,0.55)-- (2,1);
\draw [line width=0.4pt,dash pattern=on 1pt off 1pt] (0.15,1.3)-- (0,-2);
\begin{scriptsize}
\end{scriptsize}
\end{tikzpicture}
\hspace{1.5cm}
\begin{tikzpicture}[line cap=round,line join=round,>=triangle 45,x=0.9cm,y=0.9cm]
\clip(-3,-3.5) rectangle (3,3.5);
\draw (-2,-2)-- (-3,1);
\draw (-3,1)-- (2,3);
\draw (2,3)-- (3,-2);
\draw (-2,-2)-- (3,-2);
\draw [dotted] (-2,-2)-- (2,3);
\draw [dotted] (-3,1)-- (3,-2);
\draw [dash pattern=on 1pt off 1pt] (-0.57,-0.21)-- (-2.94,-0.84);
\draw [dash pattern=on 1pt off 1pt] (-0.57,-0.21)-- (-0.48,2.46);
\draw [dash pattern=on 1pt off 1pt] (-0.57,-0.21)-- (2.95,0.71);
\draw [dash pattern=on 1pt off 1pt] (-0.57,-0.21)-- (0.4,-2.46);
\draw (-3.39,0.53)-- (-2.2,0.21);
\draw (-3.17,-0.15)-- (-2.06,-0.18);
\draw (-2.94,-0.84)-- (-1.93,-0.58);
\draw (-2.68,-1.62)-- (-1.78,-1.02);
\draw (-1.5,-2.46)-- (-1.11,-1.51);
\draw (-0.57,-1.51)-- (-0.56,-2.46);
\draw (-0.01,-1.51)-- (0.4,-2.46);
\draw (1.25,-2.46)-- (0.48,-1.51);
\draw (1.06,-1.51)-- (2.25,-2.46);
\draw (3.33,-1.23)-- (2.3,-0.96);
\draw (3.2,-0.55)-- (2.2,-0.46);
\draw (3.07,0.1)-- (2.1,0.01);
\draw (2.01,0.46)-- (2.95,0.71);
\draw (2.82,1.34)-- (1.92,0.93);
\draw (1.82,1.44)-- (2.68,2.04);
\draw (2.55,2.68)-- (1.72,1.92);
\draw (1.36,3.19)-- (0.7,2.03);
\draw (0.53,2.86)-- (0.15,1.81);
\draw (-0.48,2.46)-- (-0.51,1.54);
\draw (-0.96,1.37)-- (-1.16,2.19);
\draw (-1.96,1.87)-- (-1.48,1.16);
\draw (-1.95,0.97)-- (-2.65,1.58);
\begin{scriptsize}
\end{scriptsize}
\end{tikzpicture}
\caption{\`A gauche, le billard projectif droit-sphérique obtenu à partir d'un exemple de billard $3$-réfléchissant sur la sphère, décrit dans \cite{bary_introuvable, VKNZ}. \`A droite, un exemple de billard projectif $4$-réfléchissant découvert au cours de cette thèse, voir Section \ref{section:examples_reflective_billiards}.}
\label{fig:right_spherical_intro}
\end{figure}

Dans ce chapitre, il est question d'étudier un analogue de la conjecture de Ivrii pour les billards projectifs. Une réponse immédiate peut être donnée à cette conjecture grâce à l'exemple déjà cité de billard $3$-réfléchissant sur la sphère $\SS^2$ \cite{bary_introuvable,VKNZ}. Il est en effet possible, en utilisant une projection centrale de la sphère sur un plan affine, d'interpréter ce billard comme un billard projectif, qui dès lors est lui-même $3$-réfléchissant. Cet exemple de billard projectif, appelé \textit{billard droit-sphérique} (voir la Figure \ref{fig:right_spherical_intro}), contredit tout de suite la conjecture de Ivrii pour les billards projectifs. 

On peut se demander s'il existe d'autres types de billards projectifs ayant des ouverts d'obites périodiques avec plus que trois réflexions. Cette thèse présente des exemples de billards projectifs dans des polygones qui sont $k$-réfléchissants pour le choix arbitraire d'un entier pair $k$ (\textit{cf} Section \ref{section:examples_reflective_billiards} et \cite{fierobe1}). Le caractère $k$-réfléchissant de ces billards vient de leur symétrie, symétrie du polygone ou du champ de droites projectives. Bien qu'ayant cherché, nous n'avons pu trouver des exemples "évidents" de billards projectifs $k$-réfléchissants avec $k$ impair, en dehors des billards droit-sphériques. On peut donc soulever la question de l'existence de billards projectifs $k$-réfléchissants dans les polygones, avec $k$ impair supérieur ou égal à $5$. Peut-être que la réponse à cette question pourrait s'inspirer de \cite{glut2}, qui montre que la conjecture de Ivrii pour des orbites de période impaire est vérifiée dans une classe assez générale de billards de bord algébrique par morceaux.

Ces exemples suggèrent donc de classifier les billards projectifs possédant des ensembles ouverts ou de mesure non-nulle d'orbites périodiques. L'avantage de cette démarche est de comprendre la conjecture de Ivrii pour d'autres billards. On pourra en tout premier lieu noter que l'existence d'un billard projectif $k$-réfléchissant fournit de nombreux exemples de billards projectifs ayant un ensemble de mesure non-nulle d'orbites $k$-périodiques par la construction suivante : étant donné un billard projectif $k$-réfléchissant ayant un ouvert $U$ d'orbites périodiques, tout billard qui coïncide avec le précédent sur un ensemble de Cantor de mesure non-nulle inclus dans $U$ possède un ensemble de mesure non-nulle d'orbites $k$-périodiques. Ainsi quand un billard $k$-réfléchissant existe, on pourra classifier uniquement les billards $k$-réfléchissants pour comprendre les obstructions à la conjecture de Ivrii. Cette thèse s'intéresse notamment au cas particulier des billards projectifs ayant des ensembles ouverts ou de mesure non-nulle d'orbites de période $3$. Elle prouve la classification suivante de ces billards en Section~\ref{section_3_reflective_proj_billiards}:

\enonce{Théorème}{\textbf{1)} Les seuls billards projectifs $3$-réfléchissants de $\RR^2$ de bord $\class^{\infty}$ par morceaux sont les billards droit-sphériques.\\
\textbf{2)} Si $d\geq 3$, il n'y a pas de billards projectifs dans $\RR^d$ de bord $\class^{\infty}$ par morceaux possédant un ensemble de mesure non-nulle d'orbites $3$-périodiques.}

La preuve de ce théorème est très largement inspirée de \cite{glut, glutkud2} et se décompose en deux étapes: il est d'abord question de traiter le résultat pour une version complexe des billards projectifs $3$-réfléchissants de bord analytique par morceaux, puis de l'élargir aux bords $\class^{\infty}$ en utilisant les systèmes Pfaffiens. Cette dernière étape est l'objet de la Section~\ref{section_pfaffian_systems}, dans laquelle sont introduits et étudiés des systèmes Pfaffiens relatifs aux billards projectifs et Euclidiens. 

L'utilité des systèmes Pfaffiens vient d'une idée de Barychnikov et Zharnitsky \cite{bary,bary2} d'associer un billard classique $k$-réfléchissant à une surface intégrale d'une certaine distribution, appelée \textit{distribution de Birkhoff}: pour les billards dans le plan, la distribution de Birkhoff est la distribution qui associe à un polygone non-dégénéré à $k$ côtés le produit cartésien de ses \textit{bissectrices extérieures} (c'est-à-dire les droites qui coupent en deux les deux angles extérieurs opposés formés par les droites supportant deux côtés consécutifs du polygone). Elle vérifie que si une surface intégrale de dimension $2$ de cette distribution est telle que la projection sur chaque sommet est une courbe lisse, alors ces courbes lisses forment des morceaux du bord d'un même billard $k$-réfléchissant. En effet, tout point de la surface intégrale est un polygone dont les bissectrices extérieures sont tangentes aux bords du billard, par définition de la distribution, et donc est une orbite de période $k$. Un système Pfaffien est alors un objet qui résume la donnée d'une distribution, de la dimension de ses variétés intégrales recherchées, et de conditions dites de transversalité, sur lequel peuvent être effectuées certaines opérations de \textit{prolongement} dans le but de trouver des surfaces intégrales. L'idée de Barychnikov et Zharnitsky a été reprise dans \cite{glutkud2}, où est conjecturé (Conjecture 5) l'énoncé suivant:

\enonce{Conjecture de Kudryashov}{Soient $k\geq 3$ et $d\geq 2$ deux entiers. Il existe un entier $r\geq 2$, dépendant uniquement de $k$ et $d$, tel que l'existence dans $\RR^d$ d'un billard de bord $\class^r$ par morceaux possédant un ensemble de mesure non-nulle d'orbites $k$-périodiques entraine l'existence d'un billard analytique par morceaux qui est $k$-réfléchissant.} 

Cette conjecture peut être résumée en disant que \textit{si la conjecture de Ivrii est fausse pour les billards de bord $\class^r$ par morceaux, alors il existe un billard analytique par morceaux $k$-réfléchissant}. Certains arguments présentés dans \cite{glutkud2} et dispersés dans l'article permettent de prouver un cas plus simple de cette conjecture en prenant $r=\infty$, mais ce résultat n'est malheureusement pas énoncé dans l'article. Comme il mérite d'être explicitement formulé, nous en donnons une preuve en Section~\ref{section_pfaffian_systems}, et dont l'essentiel des arguments provient de \cite{glutkud2}.  

\enonce{Théorème}{La conjecture de Kudryashov est valable pour $r=\infty$.} 

Nous prouvons de plus que si un billard $k$-réfléchissant de bord $\class^{\infty}$ par morceaux existe, alors pour tout entier $r\geq 1$ son bord peut être approché par des $r$-jets de billards $k$-réfléchissants de bord analytique par morceaux. Nous élargissons alors aussi au cas des billards projectifs la preuve de la conjecture de Kudryashov avec $r=\infty$  (\textit{cf} Section \ref{section_pfaffian_systems}), en prouvant le résultat suivant:

\enonce{Théorème}{S'il existe un billard projectif de bord $\class^{\infty}$ par morceaux (avec un champ de droites transverses $\class^{\infty}$ par morceaux) possédant un ensemble de mesure non-nulle d'orbites $k$-périodiques, alors il existe un billard projectif analytique $k$-réfléchissant.}

Ainsi ces arguments peuvent fournir des outils intéressants pour la résolution éventuelle de la conjecture de Ivrii~: se ramener aux cas des billards $k$-réfléchissants de bord analytique par morceaux ou bien étudier ces mêmes billards dans un cadre projectif. Généraliser peut parfois permettre de simplifier.

\section*{Perspectives de recherche}

Pour récapituler, le travail accompli pendant cette thèse a permis de mieux comprendre les billards projectifs ayant des ensembles de mesure non nulle d'orbites périodiques, de les classifier lorsqu'il s'agit en particulier des orbites triangulaires, de mettre en évidence des caustiques dites complexes du billard sur une conique complexifiée, de proposer des structures projectives sur des coniques et quadriques de sorte que ces dernières admettent des caustiques, et d'étendre un résultat de Berger pour les caustiques de billards projectifs en dimension au moins $3$ qui s'applique à la classification des billards pseudo-Euclidiens ayant des caustiques. Mais l'étude réalisée dans cette thèse n'est pas terminée et soulève peut-être plus de questions qu'elle n'apporte de réponses...

Le problème des billards projectifs admettant des caustiques en dimension $d\geq 3$ n'est que très partiellement résolu : certes un argument clé de Berger a pu être étendu à cette classe de billards, mais aucun résultat général similaire à celui de Berger n'a pu être prouvé, à part pour le cas très particulier des espaces pseudo-Euclidiens. Il serait intéressant de le généraliser à une classe plus vaste de billards projectifs, par exemples aux billards projectifs ayant un champ dit \textit{exact} de droites transverses \cite{taba_projectif_ball}. Peut-on avancer une conjecture ? Peut-être que les seules caustiques possibles d'un billard projectif en dimension $d\geq 3$ sont les quadriques. Je serais très curieux de connaître le résultat.

La conjecture de Ivrii est un problème majeur de théorie des billards. Sans chercher à en donner une réponse définitive, il pourrait être intéressant d'étudier des classes simples de billards projectifs $k$-réfléchissants. On pourrait par exemple essayer de savoir s'il existe des billards projectifs $k$-réfléchissants dans des polygones avec $k\geq 5$ impair. Notre recherche n'a en effet pas permis d'en trouver. On peut plus généralement se demander si les exemples de billards $k$-réfléchissants que nous présentons en Section \ref{section:examples_reflective_billiards} sont les seuls billards projectifs $k$-réfléchissants dans des polygones. Enfin il serait à envisager de comprendre si les arguments de classification des billards projectifs $3$-réfléchissants avancés par \cite{glut, glutkud2} et repris dans le chapitre \ref{chapter_ivrii} peuvent être synthétisés et généralisés à un nombre général de réflexions. 

\chapter*{Introduction in English}
\addcontentsline{toc}{chapter}{\protect\numberline{}Introduction in English}%
	
A billiard can be described as a dynamical system describing the trajectory of an infinitely small object without mass moving in a homogeneous domain bounded by a reflective boundary, like the trajectory of a ray of light inside a room covered by mirrors or of a particle. As stated by Valerii V. Kozlov et Dmitrii V. Treshchëv \cite{treshchev}: \go Starting with the works of G. D. Birkhoff, billiards have been a popular topic of investigation where various subjects of ergodic theory, Morse theory, KAM theory, etc. are intertwined. On the other hand, billiard systems are further remarkable in that they arise naturally in a number of important problems of mechanics and physics (vibro-impact systems, the diffraction of shortwaves, etc.). \gf The present manuscript investigates this field of research and present modest results about billiards.

The dynamic of the billiard trajectory is induced by the two following statements: 1) it moves along \textit{straight lines} inside the domain 2) and it is reflected on the boundary following the usual law of optics: \textit{angle before reflection = angle after reflection}. There are different ways to model statements 1) and 2), and the most common one consists of considering that the domain is inside a complete Riemannian manifold: the straight lines have to be understood as geodesics and the angles are defined by the metric. We can therefore study billiards in the usual plane, the space, on a hyperboloïd or on a sphere, when for example we study the movement of a small object inside a wide domain on the surface of a planet for which the planet's curvature cannot be neglected. However there are other models of billiards than this so-called \textit{classical} model, such as pseudo-Euclidean billiards, complex billiards, outer billiards or wire billiards. In this manuscript, we focus our attention to the so-called \textit{projective} billiards and \textit{complex} billiards. These billiards are linked with the classical billiard, as it will be shown in this thesis.

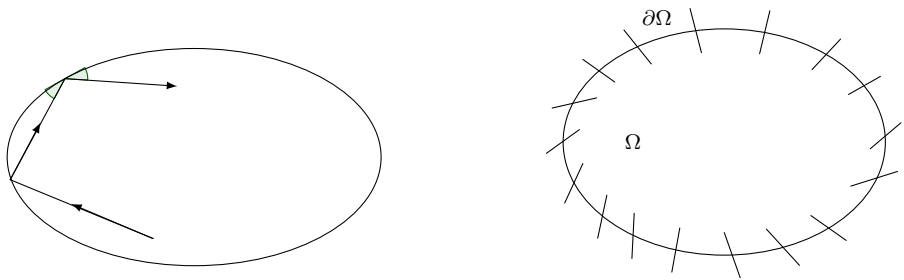
\begin{figure}[!h]
\centering
\definecolor{qqwuqq}{rgb}{0,0.39,0}
\begin{tikzpicture}[line cap=round,line join=round,>=triangle 45,x=2.0cm,y=2.0cm]
\clip(-1.4,-0.8) rectangle (1.4,0.8);
\draw [shift={(-0.85,0.52)},color=qqwuqq,fill=qqwuqq,fill opacity=0.1] (0,0) -- (-150.66:0.15) arc (-150.66:-117.76:0.15) -- cycle;
\draw [shift={(-0.85,0.52)},color=qqwuqq,fill=qqwuqq,fill opacity=0.1] (0,0) -- (-3.56:0.15) arc (-3.56:29.34:0.15) -- cycle;
\draw [rotate around={0:(0,0)}] (0,0) ellipse (2.46cm and 1.44cm);
\draw (-0.27,-0.54)-- (-1.21,-0.15);
\draw (-1.21,-0.15)-- (-0.85,0.52);
\draw [-latex] (-0.27,-0.54) -- (-0.81,-0.31);
\draw [-latex] (-1.21,-0.15) -- (-1.01,0.23);
\draw [-latex] (-0.85,0.52) -- (-0.11,0.47);
\begin{scriptsize}
\end{scriptsize}
\end{tikzpicture}
\hspace{1.5cm}
\begin{tikzpicture}[line cap=round,line join=round,>=triangle 45,x=1.5cm,y=1.5cm]
\clip(-1.6,-1.2) rectangle (1.6,1.2);
\draw [rotate around={0:(0,0)}] (0,0) ellipse (2.12cm and 1.5cm);
\draw (-1.56,-0.1)-- (-1.27,0.11);
\draw (-1.51,0.3)-- (-1.12,0.39);
\draw (-1.24,0.73)-- (-0.96,0.53);
\draw (-0.88,0.99)-- (-0.64,0.69);
\draw (-0.3,1.18)-- (-0.2,0.79);
\draw (0.4,1.16)-- (0.31,0.78);
\draw (1.02,0.91)-- (0.77,0.64);
\draw (1.37,0.58)-- (1.09,0.42);
\draw (1.55,0.17)-- (1.28,-0.06);
\draw (1.07,-0.87)-- (0.76,-0.65);
\draw (0.66,-1.09)-- (0.37,-0.77);
\draw (0.14,-1.2)-- (0,-0.8);
\draw (-0.46,-1.15)-- (-0.39,-0.76);
\draw (-0.81,-1.03)-- (-0.79,-0.63);
\draw (-1.11,-0.85)-- (-1.03,-0.47);
\draw (-1.4,-0.54)-- (-1.24,-0.19);
\draw (1.11,-0.39)-- (1.52,-0.26);
\begin{scriptsize}
\draw[color=black] (-0.59,1.1) node {$\partial\Omega$};
\draw[color=black] (-0.8,0) node {$\Omega$};
\end{scriptsize}
\end{tikzpicture}
\caption{On the left, a ray of light reflected on the boundary of a reflective domain. On the right, a projective billiard with its field of projective transverse lines.}
\label{figure:billiard_intro}
\end{figure}

\textit{Complex billiards} are a natural generalization of the classical billiards of the Euclidean plane $\RR^2$ to its complexification $\CC^2$. They were introduced and studied by Glutsyuk \cite{glut,glut1,glut2} to solve Ivrii's conjecture for $4$ reflections, the commuting billiard conjecture in dimension $2$, or Plakhov's invisibility conjecture (planar case with $4$ reflections). Combined to Pfaffian systems, complex billiards can be used to apply methods of complex analytic geometry to problems of standard (real) geometry. These points will be discussed in more details below.

\textit{Projective billiards} were introduced by Tabachnikov \cite{taba_projectif, taba_projectif_ball} as a generalization of classical billiards of the Euclidean space. A projective billiard is a bounded domain of a Euclidean space whose boundary is endowed with a field of transverse lines, called \textit{projective lines}. A trajectory is then reflected at a point on the boundary by a specific law of reflection depending on the projective line at the point of impact. When the latter projective line is orthogonal to the boundary, the reflection of the trajectory is the same as the usual law of optics. This statement is still valid for other billiards, like billiards in pseudo-Euclidean manifolds or in metrics projectively equivalent to the Euclidean one (which are metrics whose geodesics are supported by lines). Therefore, the model of projective billiards contain other models of billiards.

In the classical model of billiard inside a domain $\Omega$ bounded by a smooth boundary, the different trajectories can be mathematically described by two objects. The first one is the \textit{phase space} which is defined as the set of oriented geodesics between two points of reflection. It can be described as the set of pairs $(p,v)$ where $p$ is a point of the boundary $\partial\Omega$ and $v$ is a unit vector with origin at $p$, pointing inside $\Omega$ and representing the direction of the corresponding geodesic. In dimension $2$, $v$ can be replaced by the angle $\theta\in[0,\pi]$ it makes with the tangent line $T_p\partial\Omega$. The dimension of the phase space is $2$ for billiards in the plane, and $2(d-1)$ for billiards in a space of dimension $d$. The second object describing a billiard is the billiard map: it is a map associating to an element $(p,v)$ of the phase space representing a trajectory moving from $p$ in the direction given by $v$ the element $(q,v)$ where $q$ is the next point of impact of the trajectory and $w$ is the directing vector of the trajectory after reflection. Both objects have similar definitions for other billiard types.

\section*{Ivrii's conjecture}

One of the main issues of billiard theory is the study of \textit{periodic orbits}, which are trajectories repeating themselves after a finite number of reflections. Ivrii \cite{ivrii} showed in 1980 that the study of periodic orbits has an application in a famous problem which was summarized by Kac \cite{kac} in one question: \textit{Can one hear the shape of a drum ?} The problem is about to understand if the eigenvalues of the Laplacien with Dirichlet initial condtions in a bounded domain $\Omega\subset\RR^d$ determine completely the shape of $\Omega$. These eigenvalues are defined as the real numbers $\lambda\in\RR$ for which the system
\begin{equation}
\label{equation:probleme_dirichlet_en}
\left\{\begin{array}{l}
\Delta u +\lambda u=0\\
\restreint{u}{\partial\Omega}=0
\end{array}\right.
\end{equation}
has non-trivial solutions $u$. They can be interpreted physically as different vibration modes of a shape given by $\Omega$. Kac's question was answered negatively since examples of distinct shapes were given in which the corresponding Dirichlet problems \eqref{equation:probleme_dirichlet_en} have the same eigenvalues. However the question of recovering data about $\Omega$ from these eigenvalues is still investigated. Weyl \cite{weyl} showed that we can \textit{hear the volume\footnote{\go The first pertinent result is that one can hear the area of $\Omega$ \gf, \cite{kac}} of $\Omega$}, meaning that we can recover the volume of $\Omega$ from Dirichlet eigenvalues. Indeed, the eigenvalues of Dirichlet problem can be enumerated into a sequence $(\lambda_n)_n$ of real numbers such that $0\leq\lambda_1\leq\lambda_2\leq\ldots\leq\lambda_n\leq\ldots$ and $\lambda_n\to+\infty$. If we denote by $N(\lambda)$ the number of eigenvalues less or equal to $\lambda$, then Weyl showed that $N(\lambda)\sim(2\pi)^{-d}v_d\text{vol}({\Omega})\lambda^{d/2}$, where $v_d$ denotes the volume of the unit Euclidean sphere in $\RR^d$. He also conjectured the second asymptotic term
\begin{equation}
\label{conjecture:weyl_en}
N(\lambda)=(2\pi)^{-d}v_d\text{vol}({\Omega})\lambda^{d/2}-\frac{1}{4(2\pi)^{d-1}}\text{area}(\partial\Omega)\lambda^{(d-1)/2}+o(\lambda^{(d-1)/2}).
\end{equation}
This conjecture is not proven yet although many results exist and confirm Weyl's conjecture. One of them is a result due to Ivrii \cite{ivrii} who proved that \eqref{conjecture:weyl_en} is satisfied under the assumption that the billiard inside $\Omega$ has a \textit{few} periodic orbits, meaning that the set of parameters in the phase space corresponding to periodic orbits has zero measure in $\Omega$. A famous conjecture was stated following this result:

\enonce{Ivrii's conjecture}{Given a bounded domain in the Euclidean space with sufficiently smooth boundary, its set of periodic orbits has zero measure.}

This conjecture still holds and is more difficult than it was expected at the beginning. Particular cases of billiards with a set of positive measure of periodic orbits are given by the so-called $k$-reflective billiards: billiards having open subsets of periodic orbits of period $k$, more precisely having open subsets in its phase space of parameters $(p,v)$ corresponding to periodic orbits. The existence of a $k$-reflective billiard is still unknown, but could lead to a rather curious construction: a room whose walls are covered by mirrors and such that there is a place in the room where any observer can still see himself from behind, even by moving or turning a little round. 

There is still no definitive answer to Ivrii's conjecture, even for $k$-reflective billiards with any integer $k$. Many partial results however already exist. Petkov and Stojanov \cite{petkovstojanov} proved it for generic billiards: the set of all domains in $\RR^d$ with $\class^{\infty}$-smooth boundary having a finite number of periodic orbits of period $k$ for all $k$ contains a residual set (a countable intersection of open dense subsets). Another answer was given by Vasiliev \cite{vasiliev} who proved the conjecture for a convex domain with analytic boundary. Rychlik \cite{rychlik} and then Stojanov \cite{stojanov} proved that the set of periodic orbits of period $3$ has zero measure in any billiard of the Euclidean plane with $\class^3$-smooth boundary. Vorobets \cite{vorobets} extended this result to billiards in any dimension. Later, Wojtkowski \cite{wojtkowski}, and then Baryshnikov and Zharnitsky \cite{bary} gave new proofs of this result. More recently, Glutsyuk and Kudryashov \cite{glutkud2} proved the conjecture for periodic orbits of period $4$ in planar billiards with $\class^4$-smooth boundary. Thus in the Euclidean case, Ivrii's conjecture remains unproved for any period and any regularity of the boundary (even for billiards with piecewise-analytic boundary).  

\begin{figure}[!h]
\centering
\includegraphics[scale=0.4]{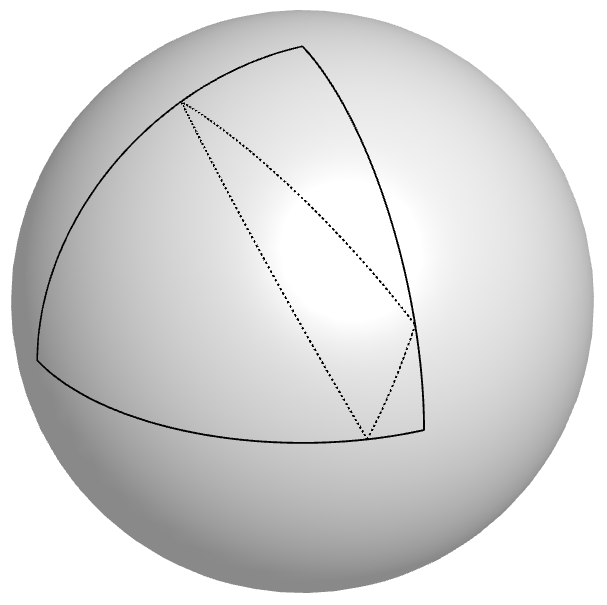}
\caption{An example of $3$-reflective billiard on the sphere presented by Barychnikov. The outer triangle is boundary of the billiard, the interior triangle in dotted lines is an orbit. Two vertices of the orbit can be moved arbitrarily without changing its periodicity.}
\label{figure:sphere_bary_en}
\end{figure}

Ivrii's conjecture can be stated analogously for non-Euclidean billiards, such as billiards in manifolds of constant curvature, on a sphere or on a hyperboloïd. Remarkable examples of $2$- and $3$-reflective billiards can be given on the $2$-dimensional sphere $\SS^2$ \cite{bary_introuvable,VKNZ}, which are linked with the existence of points joined by an infinite number of geodesics, see Figure \ref{figure:sphere_bary_en}. The cited articles give a classification of billiards on the unit sphere $\SS^2$ having a set of non-zero measure of periodic orbit of period $3$. They also prove that Ivrii's conjecture for $3$-periodic orbits is also true for billiards on the hyperboloïd.

\section*{Integrable billiards}

An other important issue of billiard theory is the study of the so-called \textit{integrable billiards}. A billiard $\Omega$ of the plane is said to be globally integrable if its phase space is foliated by smooth closed curves invariant by the billiard map. $\Omega$ is said to be locally integrable if such a foliation exists only in neighborhood of the curve $\{\theta=0\}$ in the phase space. This property is strongly linked with the existence of caustics corresponding to these invariant curves, and which can be defined independantly in all dimensions: a caustic of a billiard $\Omega$ is a hypersurface $\Gamma\subset\Omega$ such that any line tangent to $\Gamma$ and intersecting the boundary $\partial\Omega$ at $p$ is reflected into a line tangent to $\Gamma$ after reflection at $p$ on $\partial\Omega$.

An example of globally integrable billiard is the disk, since any concentric circle inside the disk is a caustic of the corresponding billiard. An ellipse is an example of a locally integrable billiard, since any billiard trajectory which do not passes between the foci of the ellipse remains tangent to a smaller confocal ellipse. Birkhoff and Poritsky asked if these examples are the only such examples of locally integrable billiard, and this question is now cited as a famous conjecture, as it is recalled in \cite{kaloshinsorrentino}.

\enonce{Birkhoff-Poritsky conjecture}{If a billiard is locally integrable, then it is an ellipse.}

Major results were discovered about this conjecture. Bialy \cite{bialy} proved that if the phase space of the billiard $\Omega$ is foliated by not null-homotopic continuous invariant closed curves, then $\partial\Omega$ is a circle. Notice that this result requires the foliation to be global and implies that the only globally integrable billiard is the circle. An algebraic proof of Birkhoff-Poritsky conjecture for planar billiards and billiards on surfaces of constant curvature was found by Bialy, Glutsyuk and Mironov \cite{bialymironov_poly, bialymironov_poly_bis, glut_integrability, glut_integrability_bis}. Kaloshin and Sorrentino \cite{kaloshinsorrentino} showed that any integrable deformation of an ellipse is an ellipse. In greater dimension, the study of billiards having caustics was ended earlier by Berger \cite{berger_caustics} who proved that if a billiard $\Omega$ in $\RR^d$, with $d\geq 3$, has a caustic, then $\partial\Omega$ is a quadric and its caustic is a confocal quadric. The assumptions of this result are weaker, and they do not require the existence of a foliation.

\section*{Results obtained in this thesis}
\addcontentsline{toc}{section}{\protect\numberline{}Results obtained in this thesis}%

This manuscript pesents different results about complex and projective billiards which some of them can also be applied to classical billiards. It is structured in three chapters: \textbf{Chapter \ref{chapter_billiards}} exposes in details both models of complex and projective billiards. \textbf{Chapter \ref{chapter_caustics}} study the existence of caustics for different billiards of both types. \textbf{Chapter \ref{chapter_ivrii}} is focused on the analogue of Ivrii's conjecture for projective billiards.

\subsection*{Details of Chapter \ref{chapter_billiards}}

This chapter presents two types of billiards studied all along this manuscript: the complex and projective billiards. We present here briefly the definitions of these billiards to understand the overviews of each chapter.

A projective billiard is a bounded domain $\Omega$ of $\RR^d$ whose boundary is smooth and endowed with a smooth field of transverse lines. This field of lines induces at each point $p\in\partial\Omega$ of the boundary a transformation of the field of oriented lines containing $p$, and which allows to construct billiards orbits: an oriented line $\ell_0$ intersecting $\Omega$ at a point $p$ is reflected by previous transformation at $p$ into a line $\ell_1$. If $\ell_1$ intersect $\partial\Omega$ in another point, this construction can be repeated, and so on.

A complex billiard is a complex curve $\gamma$ of $\CP^2$ on which we can also define a law of reflection on lines intersecting it. This construction can be realised by considering the complexification of the Euclidean metric $\dd x^2+\dd y^2$ to $\CC^2$. Given a so-called \textit{non-isotropic} complex line $L\subset\CC^2$, one can define a symetry of complex lines with respect to $L$ as the unique non-trivial affine involution preserving the latter complex quadratic form and fixing the points of the line $L$. Two complex lines $\ell,\ell'$ intersecting $\gamma$ at a point $p$ are said to be symetric for the complex reflection law if the symetry of lines with respect to the tangent lines $T_p\gamma$ sends $\ell'$ to $\ell$ or $\ell$ to $\ell'$.

\subsection*{Details of Chapter \ref{chapter_caustics}}

In this chapter, we present results related to the existence of caustics in projective and complex billiards. Section \ref{section:general_properties_on_quadrics} describes a first result on the so-called \textit{complex caustics} of an ellipse or hyperbola. We say that a conic $C'\subset\CP^2$ is a complex caustic of another conic $C\subset\CP^2$ if any line $\ell$ tangent to $C'$ is reflected into a line tangent to $C$ by the complex law of reflection in one of the intersection point of $\ell$ with $C$. Given $a,b\in\polar{\RR}$, we introduce the set $(\mathcal{C}_{\lambda})_{\lambda\in\CC}$ of conics of $\CP^2$ given by the equation
$$\mathcal{C}_{\lambda}:\frac{x^2}{a-\lambda}+\frac{y^2}{b-\lambda}=1$$
and we study the complex billiard defined by $\mathcal{C}_0$. It is known that in the case of the usual billiard on the real conic $\mathcal{C}_0$, the real conics $\mathcal{C}_{\lambda}$ are caustics. We answer the question if this is still true for the complex billiard, and which are the conics inscribed in periodic orbits. We prove the following results:

\enonce{Proposition}{Any conic $\mathcal{C}_{\lambda}$ is a complex caustic of $\mathcal{C}_0$.}

\enonce{Proposition}{Let $n\geq 3$ be an integer. There is a polynomial in $(a,b,\lambda)$, denoted by $\mathcal{B}^n_{a,b}(\lambda)$, whose complex roots in $\lambda$ corresponds to the caustics $\mathcal{C}_{\lambda}$ inscribed in periodic orbits of period $n$. For all $(a,b)$ outside a finite number of values of $a/b$, the degree in $\lambda$ of the polynomial $\mathcal{B}^n_{a,b}(\lambda)$ is $(n^2-1)/4$ if $n$ is odd, and $n^2/4-1$ if $n$ is even.}

Thus the distinct roots in $\lambda$ of $\mathcal{B}^n_{a,b}(\lambda)$ different from $a$ and $b$ corresponds to the complex caustics of $\mathcal{C}_{0}$ inscribed in such periodic orbits of period $n$. We were able to show that for a generic number of pairs $(a,b)$ (in the sense of previous result), neither $a$ nor $b$ are roots (in $\lambda$) of $\mathcal{B}^n_{a,b}(\lambda)$. It remains to understand if $\mathcal{B}^n_{a,b}(\lambda)$ has generically simple roots in $\lambda$ or not. For now, the result is still unknown, but is true for small periods. And a surprising phenomenon appears for period $3$ when $\mathcal{C}_0$ is an ellipse (and similar results have been achieved for a hyperbola or periodic orbits of period $4$):

\enonce{Proposition}{If $a,b>0$, there are exactly two complex conics confocal to $\mathcal{C}_0$ which are inscribed in periodic orbits of period $3$. They are complexified ellipses: one of them $\mathcal{C}_i$ is included in $\mathcal{C}_0$, the other one $\mathcal{C}_e$ \textit{contains} $\mathcal{C}_0$ (see Figure \ref{figure:caustics_complex_en}).}

By curiosity, we looked for specific billiards properties of these ellipse, like the possibility for $\mathcal{C}_0$ or $\mathcal{C}_i$ to be a caustic of $\mathcal{C}_e$ inscribed in periodic orbits of the classical billiard. But simulations failed to show such eventual curious result. We then show that an invariant of the real elliptic billiard known as \textit{Joachimsthal invariant} can be generalized to the complex billiard.

This thesis then presents a result related to the existence of caustics in projective billiards. Let us first note that numerous results were obtained by Tabachnikov \cite{taba_projectif_ball,taba_projectif} on the existence of area forms of the phase space invariant by the projective billiard map, and on their consequences about the integrability of the billiard. For example Corollary F of \cite{taba_projectif} states that \textit{if the projective billiard inside a circle has an invariant area form smooth up to the boundary, then the billiard is integrable.} Note also that a new proof of the integrability of the elliptic billiard in the Euclidean plane, on the sphere or on a hyperboloid was given using considerations about projective billiards (see Corollary G of \cite{taba_projectif}).

In Section \ref{section_projective_caustics}, we investigate the existence of caustics for quadrics endowed with a structure of projective billiard. Let us precise that in the following results the term quadric contains the conics. We show the following result which is a consequence of a construction contained in \cite{CKS} to generalize Poncelet theorem, but the latter does not mention the projective billiards:

\enonce{Proposition}{Let $Q_1$ and $Q_2$ be two distinct conics or quadrics. There is an open dense subset of $Q_1$ which can be endowed with a structure of projective billiard such that $Q_2$ is caustic of the corresponding projective billiard on $Q_1$.}

Given two distinct quadrics $Q_1$ and $Q_2$, we can consider the pencil of quadrics $\polar{\mathcal{F}}(Q_1,Q_2)$, which contains $Q_1$ and $Q_2$ and is defined by duality: the dual quadrics of the quadrics contained in $\polar{\mathcal{F}}(Q_1,Q_2)$ is a line containing the dual quadrics of $Q_1$ and $Q_2$ (in the space of quadrics). We can interpret $\polar{\mathcal{F}}(Q_1,Q_2)$ as a generalization of the notion of pencil of confocal quadrics. Then we prove:

\begin{figure}
\centering
\begin{tikzpicture}[line cap=round,line join=round,>=triangle 45,x=1cm,y=1cm]
\clip(-3.5,-2.7) rectangle (3.5,2.7);
\draw [rotate around={0:(0.0025086145409116416,0)},line width=1.5pt] (0.0025086145409116416,0) ellipse (1.9974913854590879cm and 1.0432433360133084cm);
\draw [rotate around={0:(0.0025086145409115726,0)},line width=1pt] (0.0025086145409115726,0) ellipse (2.9974913854590883cm and 2.4664426668897765cm);
\draw [rotate around={0:(0.002508614540915727,0)},line width=1pt] (0.002508614540915727,0) ellipse (1.709768500390961cm and 0.1472859398657943cm);
\begin{scriptsize}
\draw[color=black] (-1.2,1.1) node {$\mathcal{C}_0$};
\draw[color=black] (-1.5597148670526204,1.8) node {$\mathcal{C}_{e}$};
\draw[color=black] (-0.7111058200015711,0.31037283751140887) node {$\mathcal{C}_{i}$};
\end{scriptsize}
\end{tikzpicture}
\hspace{1cm}
\begin{tikzpicture}[line cap=round,line join=round,>=triangle 45,x=1cm,y=1cm]
\clip(-3,-2.7) rectangle (3,2.7);
\draw [rotate around={0:(0.004465809016780549,0)},line width=1.5pt] (0.004465809016780549,0) ellipse (1.9955341909832196cm and 1.0008587655128343cm);
\draw [rotate around={0:(0.004465809016780692,0)},line width=1pt] (0.004465809016780692,0) ellipse (1.784232681737871cm and 0.45060828188388125cm);
\draw [rotate around={0:(0.004465809016780525,0)},line width=1pt] (0.004465809016780525,0) ellipse (2.2943056805364748cm and 1.5110923588129084cm);
\draw [samples=50,domain=-0.99:0.99,rotate around={0:(0.0044658090167803855,0)},xshift=0.0044658090167803855cm,yshift=0cm,line width=1pt] plot ({1.6329479850417183*(1+(\x)^2)/(1-(\x)^2)},{0.5602850319501439*2*(\x)/(1-(\x)^2)});
\draw [samples=50,domain=-0.99:0.99,rotate around={0:(0.0044658090167803855,0)},xshift=0.0044658090167803855cm,yshift=0cm,line width=1pt] plot ({1.6329479850417183*(-1-(\x)^2)/(1-(\x)^2)},{0.5602850319501439*(-2)*(\x)/(1-(\x)^2)});
\begin{scriptsize}
\draw[color=black] (-0.9777216831311774,1.1) node {$\mathcal{C}_0$};
\end{scriptsize}
\end{tikzpicture}
\caption{On the left, an ellipse $\mathcal{C}_0$ with its two caustics $\mathcal{C}_i$ and $\mathcal{C}_e$ inscribed in triangular orbits. These are two complexified ellipses, one of them is included in $\mathcal{C}_0$ and the other one contains it. The graphic represents their real parts. On the right, the complex caustics of $\mathcal{C}_0$ for periodic orbits of period $4$.}
\label{figure:caustics_complex_en}
\end{figure}
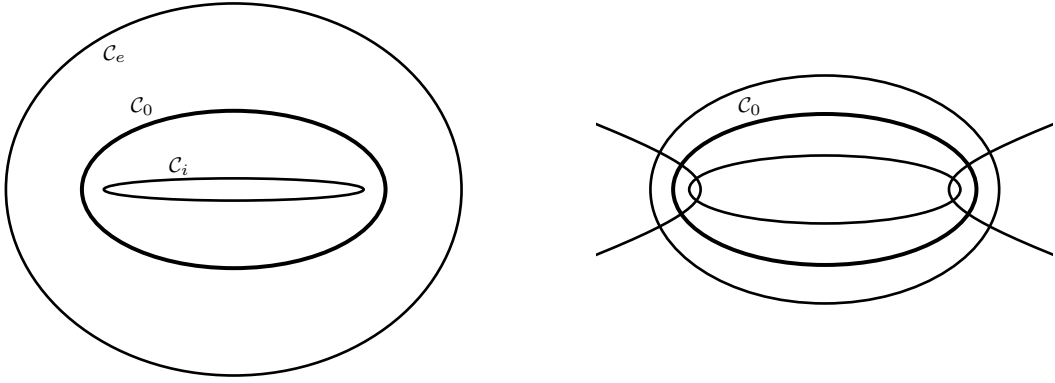

\enonce{Proposition}{The quadrics of $\polar{\mathcal{F}}(Q_1,Q_2)$ are caustics of $Q_1$ for the structure of projective billiard induced by $Q_2$ on $Q_1$. Any quadric of $\polar{\mathcal{F}}(Q_1,Q_2)$ induces the same projective structure on $Q_1$ as the one induced by $Q_2$.}

In dimension greater than $2$, the study of billiards having caustics has been ended by Berger \cite{berger_caustics}. He stated a result whose assumptions are weaker than Birkhoff-Poritsky conjecture: Berger showed that if there are hypersurfaces $S$, $U$, $V$ of $\RR^d$, with $d\geq 3$, having non-degenerate second fundamental forms and such that there is an open subset of lines tangent to $U$ and intersecting $S$ which are reflected by $S$ in lines tangent to $V$, then $S$ is a piece of quadric, and $U,V$ are pieces of one and the same confocal quadric. 

We prove at Section \ref{section_berger_property} that a key argument of Berger's proof can be generalized to projective billiards of $\RR^d$, $d\geq3$, and we apply it to generalize Berger's result to pseudo-Euclidean billiards:

\enonce{Theorem}{Let $\Omega\subset\RR^d$, $d\geq 3$, be a strictly convex pseudo-Euclidean billiard having a caustic $\Gamma$. Then $\partial\Omega$ is an ellipsoid and $\Gamma$ is a piece of quadric which is confocal for the pseudo-Euclidean metric.}

The argument of Berger we generalize can be described as follows. Let $S\subset\RR^d$ be a hypersurface and $U,V$ be as in the previous mentionned result of Berger. Any line $\ell$ of the open subset of lines tangent to $U$, intersecting $S$ at $p$ and reflected in a line $\ell'$ tangent to $V$, is such that the hyperplane tangent to $U$ containing $\ell$ and the hyperplane tangent to $V$ containing $\ell'$ intersect $T_pS$ in the same hyperplane $H$ of $T_pS$. Such hyperplane $H\subset T_pS$ is said to be \textit{permitted}. Berger's key argument states that for a fixed $p$ there are at most $d-1$ such permitted hyperplanes. We show that in the case of projective billiards, this argument is still satisfied \textit{generically} (a more precise meaning to this word will be given later):

\enonce{Proposition}{Generically at a point of reflection of a projective billiard in dimension $d\geq 3$, the number of permitted hyperplanes is at most $d-1$.}

We think that this result is applicable not only to pseudo-Euclidean billiards. Maybe it could be used at least to show that if a projective billiard has a caustic, then this caustic is a quadric. A first step would consist for example in proving it for a wider class of projective biliards containing pseudo-Euclidean billiards, and called projective billiards with exact transverse line fields, see \cite{taba_projectif_ball}.

\subsection*{Details of Chapter \ref{chapter_ivrii}}

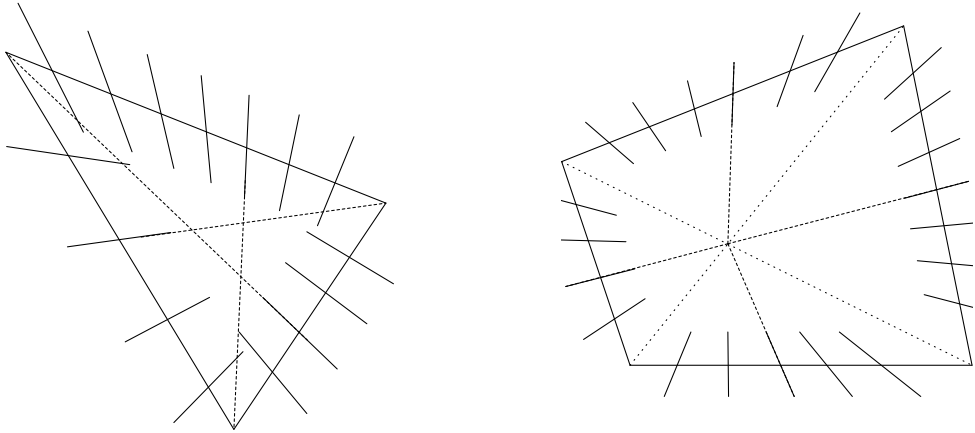
\begin{figure}[!h]
\centering
\begin{tikzpicture}[line cap=round,line join=round,>=triangle 45,x=1cm,y=1cm]
\clip(-3.4,-2.5) rectangle (2.4,4);
\draw (2,1)-- (0,-2);
\draw (0,-2)-- (-3,3);
\draw (-3,3)-- (2,1);
\draw (-2.84,3.65)-- (-1.98,1.94);
\draw (-1.92,3.28)-- (-1.34,1.68);
\draw (-1.14,2.97)-- (-0.79,1.46);
\draw (-0.43,2.69)-- (-0.3,1.27);
\draw (0.2,2.43)-- (0.14,1.09);
\draw (0.86,2.17)-- (0.6,0.9);
\draw (1.58,1.88)-- (1.1,0.7);
\draw (-2.99,1.75)-- (-1.37,1.51);
\draw (-2.19,0.42)-- (-0.83,0.61);
\draw (-1.43,-0.84)-- (-0.32,-0.25);
\draw (-0.79,-1.91)-- (0.12,-0.97);
\draw (0.07,-0.71)-- (0.96,-1.79);
\draw (1.36,-1.2)-- (0.38,-0.25);
\draw (0.68,0.21)-- (1.75,-0.6);
\draw (0.96,0.62)-- (2.1,-0.07);
\draw [line width=0.4pt,dash pattern=on 1pt off 1pt] (0.97,-0.82)-- (-3,3);
\draw [line width=0.4pt,dash pattern=on 1pt off 1pt] (-1.22,0.55)-- (2,1);
\draw [line width=0.4pt,dash pattern=on 1pt off 1pt] (0.15,1.3)-- (0,-2);
\begin{scriptsize}
\end{scriptsize}
\end{tikzpicture}
\hspace{1.5cm}
\begin{tikzpicture}[line cap=round,line join=round,>=triangle 45,x=0.9cm,y=0.9cm]
\clip(-3,-3.5) rectangle (3,3.5);
\draw (-2,-2)-- (-3,1);
\draw (-3,1)-- (2,3);
\draw (2,3)-- (3,-2);
\draw (-2,-2)-- (3,-2);
\draw [dotted] (-2,-2)-- (2,3);
\draw [dotted] (-3,1)-- (3,-2);
\draw [dash pattern=on 1pt off 1pt] (-0.57,-0.21)-- (-2.94,-0.84);
\draw [dash pattern=on 1pt off 1pt] (-0.57,-0.21)-- (-0.48,2.46);
\draw [dash pattern=on 1pt off 1pt] (-0.57,-0.21)-- (2.95,0.71);
\draw [dash pattern=on 1pt off 1pt] (-0.57,-0.21)-- (0.4,-2.46);
\draw (-3.39,0.53)-- (-2.2,0.21);
\draw (-3.17,-0.15)-- (-2.06,-0.18);
\draw (-2.94,-0.84)-- (-1.93,-0.58);
\draw (-2.68,-1.62)-- (-1.78,-1.02);
\draw (-1.5,-2.46)-- (-1.11,-1.51);
\draw (-0.57,-1.51)-- (-0.56,-2.46);
\draw (-0.01,-1.51)-- (0.4,-2.46);
\draw (1.25,-2.46)-- (0.48,-1.51);
\draw (1.06,-1.51)-- (2.25,-2.46);
\draw (3.33,-1.23)-- (2.3,-0.96);
\draw (3.2,-0.55)-- (2.2,-0.46);
\draw (3.07,0.1)-- (2.1,0.01);
\draw (2.01,0.46)-- (2.95,0.71);
\draw (2.82,1.34)-- (1.92,0.93);
\draw (1.82,1.44)-- (2.68,2.04);
\draw (2.55,2.68)-- (1.72,1.92);
\draw (1.36,3.19)-- (0.7,2.03);
\draw (0.53,2.86)-- (0.15,1.81);
\draw (-0.48,2.46)-- (-0.51,1.54);
\draw (-0.96,1.37)-- (-1.16,2.19);
\draw (-1.96,1.87)-- (-1.48,1.16);
\draw (-1.95,0.97)-- (-2.65,1.58);
\begin{scriptsize}
\end{scriptsize}
\end{tikzpicture}
\caption{On the left, the right-spherical billiard obtained from an example of $3$-reflective billiard on the sphere, as described in \cite{bary_introuvable, VKNZ}. On the right, an example of $4$-reflective projective billiard presented in this manuscript,see Section \ref{section:examples_reflective_billiards}.}
\label{fig:right_spherical_intro_en}
\end{figure}

We study in this chapter the analogue of Ivrii's conjecture for projective billiards. A first answer can be given thanks to the above mentionned example of $3$-reflective billiard on the unit sphere $\SS^2$ \cite{bary_introuvable,VKNZ}. Indeed, a central projection from the sphere onto an affine plane projects such $3$-reflective billiard into a $3$-reflective projective billiard of the plane. This example of projective billiard, called \textit{right-spherical billiard} (see Figure \ref{fig:right_spherical_intro_en}), immediately contradicts Ivrii's conjecture for projective billiards.

We can ask if there are other examples of projective billiards having open subsets of periodic orbits with more than $3$ reflections. This thesis presents examples of projective billiards inside polygons which are $k$-reflective for any choice of an arbitrary even integer $k$ (\textit{cf} Section \ref{section:examples_reflective_billiards} and \cite{fierobe1}). Their $k$-reflectivity comes from the particular symmetry of the polygons and of their projective fields of lines. We were unable to find other examples of $k$-reflective billiards with an odd $k$. We can ask the question wether there exist or not $k$-reflective billiards in polygons with an odd $k\geq 5$. Maybe the answer to this question could use a similar argument to \cite{glut2}, which prove Ivrii's conjecture for periodic orbits of odd periods inside billiards with piecewise algebraic boundary.

These examples suggest to classify the projective billiards having open subsets or subsets of non-zero measure of periodic orbits. The benefit of this method is to understand Ivrii's conjecture in other geometries. We can first note that the existence of a $k$-reflective projective billiard gives numerous examples of projective billiards having a subset of non-zero measure of $k$-periodic orbits by the following construction: given a $k$-reflective projective billiard having an open subset $U$ of $k$-periodic orbits, any billiard which coincide with the first one on a Cantor set of positive measure included in $U$ has a subset of non-zero measure of $k$-periodic orbits. Therefore we can focus on classifying $k$-reflective projective billiards only, as soon as a $k$-reflective billiard already exists. This manuscript gives a classification of billiards having open subsets of periodic orbits (in dimension $2$) and subset of non-zero measure of periodic orbits (in dimension $d \geq 3$):

\enonce{Proposition}{\textbf{1)} The only $3$-reflective projective billiard of $\RR^2$ with piecewise $\class^{\infty}$-smooth boundary is the right-spherical billiard.\\
\textbf{2)} If $d\geq 3$, there is no projective billiard in $\RR^d$ with $\class^{\infty}$-smooth boundary having a set of non-zero measure of $3$-periodic orbits.}

The proof of this theorem is widely inspired from \cite{glut, glutkud2} and can be decomposed in two steps: we first study a complex version of $3$-reflective projective billiards with piecewise analytic boundary, then we extend the result to $\class^{\infty}$-smooth boundary using the theory of Pfaffian systems. This last step is presented in Section~\ref{section_pfaffian_systems}, in which Pfaffian systems related to projective and Eulidean billiards are introduced and studied. 

Pfaffian systems are a tool based on analytic distribution, and their application to billiard theory can be attributed to Barychnikov and Zharnitsky \cite{bary,bary2}: they had the idea to associate to a $k$-reflective billiard an intergal surface of a certain distribution, called \textit{Birkhoff's distribution}. In the case of planar billiards, Birkhoff's distribution is the distribution associating to a non-degenerate $k$-sided polygon the cartesian product of its outer bisectors (which are the lines splitting in half the outer opposite angles formed by the lines supporting two consecutive sides of the polygon). Thus, if a $2$-dimensional integral surface of Birkhoff's distribution is such that its projections onto each vertex are smooth curves, then these smooth curves are on the boundary of a $k$-reflective billiard. Indeed, any point of the integral surface is a polygon whose outer bisectors are tangent to the boundary of the billiard, by definition of the distribution, hence is a $k$-periodic orbit. A Pfaffian system is then an object which contains the data of a distribution, the dimension of its integral surfaces of interest, and some transversality conditions, on which can be applied what are called \textit{prolongations} in order to find intergal surfaces. Barychnikov and Zharnitsky's idea was also used in \cite{glutkud2}, where the following conjecture (Conjecture 5) is stated:

\enonce{Kudryashov's conjecture}{Let $k\geq 3$ and $d\geq 2$ be integers. There is an integer $r\geq 2$, uniquely depending on $k$ and $d$, such that if there is a piecewise $\class^r$-smooth billiard in $\RR^d$ having a set of non-zero measure of $k$-periodic orbits, then there is a $k$-reflective billiard with piecewise analytic boundary.} 

This conjecture can be understood as follows: \textit{If Ivrii's conjecture is false for billiards with piecewise $\class^r$-smooth boundary, then there is a $k$-reflective billiard with piecewise analytic boundary}. Some arguments of \cite{glutkud2} can be used to prove the case $r=\infty$, but the corresponding result is not mentioned. In our opinion, it is a remarkable result which needs to be explicitely formulated. Hence we give a complete proof of it in Section~\ref{section_pfaffian_systems}, whose arguments comes from \cite{glutkud2}.

\enonce{Theorem}{Kudryashov's conjecture holds for $r=\infty$.} 

We also prove that if a $k$-reflective billiard with piecewise $\class^{\infty}$-smooth boundary exists, then for any integer $r\geq 1$ its boundary can be approwimated by $r$-jets of $k$-reflective billiards with piecewise analytic boundary. We further extend this proof to the class of projective billiards (\textit{cf} Section \ref{section_pfaffian_systems}):

\enonce{Theorem}{If there is a piecewise $\class^{\infty}$-smooth projective billiard (with a piecewise $\class^{\infty}$-smooth field of transverse lines) having a subset of non-zero measure of periodic orbits, then there is a piecewise analytic $k$-reflective projective billiard.}

These arguments can give interesting tools towards the possible resolution of Ivrii's conjecture, like for example studying the more simple case of $k$-reflective billiards with piecewise analytic boundary, or studying these billiards in the class of projective billiards. Generalizations could maybe lead to simplifications.

\section*{Perspectives}

To conclude, the main results obtained during this thesis helped to better understand projective billiards with sets of non-zero measure of periodic orbits, to classify them in the particular case of $3$-periodic orbits, to expose so-called \textit{complex} caustics of the elliptic billiard, to show the existence of projective billiard structures on conics and quadrics so that the latter admit caustics, and to generalize a result of Berger to projective billiards in dimension at least $3$, which was applied to classify pseudo-Euclidean billiards having caustics. Nevertheless, the study realised during this thesis is not over and raises maybe more questions than it gives answers...

The problem of projective billiards having caustics in dimension $d\geq 3$ has only partial answers: a key argument of Berger was succesfuly generalized to projective billiards, but the result of Berger was itself generalized only to a small class of projective billiards (the pseudo-Euclidean ones). It could be interesting to find a more general class of billiards in which this result can be proven to be true, for example the so-called \textit{projective billiards with exact transverse line fields} \cite{taba_projectif_ball}. We can maybe state a conjecture: possibly, if a projective billiard in dimension $d\geq 3$ has a caustic then this caustic is a quadric. I am very curious about the answer.

Ivrii's conjecture is also a major problem of billiard theory. We do not pretend to give an answer, but it could be interesting to study "simple" classes of $k$-reflective projective billiards. We can try fro example to answer the question if there are $k$-reflective projective billiards with an odd $k\geq 5$ inside polygons. We were unable to find examples of such billiards. More generally, we can investigate the question if the examples of $k$-reflective billiards presented in Section \ref{section:examples_reflective_billiards} are the only $k$-reflective projective billiards inside polygons. We can finally try to understand if the arguments given in \cite{glut, glutkud2} and also studied in Chapter \ref{chapter_ivrii} to classify $3$-reflective projective billiards can be generalized to a finite number of reflections. 


\chapter{Complex and projective billiards}
\label{chapter_billiards}

Billiards are usually defined as bounded domains $\Omega$ in complete Riemannian manifolds, on the boundary of which the geodesics can be reflected into new ones by the classical \textit{law of reflection} of physical optics. In the case when $\Omega$ is of dimension $2$, this law states that the angle with the boundary made by the geodesic before impact has to be the same as the angle with the boundary made by the reflected geodesic. In dimension at least $3$, the vectors directing the incident and reflected geodesics together with any normal vector to the boundary at the point of impact should also be contained in the same plane.

\begin{figure}[!h]
\centering
\begin{tikzpicture}[line cap=round,line join=round,>=triangle 45,x=1.5cm,y=1.5cm]
\clip(-2,-0.07) rectangle (2,2.8);
\draw [domain=0:2, samples=30] plot (\x,1+\x^2/4);
\draw [domain=-2:0, samples=50] plot (\x,1-\x^2/4);
\draw [domain=-2.66:2.57] plot(\x,{(--1-0*\x)/1});
\draw (0,1)-- (-0.57,3.02);
\draw (0,1)-- (1.73,2.23);
\draw [-latex] (0,1) -- (0.57,1.41);
\draw [-latex] (-0.57,3.02) -- (-0.11,1.37);
\begin{scriptsize}
\draw[color=black] (0,0.9) node {$p$};
\draw[color=black] (1,1.4) node {$\gamma$};
\end{scriptsize}
\end{tikzpicture}
\hspace{1cm}
\includegraphics[scale=0.3]{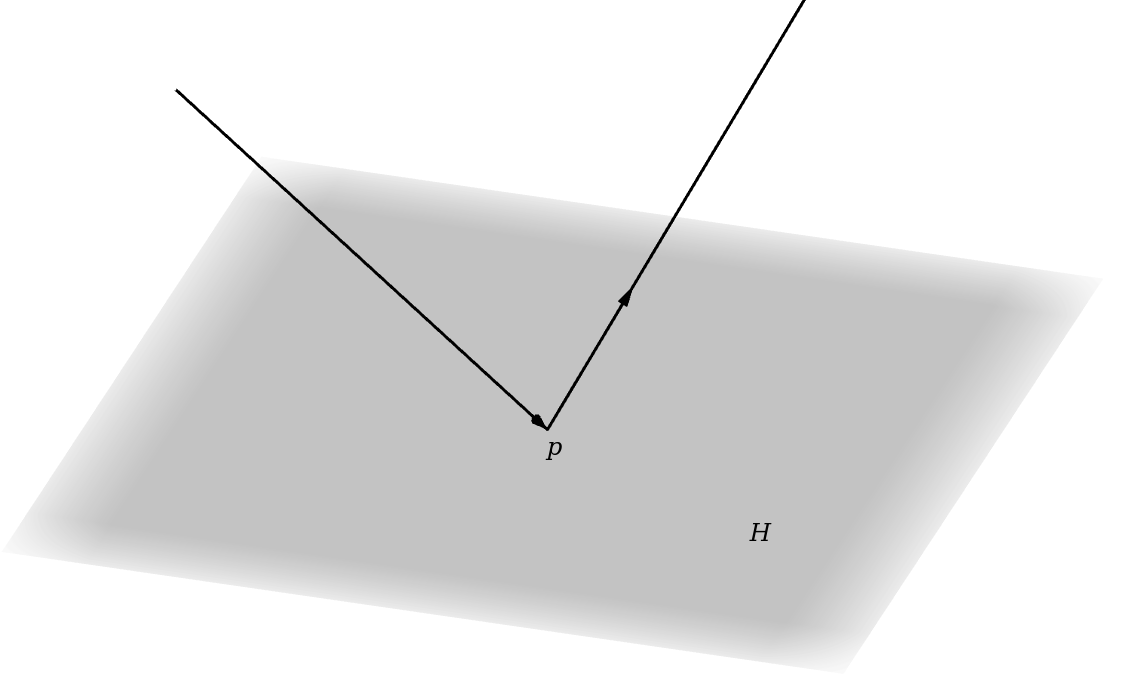}
\caption{An oriented line reflected at $p$ by a certain law of reflection on a line tangent to a curve $\gamma$ (left)/a hyperplane $H$ (right).}
\label{figure:raw_reflection2}
\end{figure}
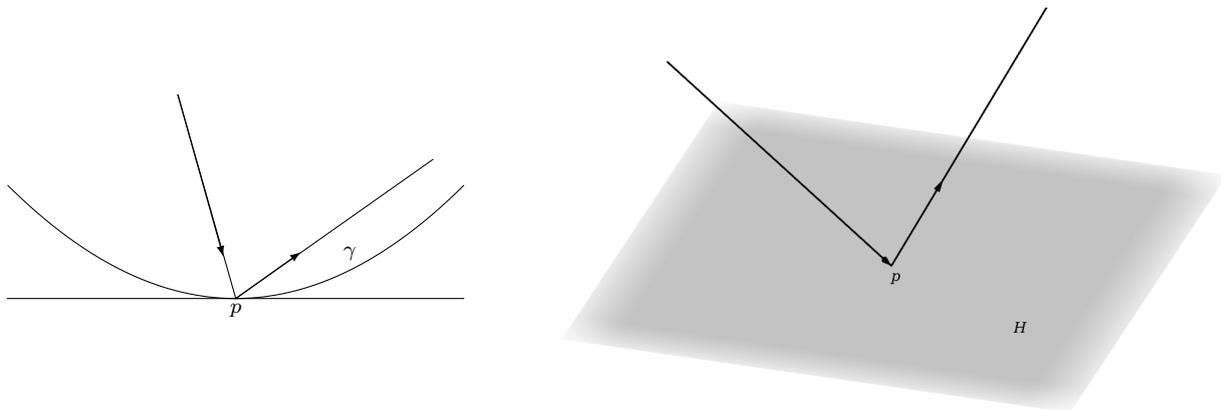

In this chapter, we define other types of reflection, or \textit{reflection laws}. Before going further into details, we would like the reader to think of them as follows: if $K$ is either the field $\RR$ or $\CC$ and $H$ is an affine hyperplane of $K^d$ (the tangent space) containing a point $p$ (the point of impact), \textit{a law of reflection at $p$ with respect to $H$ can be thought of as a non-trivial involutive map of the set of lines containing $p$ fixing the lines included in $H$.} When $K=\RR$, we can further orient the lines containing $p$ with respect to $H$, so that the image by the reflection law of an oriented line has an opposite orientation with respect to $H$ (see Figure \ref{figure:raw_reflection2}).

This chapter presents two types of billiards, the \textit{projective} and \textit{complex billiards}, defined by laws of reflections inspired from previous idea, and described in different sections. The law of reflection of projective billiards, or \textit{projective law of reflection} (see Section \ref{section:projective_billards}), is defined with help of a transverse line $L$ to $H$ at $p$. It was introduced and studied by Tabachnikov \cite{taba_projectif_ball, taba_projectif}. The law of reflection of complex billiards, or \textit{complex law of reflection} (see Section \ref{section_complex_reflection_law}), is defined in $\CC^2$ using a complexification of the Euclidean metric. It was introduced and studied by Glutsyuk \cite{glut,glut1,glut2}.

\section{Projective billiards}
	\label{section:projective_billards}
	In this section, we define the usual model of \textit{projective billiard} in $\RR^d$ as it is presented in \cite{taba_projectif_ball, taba_projectif}. This model of billiard generalizes the usual model of Euclidean billiard, but also of pseudo-Euclidean billiards and of billiards in metrics projectively equivalent to the Euclidean one (metrics in $\RR^d$ whose geodesics are contained in lines).

A \textit{projective billiard} in $\RR^d$ is a hypersurface $S$ or a collection of hypersurfaces endowed with a field of transverse lines to $S$, called \textit{field of projective lines}. For example, if $\RR^d$ is endowed with a metric or a field of non-degenerate quadratic forms, we can define a field of lines on a hypersurface $S\subset\RR^d$ as follows: for $p\in S$, define the line $L(p)$ to be the line containing $p$ and orthogonal to $T_pS$ with respect to the metric or quadratic form. It is however possible that line $L(p)$ is not transverse to $S$ at $p$ if the restriction to $T_pS$ of the field of quadratic forms is degenerate. Otherwise, $S$ has the structure of a projective billiard induced by the metric or the field of quadratic forms.

A reflection law, called \textit{projective reflection law}, can be defined on a hypersurface $S$ endowed with a field of transverse lines $L$: given an oriented line of $\RR^d$ intersecting $S$ at a certain point $p$, we define the reflected line $\ell'$ to be a line containing $p$ and satisfying a condition of harmonicity with $L(p)$ (see Definition \ref{definition:projective_reflection_law}). In the case when the projective lines $L(p)$ at $p$ is orthogonal to $T_pS$, the reflected line $\ell'$ coincides with the line reflected by the usual law of reflection (which preserves the angles of reflection in the Euclidean case).

We first recall some properties about harmonic quadruples of lines in Subsection \ref{subsection:harmonicity}, then we apply it to define projective billiards in Subsection \ref{subsection:projective_billiards}, and we finally introduce the projective billiard map in Subsection \ref{subsection:projective_billiard_map}.

\subsection{Harmonic quadruple of lines}
\label{subsection:harmonicity}

In this section, $K$ is the field $\RR$ or $\CC$. We recall some properties of the cross-ratio and harmonic quadruple of points in $\PP{1}{K}$. They can be extended to quadruple of lines containing the same point, and this will lead to the definition of projective reflection law. Most of the results on harmonic quadruples of points are very basic, and we refer the reader for example to \cite{berger_geometry} for more details.

Let $d\geq 1$ be an integer.  We denote by $\PP{d}{K}$ the $d$-dimensional projective space, which is the set of equivalence classes in $K^{d+1}\smallsetminus\{0\}$ for the relation $\sim$, defined for all $x,y\in K^{d+1}\smallsetminus\{0\}$ by $x\sim y$ if and only if there is $\lambda\in K\smallsetminus\{0\}$ such that $y=\lambda x$. For $x=(x_0,\ldots,x_d)\in K^{d+1}\smallsetminus\{0\}$, write $(x_0:\ldots:x_d)\in\PP{d}{K}$ the equivalence class of $x$ for this relation. 

\underline{Cross-ratio}. The \textit{cross-ratio} of four distinct points $p_1,p_2,p_3,p_4$ of $\PP{1}{K}$ is a well-known quantity which can be defined in many different ways. Here we adopt the definition of \cite{berger_geometry} Vol. I Chap. 6. based on the sharp $3$-reflectivity of the projective line's group of transformations:

\begin{definition}
\label{definition:cross-ratio}
The \textit{cross-ratio} of four distinct points $p_1,p_2,p_3,p_4$ of $\PP{1}{K}$ is the image $h(p_4)$ of the only projective transformation $h$ of $\PP{1}{K}$ satisfying $h(p_1)=\infty$, $h(p_2)=0$ and $h(p_3)=1$, where $\infty=(1:0)$ and $x$ stands for $(x:1)$ given any $x\in K$.
\end{definition}

The cross-ratio of four distinct points is invariant under projective transformations of $\PP{1}{K}$ (\cite{berger_geometry} Sec. 6.1.4.). We say that the quadruple $(p_1,p_2,p_3,p_4)$ is \textit{harmonic} if the cross-ratio of the corresponding points is $-1$. If we permute $p_1$ with $p_2$, or $p_3$ with $p_4$, or even $(p_1,p_2)$ with $(p_3,p_4)$, then the corresponding quadruple of points is still harmonic (\cite{berger_geometry} Prop. 6.3.1.).

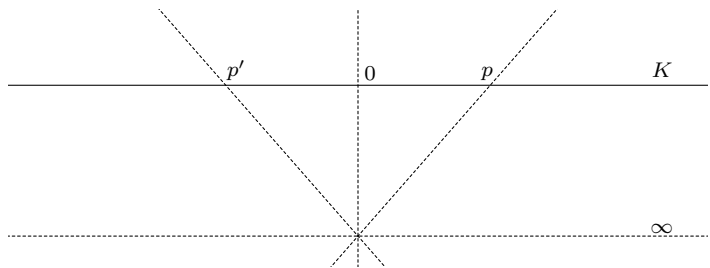
\begin{figure}[!h]
\centering
\begin{tikzpicture}[line cap=round,line join=round,>=triangle 45,x=2.0cm,y=2.0cm]
\clip(-2.3,-0.2) rectangle (2.3,1.5);
\draw [domain=-2.31:2.68] plot(\x,{(--1-0*\x)/1});
\draw [dash pattern=on 1pt off 1pt,] (0,-0.49) -- (0,2.43);
\draw [dash pattern=on 1pt off 1pt, domain=-2.31:2.68] plot(\x,{(-0-0*\x)/1});
\draw [dash pattern=on 1pt off 1pt, domain=-2.31:2.68] plot(\x,{(-0--1*\x)/0.87});
\draw [dash pattern=on 1pt off 1pt, domain=-2.31:2.68] plot(\x,{(-0--1*\x)/-0.87});
\begin{scriptsize}
\draw[color=black] (0.08,1.08) node {$0$};
\draw[color=black] (2,1.1) node {$K$};
\draw[color=black] (2,0.05) node {$\infty$};
\draw[color=black] (0.85,1.08) node {$p$};
\draw[color=black] (-0.8,1.1) node {$p'$};
\end{scriptsize}
\end{tikzpicture}
\caption{The harmonic quadruple of points $(p,p',0,\infty)$ represented 1) by points on the affine chart $K$ 2) by their equivalence classes in $\PP{1}{K}$ as dotted lines.}
\label{figure:harmonic_points}
\end{figure}

\begin{example}
\label{example:harmonic_points}
Denote by $0$ the point $(0:1)$ and by $\infty$ the point $(1:0)$. Given any point $p=(x:1)$ of $\PP{1}{K}$, the point $p'=(-x:1)$ is the only point such that the quadruple $(p,p',0,\infty)$ is harmonic (see Figure \ref{figure:harmonic_points}). Hence a quadruple of points $(p,p',0,\infty)$ is harmonic if and only if $0$ is the midpoint of $[p,p']$.
\end{example}

\underline{Harmonicity and involutive transformations.} Harmonic quadruple of points are closely related to the existence of involutive maps of the projective line $\PP{1}{K}$ (\cite{berger_geometry} Sec. 6.7.). Indeed, given two distinct points $p_3,p_4$ of $\PP{1}{K}$, there is a unique non-trivial projective involution $s$ of $\PP{1}{K}$ fixing $p_3$ and $p_4$. The map $s$ has the property that any quadruple of points of the type $(p_1,p_2,p_3,p_4)$ is harmonic if and only if $s(p_1)=p_2$.

\begin{example}
\label{example:harmonic_points_involution}
Using the same notations as in Example \ref{example:harmonic_points}, the non-trivial projective involution of $\PP{1}{K}$ fixing $0$ and $\infty$ is the map represented in the chart $\ensemble{(x:1)}{x\in K}$ as $x\mapsto-x$.
\end{example}

\underline{Space of lines}. The space of lines in $\PP{2}{K}$ is the set containing all lines of $\PP{2}{K}$. We can identify it with a $2$-dimensional projective space as follows: we see $\PP{2}{K}$ as the projectivization $\PP{}{V}$ of the space $V=K^3$. In this representation, the space of lines of $\PP{2}{K}$ can be identified with $\PP{}{\polar{V}}$, where $\polar{V}$ is the dual space of $V$: to any hyperplane $H$ of $V$ corresponds a unique set of colinear linear forms on $V$ having $H$ as a kernel. 

We can also identify it in a non-unique way with $\PP{2}{K}$ via a non-degenerate quadratic form, since the latter induces an isomorphism between $V$ and $\polar{V}$ (more details will be given in Section \ref{section:general_properties_on_quadrics}).

\underline{Space of lines containing a fixed point}. The set of lines containing a point $p\in\PP{2}{K}$ can be identified with a projective line $\PP{1}{K}$. We give two ways to state this identification, the first one being canonical, the other one being more geometric:

\textit{Identification 1.} The set of lines $\polar{p}$ containing a fixed point $p$ is a line in $\PP{}{\polar{V}}$. Indeed, if $x$ is a non-zero vector of $V$ whose equivalence class in $\PP{}{V}$ is $p$, the map $\alpha\in\polar{V}\mapsto \alpha(x)\in K$ is a non-zero linear form and its kernel is a hyperplane of $V$. Hence $\polar{p}$ is a one-dimensional projective space. 

\textit{Identification 2.} Consider a line $L$ which do not contain the point $p$. We can define a projective transformation $L\to\polar{p}$ by associating to any $q\in L$ the line $pq$. This gives a projective correspondance between the lines containing $p$ and the points on $L$.

Therefore the cross-ratio of four lines containing $p$ is well-defined in any identification of $\polar{p}$ with $\PP{1}{K}$ and doesn't depend on the identification since it is invariant by projective transformations:

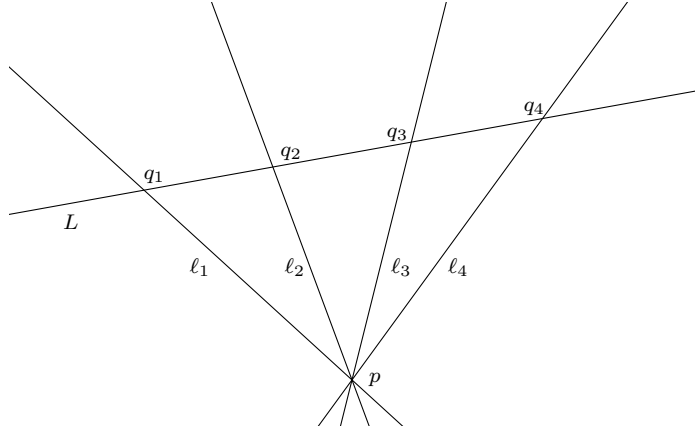
\begin{figure}[!h]
\centering
\begin{tikzpicture}[line cap=round,line join=round,>=triangle 45,x=1.0cm,y=1.0cm]
\clip(-4.5,-0.6) rectangle (4.5,5);
\draw [domain=-4.85:6.11] plot(\x,{(--22.95--1.38*\x)/7.62});
\draw [domain=-4.85:6.11] plot(\x,{(-0--2.52*\x)/-2.74});
\draw [domain=-4.85:6.11] plot(\x,{(-0--2.82*\x)/-1.04});
\draw [domain=-4.85:6.11] plot(\x,{(-0--3.15*\x)/0.78});
\draw [domain=-4.85:6.11] plot(\x,{(-0--3.47*\x)/2.51});
\begin{scriptsize}
\draw[color=black] (0.3,0) node {$p$};
\draw[color=black] (-2.6,2.7) node {$q_1$};
\draw[color=black] (-0.8,3) node {$q_2$};
\draw[color=black] (0.6,3.25) node {$q_3$};
\draw[color=black] (2.4,3.6) node {$q_4$};
\draw[color=black] (-2,1.5) node {$\ell_1$};
\draw[color=black] (-0.75,1.5) node {$\ell_2$};
\draw[color=black] (0.65,1.5) node {$\ell_3$};
\draw[color=black] (1.4,1.5) node {$\ell_4$};
\draw[color=black] (-3.7,2.1) node {$L$};
\end{scriptsize}
\end{tikzpicture}
\caption{$\ell_1,\ell_2,\ell_3,\ell_4$ form a harmonic quadruple of lines if and only if their intersection points $q_1,q_2,q_3,q_4$ with $L$ form a harmonic quadruple of points.}
\label{figure:harmonicity_lines}
\end{figure}

\begin{definition}[Harmonic quadrupe of lines in the plane]
\label{definition:harmonic_quadrupe_lines_plane}
Let $\ell_1$, $\ell_2$, $\ell_3$, $\ell_4$ be distinct lines $\ell_1$, $\ell_2$, $\ell_3$, $\ell_4$ containing a point $p\in\PP{2}{K}$. We say that \textit{the quadruple of lines $(\ell_1,\ell_2,\ell_3,\ell_4)$ is harmonic} if at least one of the following equivalent conditions is satisfied:

1) The cross-ratio of the corresponding lines is $-1$ in any identification of $\polar{p}$ with $\PP{1}{K}$.

2) The intersection points $q_1,q_2,q_3,q_4$ of $\ell_1,\ell_2,\ell_3,\ell_4$ with a line $L$ not containing $p$ form a harmonic quadruple of points (see Figure \ref{figure:harmonicity_lines}).

3) The unique non-trivial projective involution of $\polar{p}$ fixing $\ell_3$ and $\ell_4$ permutes $\ell_1$ and $\ell_2$.
\end{definition}

\begin{remark}
\label{remark:degenerate_harmonicity}
In fact condition 3) allows to extend the condition of harmonicity in the case when $\ell_1$ and $\ell_2$ are both equal to either $\ell_3$ or $\ell_4$.
\end{remark}

\begin{remark}
\label{remark:order_harmonicity}
Notice that if the quadruple of lines $(\ell_1,\ell_2,\ell_3,\ell_4)$ is harmonic, then so are the quadruples of lines obtained by permuting $\ell_1$ with $\ell_2$, or $\ell_3$ with $\ell_4$, or even $(\ell_1,\ell_2)$ with $(\ell_3,\ell_4)$. We will often use this remark.
\end{remark}

\underline{Azimuth of a line.} A computational way to work with harmonic quadruple of lines can be described by the following idea from \cite{glut}. Consider an identification of a line $L$ not containing $p$ or of $\polar{p}$ with $\PP{1}{K}=K\cup\{\infty\}$: any line $\ell$ containing $p$ can be associated with a value $z\in K\cup\{\infty\}$ called \textit{azimuth} of $\ell$, denoted by $\az(\ell)$, and defined as the corresponding coordinate of $\ell$ in $\PP{1}{K}$.

\begin{proposition}[\cite{berger_geometry} Prop. 6.7.2.]
Let $(\ell_1,\ell_2,\ell_3,\ell_4)$ be a quadruple of lines through $p$. Denote by $(z_1,z_2,z_3,z_4)$ their corresponding azimuths. The quadruple of lines is harmonic if and only if there is a non-trivial involutive projective transformation $h$ of $K\cup\{\infty\}$ fixing $z_3$, $z_4$ and permuting $z_1$ and $z_2$. The latter transformation is given for all $z\in\PP{1}{K}$ by
\begin{equation}
\label{equation:projective_transformation_harmonic}
h(z)=\frac{(z_3+z_4)z-2z_3z_4}{2z-(z_3+z_4)}.
\end{equation} 
\end{proposition}

\begin{proof}
A proof of the first statement is given in \cite{berger_geometry} Prop. 6.7.2. Formula \eqref{equation:projective_transformation_harmonic} for $h$ is not explicitely given in \cite{berger_geometry}, but the reader may check that it defines a non-trivial involutive transformation fixing $z_3$ and $z_4$.
\end{proof}

\underline{In any dimension $d\geq 2$.} We can extend statement $3)$ of Definition \ref{definition:harmonic_quadrupe_lines_plane} to lines of $\PP{d}{K}$ as follows. Let $p\in\PP{d}{K}$: the set $\polar{p}$ of lines containing the point $p$ is a projective space of dimension $d-1$ (by the same argument as for $\PP{2}{K}$). Let $H\subset\PP{d}{K}$ be a projective hyperplane and $L$ a line intersecting $H$ transversally at $p$. 

\begin{proposition}
\label{proposition:projective_law_of_reflection}
There is a unique non-trivial projective involution $s$ of $\polar{p}$, fixing $L$ and the lines included in $H$. Given any pair of lines $\ell,\ell'$ intersecting $H$ transversally at $p$, $s$ satisfies the following equivalent statements:

1) $\ell'=s(\ell)$;

2) The lines $\ell,\ell',L$ are contained in the same plane $\mathcal P$ and the quadruple of lines $(\ell,\ell',L,H\cap \mathcal P)$ is harmonic.

The involution $s$ is called the \textit{projective reflection law} with respect to $(L,H)$.
\end{proposition}

\begin{proof}
Identify $\polar{p}$ with $\PP{d-1}{K}$, so that the set of lines of $\polar{p}$ contained in $H$ is a projective hyperplane of $\PP{d-1}{K}$ and $L$ is a point of $\PP{d-1}{K}\setminus H'$: the latter are the projections in $\PP{d-1}{K}$ of a linear hyperplane $H_0\subset K^d$ and of a one-dimensional linear subspace $L_0\subset K^d$ such that $K^d=H_0\oplus L_0$. Consider the linear map acting identically on $H_0$ and on $L_0$ as $x\mapsto -x$. Then the map $s$ is obtained from it by passing to the quotients. In the same way, a linear map of $K^d$ preserving the one-dimensional subspaces of $H_0$ restricts to $H_0$ as a homothety, and the unicity of $s$ follows.

$1)\Leftrightarrow 2)$ The restriction of $s$ to any plane $\mathcal{P}$ containing $L$ is well-defined, non-trivial and involutive. Hence the equivalence between both statements is a consequence of condition $3)$ of Definition \ref{definition:harmonic_quadrupe_lines_plane}.
\end{proof}

\subsection{Line-framed hypersurfaces and projective reflection law}
\label{subsection:projective_billiards}

In this section, we introduce line-framed hypersurfaces and their reflection law, which are the formal objects used to define projective billiards. These definitions are based on the following identification: given a point $p\in\RR^d$, a line through $p$ can be seen as an element of $\mathbb{P}(T_p\RR^d)$ via the exponential map $\exp_p:T_p\RR^d\to\RR^d$. Hence we consider the following (trivial) fiber bundle $\mathbb{P}(T\RR^d)$ together with the usual projection $\pi:\mathbb{P}(T\RR^d)\to\RR^d$.

\begin{figure}[!ht]
\centering
\input{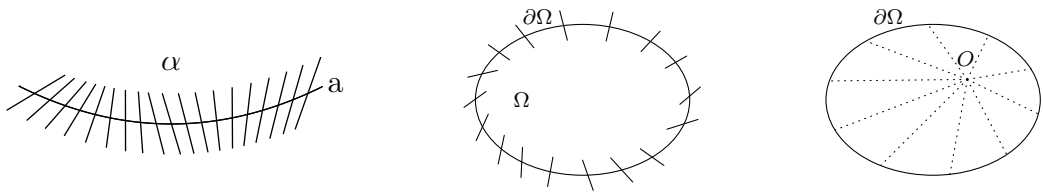}
\hspace{1cm}
\begin{tikzpicture}[line cap=round,line join=round,>=triangle 45,x=1.0cm,y=1.0cm]
\clip(-1.6,-1.2) rectangle (1.6,1.2);
\draw [rotate around={0:(0,0)}] (0,0) ellipse (1.41cm and 1cm);
\draw (-1.56,-0.1)-- (-1.27,0.11);
\draw (-1.51,0.3)-- (-1.12,0.39);
\draw (-1.24,0.73)-- (-0.96,0.53);
\draw (-0.88,0.99)-- (-0.64,0.69);
\draw (-0.3,1.18)-- (-0.2,0.79);
\draw (0.4,1.16)-- (0.31,0.78);
\draw (1.02,0.91)-- (0.77,0.64);
\draw (1.37,0.58)-- (1.09,0.42);
\draw (1.55,0.17)-- (1.28,-0.06);
\draw (1.07,-0.87)-- (0.76,-0.65);
\draw (0.66,-1.09)-- (0.37,-0.77);
\draw (0.14,-1.2)-- (0,-0.8);
\draw (-0.46,-1.15)-- (-0.39,-0.76);
\draw (-0.81,-1.03)-- (-0.79,-0.63);
\draw (-1.11,-0.85)-- (-1.03,-0.47);
\draw (-1.4,-0.54)-- (-1.24,-0.19);
\draw (1.11,-0.39)-- (1.52,-0.26);
\begin{scriptsize}
\draw[color=black] (-0.59,1.1) node {$\partial\Omega$};
\draw[color=black] (-0.8,0) node {$\Omega$};
\end{scriptsize}
\end{tikzpicture}
\hspace{1cm}
\begin{tikzpicture}[line cap=round,line join=round,>=triangle 45,x=1.0cm,y=1.0cm]
\clip(-1.6,-1.2) rectangle (1.6,1.2);
\draw [rotate around={0:(0,0)}] (0,0) ellipse (1.41cm and 1cm);
\draw [dotted](0.45,0.27)-- (-0.07,1);
\draw [dotted](0.45,0.27)-- (-0.89,0.78);
\draw [dotted](0.45,0.27)-- (-1.37,0.25);
\draw [dotted](0.45,0.27)-- (-1.29,-0.4);
\draw [dotted](0.45,0.27)-- (-0.69,-0.87);
\draw [dotted](0.45,0.27)-- (0.22,-0.99);
\draw [dotted](0.45,0.27)-- (0.97,-0.73);
\draw [dotted](0.45,0.27)-- (1.39,-0.19);
\draw [dotted](0.45,0.27)-- (1.29,0.41);
\draw [dotted](0.45,0.27)-- (0.73,0.86);
\begin{scriptsize}
\draw[color=black] (-0.59,1.1) node {$\partial\Omega$};
\draw[color=black] (0.43,0.55) node {$O$};
\end{scriptsize}
\end{tikzpicture}
\caption{\underline{left:} A line-framed curve $\alpha$ over the curve $a$. \underline{center:} A bounded domain $\Omega$ whose boundary $\partial\Omega$ is a line-framed curve. \underline{right:} The same domain $\Omega$ endowed with a so-called \textit{centrally-projective} \cite{taba_projectif} field of transverse lines (dotted lines).}
\label{fig:projective_billiard}
\end{figure}

\begin{definition} 
\label{definition:line_framed_hyp}
A {\it line-framed hypersurface} (see Figure \ref{fig:projective_billiard}) is a regularly embedded connected $(d-1)$-dimensional surface $\Sigma\subset\mathbb P(T\RR^d)$ with the following properties:\\
\enum The projection $\pi$ sends $\Sigma$ diffeomorphically to a regularly embedded hypersurface $S\subset\RR^2$, 
which will be identified with $\Sigma$ and called the {\it classical boundary} of the hypersurface $\Sigma$.\\
\enum For every $(p,L)\in\Sigma$ the line $L$ is transverse to $T_pS$. \\
We will often say that $\Sigma$ is a line framed-hypersurface \textit{over} $S=\pi(\Sigma)$, and that $L$ is the \textit{field of projective lines} of $\Sigma$. In particular, $L(p)$ is the line such that $(p,L(p))\in\Sigma$.
\end{definition}

\begin{remark}
An analogous definition can be given without supposing that $L$ is transverse to $T_pS$. In this case we say that such line-framed hypersurface has \textit{projective singularities}.
\end{remark}

\begin{remark}
\label{remark:line_framed_hyp_proj}
Line-framed hypersurface can also be defined on $\PP{}{T\PP{d}{\RR}}$ with analogue statements as in Definition \ref{definition:line_framed_hyp}.
\end{remark}

Let $\Sigma$ be a line-framed hypersurface over an hypersurface $S\subset\RR^d$. The projective reflection law on $\Sigma$ can be defined as follows:
\begin{definition}
\label{definition:projective_reflection_law}
Let $p\in S$ and $\ell,\ell'$ be oriented lines intersecting $S$ at $p$. We say that $\ell'$ is \textit{obtained from $\ell$ by the projective reflection law on $\Sigma$ at $p$} if\\
\enum the lines $\ell$, $\ell'$, $L(p)$ are contained in a plane $\mathcal{P}$\\
\enum the quadruple of lines $\ell,\ell',L(p),T_pS\cap \mathcal{P}$ is harmonic in $\mathcal{P}$;\\
\enum the orientations of $\ell$ and $\ell'$ with respect to $T_pS$ are opposite.
\end{definition}

Using previous statements, projective billiards can be defined as follows:

\begin{definition}
A \textit{projective billiard} is a domain $\Omega$ whose boundary $S=\partial\Omega$ is the classical boundary of a $\class^ 1$-smooth line-framed hypersurface $\Sigma$ together with the corresponding projective reflection law on $\Sigma$.  See Figure \ref{fig:projective_billiard}.
\end{definition}

\subsection{Projective orbits and projective billiard map}
\label{subsection:projective_billiard_map}

One can study the orbits of the projective reflection law inside bounded domain $\Omega$ whose boundary $S=\partial\Omega$ is the classical boundary of a $\class^ 1$-smooth line-framed hypersurface $\Sigma$.

\begin{definition}
A \textit{projective orbit}, or simply an \textit{orbit}, of the projective billiard $\Omega$ is a sequence of points $p_1,\ldots,p_k\in\partial\Omega$ such that for each $j=1,\ldots, k-1$\\
\enum $p_j\neq p_{j+1}$, the line $p_{j}p_{j+1}$ is oriented from $p_j$ to $p_{j+1}$;\\
\enum the interior of each segment $p_{j}p_{j+1}$ is included in $\Omega$ ;\\
\enum for $j>1$, the lines $p_{j-1}p_j$ and $p_{j}p_{j+1}$ are transverse to $S$ at $p_j$;\\
\enum for $j>1$, the line $p_{j}p_{j+1}$ is obtained from $p_{j-1}p_j$ by the projective reflection law at $p_j$.\\
The orbit is said to be \textit{$k$-periodic} if $(p_1,\ldots,p_k,p_1,p_2)$ is an orbit.
\end{definition}

Let $(p_1,p_2,p_3)$ be a projective orbit of $\Sigma$ such that the line $p_2p_3$ is transverse to $S$ at $p_3$. There is an open subset $U_{(p_1,p_2)}$ of $S\times S$ containing $(p_1,p_2)$ such that for all $(q_1,q_2)\in U_{(p_1,p_2)}$, $q_1\neq q_2$, the line $q_1q_2$ is transverse to $S$ at $q_2$ and is reflected into a line intersecting $S$ transversaly at a point $q_3$ by the projective law of reflection at $q_2$. We can define on $U_{(p_1,p_2)}$ the projective billiard map using above description as the map $\mathcal{B}:U_{(p_1,p_2)}\to S\times S$ satisfying
\begin{equation}
\label{equation:projective_billiard_map}
\mathcal{B}(q_1,q_2)=(q_2,q_3).
\end{equation}

\begin{proposition}
\label{proposition:billiard_regularity_rank}
Let $r\geq 2$ be an integer. If $\Sigma$ is $\class^r$-smooth (respectively analytic) then $\mathcal{B}$ is a $\class^{r-1}$-smooth (respectively an analytic) map of rank $2(d-1)$.
\end{proposition}

\begin{proof}
We first show that $\mathcal{B}$ is of class $\class^{r-1}$ (resp. analytic). Indeed, notice that there is a $\class^{r-1}$-smooth (resp. an analytic map) defined on the restriction of the set $\restreint{\PP{}{T\RR^d}}{S}$ which associate to $(p,\ell)$, where $p\in S$ and $\ell$ is a line containing $p$, the element $(p,\ell')$ where $\ell'$ is the line containing $p$ and obtained by the projective reflection law at $p$ defined by $\Sigma$. In fact the restriction of such map on each fiber $\{p\}\times \PP{}{T\RR^d}$ is a projective transformation depending $\class^{r-1}$-smoothly (resp. analyticaly) on $p$. 

We conclude on the regularity by proving the following result: consider a line $\ell$ intersecting a $\class^{r}$-smooth (resp. an analytic) hypersurface $S$ transversaly at a point $p$; then if another line $\ell'$ is close to $\ell$, then $\ell'$ intersects $S$ at a point $q$ close to $p$ and that the map $\ell'\mapsto q$ is of class $\class^r$ (resp. is analytic). Indeed, consider a affine hyperspace $H$ intersecting $\ell$ transversally at a point $p_1$ and $v\neq 0$ be a unit vector directing $\ell$. There is a diffeomorphism between a neighborhood of $U_{\ell}$ of lines containing $\ell$ and a neighborhood $U_{(p_1,v)}$ of $(p_1,v)$ in $H\times\ell$. Now consider an open subset $U_p$ of $p$ and a $\class^{r}$-smooth (resp. an analytic) submersion $f:U_p\to\RR$ such that $S\cap U_p=\inverse{f}(\{0\})$. The map $F:U_{(p_1,v)}\times\RR\to\RR$ defined by $F(q_1,v',t)=f(q_1+tv')$ is well-defined in a neighborhood of $(p_1,v,\tau)$ where $p=p_1+\tau v$, and is $\class^r$-smooth (resp. analytic). Its differential in $t$ at $(p_1,v,\tau)$ is $df(p)\cdot v$ and the latter is non-zero since $v$ is not in the tangent space to $S$ at $p$. The conclusion follows from the implicit function theorem.

Finally, the map $\mathcal{B}:U_{(p_1,p_2)}\to S\times S$ is a local diffeomorphism onto its image, since if $\mathcal{B}(q_1,q_2)=(q_2,q_3)$ then $\mathcal{B}(q_3,q_2)=(q_2,q_1)$ and we can easily contruct a smooth inverse map for $\mathcal{B}$.
\end{proof}

\subsection{Projective billiards induced by a metric}
\label{subsection:projective_and_others}

As explained in the introductive section, other types of billiards such as the usual billiards, billiards in metrics projectively equivalent to the Euclidean metric or billiards in pseudo-Euclidean spaces can be defined as specific projective billiards. In this section, we recall briefly these different types of billiards and give an explanation on why they can be seen as projective billiards. We first define different metrics on $\RR^d$:

\underline{Euclidean metric.} It is the canonical Riemannian metric on $\RR^d$: $\sum_{j=1}^d \dd x_j^2$ on $\RR^d$.

\begin{figure}
\centering
\includegraphics[scale=0.4]{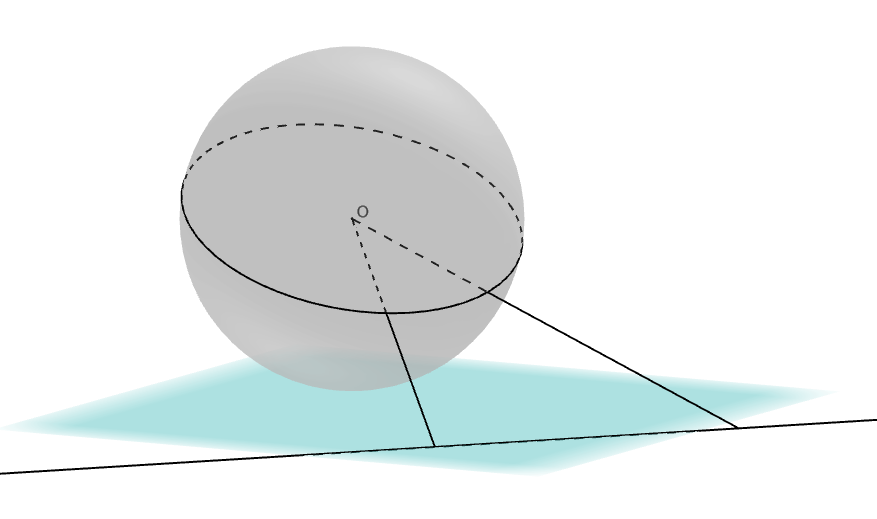}
\caption{A great circle of the sphere projected onto an affine horizontal plane.}
\label{figure:sphere_projection_geodesics}
\end{figure}

\underline{Metrics projectively equivalent to the Euclidean metric.} They are Riemannian metrics in $\RR^d$ whose geodesics are lines. A theorem of Beltrami \cite{beltrami, matveev} improved in all dimensions by Lipschitz and Schur implies that such metrics have constant sectional curvature. We describe two famous examples of such metrics (which can also be found in \cite{taba_projectif}):

\textit{Sphere.} Consider the upper half open hemisphere $H_N$ of the unit sphere $\SS^2$ of center $O$, given by the equations $x^2+y^2+z^2=1$ and $z>0$ in $\RR^3$. Any point $p\in H_N$ can be mapped to a unique point $q$ of the plane $P\subset\RR^3$ given by equation $z=-1$: $q$ is defined to be the intersection point of the line $Op$ with $P$. This defines a diffeomorphism $\varphi:H_N\to P$. The geodesics of $\SS^2$ for the usual spherical metric $g_{\SS^2}$ are contained in great circles, which are the intersection of $\SS^2$ with a plane containing $O$. Therefore, their image by $\varphi$ are lines of $P$, see Figure \ref{figure:sphere_projection_geodesics}. Hence the geodesics of $P$ for the pushforwarded metric $\varphi_{\ast}g_{\SS^2}$ are lines.

\textit{Hyperboloïd.} Consider the upper sheet ($x>0$) of the hyperboloïd of equation $x^2-y^2-z^2=1$, denoted by $\HH^2$. The usual Minkovski metric $g_{\HH^2}$ on $\HH^2$ is the restriction of $\dd y^2+\dd z^2-\dd x^2$ to the tangent planes of $\HH^2$. The geodesics of $\HH^2$ are the intersections of $\HH^2$ with a plane containing the origin $O$. Hence the same construction can be applied to push forward the metric $g_{\HH^2}$ on a plane where it is projectively equivalent to the Euclidean metric.

\underline{Pseudo-Euclidean metrics.}\cite{DragRad_minkowski2, DragRad_minkowski, khesin_taba} A pseudo-Euclidean space of signature $(k,\ell)$, with $k+\ell=d$, is the space $\RR^d$ endowed with the non-degenerate bilinear form $\scalar{\cdot}{\cdot}$ defined for all $x,y\in\RR^d$ by 
\begin{equation}
\label{equation:pseudo_euclidean_metric}
\scalar{x}{y}=\sum_{j=1}^k x_jy_j-\sum_{j=k+1}^d x_jy_j.
\end{equation}

Now consider the following situation which gives the definition on how a line intersecting a hyperplane is reflected in these metrics. Let $H\subset\RR^d$ be a hyperplane, $q$ a non-degenerate quadratic form on $\RR^d$ (for example one of the previous defined metrics) such that the $q$-orthogonal space to $H$, $\orth{H}$, which is one-dimensional, is not included in $H$ (for example when $q$ is positive-definite). Any vector $v\in\RR^d$ has a unique decomposition $v=h+n$ where $h\in H$ and $n\in\orth{H}$, and can be associated to the vector $s(v)=h-n$. The map $s$ is linear and induces a non-trivial involution on the set of lines containing the origin $O$ which fixes the line $\orth{H}$ and any line included in $H$. 

\begin{proposition}
\label{proposition:usual_is_projective}
The map $s$, called the \textit{usual law of reflection} in the metric $q$, coincides with the projective law of reflection with respect to $(L,H)$ (Proposition \ref{proposition:projective_law_of_reflection}).
\end{proposition}

\begin{proof}
Both maps satisfy the same properties, hence coincide by Proposition \ref{proposition:projective_law_of_reflection}.
\end{proof}

Therefore, if $S\subset\RR^d$ is a smooth hypersurface, and $g$ is one of previous metrics (a pseudo-Euclidean metric, the Euclidean metric or a projectively equivalent one), then we can define at each point $p\in S$ the $g$-orthogonal line $L(p)$ to $T_pS$ containing $p$. If at each point $p$, the line $L(p)$ is transverse to $T_pS$, then $g$ induces a line-framed hypersurface over $S$ denoted by
$$\linep{S}{g} = \ensemble{(p,L)\in\mathbb{P}(T\RR^d)}{L=\left(T_pS\right)^{\bot_g}}.$$
By Proposition \ref{proposition:usual_is_projective}, any orbit in $S$ for the usual reflection law in the metric $g$ is an orbit of the corresponding projective billiard.

\section{Complex billiards}
	\label{section_complex_reflection_law}
	
In this section, we present a natural generalization of the usual reflection law in the Euclidean plane to $\CC^2$ and also $\CP^2$: \textit{the complex reflection law}. It was introduced, together with complex planar billiards, by Glutsyuk in \cite{glut1} and \cite{glut2}. See also \cite{glut} where they were applied to solve the two-dimensional Tabachnikov's Commuting Billiard conjecture and a particular case of two-dimensional Plakhov's Invisibility conjecture with four reflections. 
	
\subsection{Complex reflection law}

We denote by $\CP^2$ the set $\PP{2}{\CC}$ defined at Section \ref{subsection:harmonicity}. Any element of $\CP^2$ can be written as a triple $(x:y:z)$, with $(x,y,z)\in \CC^3\smallsetminus\{0\}$. By construction $(tx:ty:tz)=(x:y:z)$ for any complex number $t\neq 0$. In this set of coordinates, the complex projective space is the disjoint union
$$\CP^2=U_z\cup L_{\infty}$$
of the so-called \textit{standard} open subset $U_z=\ensemble{(x:y:1)}{(x,y)\in\CC^2}$ and the line $L_{\infty}=\ensemble{(x:y:0)}{(x,y)\in\CC^2\smallsetminus\{0\}}$ called \textit{line at infinity}. The map $(x,y)\in\CC^2\mapsto(x:y:1)\in U_z$ is an analytic chart mapping $\CC^2$ to $U_z$. Hence we can consider the pushforward of the non-degenerate quadratic form $q=dx^2+dy^2$ defined on $\restreint{T\CP^2}{U_z}:= TU_z$.

\begin{definition}
\label{definition:isotropic_line}
A line of $\CP^2$ is said to be \textit{isotropic} if it contains either the point $I=(1:i:0)$ or the point $J=(1:-i:0)$, and \textit{non-isotropic} if not. Notice that the line at infinity is isotropic.
\end{definition}

In the case of a non-isotropic line $L\subset\CP^2$, we can define a complex $q$-isometric involution of the space $U_z\simeq \CC^2$ fixing the points of $L$. This involution can be constructed by considering the projective transformations preserving $L$ and its $q$-orthogonal lines. This involution induces a symmetry on lines of $U_z$, and can be extended to all lines of $\CP^2$ by sending $L_{\infty}$ to $L_{\infty}$. In the case of an isotropic line, this contruction fails since $q$-orthogonal lines to $L$ are its parallel lines.

\begin{definition}[\cite{glut1}, definition 2.1]
The \textit{symmetry} with respect to a line $L\neq L_{\infty}$ is defined as follows:\\
\enum Case 1: $L$ is non-isotropic. The \textit{symmetry acting on $\CC^2$} is the unique non-trivial complex $q$-isometric involution fixing the points of the line $L$. It induces the same symmetry acting on lines.\\
\enum Case 2: $L$ is isotropic. We define the \textit{symmetry of lines through a point $p\in L\cap U_z$}: two lines $\ell$ and $\ell'$ which contain $p$ are called symmetric if there are sequences $(L_n)_n$, $(\ell_n)_n$, $(\ell'_n)_n$ of lines through points $p_n$ so that $L_n$ is non-isotropic, $\ell_n$ and $\ell'_n$ are symmetric with respect to $L_n$, $\ell_n \to \ell$, $\ell'_n \to \ell'$, $L_n \to L$ and $p_n \to p$.
\end{definition}

We recall now lemma $2.3$ \cite{glut1} which gives an idea of this notion of symmetry in the case of an isotropic line through a finite point.

\begin{lemma}[\cite{glut1}, lemma 2.3]
\label{lemma_glutsyuk}
If $L$ is an isotropic line through a point $p\in U_z$ and $\ell$, $\ell'$ are two lines which contain $x$, then $\ell$ and $\ell'$ are symmetric with respect to $L$ if and only if either $\ell=L$, or $\ell'=L$. 
\end{lemma}

\subsection{Complex orbits}

Let $\gamma\subset\CP^2$ be a complex curve of $\CP^2$ (that is smooth at each point).

\begin{definition}[\cite{glut1}]
\label{definition:complex_ndegenerate_orbit}
A \textit{non-degenerate orbit} on $\gamma$ is a finite sequence $(p_1,\ldots,p_k)\in\gamma^k$ such~that\\
\enum $p_j\neq p_{j+1}$ for each $j\in\{1,\ldots,k-1\}$;\\
\enum $T_{p_j}\gamma$ is not isotropic for each $j\in\{1,\ldots,k\}$;\\
\enum the lines $p_{j-1}p_j$ and $p_jp_{j+1}$ are symmetric with respect to the tangent line $T_{p_j}\gamma$ for each $j\in\{2,\ldots,k-1\}$.\\
The \textit{side} of an orbit is one of the lines $p_jp_{j+1}$.
A \textit{non-degenerate $k$-periodic orbit} is a non-degenerate orbit $(p_1,\ldots,p_{k})\in\gamma^k$ such that $(p_1,\ldots,p_{k},p_1,p_2)$ is a non-degenerate orbit.
\end{definition}

By opposition we can define degenerate orbits as follows:

\begin{definition}[\cite{glut1}]
A \textit{degenerate orbit} (resp. a \textit{degenerate $k$-periodic orbit}) on $\gamma$ is a set of points $(p_1,\ldots,p_k)\in\gamma^k$ which is the limit of non-degenerate orbits (resp. non-degenerate $k$-periodic orbits) and is not a non-degenerate orbit (resp. non-degenerate $k$-periodic orbits).
\end{definition}

We can also define the \textit{side} of a degenerate orbit as the limit of the sides of non-degenerate orbits converging to it. In the case when $p_j=p_{j+1}$, we can naturally define the side $p_jp_{j+1}$ as the tangent line $T_{p_j}\gamma$.

\section{Proof by complexification: circumcenters of triangular orbits}
	\label{section_circum_centers}
	
	\newcommand{\ellipse}{\mathcal{E}}

In this section, we present a published result \cite{fierobe_circumcenters} on the circumcenters of triangular orbits in an elliptic billiard, which is of great interest for us since its proof uses complex billiards to solve a problem of real geometry. More precisely, we are interested in the usual billiard inside an ellipse and its $3$-periodic or \textit{triangular} orbits. We show that  the set of all circumcenters to these orbits is an ellipse. The proof of this result is based on the complexification of the problem and on the use of the complex reflection law introduced at Section \ref{section_complex_reflection_law}.

\begin{theorem}
	\label{theorem:circumcenters}
	The set $\mathcal{C}$ of the circumcenters of all triangular orbits of the billiard within an ellipse is also an ellipse.
\end{theorem}

\begin{remark}
\label{rmk_circle}
Theorem \ref{theorem:circumcenters} is obvious in the particular case where the ellipse is a circle, because then the set of circumcenters is reduced to a single point. Thus, \textit{we will assume that the ellipse is not a circle}.
\end{remark}

There are many other results similar to theorem \ref{theorem:circumcenters}. Dan Reznik discovered experimentally the same result for the incenters of triangular orbits, see the video \cite{reznik_youtube} and the github page \cite{reznik_github} written with Jair Koiller. Romaskevich (see \cite{romaskevich2014}) confirmed these observations by proving them and her proof widely inspired ours. Tabachnikov and Schwartz, in \cite{taba_centers}, proved that the loci of the centers of mass (and of an other particular point) of a $1$-parameter family of Poncelet $n$-gons in an ellipse is an ellipse homothetic to the previous one. They also mention that a similar result was proved by Zaslawski, Kosov and Muzafarov for the orthocenters (\cite{zaslawski}, reference from \cite{taba_centers}). And Garcia (see \cite{ronaldo}) uses explicit calculations to prove that the loci of circumcenters, incenters and orthocenters of triangular orbits are ellipses, and describes them precisely. His proof of the result about circumcenters was found simultaneously and independantly to us.

Before going into details, we give here a brief summary of the proof, which is inspired by \cite{romaskevich2014}, and in which we use the same complex methods. We consider a projective complexified version of $\mathcal{C}$, denoted by $\hat{\mathcal{C}}$, which turns out to be an algebraic curve as a consequence of Remmert proper mapping theorem and Chow's theorem, see \cite{GH} p. 34. Then we show that the intersection of the complex curve $\hat{\mathcal{C}}$ with the foci line of the boundary ellipse $\ellipse$ is reduced to two points, each one of them corresponding to a single triangular orbit. Further algebraic arguments on the intersection type of $\hat{\mathcal{C}}$ with the foci line of $\ellipse$ allow to conclude that it is a conic, using Bezout theorem. It's then possible to check that $\hat{\mathcal{C}}$ is an ellipse since its real part is bounded.

As explained, one considers the projective complex Zariski closure of the ellipse $\ellipse$ and a complexified version $\hat{\mathcal{C}}$ of $\mathcal{C}$. In order to define $\hat{\mathcal{C}}$ and to prove the first statement concerning the intersection with the foci line, we study an extension of the reflection law and of the triangular orbits to complex domain, as in \cite{romaskevich2014}, and we use some of the results contained in the latter article such as Proposition \ref{lemma_olga}. 

\textbf{Section \ref{sec_general_reflection}} is devoted to the complex reflection law and to complex orbits in a complexified ellipse: \textbf{Subsection \ref{sec_cmplxconic}} recalls some results about complexified conics; we further define what is a triangular complex orbit in \textbf{Subsection \ref{sec_orbits}}; then, in \textbf{Section \ref{sec_circumcircles}} we introduce the definition and we study properties of complex circumscribed circles to such orbits: Proposition \ref{prop_orbit_form} is the main result of this section. Finally, \textbf{Section \ref{sec_proof}} is devoted to the proof of Theorem \ref{theorem:circumcenters}, using previous results.

\subsection{Complex triangular obits on an ellipse}
\label{sec_general_reflection}

In this section, we recall some results about complexified conics and we study results about triangular orbits of the complexified ellipse $\ellipse$.

\subsubsection{Preliminary results on complexified conics}
\label{sec_cmplxconic}

We define a \textit{complexified conic} as the algebraic closure of a real conic in $\RR^2$: an ellipse, a hyperbola or a parabola. We recall that an ellipse cuts the line at infinity in two distinct points with strictly complex coordinates, a hyperbola in two distinct points with real coordinates, and a parabola is tangent to the line at infinity. The following results on conics are well-known and can be found in \cite{berger_geometry,klein26}.

\begin{proposition}[\cite{berger_geometry} subsection 17.4.2.1]
\label{circle_cycpoints}
A conic is a circle if and only if some of the points $I$ or $J$ belong to it. Furthermore, if a conic is a circle, then both $I$ and $J$ belong to it.
\end{proposition}

In fact, a circle has two isotropic tangent lines intersecting at its center (see the following propositions).

\begin{proposition}[\cite{berger_geometry} subsection 17.4.3.1]
\label{tangent_cycpoints}
A focus $f$ of a conic lies in the intersection of two isotropic tangent lines to the conic.
\end{proposition}

\begin{proposition}[\cite{klein26}, p.~179]
\label{confocal_conics}
Two complexified confocal ellipses have the same tangent isotropic lines, which are four isotropic lines taken with multiplicities: one pair intersecting at a focus, and the other one - at the other focus.
\end{proposition}

This brings us to the following redefinition of the foci:

\begin{definition}[\cite{berger_geometry} subsection 17.4.3.2]
The \textit{complex foci} of an ellipse are the intersection points of its isotropic tangent lines.
\end{definition}

\begin{remark}
The complex projective closure of a real ellipse has four complex foci, including two real ones.
\end{remark}

\begin{corollary}
\label{istropy_foci}
A conic has at most four dinstinct finite isotropic tangent lines, each two of them intersecting either at a focus, or at an isotropic point at infinity.
\end{corollary}




\subsubsection{Triangular orbits}
\label{sec_orbits}

Let $\ellipse\subset\CP^2$ be a complexified ellipse which is not a circle.

\begin{definition}
A \textit{non-degenerate triangular orbit} is a non-degenerate  $3$-periodic orbit (see Definition \ref{definition:complex_ndegenerate_orbit}).
\end{definition}

\begin{remark}
The vertices of a non-degenerate orbit are not collinear since a line intersects the ellipse in at most two points.
\end{remark}

\begin{remark}
\label{rem_non_isotropic_sides}
As explained in \cite{glut2}, the reflection with respect to a non-isotropic line permutes the isotropic directions $I$ and $J$. This argument implies that a non-degenerate triangular orbit has no isotropic side.
\end{remark}

\begin{proposition}[\cite{romaskevich2014}, lemma 3.4]
\label{lemma_olga}
A degenerate triangular orbit of $\ellipse$ has an isotropic side $A$ which is tangent to $\ellipse$, and two coinciding non-isotropic sides $B$.
\end{proposition}

During the proof, it will be convenient to distinguish two types of orbits : the ones with no points at infinity, and the others, with at least one point at infinity:

\begin{definition}
An \textit{infinite triangular orbit} on  $\ellipse$ is an orbit which has at least one vertex on the line at infinity. The orbits with only finite vertices are called \textit{finite orbits}.
\end{definition}

\begin{proposition}
\label{lemma_infinite_orbit}
An infinite triangular orbit is non-degenerate, and has exactly one vertex at infinity.
\end{proposition}

\begin{proof}
First note that the results recalled in Subsection \ref{sec_cmplxconic} imply that a tangent line of the ellipse $\ellipse$ at a point on $L_{\infty}$ cannot be isotropic.

Suppose two vertices, $\alpha,\beta$, of the orbit are at infinity. Then, $\alpha \beta$ is the line at infinity. But the tangent $T_{\beta}$ to the ellipse $\ellipse$ in $\beta$ is not isotropic, and the line at infinity reflects to itself through the reflection by $T_{\beta}$. Hence, the orbit is $\{\alpha,\beta\} = L_{\infty} \cap \ellipse$, which should be a degenerate orbit. But it cannot be a degenerate orbit by Proposition \ref{lemma_olga} since the tangent lines to its vertices $\alpha, \beta$ are not isotropic. Thus, only one vertex lies at infinity.

Therefore, if it is a degenerate orbit, it has two vertices, $\alpha,\beta$, corresponding by Proposition \ref{lemma_olga} to two sides, $A$ which is isotropic and tangent to the ellipse in $\alpha$, and $B$ which is a line containing $\alpha$ and $\beta$. Since the tangency points of isotropic tangent lines are finite, $\alpha$ is finite. Thus $\beta$ is infinite (because the orbit is supposed infinite). Then $B$ and the tangent line $T_{\beta}\ellipse$ to the ellipse in $\beta$ are collinear (since they have the same intersection point at infinity). But both are stable by the complex reflection by $T_{\beta}$, hence $T_{\beta}\ellipse = B$ which is impossible since $B$ is not tangent to the ellipse.
\end{proof}

\subsection{Circumcircles and circumcenters of complex orbits}
\label{sec_circumcircles}

Here we present the last part of the required definitions, which concerns the complex circles circumscribed to triangular orbits. This part is different from the previous one, because here the considered conics are complex and not necessarily complexified versions of real conics.

\begin{definition}
A \textit{complex circle} is a regular complex conic passing through both isotropic points at infinity, $I$ and $J$. Its \textit{center} is the intersection point of its tangent lines at $I$ and $J$.
\end{definition}

\begin{proposition}
\label{prop_circle_reg_orbit}
For a non-degenerate finite orbit, there is a unique complex circle passing through the vertices of the orbit and both isotropic points at infinity. It is called the circumscribed circle or circumcircle to the non-degenerate orbit.
\end{proposition}

\begin{proof}
Denote by $\alpha, \beta, \gamma$ the vertices of the orbit. We have to prove that no three points of $\alpha, \beta, \gamma, I, J$ are collinear. Indeed, as no vertices are on the line at infinity, we only need to study two different cases:
1) $\alpha, \beta, \gamma$ are not collinear because they are distinct and they lie on the ellipse which has at most two intersection points with any line. 2) $\alpha, \beta, I$ are not collinear or else the line $\alpha \beta$ would be isotropic. But this is impossible for a non-degenerate triangular orbit by Remark \ref{rem_non_isotropic_sides}. We then exclude all other possible combinations of two vertices of the orbit with $I$ or $J$, using the same arguments.
\end{proof}

Let us extend this definition to degenerate orbits.

\begin{definition}
\label{def_circle_deg_orbit}
Let $T$ be a degenerate or infinite orbit. A \textit{circumscribed circle} of $T$ is the limit (in the space of conics) of a converging sequence of circumscribed circles of non-degenerate finite orbits converging to $T$. If a sequence of complex circles converges to a conic so that their centers converge to a point $c\in\CP^2$, then $c$ is called a center of the limit conic. A \textit{circumcenter} of $T$ is a center of its circumscribed circle.
\end{definition}

\begin{remark}
\textit{A priori}, a limit conic $\mathcal{K}$ may have several centers in the sense of this definition. Indeed, $c$ depends on the choice of the sequence of circles converging to $\mathcal{K}$. See Case $4$ of Proposition \ref{prop_types_circle} and its proof for more details.
\end{remark}

Even if they are called \textit{circles}, the circumscribed circles to a degenerate or infinite orbit can degenerate into pairs of lines, as described below.

\begin{proposition}
\label{prop_types_circle}
The limit of a converging sequence of complex circles is one of the following:
\begin{enumerate}
	\item a regular circle ;
	\item a pair of isotropic  non-parallel finite lines ; the corresponding center lies on their intersection ; 
	\item the infinite line and a finite line $d$ ; the center $c$ lies on the line at infinity and represents a direction which is orthogonal to $d$ ; 
	\item the line at infinity taken twice : its center can be an arbitrary point in $\CP^2$.
\end{enumerate}
\end{proposition}

\begin{proof}
The equation of a regular circle $\mathcal{D}$ is of the form
$$a(x^2+y^2)+pxz+qyz+rz^2 = 0$$
where $a, p,q,r\in\CC$, $a\neq 0$ and $4a r \neq p^2+ q^2$. Both isotropic tangent lines to $\mathcal{D}$ have equations $2a(x \pm iy) +(p\pm i q)z = 0$, whose intersection is $c = (p:q:-2a)$, which is the center of $\mathcal{D}$ by definition.

If we take a limit of regular circles, the equation of the limit circle is of the same type, that is
$$a(x^2+y^2)+pxz+qyz+rz^2 = 0$$
but maybe with $a=0$ or $4a r = p^2+ q^2$. And the center $c$ is still of coordinates $(q:p:-2a)$. 

If $a=0$, the limit circle is the union of the line at infinity ($z=0$) and the line $d$ of equation $px+qy+rz = 0$. The line $d$ is finite if and only if $(p,q)\neq 0$, and in this case it has direction $(q,-p)$. Since $c = (p:q:0)$, the direction represented by $c$ is orthogonal to $d$. If $d$ is infinite, the limit circle is the (double) line at infinity. Note that in this case the center can be an arbitrary point.

If $a\neq 0$, but $4a r = p^2+ q^2$, the equation of the limit circle becomes
$$\left(x+\frac{p}{2a}z\right)^2+\left(y+\frac{q}{2a}z\right)^2=0$$
which is the equation of two isotropic non collinear lines intersecting at the point 
$(-\frac{p}{2a}:-\frac{q}{2a}:1) = (p:q:-2a) = c$. If $a\neq 0$ and $4a r \neq p^2+ q^2$, the limit circle is regular.
\end{proof}

Now let us find which triangular orbits have their center on the line of real foci of $\ellipse$. 

\begin{proposition}
\label{prop_orbit_form}
Suppose that $T$ is a complex triangular orbit whose circumcenter lies on the real foci line. Then $T$ is finite, non-degenerate, symmetric with respect to the real foci line of $\ellipse$, and has a vertex on it.
\end{proposition}

\begin{proof}
Let $T$ be a triangular orbit with a circumscribed circle $C$ having a center $c$ on the real foci line of $\ellipse$.

\textbf{First case : Suppose $T$ is finite and non-degenerate.} We follow the arguments of Romaskevich \cite{romaskevich2014} who treated the similar case for incenters. Indeed, at least two vertices should lie outside the foci line. If the line through them is not orthogonal to the foci line, then this pair of vertices together with their symmetric points and the remaining third vertex in $T$ are five distinct points contained in the intersection $\ellipse\cap C$. This is impossible, since
$\ellipse$ is not a circle. Finally, the remaining vertex has to be on the foci line, or else we could find two distinct orbits sharing a common side, which is impossible by definition of the reflection law with respect to non-isotropic lines.

\textbf{Second case : Suppose $T$ is infinite.} Then the line at infinity cuts $C$ in three distinct points, hence $C$ is degenerate. By Proposition \ref{prop_types_circle}, $C$ contains the line at infinity. Since $T$ has only one infinite vertex $\alpha$ by Proposition \ref{lemma_infinite_orbit}, and two other finite vertices $\beta,\gamma$, the other line $d\subset C$ is not the line at infinity. Again by Proposition \ref{prop_types_circle}, the center is infinite and represents the orthogonal direction to $d$. Since it is on the real foci line, the latter is orthogonal to $d$. Thus $d$ intersects the infinity line at the same point as the line orthogonal to the foci line. This point does not lie in $\ellipse$, and in particular, $d$ does not contain $\alpha$. Hence, we have $d=\beta\gamma$ is a side of $T$, $\alpha\notin d$ and by the same symmetry argument as in the first case $\alpha$ should belong to the real foci line. But this is impossible since the latter intersects $\ellipse$ in only two finite points.

\textbf{Last case : Suppose $T$ is degenerate.} Then $C$ cannot be a regular circle, otherwise the latter would be tangent to $\ellipse$ in a point of isotropic tangency (by Proposition \ref{lemma_olga}): this would imply that this point of isotropic tangency is $I$ or $J$, which is impossible since they do not belong to $\ellipse$, assumed not to be a circle.

The circumcircle $C$ cannot be the union of the line at infinity and another line $d$. Otherwise, by the same arguments as in the second case, this line would be othogonal to the real foci line. Since $T$ is finite (Proposition \ref{lemma_infinite_orbit}), $d$ goes through its both vertices, implying that they are symmetric with respect to the foci line. Therefore, both vertices are points of isotropic tangency but this cannot happen for a degenerate triangular orbit.

Finally suppose $C$ is  the union of two isotropic lines having different directions.

\begin{lemma}
\label{lemma_isotropic_circles}
Let $C_n$ be a sequence of circles containing two distinct points $M_n$ and $N_n$ of $\ellipse$ converging to the same finite point $\alpha$. Suppose $C_n$ has a center $c_n$ converging to a finite point $c\neq\alpha$. Then the line $c\alpha$ is orthogonal to the line $T_{\alpha}\ellipse$.
\end{lemma}

\begin{proof}
The tangent line to $C_n$ at $M_n$ is orthogonal to the line $M_nc_n$ hence the same is true for their limits. The limit of $T_{M_n}C_n$ is obviously the limit of the line $M_nN_n$. Since $M_n$ and $N_n$ are on $\ellipse$, the line $M_nN_n$ also converges to the tangent line $T_{\alpha}\ellipse$. Hence $T_{\alpha}\ellipse$ is orthogonal to $\alpha c$.
\end{proof}

Thus if $\alpha$ is a vertex of isotropic tangency of the orbit, Lemma \ref{lemma_isotropic_circles} implies that $\alpha c$ is orthogonal to $T_{\alpha}\ellipse$, hence $\alpha c = T_{\alpha}\ellipse$ since the latter is isotropic. Recall that $\alpha$ does not lie in the real foci line. Since both isotropic lines constituing the circle go through $c$, one of them is $T_{\alpha}\ellipse$. Hence, they are both tangent to $\ellipse$ by symmetry with respect to the real foci line. Thus the other vertex of $T$ is a point of isotropic tangency of $\ellipse$, which is not possible by the previous arguments (such an orbit is not closed).
\end{proof}

\subsection{Proof of Theorem \ref{theorem:circumcenters}}
\label{sec_proof}

We reall that $\ellipse$ is a complexified ellipse, which we will identify with $\CP^1$. As stated in \cite{poncelet}, the $3$-periodic real orbit are tangent to a smaller confocal ellipse, whose complexification is denoted by $\gamma$.

Consider the Zariski closure $\mathcal{T}$ of the set of real triangular orbits (which are circumscribed about $\gamma$). Let $\mathcal{T}_3$ denote the set of triangles with vertices in $\ellipse$ that are circumscribed about $\gamma$. It is a Zariski closed subset of $\ellipse^3\simeq\left(\CP^1\right)^3$ that contains the real orbits and can be identified with the set of pairs $(A,L)$, where $A$ is a point of the complexified ellipse $\ellipse$ and $L$ is a line through $A$ that is tangent to $\gamma$. The set of the above  pairs $(A,L)$ is identified with an elliptic curve, and each pair extends to a circumscribed triangle as above, see the complex Poncelet Theorem and its proof in \cite{flatto} for more details. Hence $\mathcal{T}_3$ is an irreducible algebraic curve. Each triangle in $\mathcal{T}$ is circumscribed about $\gamma$, by definition and since this is true for the real triangular orbits and the tangency condition of the edges with $\gamma$ is algebraic. Thus $\mathcal{T}\subset\mathcal{T}_3$. Hence $\mathcal{T}=\mathcal{T}_3$, by definition and since the curve of real triangular orbits (which is contained in $\mathcal{T}$) is Zariski dense in $\mathcal{T}_3$ (irreducibility). Now the set $\hat{\mathcal{T}}\subset\mathcal{T}$ of complex non-degenerate triangular orbits circumscribed about the Poncelet ellipse $\gamma$ is a subset of $\mathcal{T}_3 = \mathcal{T}$, Zariski open in $\mathcal{T}$ (because $\mathcal{T}\setminus\hat{\mathcal{T}}$ is defined by polynomial equations). Note that $\mathcal{T}\setminus\hat{\mathcal{T}}$ is finite (since it is a proper Zariski closed subset of an algebraic curve $\mathcal{T}$), and $\hat{\mathcal{T}}$ is dense in $\mathcal{T}$ for the usual topology. Thus the analytic map $\phi : \hat{\mathcal{T}} \to \CP^2$ which assigns to a non-degenerate orbit its circumcenter can be extended to a holomorphic map $\mathcal{T}\to\CP^2$, being a rational map. And by Remmert proper mapping theorem (see \cite{GH}), its image denoted by $\hat{\mathcal{C}}$ is an irreducible analytic curve of $\CP^2$, hence it is an irreducible algebraic curve by Chow theorem (see \cite{GH}). 


Let us show that $\hat{\mathcal{C}}$ is a conic, using Bezout theorem and studying its intersection with the real foci line of $\ellipse$. In fact, we already know two distinct points lying on this intersection: the circumcenters $c_1$ and $c_2$ of both triangular real orbits $T_1$ and $T_2$ circumscribed about Poncelet's ellipse $\gamma$ and having a vertex on the foci line.

\begin{lemma}
\label{lemma_nb_intersection}
The foci line of the ellipse intersects $\hat{\mathcal{C}}$ in only $c_1$ and $c_2$ which are distinct, and for each $i$ the only triangular orbit of $\mathcal{T}$ having $c_i$ as a circumcenter is $T_i$.
\end{lemma}

\begin{proof}
Take a point $c$ of $\hat{\mathcal{C}}$ lying on the foci line. Then by Proposition \ref{prop_orbit_form}, an orbit of center $c$ is finite, non-degenerate, and has a vertex on the foci line. If this orbit is in $\mathcal{T}$, it is circumscribed about $\gamma$. One of its vertices lies on the foci line, hence coincides with a vertex of some $T_i$. Hence it is $T_1$ or $T_2$, otherwise we could find a number strictly greater than two of tangent lines to $\gamma$ containing a vertex of $\ellipse$. Furthermore, if $c_1=c_2$, the circumcircle of $T_1$ would be the same as the one of $T_2$ by symmetry, and $\ellipse$ would share six dictinct points with the former, which is impossible. The result follows.
\end{proof}

\begin{theorem}
\label{lemma_intersection}
The set $\hat{\mathcal{C}} \subset \CP^2$ is a complexified ellipse.
\end{theorem}

\begin{proof}
Let us show that $c_1$ is a regular point of $\hat{\mathcal{C}}$, and that the latter intersects the foci line transversally. Fix an order on the vertices of $T_1$ and consider the germ $(\mathcal{T},T_1)$. The latter is irreducible (because parametrized by $\gamma$), hence the germ $(V,c_1)\subset(\hat{\mathcal{C}},c_1)$ defined as $\phi(\mathcal{T},T_1)$ is also irreducible. By Lemma \ref{lemma_nb_intersection}, any other irreducible component $V'$ of $(\hat{\mathcal{C}},c_1)$ is parametrized locally by $\phi$ and a germ $(\mathcal{T}, T_1')$, where $T_1'$ is obtained from $T_1$ by a permutation of its vertices. Thus $V'=V$ since $\phi$ doesn't change by permutation of the vertices of the orbits: $(\hat{\mathcal{C}},c_1)$ is irreducible. 

We fix a local biholomorphic parametrization $P(t)$ of the complexified ellipse $\ellipse$, so that $P_0 = P(0)$ is a vertex of the real ellipse $\ellipse$ that is also a vertex of the real triangular orbit $T_1$. This gives local parametrizations of the orbits $T(P)$ whose first vertex is $P$ and of their circumcenters $c(t) = \phi(T(P(t)))$. We restrict $P$ to the curve $P(t)$ parametrizing the real points of $\ellipse$. We can suppose that $P(t)$ and $P(-t)$ are symmetric with respect to $\mathcal{F}$. Write $r(t) = |P(t)c(t)|$ for the radius of the circumscribed circle to $T(t)$. Thus we have $c(0)=\phi(T_1)=c_1$, and we need to show that $c'(0) \neq 0$ and that $c'(0)$ has not the same direction as the line of real foci of $\ellipse$. 

First, we have $r(t) = r(-t)$ by symmetry, and $r$ is smooth around $0$ since $P(0) \neq c(0)$. Thus, $r'(0) = 0$. This implies that the vector $c'(0)-P'(0)$ is orthogonal to the line $c(0)P(0)$, which is the real foci line by definition. But $P'(0)$ is already orthogonal to the foci line (being a vector tangent to $\ellipse$ at its vertex $P_0$), hence the same hold for $c'(0)$. It's then enough to show that $c'(0) \neq 0$.

Suppose the contrary, i.e. $c'(0) = 0$. We use again $r'(0) = 0$. If we denote by $Q(t)$ one of the other vertices of $T(t)$ and $Q_0=Q(0)$, then since also $r(t) = |Q(t)c(t)|$, the equality $r'(0) = 0$ gives that the line $Q_0c_1$ is orthogonal to $T_{Q_0}\ellipse$. It means that the circumscribed circle $\mathcal{D}$ to $T_1$ has the same tangent line in $Q_0$ as $\ellipse$. Since this is also true in $P_0$ and in the third point of $T_1$ (same proof), we get that $\ellipse$ and $\mathcal{D}$ have three common points with the same tangent lines, which means that $\ellipse$ is a circle. But this case was excluded at the beginning (remark \ref{rmk_circle}).

Hence $c'(0) \neq 0$ and $c'(0)$ is orthogonal to the line of real foci. The proof is the same for $c_2$. Hence by Bezout theorem, $\hat{\mathcal{C}}$ is a complexified conic. Since its real part is bounded, it is a complexified ellipse.
\end{proof}

\chapter{On the existence of caustics}
	\label{chapter_caustics}
	
This chapter is devoted to the study of \textit{caustics} in complex billiards and projective billiards.

In the classical model of billiard, a \textit{caustic} of a billiard $\Omega$ is a hypersurface $\mathcal{C}$ inside $\Omega$ such that any oriented line tangent to $\mathcal{C}$ and intersecting $\partial\Omega$ transversally is reflected on $\partial\Omega$ into a line tangent to $\mathcal{C}$. This implies that any iterated reflections of a line tangent to $\mathcal{C}$ will produce tangent lines to $\mathcal{C}$.

\textit{Caustics of projective billiards} can be defined similarly:

\begin{definition}
\label{definition:caustic_proj_billiard}
Let $\Sigma\subset\mathbb{P}(T\RR^d)$ be a line-framed hypersurface over a hypersurface $S\subset\RR^d$. A \textit{caustic} of $\Sigma$ is a hypersurface $\Gamma\subset\RR^d$ such that any line $\ell\subset\RR^d$ tangent to $\Gamma$ and interecting $S$ transversally at a point $p$, is reflected into a line tangent to $\Gamma$ by the projective law of reflection at $p$.
\end{definition}

\textit{Caustics of complex billiards} (or \textit{complex caustics}) are difficult to define in the general case since there is no possible orientation of lines. For our purpose to work on conics such definition is more simple, since any line of $\CP^2$ is either tangent to a fixed conic or intersects it in exactly two distinct points. Therefore, complex caustics can be defined as follows in this specific case: let $C,C'\subset\CP^2$ be two distinct conics. We say that $C'$ is a \textit{complex caustic} of $C$ if for any line $\ell$ tangent to $C'$ and $p$ a point of intersection of the line with $C$, the line reflected from $\ell$ by the complex law of reflection at $p$ on $C$ is also tangent to $C'$.

This chapter is structured as follows. Basic results about conics and quadrics are first recalled at Section \ref{section:general_properties_on_quadrics}. Then we present results about complex caustics of the billiard on a complexified ellipse or hyperbola at Section \ref{section_cmplx_caustics}. It is followed by Section \ref{section_projective_caustics} which explains that given a certain pencil of conics or quadrics and any fixed conic or quadric $Q$ of this pencil, $Q$ can be endowed with a structure of projective billiard such that any element of the pencil is a caustic for $Q$. Finally, an argument of Berger \cite{berger_caustics} will be generalized to projective billiards at Section \ref{section_berger_property}, and applied in the case of pseudo-Euclidean billiards to show that if a pseudo-Euclidean billiard has a caustic, then it is itself a quadric.

\section{General properties of quadrics}
	\label{section:general_properties_on_quadrics}
	In this section we describe general properties on conics which we are going to use all along Chapter \ref{chapter_caustics}. They are very classic and can be found in \cite{berger_geometry}, Vol. II, Chap. 13 to 17. 

Let $K$ be the field $\RR$ or $\CC$, $d\geq 1$ an integer and $\pi:K^{d+1}\smallsetminus\{0\}\to\PP{d}{K}$ the natural projection. 

\begin{definition}
A \textit{quadric} $Q$ of $K$ is defined as the image by $\pi$ of sets of the form 
$$Z_q=\ensemble{x\in K^{d+1}\smallsetminus\{0\}}{q(x)=0}$$
where $q$ is a non-zero quadratic form over $K^{d+1}$. The quadric $Q$ is said to be \textit{non-degenerate} if $q$ is non-degenerate, and \textit{non-empty} if $Q\neq\emptyset$. In the specific case when $d=2$, we can also say that $Q$ is a \textit{conic} (in this study, a conic is a quadric).
\end{definition}

The space $\mathcal{Q}(K^{d+1})$ of quadratic forms over $K^{d+1}$ is a vector space such that two non-zero colinear quadratic forms define the same quadric. The converse is false with $K=\RR$, by considering for example the quadratic forms on $\RR^2$ defined by $q_1(x,y)=x^2+y^2$ and $q_2(x,y)=x^2+2y^2$. But in the case when $K=\CC$ the converse is true and is part of a more general theorem on algebraic curves:

\begin{theorem}[Nullstellensatz for quadrics, see \cite{berger_geometry}]
The map $[q]\in\mathbb{P}\mathcal{Q}(\CC^{d+1})\mapsto\pi(Z_q)\subset\PP{d}{\CC}$ is a one-to-one correspondance between equivalence classes $[q]$ of quadratic forms $q$ over $\CC^{d+1}$ and quadrics of $\PP{d}{\CC}$.
\end{theorem}

\begin{figure}[!t]
\centering
\includegraphics[scale=0.3]{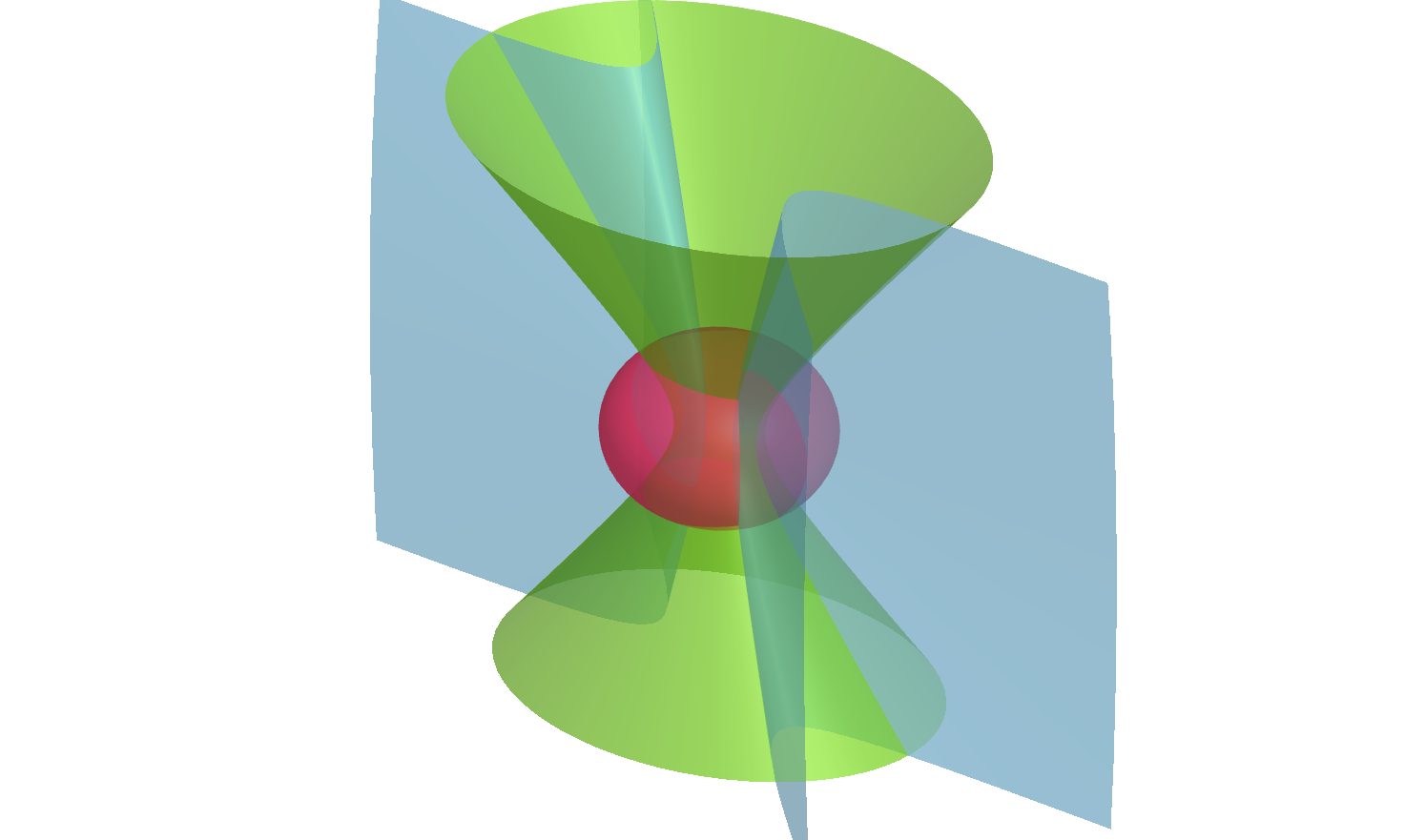}
\caption{Different types of confocal quadrics in dimension $d=3$ depending on the choice of~$\lambda$.}
\label{figure:confocal_quadrics}
\end{figure}

\begin{figure}[!h]
\centering
\includegraphics[scale=0.3]{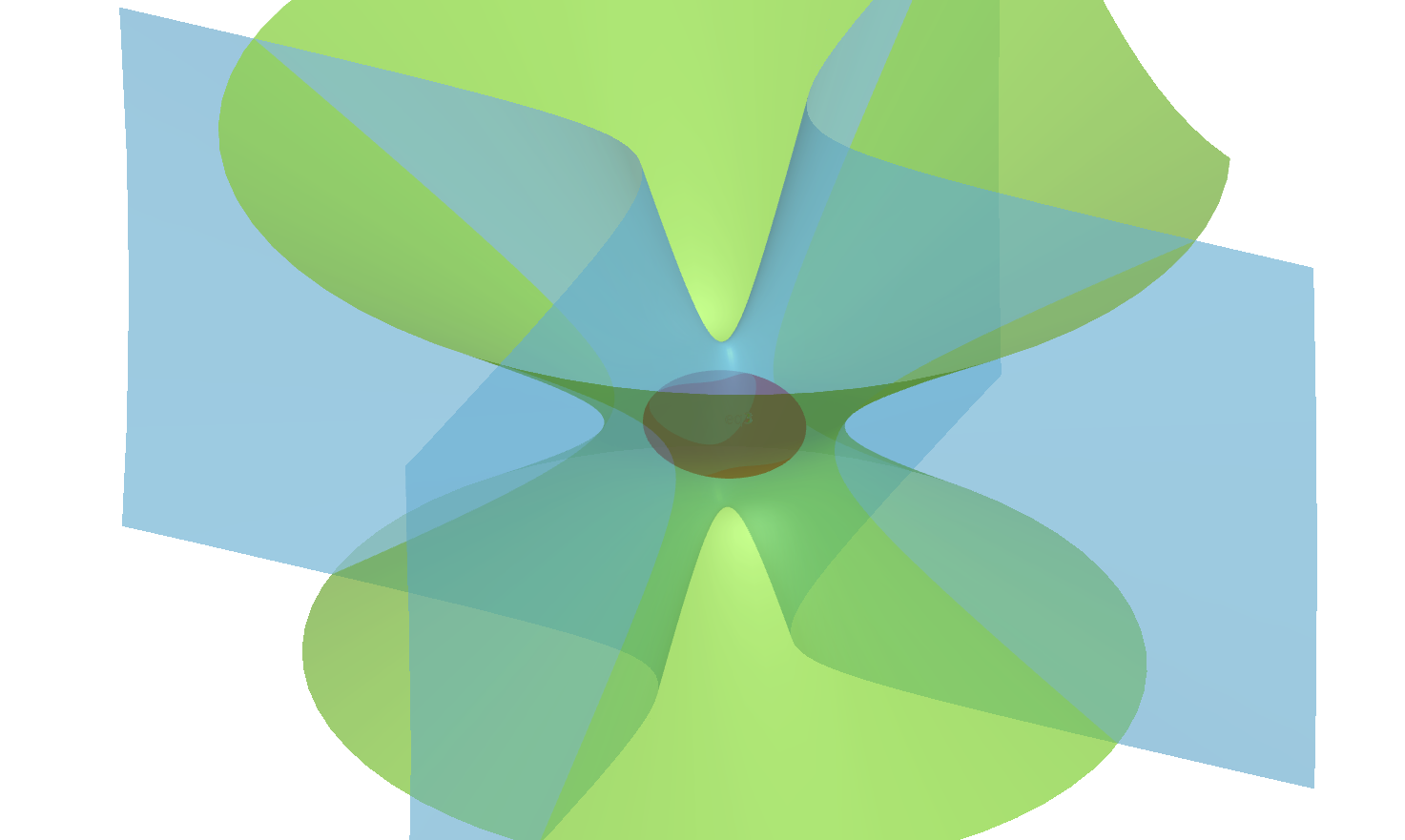}
\caption{Different types of pseudo-confocal quadrics in dimension $d=3$ depending on the choice of~$\lambda$.}
\label{figure:pseudo_confocal_quadrics}
\end{figure}

\begin{example}
Given an integer $k\in\{0,\ldots,d-1\}$, and real numbers $a_0<a_1<\ldots<a_d$, consider the family of quadrics $Q^k:=(Q_{\lambda}^k)_{\lambda\neq a_j}$ of $\RP^d$ given by the equation
\begin{equation}
\label{equation:confocal_quadrics_minkowski}
Q_{\lambda}^k: \sum_{j=0}^k \frac{x_j^2}{a_j-\lambda}+\sum_{j=k+1}^{d-1} \frac{x_j^2}{a_j+\lambda}=x_d^2.
\end{equation}
The quadric $Q^k_0$ is the same for all $k$ and contained in all families $Q^k$. Any non-degenerate quadric can be described by an equation of this form by an appropriate orthogonal change of coordinates. 

When $k=d-1$, this family is the standard family of confocal quadrics, see Figure \ref{figure:confocal_quadrics}. In the case when $k<d-1$, this family is considered to be the family of confocal quadrics for a pseudo-Euclidean metric (Figure \ref{figure:pseudo_confocal_quadrics}, see for example \cite{khesin_taba} or \cite{DragRad_minkowski}) defined as the following non-degenerate quadratic form of $\RR^{d}$
$$\sum_{j=0}^k x_j^2-\sum_{j=k+1}^{d-1} x_j^2.$$
More details on pseudo-Euclidean metrics and one these \textit{pencils of pseudo-confocal quadrics} will be given in Section \ref{section_berger_property}.
\end{example}

\subsection{Polarity with respect to a quadratic form}
\label{subsection:polarity}

In this section we recall some very basic and well-know facts about polarity. We refer the reader to \cite{berger_geometry} Chap. 14, \cite{DragRad} Chap. 4, or \cite{izmestiev} Sec. 2 for more details.

Let $K$ be the field $\RR$ or $\CC$, $d\geq 1$ an integer, $V$ the space $K^{d+1}$, and $\pi:K^{d+1}\smallsetminus\{0\}\to\PP{d}{K}$ the natural projection: if $x\in V$, its equivalence classe in $\PP{d}{K}$ is $\pi(x)$. Given a non-trivial vector subspace $H\subset V$, we define $\PP{}{H}$ to be the set of equivalence classes of non-zero vectors contained in $H$, \textit{ie} $\PP{}{H}=\pi(H\smallsetminus\{0\})$.

\begin{definition}
A \textit{polarity} is the choice of a non-degenerate quadratic form over $V$. When we speak about polarity without an explicit choice of a quadratic form, this implicitely refers to the polarity with respect to the quadratic form 
\begin{equation}
\label{equation:ref_quadratic form}
\mathcal{Q}_0=\sum_{j=0}^d x_j^2
\end{equation}
which is called sometimes \textit{absolute polarity} (see \cite{izmestiev} for this terminology). When $q$ defines a non-empty quadric $Q$, we can also speak of polarity \textit{with respect to $Q$}.
\end{definition}

This choice of quadratic form has concrete geometric consequences, described in what follows.

\underline{An isomorphism.} Given a non-degenerate quadratic form $q$ over $V$, one can consider the isomorphism from $V$ to its dual space $V^{\ast}$ defined by 
$$x\in V\mapsto q(x,\cdot)\in V^{\ast}.$$
It induces a projective isomorphism $\mathcal{I}_{q}:\mathbb{P}(V)\to\mathbb{P}(V^{\ast})$ which only depends on the equivalence class of $q$ in $\mathbb{P}\mathcal{Q}(V)$.

\underline{A bijection between points and projective hyperplanes.} The projective hyperplanes of $\PP{}{V}$ can be identified with $\PP{}{\polar{V}}$, by associating to any non-zero linear form $\alpha$ on $V$ its kernel $\ker\alpha\subset V$. Therefore the map $\mathcal{I}_{q}$ induces an explicit realization of this identification via the map
$$\pi(x)\in\PP{}{V}\mapsto\PP{}{\ker q(x,\cdot)}\subset\PP{}{V}.$$
More generally, we can define a bijective correspondance between $k$-dimensional and $(d-1-k)$-dimensional projective subspaces of $\mathbb{P}(V)$ via the map
$$\mathbb{P}(W)\mapsto\mathbb{P}(W^{\bot_q})$$
where $W^{\bot_q}$ defines the $q$-orthogonal vector subspace of $W$ in $V$.

\begin{definition}
\label{definition:polarity}
Given a projective space $H=\PP{}{W}\subset\PP{}{V}$, we call \textit{polar space of $H$ with respect to $q$} the projective space $\mathbb{P}(W^{\bot_q})$. The \textit{polar space of a point with respect to $q$} is a projective hyperplane. The polar space of a hyperplane is a point, also called \textit{pole} with respect to $q$ of this hyperplane.
\end{definition}

\underline{Dual hypersurfaces/curves.} Let $\Gamma$ be a $\class^1$-smooth hypersurface of $\mathbb{P}(V)$. The projective hyperplane containing $p$ and tangent to $\Gamma$ is a projective space $T_p\Gamma=\PP{}{W}$, and we can consider its pole with respect to $q$: the latter is the point $u_p=\pi(x)$ such that $x$ is $q$-orthogonal to $W$. The collection of all $u_p$ when $p$ describes $\Gamma$ is called \textit{dual} of $\Gamma$ with respect to $q$  and denoted by $\Gamma^{\ast}$.

\underline{Dual quadric.} Let $q_1$ be another non-degenerate quadratic form. There is a $(d+1)\times(d+1)$ invertible matrix $M$ with coefficients in $K$ such that for all $x,y\in V$, $q_1(x,y)=q(Mx,y)$. We define the \textit{dual} of $q_1$ with respect to $q$ as the quadratic form $\polar{q_1}$ over $V$ satisfying for all $x,y\in V$ the equality $\polar{q_1}(x,y)=q(\inverse{M}x,y)$. It is well-known that the dual with respect to $q$ of a quadric $Q_1$ defined by the quadratic form $q_1$ is a quadric defined by the dual $\polar{q_1}$ of $q_1$ with respect to $q$, see for example \cite{izmestiev}:

\begin{proposition}
\label{proposition:dual_quadrics_quadratic_forms}
The dual of a non-empty non-degenerate quadric $Q_1$ defined by the quadratic form $q_1$ is the quadric defined by the quadratic form $\polar{q_1}$.
\end{proposition}

\subsection{Pencil of quadrics}

We recall that $\mathcal{Q}(K^{d+1})$ is the set of quadratic forms over $K^{d+1}$. In this subsection we will abusively write quadrics for quadratic forms.

\begin{definition}[see \cite{berger_geometry}, Ch. 14.1]
A \textit{pencil of quadrics} is a line in $\mathbb{P}\mathcal{Q}(K^{d+1})$. It can be equivalently defined as a set of quadratic forms of the type
$$\mathcal{F}(q_1,q_2):=\ensemble{\lambda q_1+\mu q_2}{(\lambda,\mu)\in K^2\smallsetminus\{0\}}$$
where $q_1,q_2$ are quadratic forms with distinct equivalence classes in $\mathbb{P}\mathcal{Q}(K^{d+1})$ (non-colinear). We say that the pencil $\mathcal{F}(q_1,q_2)$ is \textit{non-degenerate} if it contains at least one non-degenerate quadratic form.
\end{definition}

Since the map $(\lambda,\mu)\mapsto\det(\lambda q_1+\mu q_2)$ is a homogeneous polynomial of degree at most $d+1$ over $K$, a pencil of quadrics contains either only degenerate quadratic forms, or a finite number less than $d+1$ of degenerate quadratic forms.

Let us consider the absolute polarity, that is the polarity with respect to the quadratic form $\mathcal{Q}_0=\sum_j x_j^2$ (see Subsection \ref{subsection:polarity}). Given a non-degenerate quadratic form $q$, we can define its dual $q^{\ast}$.

\begin{definition}
The \textit{dual pencil of quadrics} associated to a non-degenerate pencil of quadrics $\mathcal{F}(q_1,q_2)$, is the set $\polar{\mathcal{F}(q_1,q_2)}$ of duals $\polar{q}$ of non-degenerate quadratic forms $q$ contained in $\mathcal{F}(q_1,q_2):$
$$\polar{\mathcal{F}(q_1,q_2)} = \ensemble{\polar{q}}{q\in \mathcal{F}(q_1,q_2), q \text{ is non-degenerate}}=\ensemble{\polar{(\lambda q_1+\mu q_2)}}{\det(\lambda q_1+\mu q_2)\neq 0}.$$
The \textit{pencil of $(q_1,q_2)$-confocal quadrics} is the set $\polar{\mathcal{F}(\polar{q_1},\polar{q_2})}$ which contains $q_1$ and $q_2$. 
\end{definition}

\begin{remark}
The \textit{pencil of $(q_1,q_2)$-confocal quadrics} $\polar{\mathcal{F}(\polar{q_1},\polar{q_2})}$ contains $q_1$ and $q_2$ by involutivity of polarity operations.
\end{remark}

\begin{example}
Consider two confocal conics $C$ and $D$ of $\RP^2$ defined by the quadratic forms $q_C(x,y,z)=\frac{x^2}{a}+\frac{y^2}{b}-z^2$ and $q_D(x,y,z)=\frac{x^2}{a-\lambda}+\frac{y^2}{b-\lambda}-z^2$ with $a,b,\lambda\in\RR$, $\lambda\notin\{a,b\}$. By Proposition \ref{proposition:dual_quadrics_quadratic_forms}, their dual conics are defined by their dual quadratic forms $\polar{q_C}=ax^2+by^2-z^2$ and $\polar{q_D}=(a-\lambda)x^2+(b-\lambda)y^2-z^2$. Hence $\polar{q_D}$ belong to the pencil of quadrics $\mathcal{F}(\polar{q_C},\polar{q_D})=\mathcal{F}(\polar{q_C},q_{\text Eucl})$ where $q_{\text Eucl}$ is the degenerate quadratic form $q_{\text Eucl}(x,y)=x^2+y^2$. 

More generally, confocal conics or quadrics (for the usual meaning) can be defined as quadrics of a pencil of the form $\polar{\mathcal{F}(\polar{q_1},\polar{q_2})}$, which contains the quadratic form defining the Euclidean metric after an eventual change of coordinates (see \cite{izmestiev}). This explains our terminology.
\end{example}

\subsection{Theorems of Poncelet and Cayley}

The theorems of Poncelet and Cayley are remarkable results on conics and many different versions of these theorems exists, see for example \cite{berger_geometry}, Section 16.6, but also \cite{CKS,DragRad,DragRad_bicent,GHponc}, and \cite{poncelet} for the original statement. Here we present the version of \cite{flatto,GHponc} since we consider conics of $\CP^2$.

Let $C,D$ be two conics in $\CP^2$. We say that $C,D$ are \textit{in general position} if their intersection consists of four distinct points. The statements of Poncelet's and Cayley's theorems are about polygons inscribed in $C$ and circumscribed about $D$: an $n$-sided \textit{polygon} is an ordered set $P=(p_1,\ldots,p_n)$ of distinct points of $\CP^2$ called the \textit{vertices} of $P$. An $n$-sided polygon $P$ is said to be \textit{inscribed in $C$} if $p_j\in C$ for all $j$ and \textit{circumscribed about} $D$ if for all $j$, the two tangent lines to $D$ containing $p_j$ are $p_{j-1}p_j$ and $p_jp_{j+1}$ (where the indices $j-1$ and $j+1$ are seen modulo $n$).

\begin{theorem}[Poncelet, \cite{GHponc} p. 3]
\label{theorem:poncelet}
Let $C,D$ be two conics of $\CP^2$ in general position. Suppose that there is an $n$-sided polygon inscribed in $C$ and circumscribed about $D$. Then for any point $p\in C$ there is an $n$-sided polygon inscribed in $C$ and circumscribed about $D$ having $p$ as a vertex.
\end{theorem}

The natural question which arises is about the existence of such $n$-sided polygons. The answer is given by Cayley's theorem.

\begin{theorem}[Cayley, \cite{GHponc} p. 4]
\label{theorem:cayley}
Let $C,D$ be two conics of $\CP^2$ in general position. Let $Q_C$, $Q_D$ be two quadratic forms defining respectively $C$ and $D$. Consider an analytic branch of $t\mapsto\sqrt{\det(tQ_C+Q_D)}$ defined in a neighborhood of $0$ and denote its analytic expansion at $0$ by 
$$\sqrt{\det(tQ_C+Q_D)} = A_0+A_1t+A_2t^2\ldots$$
Then there is an $n$-sided polygon inscribed in $C$ and circumscribed about $D$ if and only if
$$\left|\begin{matrix}
A_2&\ldots&A_{m+1}\\
\vdots & \ddots & \vdots\\
A_{m+1}&\ldots&A_{2m}
\end{matrix}\right|=0, \qquad \text{when }n\text{ is odd, with } m=\frac{n-1}{2},$$
or
$$\left|\begin{matrix}
A_3&\ldots&A_{m+1}\\
\vdots & \ddots & \vdots\\
A_{m+1}&\ldots&A_{2m-1}
\end{matrix}\right|=0 \qquad \text{when }n\text{ is even, with } m=\frac{n}{2}.$$
\end{theorem}

\section{Complex caustics of complexified conics}
	\label{section_cmplx_caustics}
	
We present in this section what can be considered as a complexified version of the result stating that given a conic $C$ of the Euclidean plane, any confocal conic $C'$ to $C$ is a caustic of the billiard on $C$. The results presented in this section can also be found in \cite{fierobe_caustics}.

\begin{definition}
\label{definition:complex_caustic}
Let $C\subset\CP^2$ be a conic. Given another conic $C'\subset \CP^2$, we say that $C'$ is a \textit{complex caustic} of $C$ if any line tangent to $C'$ and intersecting $C$ at a certain point $p$ is reflected into a line tangent to $C'$ by the complex reflection law at $p$.
\end{definition}

Let $C,C'$ be conics such that $C'$ a complex caustic of $C$. Suppose we have $n$ distinct points $p_1,\ldots,p_n$ on $C$. Definition \ref{definition:complex_caustic} implies that the following statements are equivalent:\\
\enum $(p_1,\ldots,p_n)$ is a piece of non-degenerate orbit of $C$ (see Definition \ref{definition:complex_ndegenerate_orbit}) such that $p_jp_{j+1}$ is tangent to $C'$ for a certain $j<n-1$;\\
\enum for each $j\in\{2,\ldots,n-1\}$, the tangent lines to $C'$ containing $p_j$ are exactly the lines $p_{j-1}p_{j}$ and $p_{j}p_{j+1}$.

In the case of $n$-periodic orbits, if $C'$ is a caustic of $C$, the $n$-periodic orbits of $C$ are the same as the $n$-sided polygons circumscribed about $C'$. Hence Poncelet's theorem (see theorem \ref{theorem:poncelet}) implies that if an orbit circumscribed about some caustic $C'$ is $n$-periodic, then all orbits circumscribed about $C'$ are $n$-periodic:

\begin{proposition}
\label{prop:poncelet_for_billiard}
Let $C,C'$ be conics in general position such that $C'$ is a complex caustic of $C$. Suppose that there is an $n$-periodic orbit of $C$ circumscribed about $C'$ (as an $n$-sided polygon). Then any billiard orbit of $C$ circumscribed about $C'$ is $n$-periodic.
\end{proposition}

This induces the following definition:

\begin{definition}
Given two conics $C,C'$ in general position, we say that $C'$ is an \textit{$n$-caustic} of $C$ if $C'$ is a caustic of $C$ about which an $n$-periodic orbit of $C$ is circumscribed.
\end{definition}

\subsection{Confocal conics are complex caustics}
\label{subsection_conics_are_cautics}

In the following we show that given the complexification $\mathcal C$ of a real conic, its confocal conics are caustics. Suppose that we are given a set of coordinates $(x:y:z)$ on $\CP^2$ such that $\mathcal C$ is defined by the following equation in the affine chart $U_z=\{z=1\}$
\begin{equation}
\label{eq:reference_conic}
\mathcal C: \frac{x^2}{a}+\frac{y^2}{b}=1
\end{equation}
where $x,y\in\CC$ and $a,b\in\RR^{\ast}$. The confocal conics $\mathcal C_{\lambda}$ to $\mathcal{C}$ are given by the following family of equations depending on a $\lambda\in\CC$ different from $a$ or $b$:
\begin{equation}
\label{eq:confocal_conics}
\mathcal C_{\lambda}:\frac{x^2}{a-\lambda}+\frac{y^2}{b-\lambda}=1.
\end{equation}

\begin{remark}
\label{remark:real_elliptic_billiard} In the case of the real elliptic billiard, that is when $a,b$ are positive and we study the usual billiard inside the ellipse $\mathcal{C}$, it is well-known (see \cite{taba_book} Chapt. 4) that the real conics given by Equation \eqref{eq:confocal_conics} with $0<\lambda<a$ and $\lambda\neq b$ are caustics in the usual meaning. Let $F_1,F_2$ be the two foci of the ellipse $\mathcal{C}$. Given an orbit of the elliptic billiard, we distinguish between three disjoint situations:\\
1) If the orbit has an edge containing a focus, then all its edges alternatively contain one of both foci.\\
2) If the orbit has an edge intersecting the interior of the segment $F_1F_2$, then all its edges intersect the interior of $F_1F_2$ and remain tangent to the same hyperbola $\mathcal{C}_{\lambda}$ with $b<\lambda<a$.\\
3) If the orbit has an edge which does not intersect $F_1F_2$, then all its edges do not intersect $F_1F_2$ and remain tangent to the same smaller confocal ellipse $\mathcal{C}_{\lambda}$ with $0<\lambda<b$.
\end{remark}

\begin{proposition}
\label{prop:confocal_conics_are_caustics}
For any $\lambda\in\CC\smallsetminus\{a,b\}$, the confocal conic $\mathcal{C}_{\lambda}$ is a complex caustic of $\mathcal{C}$.
\end{proposition}

\begin{proof}
First notice that given $p\in\mathcal{C}$, the tangent line to $\mathcal C$ at $p$ is not the line at infinity defined by $L_{\infty}=\{z=0\}$. Therefore the complex reflection law induces a projective transformation on the set of lines containing $p$, hence on the projective line $p^{\ast}\simeq\CP^1$ defined as the polar space of $p$, and its action on $p^{\ast}$ is denoted by $q\mapsto q'$. 

For $\lambda\neq a,b$, the absolute dual conic $\mathcal{C}_{\lambda}^{\ast}$ of $\mathcal{C}_{\lambda}$ is given by the equation 
$(a-\lambda)x^2+(b-\lambda)y^2=1$
and thus is also defined for $\lambda = a$ or $b$ (as a degenerate conic). Hence we can consider the set 
$$V=\ensemble{(p,q,\lambda)\in \mathcal C\times\CP^2\times \CC}{q\in p^{\ast}\cap\mathcal{C}_{\lambda}^{\ast}}$$
which is an algebraic subset of $\CP^2\times \CC$ since it is given by polynomial equations. Let $V_0$ be the algebraic subset of $V$ containing the elements $(p,q,\lambda)\in V$ such that $(p,q',\lambda)\in V$.

If $\lambda$ is a real number, denote by $\mathcal C_{\lambda}^{\RR}$ the points of $\mathcal{C}_{\lambda}$ with real coordinates which can be considered as the conic of $\RP^2$ defined by Equation \eqref{eq:confocal_conics}. If $p$ is a point on $\mathcal{C}^{\RR}$ and $\lambda$ is a real number different from $a$ or $b$, we know that the line of $\RP^2$ containing $p$ and tangent to $\mathcal{C}_{\lambda}^{\RR}$ is reflected into a line tangent to $\mathcal{C}_{\lambda}^{\RR}$ by the usual reflection law at $p$ (see Remark \ref{remark:real_elliptic_billiard}). Hence the same holds for the complexification of these objects since the same equations are satisfied.

Hence the map $s:\mathcal C\times\CP^2\times \CC\to\mathcal C\times \CC$ defined by $(p,q,\lambda)\mapsto(p,\lambda)$ is such that $s(V_0)$ contains $\mathcal{C}^{\RR}\times\left(\RR\smallsetminus\{a,b\}\right)$. Now since $\mathcal C$ can be identified with $\CP^1$, $\mathcal C\times\CP^2\times \CC$ and $\mathcal C\times \CC$ are projective spaces, and therefore $s(V_0)$ is an algebraic subset of $\mathcal{C}\times\CC$. From the identification $\mathcal{C}^{\RR}\simeq\RP^1$ we get that $\mathcal{C}^{\RR}\times\left(\RR\smallsetminus\{a,b\}\right)$ is Zariski-dense. Hence $s(V_0)=\mathcal{C}\times\CC$: this means that if $(p,\lambda)\in\mathcal{C}\times\CC$ with $\lambda\neq a$ or $b$, there is a $q\in p^{\ast}\cap\mathcal{C}_{\lambda}^{\ast}$ for which $(p,q,\lambda)\in V_0$, and by construction we have exactly $p^{\ast}\cap\mathcal{C}_{\lambda}^{\ast}=\{q,q'\}$. Therefore both lines tangent to $\mathcal{C}_{\lambda}$ and containing $p$ are reflected into each other by the complex reflection law at $p$.
\end{proof}

\subsection{Number of complex confocal $n$-caustics}

Given an integer $n\geq 2$, and a real conic $\mathcal C$, we would like to study the $n$-caustics of the complex billiard $\mathcal C$. Caustics of $n$-periodic orbits of the real elliptic billiard are such that their complexifications are $n$-caustics of the corresponding complex billiard by definition. We will show that other complex $n$-caustics can appear.

In the case when $n=3$ and $\mathcal{C}$ is an ellipse, it is well-known that the usual $3$-periodic orbits of $\mathcal{C}$ are all circumscribed about exactly one smaller confocal ellipse $\gamma_3$. Therefore, the complexification of $\gamma_3$ is a $3$-caustic of $\mathcal{C}$. We can ask if it is the only $3$-caustic confocal to the complexification of $\mathcal{C}$. The answer is no, since there is another complexified ellipse confocal to $\mathcal{C}$ which is a $3$-caustic as it will be shown in this subsection. Interestingly, this caustic is bigger than $\mathcal{C}$.

\begin{remark} In the case of the real elliptic billiard, we can associate to a caustic of a periodic orbit an invariant quantity called \textit{rotation number} which is an integer. It can be defined as follows (see \cite{taba_book} Chapt. 6). Parametrize the ellipse $\mathcal{C}$ by $\SS^1=\RR/\ZZ$. For a periodic orbit given by parameters $(x_1,\ldots,x_n)\in\left(\SS^1\right)^n$, consider $t_1,\ldots,t_n\in(0,1)$ such that for each $k$ modulo $n$, the class of $t_k$ in $\RR/\ZZ$ is $x_{k+1}-x_k$. Since the orbit is closed, the quantity $\rho=t_1+\ldots+t_n$ is an integer called rotation number of the orbit. Now since this quantity depends continuously on the orbit, it is the same for all periodic orbits circumscribed about the same caustic. As a consequence, periodic orbits with different rotation numbers are circumscribed about distinct caustics. Birkhoff's theorem (see \cite{taba_book} Chapt. 6 or \cite{treshchev} Chapt. II) states, for all $n\geq 2$ and $\rho\leq\lfloor(n-1)/2\rfloor$ coprime with $n$, the existence of $n$-periodic orbits, hence the existence of caustics with rotation number $\rho$ in the elliptic case.
\end{remark}

\begin{remark} If a complex $n$-caustic $\mathcal{C}_{\lambda}$ of $\mathcal{C}$ is inscribed in a periodic orbit with all its vertices having real coordinates then $\lambda$ is a real number comprised between $0$ and $a$ (see Remark \ref{remark:real_elliptic_billiard}). Hence if $\lambda$ is a complex number outside $[0,a]$ corresponding to an $n$-caustic $\mathcal{C}_{\lambda}$, then the periodic orbits circumscribed about $\mathcal{C}_{\lambda}$ have at least one point with a strictly complex coordinate. They corresponds to either complexified bigger confocal ellipses (case $\lambda\in\RR^{-}$), or to what will be called \textit{strictly complex confocal conics} (case $\lambda\in\CC\smallsetminus]-\infty,a]$). As it will be shown, the case $n=4$ provides examples of $4$-caustics of each of the above described types.
\end{remark}

\subsubsection{Counting $n$-caustics using Cayley's determinant}

Let $\mathcal{C}$ be the conic given by Equation \eqref{eq:reference_conic}, and $\mathcal{C}_{\lambda}$ the family of its confocal conics given by Equation \eqref{eq:reference_conic}. Fix an integer $n\geq 3$: we study the number $\mathcal{N}_{a,b}(n)$ of confocal complex $n$-caustics of $\mathcal{C}$. 

As stated in Proposition \ref{prop:confocal_conics_are_caustics}, each $\mathcal{C}_{\lambda}$ is a caustic of $\mathcal{C}$ (with $\lambda\neq a,b$). For $\lambda\neq 0$, $\mathcal{C}_{\lambda}$ and $\mathcal{C}$ are in general position, hence we can study the complex numbers $\lambda$ for which $\mathcal{C}_{\lambda}$ is an $n$-caustic of~ $\mathcal{C}$.

In this subsection we prove the following results:

\begin{proposition}
\label{prop:cmplx_caustics_finite_confocal}
Let $n\geq3$. There is a polynomial $\mathcal{B}^n_{a,b}(\lambda)$ such that  $\lambda\notin\{a,b\}$ is a root of $\mathcal{B}^n_{a,b}(\lambda)$ if and only if $\class_{\lambda}$ is an $n$-caustic of $\mathcal{C}$. 

The degree of $\mathcal{B}^n_{a,b}$ satisfies
$$\deg\mathcal{B}^n_{a,b}\leq\left\{\begin{matrix}
\frac{n^2-1}{4}&\text{if } n \text{ is odd}\\
\frac{n^2}{4}-1&\text{if } n \text{ is even.}\\
\end{matrix}\right.$$
If $\mathcal{B}^n_{a,b}$ has only simple roots distinct from $a$ and $b$ then $\mathcal{N}_{a,b}(n)=\deg \mathcal{B}^n_{a,b}$.
\end{proposition}

\begin{proposition}
\label{prop:cmplx_caustics_degree_polynomial}
There exist $r_1,\ldots,r_p\in\RR$ such that for all $(a,b)$ with $a/b\notin\{r_1,\ldots,r_p\}$, we have 
$$\deg\mathcal{B}^n_{a,b}=\left\{\begin{matrix}
\frac{n^2-1}{4}&\text{if } n \text{ is odd},\\
\frac{n^2}{4}-1&\text{if } n \text{ is even.}\\
\end{matrix}\right.$$
\end{proposition}

\begin{proposition}
\label{prop:cmplx_caustics_forbidden_roots}
There exist $r'_1,\ldots,r'_q\in\RR$ such that for all $(a,b)$ with $a/b\notin\{r'_1,\ldots,r'_q\}$, $a$ and $b$ are not roots of $\mathcal{B}^n_{a,b}$.
\end{proposition}

\begin{remark}
We show in Proposition \ref{prop:circle_polynomial} that $p\geq 1$ by studying the case of the circle, more precisely that $1$ belongs to the collection of $\{r_1,\ldots,r_p\}$.
\end{remark}

\newcommand{\cat}{\text{Cat}}

\begin{proof}[Proof of Proposition \ref{prop:cmplx_caustics_finite_confocal}]
Suppose first that $n=2m+1$ is odd and fix a $\lambda\neq a,b$. As explained, $\mathcal{C}_{\lambda}$ is an $n$-caustic if and only if on can find an $n$-sided polygon inscribed in $\mathcal{C}$ and circumscribed about $\mathcal{C}$. Hence we apply Cayley's theorem (see Theorem \ref{theorem:cayley}): there is such a polygon if and only if the determinant
$$\mathcal{A}^n(\lambda) = \left|\begin{matrix}
A_2(\lambda)&\ldots&A_{m+1}(\lambda)\\
\vdots & \ddots & \vdots\\
A_{m+1}(\lambda)&\ldots&A_{2m}(\lambda)
\end{matrix}\right|$$
vanishes, where the $A_k(\lambda)$ are the coefficients in the analytic expansion of 
$$f:t\to\sqrt{\det(tQ_0+Q_{\lambda})}$$
where $Q_0$ and $Q_{\lambda}$ are quadratic forms respectively associated to $\mathcal{C}$ and to $\mathcal{C}_{\lambda}$. The quadratic form $Q_{\lambda}= (a-\lambda)^{-1}x^2+(b-\lambda)^{-1}y^2-z^2$ defines $\mathcal{C}_{\lambda}$ in $\CP^2$. Replacing $\lambda$ by $0$, we get $Q_0$. Therefore 
$$tQ_0+Q_{\lambda} = \left(\frac{t}{a}+\frac{1}{a-\lambda}\right)x^2+\left(\frac{t}{b}+\frac{1}{b-\lambda}\right)y^2-(t+1)z^2$$
hence
$$\det(tQ_0+Q_{\lambda}) =-\left(\frac{t}{a}+\frac{1}{a-\lambda}\right)\left(\frac{t}{b}+\frac{1}{b-\lambda}\right)(t+1)$$
which we factorize in
$$\det(tQ_0+Q_{\lambda}) =-\frac{1}{(a-\lambda)(b-\lambda)}\left(\frac{a-\lambda}{a}t+1\right)\left(\frac{b-\lambda}{b}t+1\right)(t+1).$$
Define the map $g:t\mapsto\sqrt{\left(\frac{a-\lambda}{a}t+1\right)\left(\frac{b-\lambda}{b}t+1\right)(t+1)}$ and write its Taylor expansion as
$$g(t) = \sum_{k=0}^{\infty} B_k(\lambda) t^.k$$
Since 
$$f(t) = \frac{ig(t)}{\sqrt{(a-\lambda)(b-\lambda)}}$$
we have
$$A_k(\lambda) =\frac{iB_k(\lambda)}{\sqrt{(a-\lambda)(b-\lambda)}}.$$
This shows that $\mathcal{A}^n(\lambda)$ is a function of $\lambda$ which vanishes at $\lambda\neq a,b$ if and only if the determinant 
$$\mathcal{B}^n(\lambda) = \left|\begin{matrix}
B_2(\lambda)&\ldots&B_{m+1}(\lambda)\\
\vdots & \ddots & \vdots\\
B_{m+1}(\lambda)&\ldots&B_{2m}(\lambda)
\end{matrix}\right|$$
also vanishes. Let us compute the $B_k$'s. Write $\sqrt{t+1} = c_0+c_1t+c_2t^2+\ldots$ where
\begin{equation}
\label{eq:c_k}
c_k = \frac{1}{k!}\left(\frac{1}{2}\right)\left(\frac{1}{2}-1\right)\ldots\left(\frac{1}{2}-k+1\right) = \frac{(-1)^{k+1}}{4^k(2k-1)}\binom{2k}{k}.
\end{equation}
Therefore for any $\beta$ we have $\sqrt{\beta t+1} = c_0+c_1\beta t+c_2\beta^2 t^2+\ldots$ Hence $B_k(\lambda)$ is given by
\begin{equation}
\label{eq:B_k}
B_k(\lambda) = \sum_{u+v+w=k} \frac{c_uc_vc_w}{a^{u}b^{v}}(a-\lambda)^u(b-\lambda)^v.
\end{equation}
Therefore each $B_k$ is a polynomial in $\lambda$ of degree at least $k$. Hence $\mathcal{B}^n(\lambda)=\mathcal{B}^n_{a,b}(\lambda)$ is a polynomial in $\lambda$ verifying: for any $\lambda\neq a,b$, $\mathcal{B}^n_{a,b}(\lambda)=0$ if and only if $\mathcal{A}^n(\lambda)=0$, which is true if and only if there exists an $n$-sided polygon inscribed in $\mathcal{C}$ and circumscribed about $\class_{\lambda}$.  The same proof also works when $n$ is even.

It remains to give an upper bound on $\deg \mathcal{B}^n_{a,b}(\lambda)$. Suppose first that $n=2m+1$ is odd. For any permutation $\sigma$ of $\{1,\ldots,m\}$ we have
$$\deg \prod_{j=1}^{m} B_{\sigma(j)+j} = \sum_{j=1}^m \deg B_{\sigma(j)+j} \leq \sum_{j=1}^m \left(\sigma(j)+j\right) = m(m+1)$$
and since $\mathcal{B}^n_{a,b}(\lambda)$ is a sum of $\pm\prod_{j=1}^{m} B_{\sigma(j)+j}$ over all $\sigma$, we have $\deg \mathcal{B}^n_{a,b}(\lambda) \leq m(m+1) = \frac{n^2-1}{4}$. If $n=2m$ is even, Cayley's determinant gives $\mathcal{B}^n_{a,b}(\lambda)=\det(B_{i+j+1})_{1\leq i,j\leq m-1}$. Hence for any permutation $\sigma$ of $\{1,\ldots,m\}$ we have
$$\deg \prod_{j=1}^{m-1} B_{\sigma(j)+j+1} = m^2-1$$
and the same argument leads to $\deg\mathcal{B}^n_{a,b}\leq m^2-1=\frac{n^2}{4}-1$.
\end{proof}

\begin{proof}[Proof of Proposition \ref{prop:cmplx_caustics_degree_polynomial}]
Suppose $n=2m+1$ is odd. By Equation \eqref{eq:B_k}, $B_k$ is of degree $\leq k$ and the coefficient in front of $\lambda^k$ is
\begin{equation}
\label{eq:domin_B_k}
d(B_k)=(-1)^k\sum_{u+v=k} \frac{c_uc_v}{a^{u}b^{v}}=\frac{1}{4^{k}}\sum_{u+v=k} \frac{1}{a^{u}b^{v}(2u-1)(2v-1)}\binom{2u}{u}\binom{2v}{v}.
\end{equation}
Fix a permutation $\sigma$ of $\{1,\ldots,m\}$. We have
$$\deg \prod_{j=1}^{m} B_{\sigma(j)+j} = \sum_{j=1}^m \deg B_{\sigma(j)+j} \leq \sum_{j=1}^m \left(\sigma(j)+j\right) = m(m+1)$$
and the coefficient in front of $\lambda^{m(m+1)}$ is $\prod_{j=1}^{m} d(B_{\sigma(j)+j})$.
Since $\mathcal{B}^n_{a,b}(\lambda)$ is a sum of $\pm\prod_{j=1}^{m} B_{\sigma(j)+j}$ over all $\sigma$, we have that $\deg \mathcal{B}^n_{a,b}(\lambda) \leq m(m+1)$, and the coefficient in front of $\lambda^{m(m+1)}$ is
$$d_n(a,b)=\left|\begin{matrix}
d(B_2)&\ldots&d(B_{m+1})\\
\vdots & \ddots & \vdots\\
d(B_{m+1})&\ldots&d(B_{2m})
\end{matrix}\right|.$$
Let us show that $d_n(a,b)\neq 0$ except for specific $(a,b)$ as described in Proposition \ref{prop:cmplx_caustics_degree_polynomial}. Note first that each $d(B_k)$ is a homogeneous polynomial in $(a^{-1},b^{-1})$ of degree $k$, and by Equation \eqref{eq:domin_B_k} the coefficient in front of $a^{-k}$ is 
$$-\frac{1}{4^{k}(2k-1)}\binom{2k}{k}=-\frac{2}{4^{k}}\cat_{k-1}$$
where $\cat_k=\frac{1}{k+1}\binom{2k}{k}$ is the $k$-th Catalan number.

Now by $m$-linearity of the determinant, $d_n(a,b)$ is also a homogeneous polynomial in $(a^{-1},b^{-1})$, and we apply the same procedure as before: for any permutation $\sigma$ of $\{1,\ldots,m\}$, we have 
$$\deg \prod_{j=1}^{m} d(B_{\sigma(j)+j}) = \sum_{j=1}^m \deg d(B_{\sigma(j)+j}) = \sum_{j=1}^m \left(\sigma(j)+j\right) = m(m+1)$$
and the coefficient in front of $a^{-k}$ is
$$\prod_{j=1}^{m} \frac{-2}{4^{j+\sigma(j)}}\cat_{j+\sigma(j)-1}=\frac{(-1)^m}{2^{m(2m+1)}}\prod_{j=1}^{m}\cat_{j+\sigma(j)-1}.$$
Since $d_n(a,b)$ is a sum of $\pm\prod_{j=1}^{m} d(B_{\sigma(j)+j})$ over all $\sigma$, we have that $\deg d_n(a,b) \leq m(m+1)$, and the coefficient in front of $a^{-m(m+1)}$ is
$$\frac{(-1)^m}{2^{m(2m+1)}}\det H_m$$
where $H_m$ is the Hankel matrix of the sequence $(\cat_{k+1})_k$ defined as
$$H_m = \left(\begin{matrix}
\cat_1&\cat_2&\cdots&\cat_m\\
\cat_2&\cat_3&&\\
\vdots&&\ddots&\vdots\\
\cat_m&&\cdots&\cat_{2m-1}
\end{matrix}\right).$$
One can show that $\det H_m=1$, see for example \cite{krattenthaler} Theorem $33$, or \cite{krattenthaler_cat} Formula $(1.2)$ for the case when $n$ is odd and Formula $(1.3)$ for the case when $n$ is even. Hence $d_n(a,b)$ is a non-zero homogeneous polynomial in $(a^{-1},b^{-1})$ and therefore there exists a finite collection of numbers $r_1,\ldots,r_p\in\RR$ such that for all $a,b>0$, we have $d_n(a,b)=0$ if and only if $a/b\in\{r_1,\ldots,r_p\}$.
\end{proof}

\begin{proof}[Proof of Proposition \ref{prop:cmplx_caustics_forbidden_roots}]
Suppose $n=2m+1$ is odd. By Equation \eqref{eq:B_k}, for $k\geq 2$,
\begin{equation}
\label{eq:domin_B_k_2}
B_k\left(-a^2\right)=\sum_{v+w=k} \frac{c_vc_w}{b^{2v}}(b^2-a^2)^v=\frac{1}{b^{2k}}\sum_{v+w=k} c_vc_wb^{2w}(b^2-a^2)^v = \frac{1}{b^{2k}}P_k(a,b)
\end{equation}
where $P_k(a,b)$ is a homogeneous polynomial in $(a,b)$ of degree $2k$. The coefficient in front of $a^{2k}$ is 
$$(-1)^kc_k=-\frac{1}{4^k(2k-1)}\binom{2k}{k}=-\frac{\cat_{k-1}}{2^{2k-1}}.$$
As in the proof of Proposition \ref{prop:cmplx_caustics_degree_polynomial}, for any permutation $\sigma$ of $\{1,\ldots,m\}$,
$$\prod_{j=1}^{m} B_{\sigma(j)+j}(-a^2)=\prod_{j=1}^{m} \frac{1}{b^{2(\sigma(j)+j)}}P_{\sigma(j)+j}(a,b)=\frac{Q_{\sigma}(a,b)}{b^{2m(m+1)}}$$
where $Q_{\sigma}(a,b)$ is a homogeneous polynomial of degree 
$$\sum_{j=1}^m \deg P_{j+\sigma(j)}=\sum_{j=1}^m 2(\sigma(j)+j) = 2m(m+1)$$
whose coefficient in front of $a^{2m(m+1)}$ is
$$\prod_{j=1}^{m}\left(-\frac{\cat_{j+\sigma(j)-1}}{2^{2(j+\sigma(j))-1}}\right)=\frac{(-1)^m}{2^{m(2m+1)}}\prod_{j=1}^{m}\cat_{j+\sigma(j)-1}.$$
As in proof of Proposition \ref{prop:cmplx_caustics_degree_polynomial}, $\mathcal{B}^n(-a^2)$ is a sum of products of the form $\pm\prod_{j=1}^{m} B_{\sigma(j)+j}(-a^2)$ hence can be written as 
$$\frac{R_n(a,b)}{b^{2m(m+1)}}$$
where $R_n(a,b)$ is the sum of $\varepsilon(\sigma)\prod_{j=1}^{m} Q_{\sigma(j)+j}(a,b)$ and $\varepsilon(\sigma)$ is the parity of $\sigma$. Thus $R_n(a,b)$ is a homogeneous polynomial of degree $2m(m+1)$ whose coefficient in front of $a^{2m(m+1)}$ is 
$$\frac{(-1)^m}{2^{m(2m+1)}}\det H_m = \frac{(-1)^m}{2^{m(2m+1)}}$$ as in the proof of Proposition \ref{prop:cmplx_caustics_degree_polynomial}. Thus $R_n(a,b)$ is a nonzero homogeneous polynomial such that
$$\mathcal{B}^n(-a^2)=\frac{R_n(a,b)}{b^{2m(m+1)}}.$$
We can do the same with $\mathcal{B}^n(-b^2)$ to obtain the same conclusion, which finishes the proof.
\end{proof}

\subsubsection{Case of the circle}

In this section we compute $\deg\mathcal{B}^n_{a,b}$ in the case of the circle ($a=b$), and show that in this case this degree is strictly less than the upper bound given in Proposition \ref{prop:cmplx_caustics_finite_confocal}.

\begin{proposition}
\label{prop:circle_polynomial}
When $a=b$ (in the case of the circle), 
$$\deg\mathcal{B}^n_{a,b}=\left\{\begin{array}{cl}
\frac{n-1}{2}&\text{if } n \text{ is odd}\\
\frac{n}{2}-1&\text{if } n \text{ is even.}\\
\end{array}\right.$$
\end{proposition}

\begin{proof}
Suppose $n=2m+1$ is odd. By Equation \ref{eq:B_k}, when $a=b=R$, for $k\geq2$
$$B_k = \sum_{w=0}^k c_{k-w}\left(1+\frac{\lambda}{a^2}\right)^w\sum_{u+v=w}c_uc_v.$$
Let us compute $\sum_{u+v=w}c_uc_v$: it is the Taylor coefficient at $t^w$ of the function $\sqrt{1+t}^2=1+t$, therefore we get that 
$$\sum_{u+v=w}c_uc_v=\left\{\begin{array}{cl}
1&\text{if } 0\leq w\leq 1\\
0&\text{if } w\geq 2.
\end{array}\right.$$
Hence $B_k = c_k+c_{k-1}x$ where $x = 1+\lambda/a^2$. Using the multilinearity of $\det$, it is not hard to see that $\mathcal{B}^n_{a,b}$ is of degree $m$ if $n$ is odd and $m-1$ if $n$ is even.
\end{proof}

\subsubsection{Explicit formulas of $\mathcal{B}^n_{a,b}$ for $n=3$ to $6$}
\label{subsubsec:cayley_explicit_formulas}

\newcommand{\X}{\textbf{X}}

We give a list of exlicit formulas of $\mathcal{B}^n_{a,b}$ for small $n$. To simplify the formulas we rather express $\tilde{\mathcal{B}}^n_{a,b}=\mu_n\mathcal{B}^n_{a,b}$ where $\mu_n$ is a non-zero real number defined by
$$\mu_n = \left\{
\begin{array}{lll}
(-1)^m2^{m(2m+1)}(ab)^{m(m+1)}&\qquad&\text{if } n=2m+1 \text{ is odd,}\\
\frac{1}{m}(-1)^{m+1}2^{(m-1)(2m+1)}(ab)^{(m-1)(m+1)}&\qquad&\text{if } n=2m \text{ is even.}
\end{array}\right.$$
We further replace its variable $\lambda$ by $\X$ for a better reading.\\

\enum Case $\boxed{n=3}$ 
$$\tilde{\mathcal{B}}^3_{a,b}=(a-b)^2\X^2+2ab(a+b)\X-3a^2b^2$$
\enum Case $\boxed{n=4}$
$$\tilde{\mathcal{B}}^4_{a,b}=(a+b)(a-b)^2\X^3-ab(a-b)^2\X^2-(ab)^2(a+b)\X+(ab)^3$$
\enum Case $\boxed{n=5}$
\begin{center}
\begin{tabular}{lll}
$\tilde{\mathcal{B}}^5_{a,b}$&$=$&$(a-b)^6\X^6 + 2ab(3a + b)(a + 3b)(a + b)(a - b)^2 \X^5$\\
&&$-(ab)^2(29a^2 + 54ab + 29b^2)(a - b)^2\X^4 + 36(ab)^3(a + b)(a - b)^2\X^3$\\
&&$ -(ab)^4(9a^2 - 34ab + 9b^2)\X^2 -10(ab)^5(a + b)\X + 5(ab)^6$
\end{tabular}
\end{center}

This list can be extended using formal calculus on a computer, but this has no interest for the present study. We rather mention that $\mathcal{B}^n_{a,b}$ has generically simple roots for small values of $n$. Moreover, on the examples of this list, the exceptional values of the pair $(a,b)$ for which the degree formula of Proposition \ref{prop:cmplx_caustics_degree_polynomial} is not satisfied are contained in the sets $a=b$ or $a=-b$. Is it always the case for all $n$ ? 

\begin{conjecture}
For all $n\geq 3$, $\mathcal{B}^n_{a,b}$ has generically simple roots. 
\end{conjecture}

Here \textit{generically} has the same meaning as in Proposition \ref{prop:cmplx_caustics_degree_polynomial}. If the conjecture is true, this would imply that the number of complex $n$-caustics is generically given by the degree of $\mathcal{B}^n_{a,b}$ as computed in Proposition \ref{prop:cmplx_caustics_degree_polynomial}.

\subsubsection{Study of complex $3$-caustics}

We study the particular case of complex caustics of $3$-periodic orbits. We will say that a complex conic is an \textit{ellipse} (respectively a \textit{hyperbola}) if its real part is an ellipse (respectively a hyperbola).

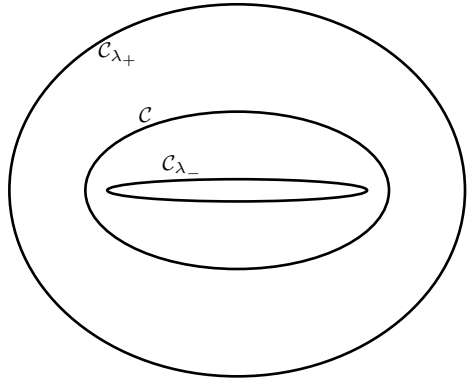
\begin{figure}[!h]
\centering
\begin{tikzpicture}[line cap=round,line join=round,>=triangle 45,x=1cm,y=1cm]
\clip(-3.5,-2.7) rectangle (3.5,2.7);
\draw [rotate around={0:(0.0025086145409116416,0)},line width=1pt] (0.0025086145409116416,0) ellipse (1.9974913854590879cm and 1.0432433360133084cm);
\draw [rotate around={0:(0.0025086145409115726,0)},line width=1pt] (0.0025086145409115726,0) ellipse (2.9974913854590883cm and 2.4664426668897765cm);
\draw [rotate around={0:(0.002508614540915727,0)},line width=1pt] (0.002508614540915727,0) ellipse (1.709768500390961cm and 0.1472859398657943cm);
\begin{scriptsize}
\draw[color=black] (-1.2,1) node {$\mathcal{C}$};
\draw[color=black] (-1.5597148670526204,1.8) node {$\mathcal{C}_{\lambda_+}$};
\draw[color=black] (-0.7111058200015711,0.31037283751140887) node {$\mathcal{C}_{\lambda_-}$};
\end{scriptsize}
\end{tikzpicture}
\caption{When $a=2$ and $b=1$, the conic $\mathcal{C}$ is an ellipse having two complexified ellipses $\mathcal{C}_{\lambda_-}$ and $\mathcal{C}_{\lambda_+}$ as complex caustics of $3$-periodic orbits.}
\label{figure:complex_3_caustics}
\end{figure}

\begin{proposition}
The complex reflection law on the billiard defined by a complexified ellipse or hyperbola $\mathcal{C}$ has exactly two $3$-caustics which are complexified conics of the same type than $\mathcal{C}$, see Figure \ref{figure:complex_3_caustics} and the following table:
\begin{center}
\begin{tabular}{|c||c||c|c|c|c|c|c|}
   \hline
   $\mathcal{C}$ & \textbf{ellipse} & \textbf{hyperbola} \\
   \hline
    Quantity of $3$-caustics & 2&2\\
   \hline
    Types of $3$-caustics & ellipse&hyperbola\\
                          & ellipse&hyperbola\\
   \hline
\end{tabular}
\end{center}
\end{proposition}

The polynomial $\mathcal{B}^3_{a,b}(\lambda)$ is computed at Subsection \ref{subsubsec:cayley_explicit_formulas}. If $a\neq b$, it has two distinct roots, $\lambda_+$ and $\lambda_-$, expressed as
$$\lambda_{\pm}=-\frac{ab}{(a-b)^2}\left(a+b\pm 2\sqrt{a^2 - ab + b^2}\right).$$
These roots are real and satisfy the following inequalities:

\enum \underline{when $\mathcal{C}$ is an \textbf{ellipse}, with $a,b>0$:}
$$a-\lambda_+>a>a-\lambda_->0\qquad\text{and}\qquad b-\lambda_+>b>b-\lambda_->0$$
Hence $\mathcal{C}_{\lambda_+}$ and $\mathcal{C}_{\lambda_-}$ are complexified ellipses, $\mathcal{C}_{\lambda_-}$ is nested in $\mathcal{C}$ which is nested in $\mathcal{C}_{\lambda_+}$.

\enum \underline{when $\mathcal{C}$ is a \textbf{hyperbola}, with $a>0>b$:}
$$a-\lambda_->a>a-\lambda_+>0\qquad\text{and}\qquad 0>b-\lambda_->b>b-\lambda_+$$
Hence $\mathcal{C}_{\lambda_+}$ and $\mathcal{C}_{\lambda_-}$ are complexified hyperbolas, and $\mathcal{C}$ is in the domain delimited by each pair of corresponding branches of $\mathcal{C}_{\lambda_+}$ and $\mathcal{C}_{\lambda_-}$.

\begin{proof}
Since the map $x\mapsto x^2-x+1$ never vanishes on $\RR$, the quantity $a^2-ab+a^2$ is always positive, hence $\lambda_+$ and $\lambda_-$ are real numbers. We can further check that $a+b+2\sqrt{a^2-ab+b^2}\geq 0$ and $a+b-2\sqrt{a^2-ab+b^2}\leq 0$ by comparing the squares of $a+b$ and of $2\sqrt{a^2-ab+b^2}$. Hence $\lambda_{+}$ and $\lambda_-$ have opposite signs. The remaining inequalities are not so difficult to prove.
\end{proof}

\subsubsection{Study of complex $4$-caustics}

We study the particular case of complex caustics of $4$-periodic orbits. We will say that a complex conic is an \textit{ellipse} (respectively a \textit{hyperbola}) if its real part is an ellipse (respectively a hyperbola).

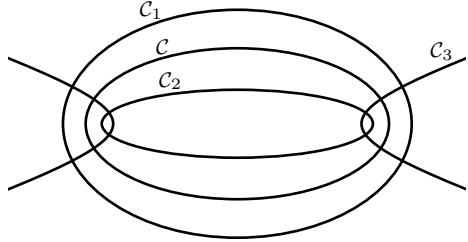
\begin{figure}[!h]
\centering
\begin{tikzpicture}[line cap=round,line join=round,>=triangle 45,x=1cm,y=1cm]
\clip(-3,-1.8) rectangle (3,1.8);
\draw [rotate around={0:(0.004465809016780549,0)},line width=1pt] (0.004465809016780549,0) ellipse (1.9955341909832196cm and 1.0008587655128343cm);
\draw [rotate around={0:(0.004465809016780692,0)},line width=1pt] (0.004465809016780692,0) ellipse (1.784232681737871cm and 0.45060828188388125cm);
\draw [rotate around={0:(0.004465809016780525,0)},line width=1pt] (0.004465809016780525,0) ellipse (2.2943056805364748cm and 1.5110923588129084cm);
\draw [samples=50,domain=-0.99:0.99,rotate around={0:(0.0044658090167803855,0)},xshift=0.0044658090167803855cm,yshift=0cm,line width=1pt] plot ({1.6329479850417183*(1+(\x)^2)/(1-(\x)^2)},{0.5602850319501439*2*(\x)/(1-(\x)^2)});
\draw [samples=50,domain=-0.99:0.99,rotate around={0:(0.0044658090167803855,0)},xshift=0.0044658090167803855cm,yshift=0cm,line width=1pt] plot ({1.6329479850417183*(-1-(\x)^2)/(1-(\x)^2)},{0.5602850319501439*(-2)*(\x)/(1-(\x)^2)});
\begin{scriptsize}
\draw[color=black] (-0.9777216831311774,1.02) node {$\mathcal{C}$};
\draw[color=black] (-0.8726664469978024,0.57) node {$\mathcal{C}_2$};
\draw[color=black] (-1.1272233653209804,1.5) node {$\mathcal{C}_1$};
\draw[color=black] (2.7,0.9599478885188634) node {$\mathcal{C}_3$};
\end{scriptsize}
\end{tikzpicture}
\caption{When $a=3$ and $b=1$, the conic $\mathcal{C}_0$ is an ellipse having two complexified ellipses $\mathcal{C}_{1}$, $\mathcal{C}_{2}$ and a complexified hyperbola $\mathcal{C}_{3}$ as complex caustics of $4$-periodic orbits.}
\label{fig_3caustics}
\end{figure}

\begin{proposition}
The $4$-caustics of the complex reflection law on the billiard defined by a complexified ellipse or hyperbola $\mathcal{C}$ are detailed in the following table:
\begin{center}
\begin{tabular}{|c||c|c|c||c|c|c|c|}
   \hline
   \tiny & \multicolumn{3}{c||}{\tiny} & \multicolumn{2}{c|}{\tiny} \\
   $\mathcal{C}$ & \multicolumn{3}{c||}{\textbf{ellipse}} & \multicolumn{2}{c|}{\textbf{hyperbola}} \\
   
   & \multicolumn{3}{c||}{\scriptsize$0<b<a$} & \multicolumn{2}{c|}{\scriptsize$b<0<a$} \\
   \hline
    $a,b$ & $a>2b$&$a=2b$&$a<2b$&$a>|b|$&$a<|b|$\\
   \hline
    Quantity of $4$-caustics & 3&2&3&3&3\\
   \hline
                          & ellipse&ellipse&ellipse      &hyperbola&hyperbola\\
    Types of $4$-caustics & ellipse&ellipse&ellipse      &hyperbola&hyperbola\\
                          & hyperbola&     &str. complex conic&ellipse&str. complex conic\\
   \hline
\end{tabular}
\end{center}
\end{proposition}

The polynomial $\mathcal{B}^4_{a,b}(\lambda)$ is computed in Subsection \ref{subsubsec:cayley_explicit_formulas}. If $a\neq b$ and $a\neq -b$, it has three distinct roots, expressed as
$$\lambda_1 = \frac{ab}{b-a}\qquad \lambda_2 = \frac{ab}{a+b}\qquad\lambda_3 = \frac{ab}{a-b}.$$
These roots satisfy the following inequalities:

\enum \underline{when $\mathcal{C}$ is an \textbf{ellipse}, with $a>b>0$:}
$$a-\lambda_1>a>a-\lambda_2>0\qquad\text{and}\qquad b-\lambda_1>b>b-\lambda_2>0$$
$$b-\lambda_3<0\qquad\text{and}\qquad a-\lambda_3\left\{\begin{array}{ccc}
>0&&\text{if }a>2b\\
=0&&\text{if }a=2b\\
<0&&\text{if }a<2b\\
\end{array}\right.$$
Hence $\mathcal{C}_{\lambda_1}$ and $\mathcal{C}_{\lambda_2}$ are always complexified ellipses, $\mathcal{C}_{\lambda_2}$ is nested in $\mathcal{C}$ which is nested in $\mathcal{C}_{\lambda_1}$. The conic $\mathcal{C}_{\lambda_3}$ is an hyperbola if $a>2b$, not defined if $a=2b$ and a complex conic if $a<2b$.

\enum \underline{when $\mathcal{C}$ is a \textbf{hyperbola}, with $a>0>b$:}\\
\begin{center}
\underline{Case $|b|<a$}
\end{center} 
$$0<a-\lambda_1<a<a-\lambda_3<a-\lambda_2\qquad\text{and}\qquad b-\lambda_1<b<b-\lambda_3<0<b-\lambda_2$$
Hence $\mathcal{C}_{\lambda_1}$ is a hyperbola, $\mathcal{C}_{\lambda_2}$ an ellipse, $\mathcal{C}_{\lambda_3}$ a hyperbola. The branches of $\mathcal{C}$ are in the domain bounded by the corresponding branches of $\mathcal{C}_{\lambda_1}$ and $\mathcal{C}_{\lambda_3}$.
\begin{center}
\underline{Case $|b|>a$}
\end{center} 
$$a-\lambda_2<0<a-\lambda_1<a<a-\lambda_3\qquad\text{and}\qquad b-\lambda_2<b-\lambda_1<b<b-\lambda_3<0$$
Hence $\mathcal{C}_{\lambda_1}$ is a hyperbola, $\mathcal{C}_{\lambda_2}$ a complex conic, $\mathcal{C}_{\lambda_3}$ a hyperbola. The branches of $\mathcal{C}$ are in the domain bounded by the corresponding branches of $\mathcal{C}_{\lambda_1}$ and $\mathcal{C}_{\lambda_3}$.

\subsection{Complex Joachimsthal invariant}
	
The first proof of Proposition \ref{prop:confocal_conics_are_caustics} which we were able to obtain was different from the one we give in Subsection \ref{subsection_conics_are_cautics}. It used a complex version of the so-called \textit{Joachimsthal invariant}, a well-known quantity in the theory of billiards on conics. This complex invariant is described in the present subsection.

The context of this subsection is as follows. We consider the conic $\mathcal{C}$  of $\RR^2$ given by equation \eqref{eq:reference_conic}:
$$\mathcal{C}:\frac{x^2}{a}+\frac{y^2}{b}=1$$
where $a,b\neq 0$ and $x,y\in\RR$.  Its complexification is given by the same equation with $x,y\in\CC$. In \cite{taba_book} Chapter $4$, Theorem $4.4$ shows that for a set of points and directions defined as successive billiard reflections on the real ellipse $\mathcal{C}$ with $a,b>0$, there is an \textit{invariant} quantity. Known as Joachimsthal invariant, it is defined by
$$\frac{xv_x}{a}+\frac{yv_y}{b}$$
where $(x,y)$ are the coordinates of a vertex of an orbit, and $v$ a unitary vector having this vertex as starting point and pointing toward the next vertex. Let us mention another reference about Joachimsthal invariant, which was given to us by the referee of our original article: see \cite{treshchev}, Chapter IV. 

In our case, we consider the complexified version of $\mathcal{C}$, and Joachimsthal invariant has to be modified to handle the complex structure. Hence from now on, we choose $a,b\neq 0$ and $\mathcal{C}$ denote the complexification of the previous defined conic. As described in Section \ref{section_complex_reflection_law}, $\CC^2$ can be endowed with the non-degenerate complex quadratic form $q(x,y)=x^2+y^2$ which vanishes on vectors of the space $\CC(1,i)\cup\CC(1,-i)$, called \textit{isotropic} vectors. We recall that we can consider complex orbits $(p_1,\ldots,p_k)$ on $\mathcal{C}$ viewed as a conic of $\CP^2$ via an embedding $\CC^2\subset\CP^2$. When an orbit has an edge $p_jp_{j+1}$ which is directed by an isotropic vector, we say that the orbit is \textit{isotropic}, and \textit{non-isotropic} otherwise. 

In the case when a point $p$ belongs to the so-called line at infinity $L_{\infty}:=\CP^2\smallsetminus\CC^2$, we say that $p$ is \textit{infinite}, and \textit{finite} otherwise. Now if $p=(x,y)\in\CC^2$ and $v=(v_x,v_y)\in\CC^2$ is not isotropic, then we define the complex Joachimsthal invariant at $(p,v)$ as the complex quantity
$$P(p,v)=\frac{1}{q(v)}\left(\frac{xv_x}{a}+\frac{yv_y}{b}\right)^2.$$

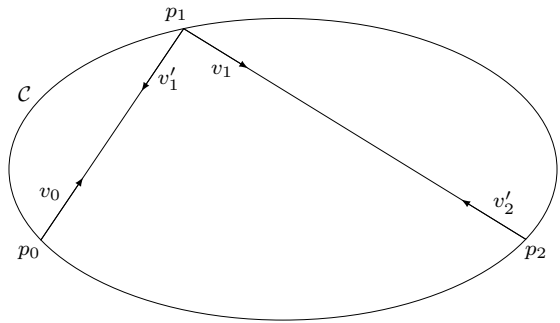
\begin{figure}[!h]
\centering

\begin{tikzpicture}[line cap=round,line join=round,>=triangle 45,x=1cm,y=1cm]
\draw [rotate around={0:(0,0)},line width=0.2pt] (0,0) ellipse (3.6055512754639882cm and 2cm);
\draw [-latex,line width=0.2pt] (-3.1850116065471195,-0.9373778593108506) -- (-2.6280785187343483,-0.10682049482641465);
\draw [-latex,line width=0.2pt] (-1.3064833284084274,1.8640816516113508) -- (-1.8634164162211988,1.033524287126915);
\draw [-latex,line width=0.2pt] (-1.3064833284084274,1.8640816516113508) -- (-0.4567448009853311,1.3368772936885167);
\draw [-latex,line width=0.2pt] (3.1938575709826025,-0.9280705901993418) -- (2.3441190435595054,-0.40086623227650797);
\draw [line width=0.2pt] (-3.1850116065471195,-0.9373778593108506)-- (-1.3064833284084274,1.8640816516113508);
\draw [line width=0.2pt] (-1.3064833284084274,1.8640816516113508)-- (3.1938575709826025,-0.9280705901993418);
\begin{scriptsize}
\draw[color=black] (-3.4,1) node {$\mathcal{C}$};
\draw[color=black] (-3.338207671594807,-1.1) node {$p_0$};
\draw[color=black] (-1.3979696107556636,2.05) node {$p_1$};
\draw[color=black] (3.3336486791209254,-1.1) node {$p_2$};
\draw[color=black] (-3.063645681853419,-0.33217314749790483) node {$v_0$};
\draw[color=black] (-1.5,1.2) node {$v_1'$};
\draw[color=black] (-0.8,1.3) node {$v_1$};
\draw[color=black] (2.9126536281841306,-0.44199794339445997) node {$v_2'$};
\end{scriptsize}
\end{tikzpicture}
\caption{In Proposition \ref{proposition:orbits_and_invariant}, we consider all quantities $P(p_0,v_0)$, $P(p_1,v_1')$, $P(p_1,v_1)$ and $P(p_2,v_2')$.}
\label{fig_proposition:orbits_and_invariant}
\end{figure}

In what follows we show that if $T=(p_0,\ldots,p_k)$ is a non-degenerate and non-isotropic orbit on $\mathcal{C}$, the quantity $P(p_j,v)$ do not depend on the choice of a finite vertex $p_j$ or of a directing vector of $p_{j-1}p_j$ (Proposition \ref{proposition:orbits_and_invariant} and Figure \ref{fig_proposition:orbits_and_invariant}). Moreover, let $\mathcal{C}_{\lambda}$ be the caustic to which $T$ remains tangent, as shown in Proposition \ref{prop:confocal_conics_are_caustics}. If we denote by $P(T)$ previous invariant quantity associated with $T$  then the caustic $\mathcal{C}_{\lambda}$ satisfies $\lambda=abP(T)$, as shown in Proposition \ref{prop_ellipse_inscrite}.
Conversely, we show that if the quantity $P(p_j,v)$ is preserved on a polygon inscribed in $\mathcal{C}$, then the latter is an orbit (except for degenerate cases of polygons), see Lemma \ref{lemma:invariant_implies_orbit} and Lemma \ref{lemma:invariant_implies_orbit_2}.

\begin{proposition}
\label{proposition:orbits_and_invariant}
Let $T= (p_0,p_1,p_2)$ be a non-degenerate and non-isotropic orbit on $\mathcal{C}$ with $p_0$ finite. Then the quantity $P(p_j,v)$
do not depend on the choice of a \textbf{finite} vertex $p_j$ of $T$ or of a directing vector of $p_{j-1}p_j$ or $p_jp_{j+1}$ (see Fig. \ref{fig_proposition:orbits_and_invariant}). 
\end{proposition}

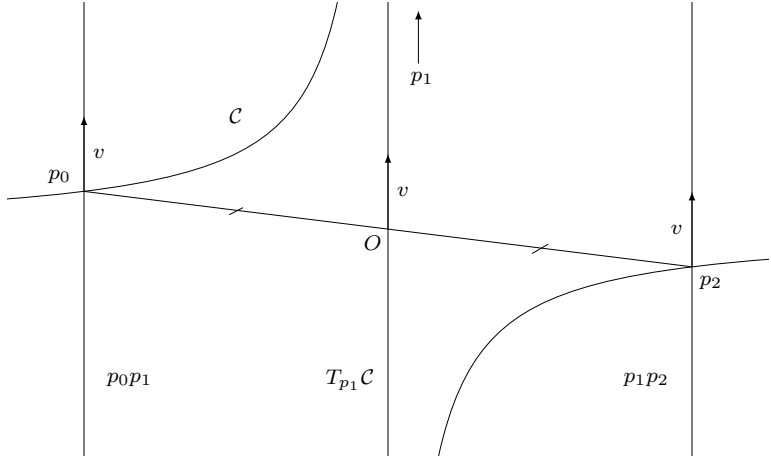
\begin{figure}[!h]
\centering
\begin{tikzpicture}[line cap=round,line join=round,>=triangle 45,x=1cm,y=1cm]
\clip(-5,-3) rectangle (5,3);
\draw [samples=50,domain=-0.99:0.99,rotate around={-45:(0,0)},xshift=0cm,yshift=0cm,line width=0.2pt] plot ({2*(1+(\x)^2)/(1-(\x)^2)},{2*2*(\x)/(1-(\x)^2)});
\draw [samples=50,domain=-0.99:0.99,rotate around={-45:(0,0)},xshift=0cm,yshift=0cm,line width=0.2pt] plot ({2*(-1-(\x)^2)/(1-(\x)^2)},{2*(-2)*(\x)/(1-(\x)^2)});
\draw [line width=0.2pt] (-4,-5.750207336486205) -- (-4,5.688558131117167);
\draw [line width=0.2pt] (0,-5.750207336486205) -- (0,5.688558131117167);
\draw [line width=0.2pt] (4,-5.750207336486205) -- (4,5.688558131117167);
\draw [line width=0.2pt] (-4,0.5)-- (4,-0.5);
\draw [line width=0.2pt] (-2.085893357114387,0.19163967179694663)-- (-1.9119521473062633,0.282391607349011);
\draw [line width=0.2pt] (1.8996291458804462,-0.322621296331418)-- (2.1038210008725913,-0.20161871559533223);
\draw [-latex,line width=0.2pt] (0.4,2.2) -- (0.4,2.9);
\draw [-latex,line width=0.2pt] (0,0) -- (0,1);
\draw [-latex,line width=0.2pt] (-4,0.5) -- (-4,1.5);
\draw [-latex,line width=0.2pt] (4,-0.5) -- (4,0.5);
\begin{scriptsize}
\draw[color=black] (-4.32034094759505,0.70) node {$p_0$};
\draw[color=black] (0.45,2) node {$p_1$};
\draw[color=black] (4.25,-0.7126491647002291) node {$p_2$};
\draw[color=black] (-0.2,-0.17438473151438424) node {$O$};
\draw[color=black] (-2,1.5) node {$\mathcal{C}$};
\draw[color=black] (-3.8,1) node {$v$};
\draw[color=black] (3.8,0) node {$v$};
\draw[color=black] (0.2,0.5) node {$v$};
\draw[color=black] (-3.4,-2) node {$p_0p_1$};
\draw[color=black] (3.4,-2) node {$p_1p_2$};
\draw[color=black] (-0.5,-2) node {$T_{p_1}\mathcal{C}$};
\end{scriptsize}
\end{tikzpicture}
\caption{An orbit $(p_0,p_1,p_2)$ on $\mathcal{C}$ with $p_1$ infinite as in the proof of Proposition \ref{proposition:orbits_and_invariant}. The points $p_0$ and $p_2$ are symmetric across $O$, hence $p_2=-p_0$ and $P(p_0,v)=P(p_2,v)$. Here $\mathcal{C}$ is represented as an hyperbola which allows us to view the tangent line at the infinity point $p_1$ as the vertical asymptote.}
\label{fig_proof_proposition:orbits_and_invariant}
\end{figure}

\begin{proof}
As explained in \cite{glut2}, the reflection with respect to a non-isotropic line permutes the isotropic directions $v_I=(1,i)$ and $v_J=(1,-i)$. Hence in our case, $q(v)\neq0$ for all $v$ taken like in the proposition we want to prove.

\textbf{First case:} If $p_0$ and $p_1$ are finite, write $p_0 = (x_0,y_0)$, $p_1=(x_1,y_1)$. Take $v_0$ a vector such that $q(v_0)=1$ and directing $p_0p_1$, and $v_1$ vector such that $q(v_1)=1$ and directing $p_1p_2$. Define the matrix
$$A=\left(\begin{matrix}1/a&0\\0&1/b \end{matrix}\right).$$ 
Then since $\transp p_j Ap_j=1$ and since $A$ is symmetric, we get
$$\transp{(p_1-p_0)}A(p_1+p_0) = \transp p_1Ap_0-\transp p_0Ap_1=0.$$
Since $v_0$ is collinear to $p_1-p_0$ we have further $\transp v_0A(p_1+p_0)=0$, thus
\begin{equation}\label{equ_first_invariant}
\transp v_0Ap_1=-\transp v_0Ap_0.
\end{equation}
But since $p_0p_1$ and $p_1p_2$ are symmetric with respect to the tangent line of $\mathcal{C}$ at $p_1$, which is also the orthogonal line to $Ap_1$ (the gradient in $p_1$ of the bilinear form defining $\mathcal{C}$), we only have two possibilities : either $v_0+v_1$ or $v_0-v_1$ is orthogonal to $Ap_1$ as we see by decomposing both $v_0$ and $v_1$ in normal and tangential components. Hence
$$\transp{(v_0+v_1)}Ap_1=0\qquad\text{or}\qquad\transp{(v_0-v_1)}Ap_1=0.$$
In both cases we get
$$\left(\transp v_0Ap_1\right)^2=\left(\transp v_1Ap_1\right)^2$$
and using equality \eqref{equ_first_invariant}, we get
\begin{equation}\label{equ_second_invariant}
\left(\transp v_0Ap_0\right)^2=\left(\transp v_1Ap_1\right)^2
\end{equation}
which proves Proposition \ref{proposition:orbits_and_invariant} for unitary vectors. For general vectors, it is enough to divide them by a square root of $q(v)$, which explains the factor $1/q(v)$ appearing in the formula of $P(p,v)$.

\textbf{Second case:} If $p_0$ is finite and $p_1$ infinite (see Fig. \ref{fig_proof_proposition:orbits_and_invariant}), then $p_2$ is finite. Indeed, $p_0p_1$ is not the line at infinity and $T_{p_1}\mathcal{C}$ is not isotropic. Hence the line symmetric to $p_0p_1$ with respect to $T_{p_1}\mathcal{C}$ is finite and parallel to $p_0p_1$ and to $T_{p_1}\mathcal{C}$ (the three lines intersects at the same infinite point). Thus the other point of intersection $p_2$ of the latter symmetric line with $\mathcal{C}$ has to be finite. If we consider $v$ a non-zero vector directing the lines $p_0p_1$, $p_1p_2$ and $T_{p_1}\mathcal{C}$, we need to prove that 
$$P(p_0,v) = P(p_2,v).$$
But $p_2 =-p_0$ since $T_{p_1}\mathcal{C}$ goes through the origin $O=(0,0)$ (by property of a tangent line at an infinite point of $\mathcal{C}$) and the ellipse $\mathcal{C}$ is symmetric across $O$ (see Fig. \ref{fig_proof_proposition:orbits_and_invariant}). This implies that $P(p_0,v) = P(p_2,v)$.
\end{proof}

\begin{corollary}
Let $T= (p_0,\ldots,p_n)$ be a non-degenerate and non-isotropic orbit on $\mathcal{C}$. Then the quantity $P(p_j,v)$ defined as before do not depend on the choice of a finite vertex $p_j$ or on $v$, a directing vector of $p_{j-1}p_j$ or $p_jp_{j+1}$.
Thus we can write $P(p_j,v) = P(T)$.
\end{corollary}

Here we prove that the invariant $P(T)$ is linked with the caustic of the orbit $T$. We first recall a result based on duality of conics, which can be deduced from Subsection \ref{subsection:polarity}. We consider coordinates on $\CP^2$ such that any point of $\CP^2$ can be denoted by $(x:y:z)$, where $(x,y,z)\in\CC^3\smallsetminus\{0\}$.

\begin{lemma}\label{lemma_conique_duale}
let $C$ be a conic in $\CP^2$ given by the equation $\transp p A p=0$ where $A$ is a $3\times3$ symmetric invertible matrix, $p=(x:y:z)$, and $v=(\alpha,\beta,\gamma)\in\CC^3$ defining the line $\ell_v$ of equation $\alpha x+\beta y+\gamma z=0$. Then $\ell_v$ is tangent to $C$ if and only if $\transp v \inverse{A} v = 0.$
\end{lemma}

\begin{figure}[!h]
\centering
\begin{tikzpicture}[line cap=round,line join=round,>=triangle 45,x=1cm,y=1cm]
\clip(-5,-2.1) rectangle (5,2.1);
\draw [rotate around={0:(0,0)},line width=0.2pt] (0,0) ellipse (2.8284271247461907cm and 2cm);
\draw [rotate around={0:(0,0)},line width=0.2pt] (0,0) ellipse (2.57838050536101cm and 1.627281791954208cm);
\draw [line width=0.2pt] (-2.822840787497323,-0.12563814795175854)-- (-1.869312067799702,1.500945101124356);
\draw [line width=0.2pt] (-1.869312067799702,1.500945101124356)-- (1.393051102240205,1.740604582687773);
\draw [line width=0.2pt] (1.393051102240205,1.740604582687773)-- (2.818498874725576,0.1674277354124261);
\draw [line width=0.2pt] (2.818498874725576,0.1674277354124261)-- (1.928469151211218,-1.463045920814977);
\begin{scriptsize}
\draw [fill=black] (-2,0) circle (1pt);
\draw[color=black] (-1.9906447176999233,0.4527745617943184) node {$F_1$};
\draw [fill=black] (2,0) circle (1pt);
\draw[color=black] (2.030966426407765,0.4773715718500229) node {$F_2$};
\draw[color=black] (-2.1013312629505934,1.7) node {$p_1$};
\draw[color=black] (1.4898322051822654,1.9) node {$p_2$};
\draw[color=black] (3.2,0.2) node {$p_3$};
\draw[color=black] (-3.2,-0.1) node {$p_0$};
\draw[color=black] (2.030966426407765,-1.7) node {$M_4$};
\draw[color=black] (-1,-1.3) node {$\class_{\lambda}$};
\draw[color=black] (0,-1.85) node {$\mathcal{C}$};
\end{scriptsize}
\end{tikzpicture}
\caption{The confocal caustic $\class_{\lambda}$ inscribed in a piece of billiard trajectory.}
\label{fig_thm_caustic}
\end{figure}
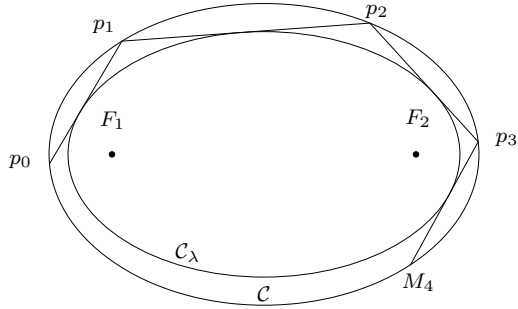

\begin{proposition}
\label{prop_ellipse_inscrite}
Let $T$ be a non-degenerate non-isotropic orbit of $\mathcal{C}$ tangent to a complex conic $\class_{\lambda}$, with $\lambda\in\CC$ different from $a$ and $b$, given by Equation \ref{eq:confocal_conics}. Then $\lambda = abP(T)$.
\end{proposition}

\begin{proof} Consider a set of coordinates of $\CP^2$ such that the conic $\mathcal{C}_{\lambda}$ is given by the equation $\transp p\inverse{B_{\lambda}} p=0$, where 
$$B_{\lambda}=\left(\begin{matrix}a-\lambda&0&0\\0&b-\lambda&0\\0&0&-1 \end{matrix}\right).$$ 
Write $T=(p_0,\ldots,p_n)$. Since the orbit is non-isotropic, two consecutive sides $p_{j-1}p_j$ and $p_jp_{j+1}$ cannot be the line at infinity. Hence we suppose without loss of generality that $p_0$ is finite. Then the line $p_0p_1$ is defined in $\CP^2$ by the equation $v_yx-v_xy+(v_xy_0-v_yx_0)z=0$, where $p_0=(x_0,y_0)$ and $v=(v_x,v_y)$ is a directing vector of $p_0p_1$ in $\CC^2$. Hence we have $p_0p_1 = \ell_w$ (in the notations of Lemma \ref{lemma_conique_duale}) where $w = (v_y,-v_x,v_xy_0-v_yx_0)$. It allows us to compute
$$\transp w B_{\lambda}w = (a-\lambda)v_y^2+(b-\lambda)v_x^2-\left(v_xy_0-v_yx_0\right)^2$$
which can be rearranged as $\transp w B_{\lambda}w=-\lambda q(v)+(a-x_0^2)v_y^2+(b-y_0^2)v_x^2+2v_xv_yx_0y_0.$
Using the fact that $p_0$ lies on $\mathcal{C}$ gives 
$$a-x_0^2=\frac{a}{b}y_0^2, \qquad b-y_0^2=\frac{b}{a}x_0^2$$
which implies that
$$\transp w B_{\lambda}w = -\lambda q(v)+ab\left(\frac{x_0v_x}{a^2}+\frac{y_0v_y}{b^2}\right)^2=-q(v)\left(\lambda-abP(M_0,v)\right).$$
Since $p_0p_1$ is tangent to $\mathcal{C}_{\lambda}$, $\transp w B_{\lambda}w=0$ and we get the result.
\end{proof}

In what follows, we show that the invariant property implies a billiard reflection property.

\begin{lemma}
\label{lemma:invariant_implies_orbit}
Let $p$ be a finite point on $\mathcal{C}$ such that the line $T_{p}\mathcal{C}$ is non-isotropic. Let $\ell_1$, $\ell_2$ two lines containing $p$ and directed by non-isitropic vectors $v_1$, $v_2$. If
\begin{equation}\label{equ_invariant_points}
P(p,v_1) = P(p,v_2)
\end{equation}
then one of the following cases holds:\\
1) $\ell_1=\ell_2$;\\
2) $\ell_1$ and $\ell_2$ are symmetric with respect to $T_{p}\mathcal{C}$.
\end{lemma}

\begin{proof}
We can suppose $q(v_1)=q(v_2)=1$. By Equality (\ref{equ_invariant_points}), we have $\transp v_1Ap = \pm \transp v_2Ap$ hence $\transp{(v_2\pm v_1)}Ap = 0.$ Thus we get that $v_1+v_2$ or $v_2-v_1$ is orthogonal to $Ap$ which is orthogonal to the tangent line of $\mathcal{C}$ at $p$. Hence $v_1+v_2$ or $v_1-v_2$ is tangent to $\mathcal{C}$ at $p$. This implies that one of these vectors is fixed by the complex reflection with respect to $T_p\mathcal{C}$. 

This means that the components of the $v_j$'s along the direction of $\orth{T_{p}\mathcal{C}}$ are the same or have opposite signs. Since the $v_j$'s are unit vectors, their components along the direction of $T_{p}\mathcal{C}$ are also the same or have opposite signs. Hence we have only three possibilities: a)$\,v_1$ and $v_2$ are symmetric with respect to $T_{M}\mathcal{C}$, b)$\,v_1$ and $v_2$ are symmetric with respect to $\orth{T_{M}\mathcal{C}}$, c)$\,v_2=\pm v_1$. All these cases imply the result.
\end{proof}

\begin{lemma}
\label{lemma:invariant_implies_orbit_2}
Suppose that $\mathcal{C}$ is not a circle (\textit{ie} $a\neq b$). Let $p_0,p_1,p_2$ be points on $\mathcal{C}$ such that $p_0,p_2$ are finite and $p_1$ infinite. Let $v_j$ be a vector directing the line $p_1p_j$, $j=0,2$. If
\begin{equation}\label{equ_invariant_points_2}
P(p_0,v_0) = P(p_2,v_2)
\end{equation}
then one of the following cases holds:\\
1) $p_0=p_2$;\\
2) $p_0p_1$ and $p_1p_2$ are symmetric with respect to $T_{p_1}\mathcal{C}$.
\end{lemma}

\begin{proof}
Since the three lines $p_0p_1$, $p_1p_2$ and $T_{p_1}\mathcal{C}$ contain the same infinite point $p_1$, they are parallel, and therefore directed by the same vector $v=v_0=v_2$. The vector $v$ cannot be isotropic, since an ellipse having an isotropic tangent line at a infinite point is a circle (it is recalled at Section \ref{section_circum_centers} or it can be shown independantly by computations). Suppose that $q(v)=1$. As before, Equation \eqref{equ_invariant_points_2} implies that $p_2-p_1$ is either colinear or orthogonal to $v$. Both cases gives the result.
\end{proof}

 \section{Caustics of quadrics endowed with a structure of projective billiard}
	\label{section_projective_caustics}
	In this section, we use an idea found in \cite{CKS} Sec. III: it appeared to us as a construction of a field of projective lines on a quadric using another quadric, also projective billiards are not mentioned in the corresponding paper. We first describe how to construct such field, and then we study its properties related to caustics. All the results taken separately are well-known, the only interest of our work is to gather them and to interpret them as results on projective billiards.

Let $Q_1,Q_2$ be two non-empty and non-degenerate distinct quadrics. Consider a point $p\in Q_1$ such that $Q_2$ is not tangent to $Q_1$ at $p$. Let $u$ be the pole of $T_pQ_1$ with respect to $Q_2$ (see Definition \ref{definition:polarity}). Since $Q_2$ is not tangent to $Q_1$ at $p$, we have $u\neq p$ and we can define the line
$$L_{Q_2}(p)=pu.$$

\begin{lemma}
\label{lemma:projquadrics_transverseline}
The line $L_{Q_2}(p)$ is tangent to $Q_1$ at $p$ if and only if $T_p Q_1$ is tangent to $Q_2$.
\end{lemma}

\begin{proof}
$L_{Q_2}(p)$ is tangent to $Q_1$ at $p$ if and only if $u\in T_p Q_1$ and the proof follows from the definition of polar spaces as projections of specific orthogonal spaces, see Section \ref{section:general_properties_on_quadrics}.
\end{proof}

The set of points $p\in Q_1$ such that $L_{Q_2}(p)$ is defined corresponds to the set of such $p$ for which $Q_2$ is not tangent to $Q_1$ at $p$. Since $Q_1\neq Q_2$, it is a dense open subset of $Q_1$. By Lemma \ref{lemma:projquadrics_transverseline}, the set $U$ of points $p\in U$ such that $L_{Q_2}(p)$ is transverse to $Q_1$ is also open and dense in $Q_1$ (the complementary set of a strict algebraic subset).

\begin{definition}
We denote by $\linep{Q_1}{Q_2}$ the line-framed hypersurface over $U$ defined by 
$$\linep{Q_1}{Q_2}=\ensemble{(p,L_{Q_2}(p))}{p\in U}.$$
If $q_2$ is a quadratic form defining $Q_2$, we can also write $\linep{Q_1}{q_2}$.
\end{definition}

\begin{proposition}
\label{proposition:quadrics_as_projective_caustics}
The quadric $Q_2$ is a caustic of the line-framed hypersurface $\linep{Q_1}{Q_2}$ over $Q_1$.
\end{proposition}

\begin{proof}
Let $p\in Q_1\smallsetminus Q_2$ such that $T_pQ_1$ is not tangent to $Q_2$. Let $u$ be the pole of $T_pQ_1$ with respect to $Q_2$. Note that the line $L_{Q_2}(p)=pu$ is not tangent to $Q_2$, since otherwise we would have $pu\subset T_pQ_1$ by a polarity argument and contradicting Lemma \ref{lemma:projquadrics_transverseline}.

Let $\ell$ be a line tangent to $Q_2$ and intersecting $Q_1$ at $p$ transversally. Since $L_{Q_2}(p)$ is not tangent to $Q_2$, the lines $\ell$ and $L_{Q_2}(p)$ are distinct: one can consider the unique $2$-dimensional plane $P$ containing both lines $\ell$ and $L_{Q_2}(p)$. The plane $P$ intersect $Q_2$ transversally: otherwise if $P$ is tangent to $Q_2$ one get that $\ell\subset T_pQ_1$ by a polarity argument, which contradicts the transversality of $\ell$ with $Q_1$. Therefore the intersection $C:=Q_2\cap P$ is a non-degenerate non-empty conic of $P$.

Let $\ell'\subset P$ be the other tangent line to $C$ in $P$ containing $p$ (since such a tangent line $\ell$ already exists, there are exactly two distinct such tangent lines). To conclude the proof, we show that the lines $\ell$, $\ell'$, $pu$ and $T_pQ_1\cap P$ form a harmonic set of lines. Denote by $T$ the line $T_pQ_1\cap P$ and by $z,z'$ the respective tangency points of $\ell$ and $\ell'$. Now consider the polarity in the plane $P$ with respect to the conic $C$: the polar line $p^{\ast}$ to $p$ contains $z$ and $z'$, hence $p^{\ast}=zz'$. Then since $u$ is the pole in $P$ of $T=T_pQ_1\cap P$, we have that $u\in zz'$. Now consider the map $s$ from $p^{\ast}$ to $p^{\ast}$ such that the image $s(x)$ of a point $x\in p^{\ast}$ is the pole of the line $px$. By construction, $s$ fixes the points $z$ and $z'$, and permutes $u$ with the intersection point of $T$ with $p^{\ast}$. Therefore these four points are harmonic, and so are the lines $\ell$, $\ell'$, $pu$ and $T$ which concludes the proof.
\end{proof}

Let $q_1$, $q_2$ be non-degenerate quadratic forms defining $Q_1$ and $Q_2$. We can consider their respective dual quadratic forms $\polar{q_1}$ and $\polar{q_2}$, together with the non-degenerate pencil of quadrics $\mathcal{F}(\polar{q_1},\polar{q_2})$. Notice that $q_1$ and $q_2$ belongs to the corresponding dual pencil of quadric $\polar{\mathcal{F}(\polar{q_1},\polar{q_2})}$ which we called \textit{pencil of $(q_1,q_2)$-confocal quadrics} (see Section \ref{section:general_properties_on_quadrics}). We prove in fact that all quadratic forms of $\polar{\mathcal{F}(\polar{q_1},\polar{q_2})}$ define the same field of projective lines over $Q_1$:

\begin{proposition}
\label{proposition:confocal_projective_caustics}
Let $q_3\in\polar{\mathcal{F}(\polar{q_1},\polar{q_2})}$ not colinear to $q_1$. Then $\linep{Q_1}{q_3}=\linep{Q_1}{q_2}$.
\end{proposition}

\begin{corollary}
Any quadric $Q\neq Q_1$ defined by a quadratic form $q\in\polar{\mathcal{F}(\polar{q_1},\polar{q_2})}$ is a caustic of $\linep{Q_1}{q_2}$.
\end{corollary}

\begin{proof}[Proof of Proposition \ref{proposition:confocal_projective_caustics}]
Fix $p\in Q_1$ and $q_3\in\polar{\mathcal{F}(\polar{q_1},\polar{q_2})}$. We want to show that the pole $u_3$ of $T_pQ_1$ with respect to $q_3$ and the pole $u_2$ of $T_pQ_1$ with respect to $q_2$ are on the same line.

By assumption, one can find $\lambda,\mu\in\RR$ such that 
\begin{equation}
\label{equation:pencil_quadrics_1}
\polar{q_3}=\lambda \polar{q_1}+\mu\polar{q_2}
\end{equation}
and $\mu\neq 0$. Denote by $[x]$ the equivalence class in $\RP^d$ of an element $x\in\RR^{d+1}$. For $j=1,2,3$, choose an $x_j\in\RR^{d+1}\smallsetminus\{0\}$ such that the dual of $T_pQ_1$ with respect to $q_j$ is $[x_j]$: we have by construction $p=[x_1]$, $u_2=[x_2]$ and $u_3=[x_3]$. Further denote by $M_j$ the $(d+1)\times(d+1)$ invertible matrix such that for all $x,y\in\RR^{d+1}$ we can write $q_j(x,y)=\mathcal{Q}_0(M_jx,y)$, where $\mathcal{Q}_0$ is the quadratic form $\sum_k x_k^2$. Equation \ref{equation:pencil_quadrics_1} can be rewritten as $\inverse{M_3}=\lambda \inverse{M_1}+\mu\inverse{M_2}$. 

Let $V\subset\RR^{d+1}$ be the hyperplane such that $T_pQ_1=\mathbb{P}(V)$. For $j=1,2,3$, the vector $x_j$ is $q_j$-orthogonal to $V$, hence $M_jx_j$ is $\mathcal{Q}_0$-orthogonal to $V$. In particular, we can find non-zero $\nu_2,\nu_3\in\RR$ such that $M_1x_1 = \nu_2 M_2x_2$ and $M_3x_3 = \nu_3 M_1x_1$. Hence 
\begin{equation}
\label{equation:projective_line}
x_3 = \nu_3 \inverse{M_3}M_1x_1=\nu_3\left(\lambda x_1+\mu\inverse{M_2}M_1x_1\right)=\nu_3\left(\lambda x_1+\mu\nu_2x_2\right)
\end{equation}
It follows from this equation that $u_3$ is on the line containing $p=[x_1]$ and $u_2=[x_2]$.
\end{proof}

Let $q_d$ be a degenerate quadratic form over $\RR^{d+1}$ of rank $d$. The kernel $\ker q_d$ of $q_d$ has dimension $1$ and is generated by a non-zero vector $x_d$. Thus given a hyperplane $V_0$ of $\RR^{d+1}$ transverse to $\ker q_d$, the restriction of the form $q_d$ to $V_0$ is non-degenerate. Consider the affine subspace $V=x_d+V_0\subset\RR^{d+1}$. Its tautological projection $\PP{}{V}$ is an affine chart identified with $V_0\simeq\RR^d$ by $x\in V_0\mapsto [x_d+x]\in\mathbb{P}(V)$, where $[y]$ denotes the equivalence class of $y$ in $\RP^{d}$. Hence we deal with $q_d$ as a non-degenerate quadratic form on $\mathbb{P}(V)\simeq V_0$. We can define its dual $\polar{q_d}$ with respect to the restriction of $\mathcal{Q}_0$ to $V_0$ (the latter is defined as a quadratic form on $V_0\simeq \PP{}{V}$).

In what follows, we take $V_0\subset\RR^{d+1}$ to be the $\mathcal{Q}_0$-orthogonal hyperplane to $\ker q_d$, where $\mathcal{Q}_0=\sum_j x_j^2$ (see Equation \ref{equation:ref_quadratic form}). 

\begin{proposition}
\label{proposition:orthogonality_projective_lines}
Let $q_d$ be a degenerate quadratic form over $\RR^{d+1}$ of rank $d$ contained in the pencil of quadrics $\polar{\mathcal{F}(\polar{q_1},\polar{q_2})}$. Let $V_0$ be the $\mathcal{Q}_0$-orthogonal hyperplane to $\ker q_d$ and $V\subset\RR^{d+1}$ an affine space parallel to $V_0$. Then given $p\in Q_1\cap\mathbb{P}(V)$, the intersections $T_pQ_1\cap\mathbb{P}(V)$ and $L_{Q_2}(p)\cap\mathbb{P}(V)$ are $\polar{q_d}$-orthogonal.
\end{proposition}

\begin{proof}
Let $p\in Q_1\cap\mathbb{P}(V)$ and $u$ be the pole of $T_pQ_1$ with respect to $q_2$. We want to show that the line $pu$ and the hyperplane $T_pQ_1$ are $\polar{q_d}$-orthogonal when intersected with $\mathbb{P}(V)$. Write $p=[X_1]$ and $u=[X_2]$.

We first find the $\polar{q_d}$-orthogonal line to $T_pQ_1$ through $p$ in $\mathbb{P}(V)$. Let $T_0\subset V_0$ be the hyperplane such that $T_pQ_1=\ensemble{[x+\lambda x_d]}{x\in T_0,\lambda\in\RR}$. Let $y\in V_0$ be such that $y$ is $\polar{q_d}$-orthogonal to $T_0$. Let $f_d:V_0\to V_0$ be the invertible linear map such that $q_d(x,x')=\mathcal{Q}_0(f_d(x),x')$ for all $x,x'\in V_0$, and $f_1:\RR^{d+1}\to\RR^{d+1}$ be the invertible linear map defined for all $x,x'\in\RR^{d+1}$ by $q_1(x,x')=\mathcal{Q}_0(f_1(x),x')$. By construction, for all $x\in T_0$ we have
$$0=\polar{q_d}(y,x)=\mathcal{Q}_0(\inverse{f_d}(y),x)$$
Writing $f_1(X_1) = \pi_{V_0}f_1(X_1)+\lambda x_d$, where $\pi_{V_0}f_1(X_1)\in V_0$ , $\lambda\in\RR$, we have for the same $x\in T_0$
$$\mathcal{Q}_0(\pi_{V_0}f_1(X_1),x)=\mathcal{Q}_0(f_1(X_1),x)=q_1(X_1,x)=0$$
since $V_0$ and $\ker q_d$ are $\mathcal{Q}_0$-orthogonal. We deduce that $\inverse{f_d}(y)$ and $\pi_{V_0}f_1(X_1)$ are colinear, hence $y$ and $f_d(\pi_{V_0}f_1(X_1))$ are colinear. Hence the $\polar{q_d}$-orthogonal space to $T_0$ is $\RR f_d(\pi_{V_0}f_1(X_1))$.

Now let us show that $f_d(\pi_{V_0}f_1(X_1))\in V_0$ considered as a point in $\PP{}{V}$ (after the above identification $V\simeq V_0\simeq \PP{}{V}$) lies in the line $pu$. This will conclude the proof. Let $f_2:\RR^{d+1}\to\RR^{d+1}$ be the invertible linear map defined for all $x,x'\in\RR^{d+1}$ by $q_2(x,x')=\mathcal{Q}_0(f_2(x),x')$. The deinition of pole and tangent line to $Q_1$ allows us to write $[f_1(X_1)]=[f_2(X_2)]$, hence $u=[X_2]=[\inverse{f_2}\circ f_1(X_1)]$. And since $q_d\in\mathcal{F}(\polar{q_1},\polar{q_2})$, we can write $\alpha_2\inverse{f_2}=f_d+\alpha_1\inverse{f_1}$ where $\alpha_1,\alpha_2\in\polar{\RR}$ and $f_d$ has been extended to the whole space $\RR^{d+1}$ by $f_d(x_d)=0$. Therefore, $u=[\inverse{f_2}\circ f_1(X_1)]=[f_d(f_1(X_1))+\alpha_1 X_1]$. This shows that $f_d(f_1(X_1))=f_d(\pi_{V_0}f_1(X_1))\in pu$.
\end{proof}

We apply the previous result to families of (pseudo-)confocal quadrics as it was described in Section \ref{section:general_properties_on_quadrics}. Fix an integer $k\in\{0,\ldots,d-1\}$, real numbers $a_0<a_1<\ldots<a_d$ and consider the family of quadrics $Q^k:=(Q_{\lambda}^k)_{\lambda\neq a_j}$ of $\RP^d$ given by the equation
\begin{equation}
\label{equation:confocal_quadrics_minkowski}
Q_{\lambda}^k: \sum_{j=0}^k \frac{x_j^2}{a_j-\lambda}+\sum_{j=k+1}^{d-1} \frac{x_j^2}{a_j+\lambda}=x_d^2.
\end{equation}
We think of this family as a family of confocal quadratics for a certain pseudo-Euclidean metric (an Euclidean metric when $k=d-1$) which is the degenerate quadratic form of $\RR^{d+1}$ defined by
$$q_d^k(x)=\sum_{j=0}^k x_j^2-\sum_{j=k+1}^{d-1} x_j^2.$$
The restriction of $q_d^k$ to the affine chart $\{x_d=1\}\simeq\RR^d$ is a non-degenerate quadratic form on $\RR^d$. Note that given a non-degenerate quadratic form $q_1$ defining a non-degenerate quadric of $Q^k$, the quadrics of $Q^k$ are defined by the pencil $\polar{\mathcal{F}(\polar{q_1},q_d^k)}$. By Propositions \ref{proposition:quadrics_as_projective_caustics} and \ref{proposition:confocal_projective_caustics}, all quadrics $Q_{\lambda}^k$ with $\lambda\neq 0$ define the same field of projective lines $L_k(p)$ on $Q:=Q_{0}^k$ for which they are caustics. 

\begin{proposition}
The line $L_k(p)$ of the field of projective lines defined on $Q:=Q_{0}^k$ by its confocal quadrics is $q_d^k$-orthogonal to $T_pQ$ in the affine chart $\{x_d=1\}\simeq\RR^d$.
\end{proposition}

\begin{proof}
This is a direct consequence of Proposition \ref{proposition:orthogonality_projective_lines}: by definition, the quadratic forms defining each $Q_{\lambda}^k$ form a pencil of quadrics containing $q_d^k$. Moreover, the $\mathcal{Q}_0$-orthogonal space to $\ker q_d^k$ is the vector space $V_0=\{x_d=0\}$. The restriction to $V_0$ of $q_d^k$ is non-degenerate and its dual $\polar{q_d^k}$ with respect to the restriction $\restreint{\mathcal{Q}_0}{V_0}$ is $q_d^k$ itself.
\end{proof}

\section{On Berger property \textit{Only quadrics have caustics}}
	\label{section_berger_property}

In this section, we are interested in projective billiards in dimension $d\geq 3$ having caustics. We try to study a generalization to projective billiards of a fundamental result discovered by Berger in \cite{berger_caustics} and which can be stated as follows:

\begin{theorem}[Berger, \cite{berger_caustics}]
Let $d\geq 3$ and $S$, $U$, $V$ be open subsets of $\class^2$-smooth hypersurfaces in $\RR^d$ with non-degenerate second fundamental forms. Suppose that there is an open subset of lines tangent to $U$ and intersecting $S$ transversally which are reflected into lines tangent to $V$. Then $S$ is a piece of quadric and $U,V$ are pieces of one and the same quadric confocal to $S$.
\end{theorem}

Glutsyuk \cite{glut_caustics} extended Berger's result to space forms of non-zero constant curvature, that is to the Euclidean unit sphere $\mathbb{S}^d$ and the hyperbolic space $\mathbb{H}^d$ with $d\geq 3$. In the corresponding paper \cite{glut_caustics}, this result is also used to prove the Commuting Billiards Conjecture in dimension $d\geq 3$. This conjecture was stated by Tabachnikov, see \cite{taba_commut, taba_book}, and was also proved by Glutsyuk in dimension $2$, see \cite{glut}. We can now call it Commuting Billiards Theorem, and the latter can be stated as follows: \textit{consider two nested billiards of $\RR^d$ (respectively $\RR^2$) with $\class^2$-smooth (respectively piecewise $\class^4$-smooth) strictly convex boundaries. Each one of their billiard maps acts on the set of oriented lines intersecting them. If these maps commute then the billiards are confocal ellipsoids (respectively ellipses).}

This section is structured as follows. We first extend in Subsection \ref{subsection:berger_argument} a key argument of Berger's proof to projective billiards. In Subsection \ref{subsection:berger_hyperplanes_distribution}, we define a distribution of hyperplanes related to cones of lines tangent to possible caustics. In Subsection \ref{subsection:berger_pseudo_euclidean}, we apply results found in the two first subsections to show that if a convex pseudo-Euclidean billiard has a caustic, then it is a quadric and its caustic is a confocal quadric for the pseudo-Euclidean metric.

\subsection{Berger's key argument for projective billiards} 
\label{subsection:berger_argument}

We first extend a result based on an observation made by Berger \cite{berger_caustics} in the case of usual billiards. This observation can be stated as follows for projective billiards: consider hypersurfaces $S,U,V$ of $\RR^d$ and a line-framed hypersurface $\Sigma$ over $S$. Suppose that there is an open subset of lines $\ell$ tangent to $U$ and intersecting $S$ which are reflected into lines $\ell'$ tangent to $V$ by the projective law of reflection on $\Sigma$. If we consider three non-colinear points $A\in U$, $B\in S$, $C\in V$ such that the above property is satisfied with $\ell=AB$ and $\ell'=BC$, then by symmetry \textit{the intersection $H=T_A U\cap T_BS$ coincides with $T_C V\cap T_BS$}. 

This observation leads to a result of finiteness on such hyperplanes $H$ of $T_BS$ which we are going to detail below. Let $\Sigma$ be a line-framed hypersurface over a hypersurface $S\subset\RR^d$, $(B,L)\in\Sigma$ and $\xi\in T_BS$ a non-zero vector.

\begin{definition}
\label{definition:H_permitted}
Let $H$ be a hyperplane of $T_BS$. $H$ is said to be \textit{permitted by $\xi$} if $\xi\notin H$ and for any $\class^{1}$-smooth germ of curve $(B(t),L(t))\in\Sigma$ with $B(0)=B$ and $L(0)=L$, and any $\class^1$-smooth germ of curve $\xi(t)\in T_{B(t)}S$ with $\xi(0)=\xi$, there exist germs of $\class^1$-smooth curves $A(t)$, $C(t)$ in $\RR^d$ such that \\
\enum $A:=A(0)$, $B$, $C:=C(0)$ are not colinear and the vector $\xi$ belongs to the plane $ABC$;\\
\enum $A'(0)$ belongs to the hyperplane containing the line $AB$ and $H$;\\
\enum $C'(0)$ belongs to the hyperplane containing the line $CB$ and $H$;\\
\enum $A(t)B(t)$ is reflected into $B(t)C(t)$ by the projective reflection law on $\Sigma$ at $B(t)$.
\end{definition}

\begin{definition}
If $\ell$ is a line intersecting $S$ transversally at $B$, we will say that a hyperplane \textit{$H\subset T_BS$ is permitted by $\ell$} if it is permitted by a non-zero vector $\xi$ in the intersection of the plane containing $L$ and $\ell$ with $T_BS$.
\end{definition}

The main result of this section is the following proposition:

\begin{proposition}
\label{proposition:H_permitted}
Suppose $S$ has non-degenerate second fundamental form at $B$. Then there is a closed subset $F$ of $T_BS$ such that $T_BS\setminus F$ is dense in $T_BS$ and such that for all $\xi\in T_BS\setminus F$, the number of hyperplanes $H\subset T_BS$ permitted by $\xi$ is at most $d-1$.
\end{proposition}

\begin{remark}
In the proof we show that $F$ is the finite union of strict vector subspaces of $T_BS$.
\end{remark}

\begin{proof}[Proof of Proposition \ref{proposition:H_permitted}]
The proof is computational and we wonder if one can find a more geometrical one. We apply the same formulas as in \cite{berger_caustics}, and use a result of linear algebra to conclude. We endow $\RR^d$ with its usual Riemannian metric: we denote by $x\cdot y$ the canonical scalar product on $\RR^d$.

Let $H$ be a hyperplane of $T_BS$ permitted by a certain $\xi\in T_BS$, and let $\eta\in T_B S$ be an orthogonal vector to $H$ of norm $1$. Choose an orthonormal basis $(u_1,\ldots,u_{d-1})$ of eigenvectors of $S$'s second fundamental form at $B$, and denote by $k_1,\ldots,k_{d-1}$ the corresponding eigenvalues. Choose $\alpha_i$ such that $\xi = \cos(\alpha_i) u_i + \sin(\alpha_i) v_i$, where $v_i$ is a vector of length $1$ orthogonal to $u_i$, and write $\ell_i=\cos(\alpha_i)\in[-1,1]$. Note that $\xi = \sum_{k=1}^{d-1}\ell_ku_k$.

Fix an $i\in\{1,\ldots,d-1\}$. Let $A(t)$, $(B(t),L(t))$, $C(t)$ be as in Definition \ref{definition:H_permitted} and verifying $B'(0)=u_i$. Let $n(t)$ be a smooth family of normal vectors to $S$ at $B(t)$ of length $1$, and $\nu(t)$ a smooth family of vectors directing $L(t)$ and such that $\nu(t)\cdot n(t)=1$.

Define $u_i(t)=B'(t)$ for all $t$, extend $v_i$ into a vector field $v_i(t)$ along $B(t)$ by parallel transport, and set $\xi(t) = \cos(\alpha_i)u_i(t) +\sin(\alpha_i)v_i(t)$. In the following, given any curve $\gamma(t)$, we will write $\gamma'$ for $\gamma'(0)$.
\vspace{0.2cm}

\textbf{1)} \textit{We first express in a matrix form the fact that $A'(0)$ belongs to the hyperplane $P_A$ containing the line $AB$ and $H$, and that $C'(0)$ belongs to the hyperplane $P_C$ containing the line $AC$ and $H$. }

Here we adapt the computations of \cite{berger_caustics} to the projective case. Let $e(t) = E_1(t)\nu(t)+E_2(t)\xi(t)$ and $\overline{e}(t)= E_1(t)\nu(t)-E_2(t)\xi(t)$ be two unit vectors directing the oriented lines $A(t)B(t)$ and $B(t)C(t)$, with $E_1(t),E_2(t)\in\RR$, and having the same orientation with respect to $T_BS$. One can write 
\begin{equation}
\label{equ:AetC}
A(t) = B(t)+a(t)e(t) \qquad\text{and}\qquad C(t) = B(t)+c(t)\overline{e}(t)
\end{equation}
where $a(t),c(t)>0$. Normal vectors to $P_A$ and $P_C$ can be respectively defined by 
$$n_A= (\eta\cdot e) n - ( n\cdot e)\eta = (E_1(\eta\cdot\nu)+E_2(\eta\cdot\xi))n-E_1\eta,$$
$$n_C= (\eta\cdot \overline{e}) n - (n\cdot\overline{e})\eta = (E_1(\eta\cdot\nu)-E_2(\eta\cdot\xi))n-E_1\eta.$$
If we denote by $'$ the derivative taken in $0$, we get
\begin{equation}
\label{equ:deriveAetC}
A'\cdot n_A = 0 \qquad\text{and}\qquad C'\cdot n_C = 0.
\end{equation}
Now since $e$ is parallel to $P_A$ we have $e\cdot n_A=0$, hence by combining Equations \eqref{equ:AetC} and \eqref{equ:deriveAetC} we get
\begin{equation}
\label{equ:deriveA}
0=A'\cdot n_A=u_i\cdot n_A + a (e'\cdot n_A)
\end{equation}
Yet, as recalled in \cite{berger_caustics} we have $n' = k_iu_i$, $u_i' = -k_in$, $v_i' = 0$, hence $\xi'=-k_i\ell_i n$. Therefore replacing $e'$ by its expression in Equation \eqref{equ:deriveA} gives
\begin{equation}
\label{equ:deriveA1}
\begin{matrix}
0&=&A'\cdot n_A &=&-E_1(\eta\cdot u_i)+a\left((E_2E_1'-E_1E_2'-k_i\ell_iE_2^2+E_1E_2(n\cdot\nu'))(\eta\cdot\xi)\right.\\
&& & &\left.\qquad\qquad +(E_1^2(n\cdot\nu')-k_i\ell_iE_1E_2)(\eta\cdot\nu) -E_1^2(\eta\cdot\nu')\right)
\end{matrix}
\end{equation}
and the same with $C$ by changing $E_2$ in $-E_2$ and $a$ in $c$:
\begin{equation}
\label{equ:deriveC1}
\begin{matrix}
0&=&C'\cdot n_C &=& -E_1(\eta\cdot u_i)+c\left((E_2'E_1-E_1'E_2-k_i\ell_iE_2^2-E_1E_2(n\cdot\nu'))(\eta\cdot\xi)\right.\\
&& & &\left.\qquad\qquad +(E_1^2(n\cdot\nu')+k_i\ell_iE_1E_2)(\eta\cdot\nu) -E_1^2(\eta\cdot\nu')\right)
\end{matrix}
\end{equation}
If we substitute Equations \eqref{equ:deriveA1} and \eqref{equ:deriveC1} in $\frac{A'\cdot n_A}{a}+\frac{C'\cdot n_C}{c}=0$ we get
\begin{equation}
\label{eq:equ1}
\begin{array}{rcl}
 \left(\frac{1}{a}+\frac{1}{c}\right)E_1(\eta\cdot u_i)&=& 2E_1^2(n\cdot\nu')(\eta\cdot\nu)-2E_1^2(\eta\cdot\nu')-2k_i\ell_iE_2^2(\eta\cdot\xi)\\
&=& -2E_1^2(\eta\cdot N_i)-2k_i\ell_iE_2^2(\eta\cdot\xi)
\end{array}
\end{equation}
where $N_i = \nu'-(n\cdot\nu')\nu$. The vector $N_i$ lies in $T_BS$ since we can check that $(N_i\cdot n)=0$. Now $(n(t)\cdot\nu(t))=1$ for all $t$, hence $(n\cdot\nu')=-(n'\cdot\nu)=-k_i(u_i\cdot\nu)$, and $N_i$ can be expressed as $N_i=d\nu\cdot u_i+k_i(u_i\cdot\nu)\nu$. \textit{$N_i$ only depends on the $2$-jet of $\Sigma$ at $(B,L)$. }

Let us rewrite Equation \eqref{eq:equ1} in a matrix form. Denote by\\
\enum $\alpha$  the quantity $-(a^{-1}+c^{-1})/2E_1$;\\
\enum $V_{\xi}$ the vector given by $V_{\xi}=\left(E_2/E_1\right)^2\sum_{i=1}^{d-1}k_i\ell_i u_i$;\\
\enum $M$ the matrix of $\matspace{d-1}{\RR}$ whose lines are given by the coordinates of $N_i$ in the orthonormal basis $(u_1,\ldots,u_{d-1})$.

Then Equation \eqref{eq:equ1}, together with the assumption that $\xi\notin H$, can be rewritten as
\begin{equation}
\label{eq:principale}
\left\{\begin{array}{l}
M\eta+(\xi\cdot\eta)V_{\xi}=\alpha\eta\\
\eta\notin\xi^{\bot}
\end{array}\right.
\end{equation}
where $\eta$ and $\xi$ are considered as vectors of $\RR^{d-1}$ with coordinates given by their coordinates in the basis $(u_1,\ldots,u_{d-1})$.
\vspace{0.2cm}

\textbf{2)} \textit{For fixed $\alpha$ and $\xi$, we now study the space of solutions $\eta$ of Equation \eqref{eq:principale}.}

Given an endomorphism $f$ of $\RR^k$, we call \textit{eigenspace} of $f$ any subspace of $\RR^k$ of the form $\ker (f-\beta id)$ for a certain $\beta$, and denote by $\im f$ the set of all $f(x)$ where $x\in \RR^k$.

\begin{lemma}
\label{lemma:berger_alg_lin}
Let $\beta_1,\ldots,\beta_s$ be the real eigenvalues of $M$. If $V_{\xi}$ doesn't belong to any of the sets $\im (M-\beta_iI_{d-1})$, then the eigenspaces of the endomorphism $f_{M,\xi}$ of $\RR^{d-1}$ defined by
$$f_{M,\xi}:x\mapsto Mx+(\xi\cdot x)V_{\xi}$$
are either of dimension at most $1$ or contain only orthogonal vectors to $\xi$.
\end{lemma}

\begin{proof}
Let $\alpha$ be an eigenvalue  of $f_{M,\xi}$. Consider an eigenvector $x$ of $f_{M,\xi}$ associated to $\alpha$: it satisfies
$$(M-\alpha I_{d-1})x=-(\xi\cdot x)V_{\xi}.$$
If $\alpha\neq\beta_k$ for all $k$, then $M-\alpha I_{d-1}$ is invertible. Hence $x\in\RR (M-\alpha I_{d-1})^{-1}V_{\xi}$ and therefore the eigenspace of $f_{M,\xi}$ associated to $\alpha$ is of dimension at most $1$. Now if $\alpha=\beta_k$ for a certain $k$, since $V_{\xi}\notin\im (M-\beta_kI_{d-1})$ we necessarily have $(\xi\cdot x)=0$ and the eigenspace of $f_{M,\xi}$ associated to $\alpha$ contains only orthogonal vectors to $\xi$.
\end{proof}

We can now finish the proof of Proposition \ref{proposition:H_permitted}. The second fundamental form of $S$ at $B$ is non-degenerate, hence none of the $k_i$ equals $0$. Hence the set $F$ of vectors $\xi\in\RR^{d-1}=\sum_{i=1}^{d-1}\ell_i u_i$ such that $V_{\xi}$ belongs to $\cup_i\im (M-\beta_iI_{d-1})$ is the finite union of strict vector subspaces of $T_BS$, and it depends neither on $E_1$ nor on $E_2$; it depends only on $M$.

Thus if we suppose that $\xi\notin F$, Lemma \ref{lemma:berger_alg_lin} implies that there are one-dimensional vector subspaces $G_1,\ldots,G_s$ of $T_BS$, $s\leq d-1$, contained in the eigenspaces of $f_{M,\xi}$ such that any solution $\eta$ of Equation \eqref{eq:principale} is contained in some $G_i$. Hence any hyperplane $H$ permitted by $\xi$ is an orthogonal space in $T_BS$ to one of the $G_i$, which ends the proof of the result.
\end{proof}

\subsection{Distributions of permitted hyperplanes}
\label{subsection:berger_hyperplanes_distribution}

Let $\Sigma$ be a line-framed hypersurface over a hypersurface $S\subset\RR^d$, and $B\in S$ such that $S$ has a non-degenerate second fundamental form at $B$. In this section, we define $d-1$ distributions of hyperplanes based on Proposition \ref{proposition:H_permitted}. We use the same notations as in the statement of this proposition.

Let $\xi\in T_BS$ be a vector outside $F$ and $H\subset T_BS$ be a hyperplane permitted by $\xi$. As a consequence of Lemma \ref{lemma:berger_alg_lin} together with the implicit function theorem, there are neighborhoods $U_H$ of $H$ in the Grassmanian and $U_{\xi}$ of $\xi$ in $T_BS$, and an analytic map $\widehat{H}:U_{\xi}\to U_H$ such that for any hyperplane $H'\in U_H$, $H'$ is permitted by a vector $\xi'\in U_{\xi}$ if and only if $H'=\widehat{H}(\xi')$. The same can be done with hyperplanes permitted by a line $\ell$.

In the case when the number of hyperplanes of $T_BS$ permitted by $\ell$ is exactly $d-1$, we have defined, in a neighborhood of $\ell$, $d-1$ analytic fields of hyperplanes $H_1,\ldots,H_{d-1}$ in $T_BS$, depending on $\ell$. For each one of them, define $\tilde{H}_k(\ell)$ the affine hyperplane of $\RR^d$ containing $H_k(\ell)$ and $\ell$. The latter can be considered as a hypersurface of the set of lines $\mathcal{L}_B\simeq \RP^{d-1}$ containing $B$, and it contains $\ell$. We denote by $h_k(\ell)=T_{\ell}\tilde{H}_k(\ell)\subset T_{\ell}\mathcal{L}_B$ its tangent space at $\ell$.

\begin{definition}
\label{definition:permitted_distribution}
The $d-1$ fields of hyperplanes $h_1,\ldots,h_{d-1}$ define distributions on an open subset of $\mathcal{L}_B\simeq\RP^{d-1}$ containing $\ell$ and called the \textit{permitted distributions of $\Sigma$ at $B$}.
\end{definition}

\begin{proposition}
\label{proposition:berger_2_jet}
Suppose that the number of hyperplanes of $T_BS$ permitted by $\ell$ is exactly $d-1$, so that the $d-1$ permitted distributions of $\Sigma$ at $B$ are well-defined. Further suppose that $\Sigma'$ is a line-framed hypersurface over a hypersurface $S'\subset\RR^d$ containing $B$. If $S$ and $S'$ have the same $2$-jet at $B$ and the fields of projective lines of $\Sigma$ and $\Sigma'$ have the same $1$-jet at $B$, then the $d-1$ permitted distributions at $B$ of $\Sigma'$ are well-defined and coincide with the permitted distributions of $\Sigma$ at $B$.
\end{proposition}

\begin{proof}
This is a direct consequence of Equation \eqref{eq:principale}, which only depends on the $1$-jet of the normal vector field $n$ to $S$ at $B$ and on the $1$-jet at $B$ of a normalized directing vector field $\nu$ of the projective field of lines of $\Sigma$.
\end{proof}

In the following proposition, we suppose that $\Sigma$ has $d-1$ permitted distributions at $B$, $h_1,\ldots,h_{d-1}$, well-defined in a neighborhood of a line $\ell$ intersecting $S$ at $B$. 

\begin{proposition}
\label{proposition:berger_caustics_permitted_distribution}
Let $U,V\subset\RR^d$ be hypersurfaces with non-degenerate second fondamental forms, with $\ell$ tangent to $U$. Suppose that there is an open subset $\Omega$, containing $\ell$, of lines tangent to $U$ and intersecting $S$ which are reflected into lines tangent to $V$ by the projective law of reflection on $\Sigma$. Then for any line $\ell'$ in $\Omega$ containing $B$ and tangent to $U$ at a point $A$, the hyperplane $H=T_AU\cap T_BS$ is permitted by $\ell'$. The corresponding hyperplanes on $T\mathcal{L}_B$ coincide with one of the $d-1$ permitted distributions of $\Sigma$ at $B$, say $h_j$. The set of lines containing $B$ and tangent to $U$ is an integral surface of $h_j$.
\end{proposition}

\begin{proof}
It is easy to check that $H$ satisfies all requirements of Definition \ref{definition:H_permitted} since the pair $U,V$ is a caustic of $S$. Hence $H$ is a hyperplane of $T_BS$ permitted by $\ell'$. The corresponding hyperplane $h(\ell')\subset T_{\ell'}\mathcal{L}_B$ coincides with one of the $d-1$ permitted distributions $h_j(\ell')$ by definition. The rest of the result follows from the definition of the distribution $h_j$.
\end{proof}

\subsection{Caustics of billiards in pseudo-Euclidean spaces}
\label{subsection:berger_pseudo_euclidean}

In this section, we are interested specifically in billiards of pseudo-Euclidean spaces of dimension $d\geq 3$ having caustics. We adopt the definition of pseudo-Euclidean spaces used in \cite{DragRad_minkowski2, khesin_taba}: the pseudo-Euclidean space $E^{k,l}$ of signature $(k,l)$, $k,l\in\NN$, with $k+l=d$, is the space $\RR^d$ endowed with the non-degenerate symetric bilinear form $\scalar{\cdot}{\cdot}$ defined for $x,y\in\RR^d$ by 
\begin{equation}
\label{equation:pseudo_euclidean_metric}
\scalar{x}{y}=\sum_{j=1}^k x_jy_j-\sum_{j=k+1}^d x_jy_j.
\end{equation}
We will denote by $q_d^k$ the quadratic form associated to $\scalar{\cdot}{\cdot}$. A line $\ell\subset E^{k,l}$ directed by a non-zero vector $v$ is said to be\\
\enum \textit{space-like} if $\scalar{v}{v}>0$;\\
\enum \textit{time-like} if $\scalar{v}{v}<0$;\\
\enum \textit{light-like} if $\scalar{v}{v}=0$.

Denote the usual scalar product on $\RR^d$ by $\rscalar{x}{y}$ for all $x,y\in\RR^d$. An \textit{affine ellipsoid} is a set containing at least two points which is of the form
$$\mathcal{E}=\ensemble{x\in\RR^d}{\rscalar{Ax}{x}+\rscalar{B}{x}+C=0}$$
where $A\in\matspace{d}{\RR}$ is a positive-definite symmetric matrix, $B\in\RR^d$ is a vector, and $C\in\RR$. As noticed in \cite{DragRad_minkowski2, khesin_taba}, since $A$ is positive-definite there is a linear change of coordinates in $\RR^d$ preserving  the pseudo-Euclidean metric \eqref{equation:pseudo_euclidean_metric}, in which $\rscalar{Ax}{x}$ takes the form 
$$\rscalar{Ax}{x}=\sum_{j=1}^d\frac{x_j^2}{a_j}$$
where $a_1,\ldots,a_d>0$. Therefore, by an appropriate choice of a new origin, the ellipsoid $\mathcal{E}$ is given by an equation of the form 
\begin{equation}
\sum_{j=1}^d \frac{x_j^2}{a_j'}=1
\end{equation}
where $a_1',\ldots,a_d'>0$. Notice that the form of the pseudo-Euclidean metric \eqref{equation:pseudo_euclidean_metric} is left unchanged in this set of coordinates. A \textit{pseudo-confocal quadric} to $\mathcal{E}$ is a quadric $Q_{\lambda}$ which can be expressed in this new set of coordinates by an equation of the form
\begin{equation}
Q_{\lambda}:\sum_{j=1}^k \frac{x_j^2}{a_j'-\lambda}+\sum_{j=k+1}^d \frac{x_j^2}{a_j'+\lambda}=1.
\end{equation}
where $\lambda\in\RR$. See Figure \ref{figure:pseudo_confocal_quadrics} for a $3$-dimensional representation.

Let us recall results from \cite{DragRad_minkowski2, khesin_taba} which can be applied to a general quadric, not only an ellipsoid:

\begin{theorem}[pseudo-Euclidean version of Chasles theorem; see \cite{khesin_taba} Theorem 4.8 and \cite{DragRad_minkowski2} Theorem 2.3]
\label{theorem:chasles_pseudo}
A space- or time-like line $\ell$ intersecting $\mathcal{E}$ is tangent to $d-1$ pseudo-confocal quadrics. The tangent hyperplanes at the tangency points with $\ell$ are pairwise orthogonal.
\end{theorem}

\begin{theorem}[pseudo-Euclidean version of Jacobi-Chasles theorem; see \cite{khesin_taba} Theorem 4.9]
\label{theorem:jacobi_chasles_pseudo}
A space- or time-like billiard trajectory in $\mathcal{E}$ remains tangent to $d-1$ fixed quadrics $Q_{\lambda}$.
\end{theorem}

We deduce the following result:

\begin{theorem}
\label{theorem:caustics_pseudo_euclidean}
Let $S,U,V\subset \RR^d$ be open subsets of hypersurfaces with non-degenerate second fundamental forms, and $S$ being convex. Suppose that there is an open subset of the set of lines intersecting $S$ and tangent to $U$ which are reflected into lines tangent to $V$ by pseudo-Euclidean reflection on $S$. Then $S$ is a piece of quadric; $U,V$ are pieces of one and the same quadric pseudo-confocal to $S$.
\end{theorem}

\begin{proof}
We first begin by proving the following

\begin{lemma}
\label{lemma:hypersurface_light_like}
Consider the set $K$ of lines tangent to $U$ and intersecting $S$ transversally. Then the set of lines of $K$ which are not light-like is an open dense subset of $K$.
\end{lemma}

\begin{proof}
Denote by $K_0$ this set. It is clearly open, and it remains to show that it is dense. Suppose the contrary, then there is a light-like line $\ell\in K$ contained in an open subset $\Omega\subset K\setminus K_0$ of $K$. Let $x\in U$ be the point of tangency of $\ell$ with $U$. All lines tangent to $U$ at $x$ and sufficiently close to $\ell$ are contained in $\Omega$, hence the pseudo-Euclidean metric $\scalar{v}{v}$ vanishes on an open subset of $T_xU$, hence on $T_xU$. This implies that $T_xU$ is contained in its orthogonal space, contradiction.
\end{proof}

We can now prove Theorem \ref{theorem:caustics_pseudo_euclidean}. Consider the line-framed hypersurface $\Sigma=\linep{S}{q_d^k}$ over $S$ defined by the pseudo-Euclidean metric $q_d^k$ (see Subsection \ref{subsection:projective_and_others}). Consider a line $\ell$ of $K$ which is not light-like and intersecting $S$ at a point $B$ transversally. Using the same notations as in Proposition \ref{proposition:H_permitted}, we can suppose that any non-zero vector $\xi$ of the intersection of $T_BS$ with the plane containing $\ell$ and the pseudo-Euclidean normal line to $S$ at $B$ is not in $F$. 

Now consider an affine ellipsoid $\mathcal{E}$ tangent to $S$ at $B$, with the same principal directions. In particular, $S$ and $\mathcal{E}$ have the same $2$-jet at $B$, and therefore their field of normal lines with respect to the pseudo-Euclidean metric have the same $1$-jet at $B$. Let $\mathcal{E}'=\linep{\mathcal{E}}{q_d^k}$ be the line-framed hypersurface over $\mathcal{E}$ induced by the pseudo-Euclidean metric $q_d^k$: let us show that $\mathcal{E}'$ has $d-1$ permitted distribution at $B$ which are integrable. By Proposition \ref{proposition:berger_2_jet}, $\Sigma$ will have the same permitted distributions at $B$.

Indeed, by Theorem \ref{theorem:chasles_pseudo}, one can find $d-1$ pseudo-confocal quadrics $Q_1,\ldots,Q_{d-1}$ tangent to $\ell$, and such that the hyperplanes containing $\ell$ and tangent to each $Q_j$  are pairwise orthogonal. By Theorem \ref{theorem:jacobi_chasles_pseudo}, each $Q_j$ is a caustic of $\mathcal{E}'$, and by Proposition \ref{proposition:berger_caustics_permitted_distribution}, the intersection of the latter hyperplanes with $T_B\mathcal{E}$ are pairwise distinct and permitted by $\ell$. Therefore $\mathcal{E}'$ has $d-1$ permitted distributions $h_1, \ldots,h_{d-1}$ defined on a neighborhood of $\ell$ in $\mathcal{L}_B$ which are induced by the latter hyperplanes. Moreover, the distributions $h_1,\ldots,h_{d-1}$ are integrable: the integral manifold of $h_j$ is the set of lines tangent to one quadric among $Q_1,\ldots,Q_{d-1}$. Therefore the union of all lines contained in the integral manifold of $h_j$ is a piece of quadratic cone (a cone defined by a quadratic form). 

Yet, the set of lines containing $B$ and tangent to $U$ and $V$ is also an integral manifold of one of the permitted distributions $h_j$ of $\Sigma$ at $B$, by Proposition \ref{proposition:berger_2_jet} and Proposition \ref{proposition:berger_caustics_permitted_distribution}. Hence $U$ and $V$ are tangent to one of previously defined quadratic cones at points defining curves of $U$ and $V$. The same operation can be applied by taking different points $B$ in a small open subset of $S$. Now by an argument of \cite{berger_caustics} working on the duals of $U$ and $V$, this implies that $U$ and $V$ are pieces of one and the same quadric $Q$.

To prove that $S$ is a piece of pseudo-confocal quadric, we use the same argument as Berger in the pseudo-Euclidean case. For points $B'\in\RR^d$ close to $B$, consider the quadratic cone $C_{B'}$ tangent to $Q$ and containing $B'$. We can define a hyperplane $H_{B'}\subset\RR^d$ containing $B'$ as the only hyperplane through $B'$ close to $T_BS$ such that $C_{B'}$ is symmetric by the pseudo-Euclidean orthogonal symmetry with respect to $H_{B'}$. The induced hyperplane distribution is integrable and its integral surfaces are quadrics pseudo-confocal to $Q$ by Theorem \ref{theorem:jacobi_chasles_pseudo}. Hence $S$ is a piece of quadric pseudo-confocal to $Q$.
\end{proof}
\chapter{Billiards with open subsets of periodic orbits}
	\label{chapter_ivrii}
	This chapter is devoted to the study of periodic orbits of projective billiards and of an analogue of Ivrii's conjecture for projective billiards. 

A projective billiard $\Omega$ is said to be \textit{$k$-reflective} if we can find a $k$-periodic orbit $p=(p_1,\ldots,p_{k})$ and an open subset $U_1\times U_2\subset (\partial\Omega)^2$ containing $(p_1,p_2)$ such that for any $(q_1,q_2)\in U_1\times U_2$, $\Omega$ has a $k$-periodic orbit $q=(q_1,q_2,\ldots,q_k)$ close to $p$. The billiard $\Omega$ is said to be \textit{$k$-pseudo-reflective} if the same statement is satisfied with $U_1\times U_2$ replaced by a subset of $(\partial\Omega)^2$ of non-zero measure. Using these definitions, Ivrii's conjecture for projective billiards can be stated as the following injunction: \textit{classify the $k$-pseudo-reflective projective billiards for all integer $k\geq 3$.}

The conditions of $k$-reflectivity and $k$-pseudo-reflectivity are \textit{local} properties in the following sense: as a consequence of Proposition \ref{proposition:billiard_regularity_rank}, if $(q_1,q_2)$ is chosen sufficiently close to $(p_1,p_2)$ on $(\partial\Omega)^2$, then the corresponding orbit $q=(q_1,q_2,\ldots,q_k)$ has its vertices $q_j$ in open small neighborhoods $V_j\subset\partial\Omega$ of the vertices $p_j$ of the orbit $p=(p_1,\ldots,p_{k})$. Therefore, if we perturb $\partial\Omega$ arbitrarily outside each $V_j$, the set of periodic orbits close to $p$ will remain unchanged, and the property of $k$-reflectivity and $k$-pseudo-reflectivity is kept intact.

\begin{figure}[!h]
\centering
\definecolor{xdxdff}{rgb}{0.49,0.49,1}
\definecolor{ffffqq}{rgb}{1,1,0}
\definecolor{tttttt}{rgb}{0.2,0.2,0.2}
\begin{tikzpicture}[line cap=round,line join=round,>=triangle 45,x=0.5cm,y=0.5cm]
\clip(-8,-6) rectangle (13,7);
\fill[dash pattern=on 2pt off 2pt,color=tttttt,fill=tttttt,fill opacity=0.1] (-7.3,-4.46) -- (3.46,5.94) -- (12.08,-2.26) -- cycle;
\draw [line width=1.2pt,color=tttttt] (-7.3,-4.46)-- (3.46,5.94);
\draw [line width=1.2pt,color=tttttt] (3.46,5.94)-- (12.08,-2.26);
\draw [line width=1.2pt,color=tttttt] (12.08,-2.26)-- (-7.3,-4.46);
\draw [line width=5.2pt,color=ffffqq] (6.54,3.01)-- (-3.44,-0.73);
\draw [line width=5.2pt,color=ffffqq] (6.54,3.01)-- (2.88,-3.3);
\draw [line width=5.2pt,color=ffffqq] (-3.44,-0.73)-- (6.48,-2.9);
\draw [line width=5.2pt,color=ffffqq] (0.82,3.39)-- (9.84,-0.13);
\draw [line width=5.2pt,color=ffffqq] (9.84,-0.13)-- (6.48,-2.9);
\draw [line width=5.2pt,color=ffffqq] (0.82,3.39)-- (2.88,-3.3);
\draw (2.88,-3.3)-- (6.54,3.01);
\draw (6.54,3.01)-- (-3.44,-0.73);
\draw (-3.44,-0.73)-- (6.48,-2.9);
\draw (6.48,-2.9)-- (9.84,-0.13);
\draw (9.84,-0.13)-- (0.82,3.39);
\draw (0.83,3.39)-- (2.88,-3.3);
\draw [dash pattern=on 5pt off 5pt] (-3.44,-0.73)-- (3.48,0.5);
\draw [dash pattern=on 5pt off 5pt] (3.48,0.5)-- (0.83,3.39);
\draw [dash pattern=on 5pt off 5pt] (3.48,0.5)-- (6.54,3.01);
\draw [dash pattern=on 5pt off 5pt] (3.48,0.5)-- (9.84,-0.13);
\draw [dash pattern=on 5pt off 5pt] (3.48,0.5)-- (6.48,-2.9);
\draw [dash pattern=on 5pt off 5pt] (3.48,0.5)-- (2.88,-3.3);
\begin{scriptsize}
\fill [color=black] (3.48,0.5) circle (1.5pt);
\fill [color=black] (2.88,-3.3) circle (1.5pt);
\draw[color=black] (3.04,-3.7) node {$p_1$};
\fill [color=black] (6.54,3.01) circle (1.5pt);
\draw[color=black] (6.85,3.35) node {$p_2$};
\fill [color=black] (-3.44,-0.73) circle (1.5pt);
\draw[color=black] (-3.6,-0.4) node {$p_3$};
\fill [color=black] (6.48,-2.9) circle (1.5pt);
\draw[color=black] (6.64,-3.3) node {$p_4$};
\fill [color=black] (9.84,-0.13) circle (1.5pt);
\draw[color=black] (10.1,0.2) node {$p_5$};
\fill [color=black] (0.82,3.39) circle (1.5pt);
\draw[color=black] (0.5,3.66) node {$p_6$};
\end{scriptsize}
\end{tikzpicture}
\caption{A local projective billiard in triangle with a $6$-periodic orbit  $(p_j)_{j=1\ldots 6}$ in yellow. The fields of projective lines are represented by dotted lines.}
\label{figure:6_reflective_projective}
\end{figure}
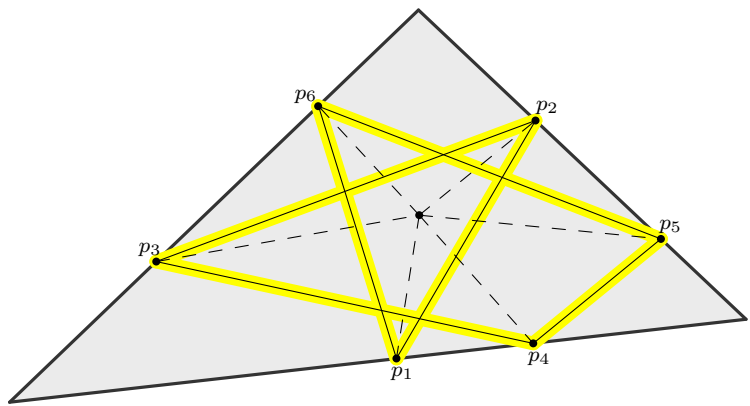

This is why we can only consider the germs of curves or hypersurfaces $(\partial\Omega,p_j)$, $j=1,\ldots,k$, and work with them to study the analogue of Ivrii's conjecture for projective billiards. We give a new local definition of projective billiards corresponding to this idea:

\begin{definition}
\label{definition:local_projective_billiard}
A \textit{local projective billiard} $\mathcal{B}$ is a collection of line-framed hypersurfaces $\alpha_1,\ldots,\alpha_k\subset\mathbb{P}(T\RR^d)$ over hypersurfaces $a_1,\ldots,a_k$ of $\RR^d$ called \textit{classical boundaries} of $\mathcal{B}$. It is said to be respectively $\class^r$-smooth (with $r=1,2,\ldots,\infty$) or analytic if all $\alpha_j$ are $\class^r$-smooth or are analytic.

An \textit{orbit} of $\mathcal{B}$ is a (finite or infinite) sequence of points $(p_j)_{j=-s\ldots t}$, with integers $s\leq t$ eventually infinite, such that for each $j$ (seen modulo $k$) \\
\enum $p_j\in a_{j\mod k}$, $p_j\neq p_{j+1}$ and the line $p_{j}p_{j+1}$ is oriented from $p_j$ to $p_{j+1}$;\\
\enum the line $p_{j}p_{j+1}$ is transverse to both $a_j$ and $a_{j+1}$;\\
\enum the line $p_{j}p_{j+1}$ is obtained from $p_{j-1}p_j$ by the projective reflection law of $\alpha_{j\mod k}$ at $p_j$.\\
If $k'$ is a multiple of $k$, an orbit $(p_j)_j$ is said to be \textit{$k'$-periodic} if $s=1$, $t=k'$ and $(p_1,\ldots,p_{k'},p_1,p_2)$ is an orbit.

A \textit{local classical billiard} (or simply \textit{local billiard}) is a local projective billiard whose line-framed hypersurfaces are induced by the Euclidean metric, \textit{i.e.} the lines of the projective fields of lines are orthogonal to the tangent hyperplanes (see Subsection \ref{subsection:projective_and_others}).

A local projective billiard $\mathcal{B}$ is said to be \textit{$k$-reflective} (respectively \textit{$k$-pseudo-reflective}) if there is a non-empty open subset (respectively a subset of non-zero measure) $U_1\times U_2\subset a_1\times a_2$ such that to any pair $(p_1,p_2)\in U_1\times U_2$ corresponds a $k$-periodic orbit of $\mathcal{B}$. The \textit{$k$-reflective} set of $\mathcal{B}$ is the set of pairs $(p_1,p_2)$ contained in open subsets $U_1\times U_2\subset a_1\times a_2$ satisfying previous property.
\end{definition}

\begin{remark}
\label{remark:billiard_analytic_pseudoref_reflective}
If an analytic local projective billiard is $k$-pseudo-reflective, then it is $k$-reflective. This result is an easy corollary of the Uniqueness Theorem for analytic extension.
\end{remark}

\begin{remark}
We will sometimes consider that the $\alpha_j$ are line-framed hypersurfaces of $\PP{}{T\RP^d}$ (see Remark \ref{remark:line_framed_hyp_proj}).
\end{remark}

This chapter is structured as follows: different examples of $k$-reflective projective billiards inside polygons are given at Section \ref{section:examples_reflective_billiards}. Section \ref{section_pfaffian_systems} introduces Pfaffian systems and applies them to the study of Ivrii's conjecture. A classification of the $3$-reflective and $3$-pseudo-reflective local projective billiards is given at Section \ref{section_3_reflective_proj_billiards}.

\section{Examples of $k$-reflective projective billiards}
	\label{section:examples_reflective_billiards}
		
In this section, we construct different types of local projective billiards in the plane which are $k$-reflective for all integer $k$ being either $3$ or an even integer. The constructed examples are projective billiards whose classical boundaries are lines, hence can be considered as projective billiards inside polygons. All the results presented here are gathered in a preprint \cite{fierobe1}. 

\begin{remark}
I apologize in advance for this remark which is not related directly to mathematical considerations. I remember the results of this section as a very pleasant moment of my thesis. In this remark I just describe a funny anecdote taking place in a train from Nizhnyi Novgorod to Novosibirsk and related to the discovery of the present results. I found some of the latter's proofs during this two days long journey, while the train was moving in the middle of beautiful empty frozen steppes. When the redaction was over, I tried to put the preprint online. But wifi was only available for short periods of time at different train stops situated midway between Nizhnyi Novgorod and Novosibirsk, so that I had to be quite dexterous and try again many times to achieve the upload. This journey appeared to me as a great moment of creativity, which probably would have been different if I decided to take the airplane instead of the train in order to join Novosibirsk from Nizhnyi Novgorod.
\end{remark}

Given distinct points $O,P,Q\in\RP^2$ not on the same line, one can consider the line $PQ$ endowed with the field of transverse lines containing $O$, denoted by $\linep{PQ}{O}$:
$$\linep{PQ}{O}=\ensemble{(p,L)\in\mathbb{P}(T\RP^2)}{p\in PQ,\, O\in L}.$$

\begin{figure}[!h]
\centering
\begin{tikzpicture}[line cap=round,line join=round,>=triangle 45,x=1.3cm,y=1.3cm]
\clip(-3.4,-2.5) rectangle (2.4,4);
\draw (2,1)-- (0,-2);
\draw (0,-2)-- (-3,3);
\draw (-3,3)-- (2,1);
\draw (-2.84,3.65)-- (-1.98,1.94);
\draw (-1.92,3.28)-- (-1.34,1.68);
\draw (-1.14,2.97)-- (-0.79,1.46);
\draw (-0.43,2.69)-- (-0.3,1.27);
\draw (0.2,2.43)-- (0.14,1.09);
\draw (0.86,2.17)-- (0.6,0.9);
\draw (1.58,1.88)-- (1.1,0.7);
\draw (-2.99,1.75)-- (-1.37,1.51);
\draw (-2.19,0.42)-- (-0.83,0.61);
\draw (-1.43,-0.84)-- (-0.32,-0.25);
\draw (-0.79,-1.91)-- (0.12,-0.97);
\draw (0.07,-0.71)-- (0.96,-1.79);
\draw (1.36,-1.2)-- (0.38,-0.25);
\draw (0.68,0.21)-- (1.75,-0.6);
\draw (0.96,0.62)-- (2.1,-0.07);
\draw [line width=0.4pt,dash pattern=on 1pt off 1pt] (0.97,-0.82)-- (-3,3);
\draw [line width=0.4pt,dash pattern=on 1pt off 1pt] (-1.22,0.55)-- (2,1);
\draw [line width=0.4pt,dash pattern=on 1pt off 1pt] (0.15,1.3)-- (0,-2);
\begin{scriptsize}
\fill [color=black] (2,1) circle (0.5pt);
\draw[color=black] (2.15,1.04) node {$P_3$};
\fill [color=black] (0,-2) circle (0.5pt);
\draw[color=black] (0.02,-2.17) node {$P_1$};
\fill [color=black] (-3,3) circle (0.5pt);
\draw[color=black] (-3.1,3.12) node {$P_2$};
\draw[color=black] (0.5,-0.1) node {$L_3$};
\draw[color=black] (-0.9,0.4) node {$L_1$};
\draw[color=black] (-0.0,1.2) node {$L_2$};
\end{scriptsize}
\end{tikzpicture}
\caption{The right-spherical billiard based at $P_1,P_2,P_3$ with each one of its fields of projective lines $L_1,L_2,L_3$}
\label{fig:right_spherical}
\end{figure}
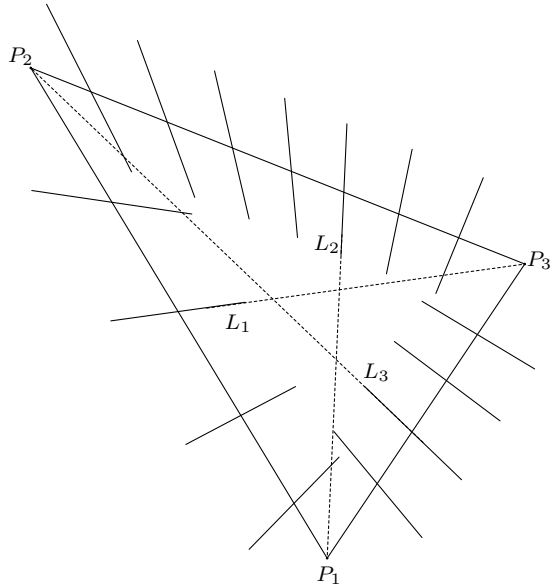

\begin{definition}
\label{definition:right_spherical_billiard}
Given three points $P_1,P_2,P_3$ not on the same line, the \textit{right-spherical billiard} based at $P_1,P_2,P_3$ is the local projective billiard $(\linep{P_1P_2}{P_3},\linep{P_2P_3}{P_1},\linep{P_3P_1}{P_2})$. See Figure \ref{fig:right_spherical}.
\end{definition}

\begin{figure}[!h]
\centering
\begin{tikzpicture}[line cap=round,line join=round,>=triangle 45,x=1.3cm,y=1.3cm]
\clip(-3.7,-2.5) rectangle (3.7,3.5);
\draw (-2,-2)-- (-3,1);
\draw (-3,1)-- (2,3);
\draw (2,3)-- (3,-2);
\draw (-2,-2)-- (3,-2);
\draw [dotted] (-2,-2)-- (2,3);
\draw [dotted] (-3,1)-- (3,-2);
\draw [dash pattern=on 1pt off 1pt] (-0.57,-0.21)-- (-2.94,-0.84);
\draw [dash pattern=on 1pt off 1pt] (-0.57,-0.21)-- (-0.48,2.46);
\draw [dash pattern=on 1pt off 1pt] (-0.57,-0.21)-- (2.95,0.71);
\draw [dash pattern=on 1pt off 1pt] (-0.57,-0.21)-- (0.4,-2.46);
\draw (-3.39,0.53)-- (-2.2,0.21);
\draw (-3.17,-0.15)-- (-2.06,-0.18);
\draw (-2.94,-0.84)-- (-1.93,-0.58);
\draw (-2.68,-1.62)-- (-1.78,-1.02);
\draw (-1.5,-2.46)-- (-1.11,-1.51);
\draw (-0.57,-1.51)-- (-0.56,-2.46);
\draw (-0.01,-1.51)-- (0.4,-2.46);
\draw (1.25,-2.46)-- (0.48,-1.51);
\draw (1.06,-1.51)-- (2.25,-2.46);
\draw (3.33,-1.23)-- (2.3,-0.96);
\draw (3.2,-0.55)-- (2.2,-0.46);
\draw (3.07,0.1)-- (2.1,0.01);
\draw (2.01,0.46)-- (2.95,0.71);
\draw (2.82,1.34)-- (1.92,0.93);
\draw (1.82,1.44)-- (2.68,2.04);
\draw (2.55,2.68)-- (1.72,1.92);
\draw (1.36,3.19)-- (0.7,2.03);
\draw (0.53,2.86)-- (0.15,1.81);
\draw (-0.48,2.46)-- (-0.51,1.54);
\draw (-0.96,1.37)-- (-1.16,2.19);
\draw (-1.96,1.87)-- (-1.48,1.16);
\draw (-1.95,0.97)-- (-2.65,1.58);
\begin{scriptsize}
\fill [color=black] (-2,-2) circle (1.5pt);
\draw[color=black] (-2.1,-2.13) node {$P_4$};
\fill [color=black] (-3,1) circle (1.5pt);
\draw[color=black] (-3.2,1.09) node {$P_1$};
\fill [color=black] (2,3) circle (1.5pt);
\draw[color=black] (2.12,3.2) node {$P_2$};
\fill [color=black] (3,-2) circle (1.5pt);
\draw[color=black] (3.19,-2.05) node {$P_3$};
\fill [color=black] (-0.57,-0.21) circle (1.5pt);
\draw[color=black] (-0.7,0.04) node {$O$};
\draw[color=black] (-1.69,-0.34) node {$L_4$};
\draw[color=black] (-0.33,1.17) node {$L_1$};
\draw[color=black] (1.88,0.26) node {$L_2$};
\draw[color=black] (-0.24,-1.35) node {$L_3$};
\end{scriptsize}
\end{tikzpicture}
\caption{The centrally-projective quadrilateral based at $O,P_1,P_2,P_3,P_4$ with each one of its fields of transverse lines $L_1,L_2,L_3,L_4$}
\label{fig:centrally_projective_quadrilateral}
\end{figure}
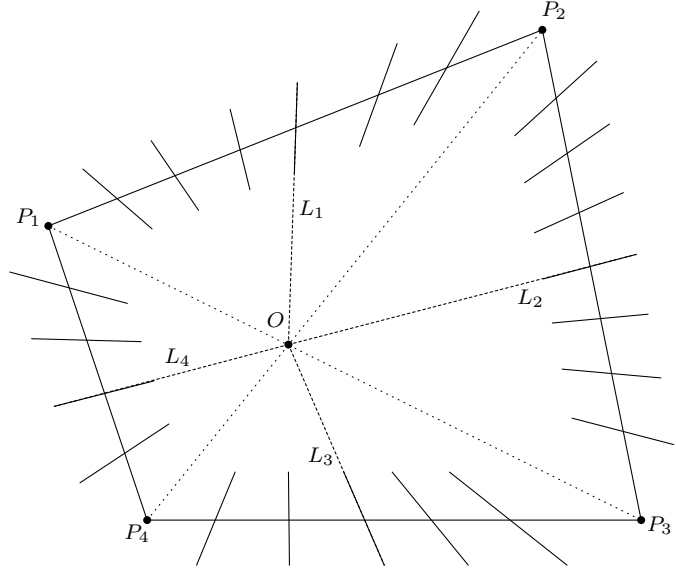

\begin{definition}
\label{definition:centrally_projective_polygon}
Let $n+1$ points $O,P_1,P_2,\ldots,P_n$ be such that for each $j$ modulo $n$, $P_j$, $P_{j+1}$ and $O$ are not on the same line. The \textit{centrally-projective polygon} based at $O,P_1,\ldots,P_n$ is the local projective billiard $(\linep{P_1P_2}{O},\linep{P_2P_3}{O},\ldots, \linep{P_nP_1}{O})$. When $n=4$ we will say \textit{quadrilateral} instead of polygon. See Figure \ref{fig:centrally_projective_quadrilateral}.
\end{definition}

We show that the right-spherical billiard is $3$-reflective and that specific centrally-projective polygons with $n$-vertices are $k$-reflective, where the different cases for $n$ and $k$ are summarized in the following proposition:

\begin{proposition}
\label{proposition:examples_projective_billiards}
Let $n,k\geq 3$ be integers. The following local projective billiards with $n$ vertices are $k$-reflective:\\
\enum $n=k=3$: the right-spherical billiard based at any points $P_1,P_2,P_3$ is $3$-reflective;\\
\enum $n=k=4$: the centrally-projective quadrilateral based at any points $O,P_1,P_2,P_3,P_4$ is $4$-reflective, where $O$ is the intersection point of $P_1P_3$ and $P_2P_4$;\\
\enum $n=k$ even: the centrally-projective polygon based at points $O,P_1,P_2,\ldots,P_n$ is $n$ reflective, where $n\geq 4$ is an even integer, $P_1,\ldots,P_n$ enumerate the vertices of a regular polygon in the clockwise (or counterclockwise) order and $O$ is its center of symmetry;\\
\enum $k=2n$: the centrally-projective polygon based at any points $O,P_1,P_2,\ldots,P_n$ is $2n$-reflective, where $n\geq 3$ is an integer.
\end{proposition}

\begin{remark}
There is another class of $4$-reflective projective billiards inside quadrilaterals which can be constructed using two right-spherical billiards and "gluing" them together. We do not give details about this construction and we refer the interested reader to \cite{fierobe1}.
\end{remark}

The next subsections will be devoted to the proof of Proposition \ref{proposition:examples_projective_billiards}. In the proof, we will consider \textit{virtual orbits}, which are the same as orbits defined in the introductive section without the statement about orientation of lines. For virtual orbits, the reflection at each point can cross the boundary. Notice also that, as for usual orbits, a side $p_jp_{j+1}$ of a virtual orbit can cross another boundary, $a_i$ with $i\neq j,j+1$, of the billiard without being reflected by it. 

In certain cases, for example when the classical boundaries $a_j$ are lines, the virtual orbit $(p_j)_j$ of $(p_1,p_2)\in a_1\times a_2$ is uniquely defined. We say that $(p_1,p_2)$ \textit{determines} the orbit $(p_j)_j$.

\begin{remark} 
\label{remark:centrally_proj_polygon}
Projective billiards in centrally-projective polygons are strongly related to dual billiards in polygons, as explained in \cite{taba_projectif}. A \textit{dual billiard} (or \textit{outer billiard}) \cite{shaidenko_vivaldi, taba_dual_billiards} is an oriented closed convex curve $\gamma$ together with a map $\varphi$ defined on the exterior of the curve as follows: given a point $p$ outside the closed domain bounded by the curve, there are two tangent lines to $\gamma$ containing $p$ which we can orient from $p$ to any point of tangency. Each of them is tangent to $\gamma$ either at a unique point, or along a segment (convexity). Choose the so-called right tangent line, which has the same orientation, as $\gamma$, at their tangency point(s). We deal only with those points $p$ for which the corresponding tangency point (denoted by $q$) is unique. In the case, when $\gamma$ is a polygon, the point $q$ is its vertex. In this case the orientation condition should be modified as follows. Turn the oriented line $pq$ around the point $q$ until it becomes tangent to $\gamma$ along a segment adjacent to $q$ (either clockwise, or counterclockwise). Then the rotated line and the latter segment should have the same orientation. Set $\varphi(p)$ to be the point obtained by reflecting $p$ with respect to $q$.

We would like to thank Sergei Tabachnikov who pointed us out the following results: in a certain class of polygons, called \textit{rational polygons}, the outer orbits are always finite (see \cite{taba_book} Chapt. 9, or \cite{shaidenko_vivaldi} for a proof). Rational polygons are polygons whose vertices lie on the affine image of a lattice. For example triangles and parallelograms are rational polygons. They have the property that the outer orbit $(p_j)_j$ of a given point $p_0$ is discrete (indeed, the vectors joining for each $j$ the points $p_j$ and $p_{j+2}$ are in a lattice, as it can be deduced from Lemma \ref{lemma:centproj_thales_corollary}). In \cite{shaidenko_vivaldi}, it is shown that orbits in rational polygons are also bounded, which proves their finiteness.

Results on dual billiards are of great interest for projective billiards endowed with a so-called centrally-projective field of lines (meaning that the transverse lines to the boundary of the billiard contain the same point $O$). Indeed, such projective billiards are conjugated by polarity with dual billiards, see \cite{taba_projectif}. Let us describe this construction for polygons: suppose we are given a centrally-projective polygon based at points $O,P_1,P_2,\ldots,P_n$ and an orbit $(p_j)_{j\in\ZZ}$ of the corresponding projective billiard. In our case, this means that each $p_j$ lies on the side $P_jP_{j+1}$ (see Definition \ref{definition:local_projective_billiard}). Consider a polarity such that the point $O$ is the pole of the line at infinity (see Section \ref{section:general_properties_on_quadrics}). For each $k$ denote by $Q_k$ the pole of the line $P_kP_{k+1}$ and by $q_k$ the pole of $p_kp_{k+1}$: then the line $Op_k$ has its pole $\omega_k$ at infinity, and the points $q_{k-1}$, $q_{k}$, $Q_k$, $\omega_k$ form a harmonic quadruple of points (since they are poles of lines in a harmonic set of lines). Therefore $Q_k$ is the midpoint of the segment $q_kq_{k+1}$ and we recover the dynamics of a dual billiard outside the polygon $Q_1\cdots Q_n$ where $(q_k)_k$ is an orbit. The dynamics in this case is more simple since by construction the midpoints of successive edges are \textit{consecutive} vertices of the polygon. Hence both projective and dual billiards are conjugated by polarity.

This link between centrally-projective polygons and dual billiards about polygons and the finiteness of orbits in rational polygons immediately implies the following result: projective billiards in centrally-projective polygons which are duals of rational polygons have only periodic orbits. It could be interesting to describe this new class of centrally-projective polygons (the centrally-projective quadrilateral of Proposition \ref{proposition:examples_projective_billiards}, case $n=k=4$, is an example of such polygon, since its associated dual billiard is a parallelogram, as it will be explained below).
\end{remark}

\subsection{$3$-reflectivity of the right-spherical billiard}

In this subsection, we denote by $P_1,P_2,P_3$ three non-colinear points of $\RP^2$ and we consider the right-spherical billiard based at $P_1,P_2,P_3$. Given an integer $j$ modulo $3$ and $p\in P_jP_{j+1}$, we denote by $L_j(p)$ the projective line at $p$ of $\linep{P_jP_{j+1}}{P_{j+2}}$, that is $L_j(p)=pP_{j+2}$, the line containing $p$ and $P_{j+2}$.

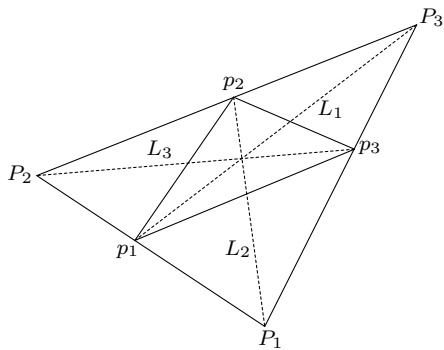
\begin{figure}[!h]
\centering
\begin{tikzpicture}[line cap=round,line join=round,>=triangle 45,x=1cm,y=1cm]
\clip(-3.5,-2.5) rectangle (3,2.5);
\draw (-3,0)-- (0,-2);
\draw (0,-2)-- (2,2);
\draw (2,2)-- (-3,0);
\draw (-1.71,-0.86)-- (-0.41,1.04);
\draw (-0.41,1.04)-- (1.18,0.35);
\draw [line width=0.4pt,dash pattern=on 1pt off 1pt] (-1.71,-0.86)-- (2,2);
\draw [line width=0.4pt,dash pattern=on 1pt off 1pt] (-0.41,1.04)-- (0,-2);
\draw [line width=0.4pt,dash pattern=on 1pt off 1pt] (1.18,0.35)-- (-3,0);
\draw (-1.71,-0.86)-- (1.18,0.35);
\begin{scriptsize}
\draw[color=black] (-3.22,0.02) node {$P_2$};
\draw[color=black] (0.08,-2.17) node {$P_1$};
\draw[color=black] (2.2,2.1) node {$P_3$};
\draw[color=black] (-1.8,-1) node {$p_1$};
\draw[color=black] (-0.4,1.22) node {$p_2$};
\draw[color=black] (1.39,0.37) node {$p_3$};
\draw[color=black] (0.87,0.88) node {$L_1$};
\draw[color=black] (-0.35,-0.97) node {$L_2$};
\draw[color=black] (-1.38,0.35) node {$L_3$};
\end{scriptsize}
\end{tikzpicture}
\caption{The right-spherical billiard based at $P_1,P_2,P_3$ and a triangular orbit $(p_1,p_2,p_3)$ obtained by reflecting any segment $p_1p_2$ two times}
\label{fig:right_spherical_orbit}
\end{figure}

\begin{proposition}
\label{prop:right_spherical_3_reflective}
Any $(p_1,p_2)\in P_1P_2\times P_2P_3$ with $p_1\neq p_2$ determines a $3$-periodic orbit of the right-spherical billiard based at $P_1,P_2,P_3$. See Figure \ref{fig:right_spherical_orbit}.
\end{proposition}

\begin{figure}[!h]
\centering
\begin{tikzpicture}[line cap=round,line join=round,>=triangle 45,x=1.4cm,y=1.4cm]
\clip(-3.5,-2.5) rectangle (3.5,2.5);
\draw (-3,0)-- (0,-1.89);
\draw (0,-1.89)-- (1.93,0.92);
\draw (1.93,0.92)-- (-3,0);
\draw (-1.71,-0.81)-- (0.6,0.67);
\draw (0.6,0.67)-- (1.06,-0.36);
\draw [line width=0.4pt,dash pattern=on 1pt off 1pt] (-1.71,-0.81)-- (1.93,0.92);
\draw [line width=0.4pt,dash pattern=on 1pt off 1pt] (0.6,0.67)-- (0,-1.89);
\draw [dotted,domain=-4.91:4.59] plot(\x,{(--0.66--1.48*\x)/2.31});
\draw [dotted,domain=-4.91:4.59] plot(\x,{(--3.66-2.81*\x)/-1.93});
\draw (-1.71,-0.81)-- (0.78,-0.77);
\begin{scriptsize}
\fill [color=black] (-3,0) circle (0.5pt);
\draw[color=black] (-3.15,0.02) node {$P_2$};
\fill [color=black] (0,-1.89) circle (2pt);
\draw[color=black] (0.15,-1.9) node {$P_1$};
\fill [color=black] (1.93,0.92) circle (2pt);
\draw[color=black] (2.12,0.95) node {$P_3$};
\fill [color=black] (-1.71,-0.81) circle (0.5pt);
\draw[color=black] (-1.74,-0.89) node {$p_1$};
\fill [color=black] (0.6,0.67) circle (0.5pt);
\draw[color=black] (0.6,0.79) node {$p_2$};
\draw [color=black] (1.06,-0.36) circle (2pt);
\draw[color=black] (1.25,-0.35) node {$p_3$};
\draw[color=black] (-0.15,-0.17) node {$L_1$};
\draw[color=black] (0.31,-1) node {$L_2$};
\fill [color=black] (2.69,2.01) circle (2pt);
\draw[color=black] (2.78,1.98) node {$A$};
\draw [color=black] (0.78,-0.77) circle (2pt);
\draw[color=black] (1,-0.79) node {$p_3'$};
\end{scriptsize}
\end{tikzpicture}
\caption{As in the proof of Proposition \ref{prop:right_spherical_3_reflective}, both quadruples of points $(p_2,A,P_1,P_3)$ and $(p_2',A,P_1,P_2)$ are harmonic, hence necessarily $p_2=p_2'$.}
\label{fig:right_spherical_proof}
\end{figure}
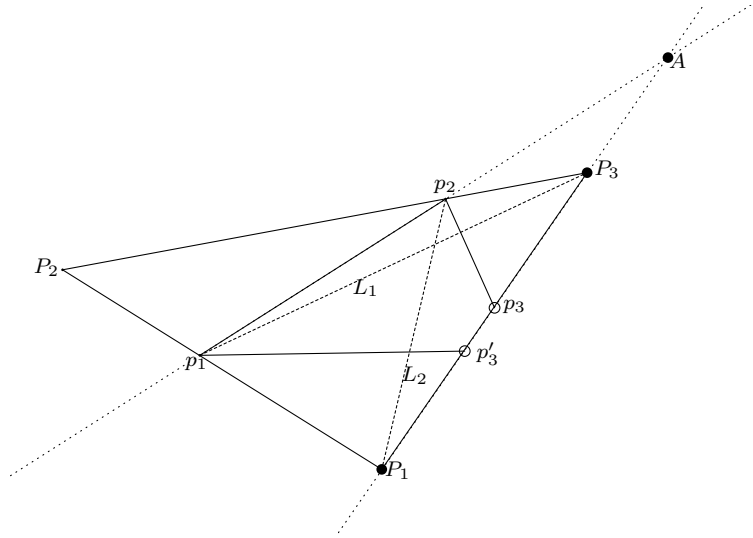

We give two proofs of Proposition \ref{prop:right_spherical_3_reflective}. The first one is based on the observation that the right-spherical billiard is obtained by the projection of a $3$-reflective billiard on the sphere $\SS^2$ to the Euclidean plane (see Section \ref{subsection:projective_and_others}). The other proof uses more intrinsic arguments about harmonic sets of lines.

\begin{proof}[Proof 1] Projective transformations do not change the cross-ratio of four distinct points on the same line, hence if four lines form a harmonic set then their images by a projective transformation also form a harmonic set. Therefore it is enough to show that at least one example of right-spherical billiard is $3$-reflective. As explained in the introductive section, a triangle on the sphere $\SS^2$ having only right angles is a $3$-reflective classical billiard on the sphere $\SS^2$ (see Figure \ref{figure:sphere_bary_en} in the introductive section). This example was found by Baryshnikov \cite{bary_introuvable,VKNZ}. Choose such triangle on the upper open hemisphere of $\SS^2$, and project it on the plane $\mathcal P=\{z=-1\}$ by a projection with respect to the center of $\SS^2$ (see Section \ref{subsection:projective_and_others} for more details). Endow $\mathcal P$ with the metric $g$ obtained by pushing forward the spherical metric with this projection: the geodesics of the Riemannian manifold $(\mathcal P, g)$ are lines. By construction, we obtain a $3$-reflective billiard whose boundary is a triangle $P_1P_2P_3$, and the $g$-normal line to $P_jP_{j+1}$ at any point $p$ is the line joining $p$ to the opposite vertex $P_{j+2}$ (see Figure \ref{fig:right_spherical_orbit}). By Proposition \ref{proposition:usual_is_projective}, two lines $\ell$ and $\ell'$ containing $p$ are symetric with respect to $P_jP_{j+1}$ in the metric $g$ if and only if the quadruple of lines $(\ell,\ell',P_jP_{j+1},pP_{j+2})$ is harmonic. Therefore, the orbits of the billiard $P_1P_2P_3$ in the metric $g$ coincide with the orbits of the right-spherical billiard based at $P_1$, $P_2$, $P_3$, and the corresponding right-spherical billiard is $3$-reflective.
\end{proof}

\begin{proof}[Proof 2]
This proof was found by Simon Allais in a talk we had about harmonicity conditions in a projective space. Let $p_3\in P_1P_3$ be such that $p_1p_2$, $p_2p_3$, $P_2P_3$, $L_2(p_2)$ are harmonic lines. Define $p_3'\in P_1P_3$ similarly: $p_1p_2$, $p_1p_3'$, $P_1P_2$, $L_1(p_1)$ are harmonic lines. Let us show that $p_3=p_3'$ (see Figure \ref{fig:right_spherical_proof}). Consider the line $P_1P_3$ and let $A$ be its point of intersection with $p_1p_2$. Let us consider harmonic quadruples of points on $P_1P_3$. By harmonicity of the previous defined lines passing through $p_2$, the quadruple of points $(A,p_3,P_3,P_1)$ is harmonic. Doing the same with the lines passing through $p_1$, the quadruple of points $(A,p_3',P_3,P_1)$ is harmonic. Hence $p_3=p_3'$ since the projective transformation defining the cross-ratio is one to one.

Now let us prove that the lines $p_2p_3$, $p_1p_3$, $P_1P_3$, $L_3(p_3)$ are harmonic lines. Consider the line $P_1P_2$: $p_2p_3$ intersects it at a certain point denoted by $B$, $p_3p_1$ at $p_1$, $P_3P_1$ at $P_1$ and $L_3(p_3)$ at $P_2$. But the quadruple of points $(B,p_1,P_1,P_2)$ is harmonic since there is a reflection law at $p_2$ whose lines intersect $P_1P_2$ exactly in those points.
\end{proof}

\subsection{$4$-reflectivity of the centrally-projective quadrilateral}

In this subsection, we denote by $P_1,P_2,P_3,P_4$ points of $\RP^2$ such that no three of them are colinear, and $O$ the intersection point of the line $P_1P_3$ with $P_2P_4$. We consider the centrally-projective quadrilateral based at $O,P_1,P_2,P_3,P_4$. Given an integer $j$ modulo $4$ and $p\in P_jP_{j+1}$, we denote by $L_j(p)$ the projective line at $p$ of $\linep{P_jP_{j+1}}{O}$, that is $L_j(p)=Op$, the line containing $p$ and $O$.

\begin{proposition}
\label{prop:centproj_quadrilateral_4_reflective}
Any $(p_1,p_2)\in P_1P_2\times P_2P_3$ with $p_1\neq p_2$ determines a $4$-periodic orbit of the centrally-projective polygon based at $O,P_1,P_2,P_3,P_4$.
\end{proposition}

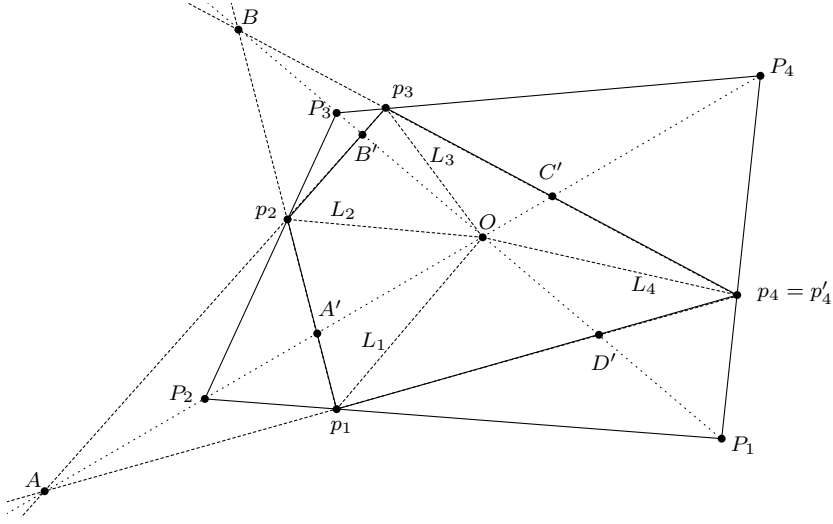
\begin{figure}[!h]
\centering
\begin{tikzpicture}[line cap=round,line join=round,>=triangle 45,x=1.7cm,y=1.7cm]
\clip(-3.8,-1.2) rectangle (3.5,2.8);
\draw (-2.27,-0.28)-- (-1.25,1.95);
\draw (-1.25,1.95)-- (2.03,2.24);
\draw (2.03,2.24)-- (1.73,-0.59);
\draw (1.73,-0.59)-- (-2.27,-0.28);
\draw [dash pattern=on 1pt off 1pt,domain=-4.36118736911263:-1.2546504345620575] plot(\x,{(--1.99--1.48*\x)/-0.38});
\draw [dash pattern=on 1pt off 1pt,domain=-4.36118736911263:-0.8682629057033406] plot(\x,{(-2.26-0.86*\x)/-0.76});
\draw [dash pattern=on 1pt off 1pt,domain=-4.36118736911263:1.850326030682974] plot(\x,{(--0.01-0.88*\x)/-3.1});
\draw [dash pattern=on 1pt off 1pt,domain=-4.36118736911263:1.850326030682974] plot(\x,{(-4.13--1.46*\x)/-2.72});
\draw [dotted,domain=-4.36118736911263:1.732822661578601] plot(\x,{(-2.64--2.55*\x)/-2.98});
\draw [dotted,domain=-4.36118736911263:2.030297933160431] plot(\x,{(-4.53-2.52*\x)/-4.3});
\draw (-0.87,1.99)-- (1.85,0.53);
\draw (1.85,0.53)-- (-1.25,-0.36);
\draw (-1.63,1.12)-- (-1.25,-0.36);
\draw (-1.63,1.12)-- (-0.87,1.99);
\draw [dash pattern=on 1pt off 1pt] (-0.12,0.98)-- (-1.63,1.12);
\draw [dash pattern=on 1pt off 1pt] (-0.12,0.98)-- (-0.87,1.99);
\draw [dash pattern=on 1pt off 1pt] (-0.12,0.98)-- (-1.25,-0.36);
\draw [dash pattern=on 1pt off 1pt] (-0.12,0.98)-- (1.85,0.53);
\begin{scriptsize}
\fill [color=black] (-2.27,-0.28) circle (1.5pt);
\draw[color=black] (-2.45,-0.24) node {$P_2$};
\fill [color=black] (-1.25,1.95) circle (1.5pt);
\draw[color=black] (-1.38,1.98) node {$P_3$};
\fill [color=black] (2.03,2.24) circle (1.5pt);
\draw[color=black] (2.2,2.31) node {$P_4$};
\fill [color=black] (1.73,-0.59) circle (1.5pt);
\draw[color=black] (1.9,-0.64) node {$P_1$};
\fill [color=black] (-0.12,0.98) circle (1.5pt);
\draw[color=black] (-0.08,1.1) node {$O$};
\fill [color=black] (-1.25,-0.36) circle (1.5pt);
\draw[color=black] (-1.21,-0.5) node {$p_1$};
\fill [color=black] (-1.63,1.12) circle (1.5pt);
\draw[color=black] (-1.79,1.17) node {$p_2$};
\fill [color=black] (-0.87,1.99) circle (1.5pt);
\draw[color=black] (-0.73,2.12) node {$p_3$};
\fill [color=black] (-3.51,-1) circle (1.5pt);
\draw[color=black] (-3.6,-0.91) node {$A$};
\fill [color=black] (-2.01,2.6) circle (1.5pt);
\draw[color=black] (-1.92,2.7) node {$B$};
\fill [color=black] (1.85,0.53) circle (1.5pt);
\draw[color=black] (2.3,0.55) node {$p_4=p_4'$};
\fill [color=black] (0.42,1.3) circle (1.5pt);
\draw[color=black] (0.41,1.5) node {$C'$};
\fill [color=black] (0.78,0.22) circle (1.5pt);
\draw[color=black] (0.82,0.01) node {$D'$};
\fill [color=black] (-1.4,0.23) circle (1.5pt);
\draw[color=black] (-1.3,0.43) node {$A'$};
\fill [color=black] (-1.05,1.78) circle (1.5pt);
\draw[color=black] (-1.02,1.64) node {$B'$};
\draw[color=black] (-1.2,1.19) node {$L_2$};
\draw[color=black] (-0.43,1.61) node {$L_3$};
\draw[color=black] (-0.96,0.17) node {$L_1$};
\draw[color=black] (1.13,0.61) node {$L_4$};
\end{scriptsize}
\end{tikzpicture}
\caption{The centrally-projective quadrilateral based at $O,P_1,P_2,P_3,P_4$ with a periodic orbit obtained by reflecting $p_1p_2$ three times. Here the notations are the same as in the proof of Proposition \ref{prop:centproj_quadrilateral_4_reflective}.}
\label{fig:centrally_projective_quadrilateral_proof}
\end{figure}

\begin{proof}
Let $p_3\in P_3P_4$ such that $p_1p_2$ is reflected into $p_2p_3$ by the reflection law at $p_2$. Let $p_4\in P_4P_1$ such that $p_2p_3$ is reflected into $p_3p_4$ by the reflection law at $p_3$. Let $p_4'\in P_4P_1$ such that $p_1p_2$ is reflected into $p_1p_4'$ by the reflection law at $p_1$. Denote by $d$ the line reflected from $p_3p_4$ by the projective reflection law at $p_4$. We have to show that $d=p_1p_4'$. 

First, let us introduce a few notations (see Figure \ref{fig:centrally_projective_quadrilateral_proof}). Consider the line $p_1p_2$; it intersects: the line $P_1P_3$ at a point $B$ and the line $P_2P_4$ at a point $A'$. Now consider the line $p_2p_3$; it intersects: the line $P_1P_3$ at a point $B'$ and the line $P_2P_4$ at a point $A$. Finally let $C'$ be the intersection point of $p_3p_4$ with $P_2P_4$ and $D'$ the intersection point of $p_1p_4'$ with $P_1P_3$.

Then, notice that by the projective law of reflection at $p_2$, the quadruple of points $(A, A',P_2,O)$ is harmonic. Since the points $P_2,A',O$ correpond to the lines $P_1P_2$, $p_1p_2$, $L_1(p_1)$, the previously defined reflected line $p_1p_4'$ needs to pass through $A$ in order to form a harmonic quadruple of lines. The same remark on the other diagonal leads to note that $p_3p_4$ passes through $B$.

Now by the reflection law at $p_3$, one observe that the quadruple of points $(A,C',O,P_4)$ is harmonic. But $P_4P_1$ passes through $P_4$, $p_4p_3$ through $C'$ and $L_4(p_4)$ through $O$. Hence $d$ needs to pass through $A$. Then, by the reflection law at $p_1$, one observe that the quadruple of points $(B,D',O,P_1)$ is harmonic. But $P_4P_1$ passes through $P_1$, $p_4p_3$ through $B$ and $L_4(p_4)$ through $O$. Hence $d$ needs to pass through $D'$.

Therefore we conclude that $d=AD'=p_1p_4'$.
\end{proof}

\begin{remark}
Another proof can be given by duality: as explained at Remark \ref{remark:centrally_proj_polygon}, we can associate a dual billiard to the centrally-projective quadrilateral $P_1P_2P_3P_4$ of Proposition \ref{fig:centrally_projective_quadrilateral_proof} by a polarity sending $O$ at infinity. Since the point $O$ is on both its diagonals $P_1P_3$ and $P_2P_4$, the dual polygon to $P_1P_2P_3P_4$ is a parallelogram $Q_1Q_2Q_3Q_4$ (hence a rational polygon). The study of the simplified dual billiard outside $Q_1Q_2Q_3Q_4$ (as described in  Remark \ref{remark:centrally_proj_polygon}) gives another proof of Proposition \ref{fig:centrally_projective_quadrilateral_proof} as a simple consequence of the famous intercept theorem in geometry.
\end{remark}

\subsection{$2m$-reflectivity of centrally-projective regular $2m$-sided polygons}

Let $n=2m\geq 4$ be an even integer, $P_1,\ldots,P_n$ be a clockwise enumeration of the vertices of a \textit{regular} polygon, and $O$ be the intersection point of its great diagonals (that is the point of intersection of the lines $P_jP_{j+k}$ where $j$ is an integer taken modulo $n$). We consider the centrally-projective polygon based at $O,P_1,\ldots,P_n$. 

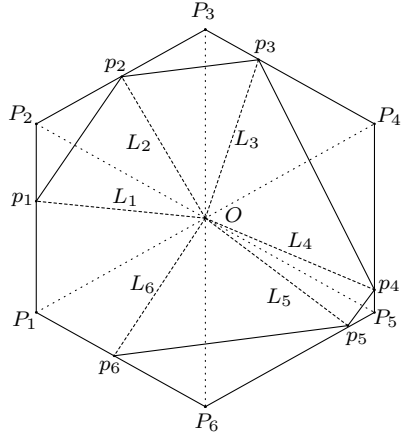
\begin{figure}[!h]
\centering
\begin{tikzpicture}[line cap=round,line join=round,>=triangle 45,x=2.5cm,y=2.5cm]
\clip(-2.27,-1.21) rectangle (1.99,1.16);
\draw (-0.89,0.5)-- (0,1);
\draw (0.89,0.5)-- (0,1);
\draw (-0.89,-0.5)-- (0,-1);
\draw (0.89,-0.5)-- (0,-1);
\draw (-0.89,0.5)-- (-0.89,-0.5);
\draw (0.89,0.5)-- (0.89,-0.5);
\draw [dotted] (0,1)-- (0,-1);
\draw [dotted] (-0.89,-0.5)-- (0.89,0.5);
\draw [dotted] (-0.89,0.5)-- (0.89,-0.5);
\draw (-0.89,0.09)-- (-0.44,0.75);
\draw (-0.44,0.75)-- (0.28,0.84);
\draw (0.28,0.84)-- (0.89,-0.38);
\draw (0.89,-0.38)-- (0.75,-0.57);
\draw (0.75,-0.57)-- (-0.48,-0.73);
\draw [dash pattern=on 1pt off 1pt] (0,0)-- (-0.44,0.75);
\draw [dash pattern=on 1pt off 1pt] (0,0)-- (0.28,0.84);
\draw [dash pattern=on 1pt off 1pt] (0,0)-- (0.89,-0.38);
\draw [dash pattern=on 1pt off 1pt] (0,0)-- (0.75,-0.57);
\draw [dash pattern=on 1pt off 1pt] (0,0)-- (-0.48,-0.73);
\draw [dash pattern=on 1pt off 1pt] (0,0)-- (-0.89,0.09);
\begin{scriptsize}
\fill [color=black] (-0.89,0.5) circle (0.5pt);
\draw[color=black] (-0.97,0.55) node {$P_2$};
\fill [color=black] (0,1) circle (0.5pt);
\draw[color=black] (-0.01,1.1) node {$P_3$};
\fill [color=black] (0.89,0.5) circle (0.5pt);
\draw[color=black] (0.97,0.54) node {$P_4$};
\fill [color=black] (0,1) circle (0.5pt);
\fill [color=black] (-0.89,-0.5) circle (0.5pt);
\draw[color=black] (-0.95,-0.54) node {$P_1$};
\fill [color=black] (0,-1) circle (0.5pt);
\draw[color=black] (0.01,-1.08) node {$P_6$};
\fill [color=black] (0.89,-0.5) circle (0.5pt);
\draw[color=black] (0.95,-0.53) node {$P_5$};
\fill [color=black] (0,-1) circle (0.5pt);
\fill [color=black] (0,0) circle (0.5pt);
\draw[color=black] (0.15,0.02) node {$O$};
\fill [color=black] (-0.89,0.09) circle (0.5pt);
\draw[color=black] (-0.96,0.1) node {$p_1$};
\fill [color=black] (-0.44,0.75) circle (0.5pt);
\draw[color=black] (-0.47,0.81) node {$p_2$};
\fill [color=black] (0.28,0.84) circle (0.5pt);
\draw[color=black] (0.32,0.91) node {$p_3$};
\fill [color=black] (0.89,-0.38) circle (0.5pt);
\draw[color=black] (0.97,-0.36) node {$p_4$};
\fill [color=black] (0.75,-0.57) circle (0.5pt);
\draw[color=black] (0.8,-0.64) node {$p_5$};
\fill [color=black] (-0.48,-0.73) circle (0.5pt);
\draw[color=black] (-0.5,-0.79) node {$p_6$};
\draw[color=black] (-0.35,0.4) node {$L_2$};
\draw[color=black] (0.22,0.42) node {$L_3$};
\draw[color=black] (0.5,-0.12) node {$L_4$};
\draw[color=black] (0.39,-0.42) node {$L_5$};
\draw[color=black] (-0.33,-0.35) node {$L_6$};
\draw[color=black] (-0.42,0.11) node {$L_1$};
\end{scriptsize}
\end{tikzpicture}
\caption{A centrally-projective regular  polygon based at $O,P_1,\ldots,P_{6}$ and a piece of trajectory after four projective reflections.}
\label{fig:centrally_projective_polygon}
\end{figure}

\begin{proposition}
\label{prop:centproj_even_polygon_reflective}
Any $(p_1,p_2)\in P_1P_2\times P_2P_3$ with $p_1\neq p_2$ determines an $n$-periodic orbit of the centrally-projective regular polygon based at $O,P_1,\ldots,P_n$. See Figure \ref{fig:centrally_projective_polygon}. 
\end{proposition}

\begin{proof}
Fix $(p_1,p_2)\in P_1P_2\times P_2P_3$ with $p_1\neq p_2$ and consider its backward and forward orbit $p=(p_j)_{j\in\ZZ}$. During the proof, all indices, except for $p_j$, will be considered modulo $n$. We first prove the following

\begin{lemma}
\label{lemma:centproj_reflective_dtes_concourantes}
Fix an integer $\ell$ and consider the great diagonal $\Delta_{\ell} =P_{\ell}P_{\ell+m}$. Then for any $r\geq 0$, the lines $p_{\ell-r-2}p_{\ell-r-1}$ and $p_{\ell+r}p_{\ell+r+1}$ intersect $\Delta_{\ell}$ at the same point. See Figure \ref{fig:centrally_projective_polygon_proof}.
\end{lemma}

\begin{figure}[!h]
\centering
\begin{tikzpicture}[line cap=round,line join=round,>=triangle 45,x=2.5cm,y=2.5cm]
\clip(-2.4,-0.65) rectangle (2.06,1.6);
\draw [line width=1.2pt] (0,0.8)-- (-0.8,0.6);
\draw [line width=1.2pt] (-0.8,0.6)-- (-1.4,0.2);
\draw [line width=1.2pt] (-1.4,0.2)-- (-1.8,-0.4);
\draw [line width=1.2pt] (0,0.8)-- (0.8,0.6);
\draw [line width=1.2pt] (0.8,0.6)-- (1.4,0.2);
\draw [line width=1.2pt] (1.4,0.2)-- (1.8,-0.4);
\draw [dash pattern=on 1pt off 1pt,domain=-1.0812456489024758:2.055306508092342] plot(\x,{(--0.62--0.3*\x)/0.72});
\draw [dash pattern=on 1pt off 1pt,domain=-1.5015993354037684:2.055306508092342] plot(\x,{(--0.57--0.36*\x)/0.42});
\draw [dash pattern=on 1pt off 1pt,domain=-2.0013299841538337:1.0812456489024758] plot(\x,{(-0.62--0.3*\x)/-0.72});
\draw [dash pattern=on 1pt off 1pt,domain=-2.0013299841538337:1.5015993354037684] plot(\x,{(-0.57--0.36*\x)/-0.42});
\draw [dotted] (0,-0.6) -- (0,1.6);
\draw (-1.5,0.05)-- (-1.08,0.41);
\draw (-1.08,0.41)-- (-0.36,0.71);
\draw (0.36,0.71)-- (1.08,0.41);
\draw (1.08,0.41)-- (1.5,0.05);
\draw (-0.36,0.71)-- (0.36,0.71);
\begin{scriptsize}
\fill [color=black] (0,0.8) circle (1.5pt);
\draw[color=black] (-0.03,0.6) node {$P_{\ell}$};
\fill [color=black] (-0.8,0.6) circle (1.5pt);
\draw[color=black] (-0.86,0.68) node {$P_{\ell-r}$};
\fill [color=black] (-1.4,0.2) circle (1.5pt);
\draw[color=black] (-1.6,0.27) node {$P_{\ell-r-1}$};
\fill [color=black] (-1.8,-0.4) circle (1.5pt);
\draw[color=black] (-2.07,-0.37) node {$P_{\ell-r-2}$};
\fill [color=black] (0.8,0.6) circle (1.5pt);
\draw[color=black] (0.83,0.67) node {$P_{\ell+r}$};
\fill [color=black] (1.4,0.2) circle (1.5pt);
\draw[color=black] (1.6,0.25) node {$P_{\ell+r+1}$};
\fill [color=black] (1.8,-0.4) circle (1.5pt);
\draw[color=black] (1.87,-0.37) node {$P_{\ell+r+2}$};
\fill [color=black] (-1.08,0.41) circle (1.5pt);
\draw[color=black] (-1.3,0.48) node {$p_{\ell-r-1}$};
\fill [color=black] (-0.36,0.71) circle (1.5pt);
\draw[color=black] (-0.37,0.83) node {$p_{\ell-r}$};
\fill [color=black] (-1.5,0.05) circle (1.5pt);
\draw[color=black] (-1.8,0.07) node {$p_{\ell-r-2}$};
\fill [color=black] (1.08,0.41) circle (1.5pt);
\fill [color=black] (0.36,0.71) circle (1.5pt);
\draw[color=black] (0.39,0.83) node {$p_{\ell+r-1}$};
\fill [color=black] (1.5,0.05) circle (1.5pt);
\draw[color=black] (1.8,0.09) node {$p_{\ell+r+1}$};
\fill [color=black] (1.08,0.41) circle (1.5pt);
\draw[color=black] (1.3,0.46) node {$p_{\ell+r}$};
\fill [color=black] (0,0) circle (0.5pt);
\draw[color=black] (0.15,0) node {$L_{\ell}$};
\end{scriptsize}
\end{tikzpicture}
\caption{As in the proof of Lemma \ref{lemma:centproj_reflective_dtes_concourantes}, since the lines $p_{\ell-r-1}p_{\ell-r}$ and $p_{\ell+r-1}p_{\ell+r}$ intersect $L_{\ell}$ at the same point, the lines  $p_{\ell-r-2}p_{\ell-r-1}$ and $p_{\ell+r}p_{\ell+r+1}$ also intersect $L_{\ell}$ at a same point.}
\label{fig:centrally_projective_polygon_proof}
\end{figure}
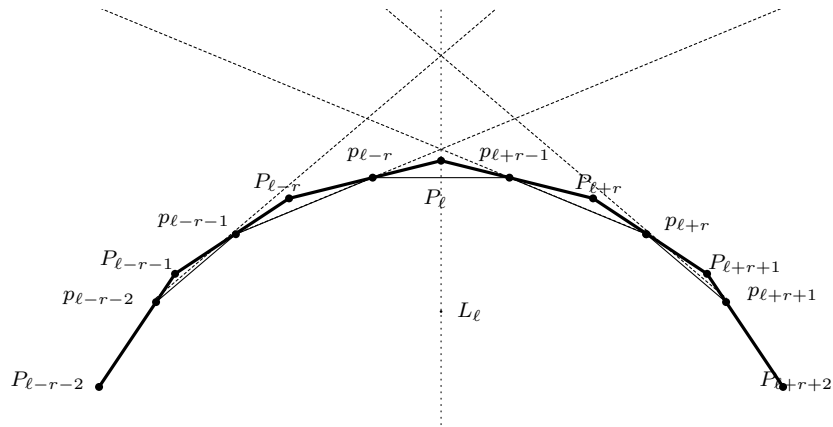

\begin{proof}
Let us prove Lemma \ref{lemma:centproj_reflective_dtes_concourantes} by induction on $r$.

\textsl{Case when $r=0$:} Fix an integer $\ell$. Let $A$ be the intersection point of $p_{\ell-2}p_{\ell-1}$ with $\Delta_{\ell}$, $A'$ the intersection point of $p_{\ell}p_{\ell+1}$ with $\Delta_{\ell}$ and $B$ the intersection point of $p_{\ell-1}p_{\ell}$ with $\Delta_{\ell}$. Consider harmonic quadruples of points on $\Delta_{\ell}$: $(A,B,P_{\ell},O)$ is harmonic by the reflection law in $p_{\ell-1}$, and $(A',B,P_{\ell},O)$ is harmonic by the reflection law in $p_{\ell+1}$. Hence $A=A'$ which concludes the proof for $r=0$.

\textsl{Inductive step:} suppose Lemma \ref{lemma:centproj_reflective_dtes_concourantes} is true for any integer $\ell$ and any $r'<r$ and let us prove it for $r$. See Figure \ref{fig:centrally_projective_polygon_proof} for a detailled drawing of the situation. Fix an integer $\ell\in\ZZ$. By assumption, we know that $p_{\ell-r-1}p_{\ell-r}$ and $p_{\ell+r-1}p_{\ell+r}$ intersect $\Delta_{\ell}$ at the same point $A$. Moreover, by symmetry of the regular polygon with respect to the line $\Delta_{\ell}$, the lines $P_{\ell-r-1}P_{\ell-r}$ and $P_{\ell+r}P_{\ell+r+1}$ intersect $\Delta_{\ell}$ at the same point. Now the intersection points of $p_{\ell-r-1}p_{\ell-r}$ with $p_{\ell+r-1}p_{\ell+r}$, of $P_{\ell-r-1}P_{\ell-r}$ with $P_{\ell+r}P_{\ell+r+1}$, and of $p_{\ell-r-1}O$ with $p_{\ell+r}O$ lie on $\Delta_{\ell}$. Hence in order to sastisfy the projective reflection law at $p_{\ell-r-1}$ and at $p_{\ell+r}$ respectively, the lines $p_{\ell-r-1}p_{\ell-r-2}$ and $p_{\ell+r}p_{\ell+r+1}$ should intersect at the same point. Hence the inductive step is over and this conclude the proof.
\end{proof}

Let us finally prove Proposition \ref{prop:centproj_even_polygon_reflective}. We have to show that $p_{0}p_1 = p_{n}p_{n+1}$. We will use Lemma \ref{lemma:centproj_reflective_dtes_concourantes}. First, by setting $\ell=m+1$ and $r=m-1$, we conclude that the lines
$p_0p_{1}$ and $p_{n}p_{n+1}$ intersect $\Delta_{m+1}=\Delta_1$ at the same point denoted by $A$. Then, by setting $\ell=m+2$ and $r=m-2$ we get that the lines $p_2p_3$ and $p_{n}p_{n+1}$ intersect $\Delta_{m+2}=\Delta_2$ at the same point denoted by $B$. Now it is also true that $p_0p_{1}$ intersects $\Delta_2$ at $B$, by setting $\ell=2$ and $r=0$ in Lemma \ref{lemma:centproj_reflective_dtes_concourantes}. Hence we have shown that $p_{n}p_{n+1}=AB=p_0p_{1}$ which concludes the proof.
\end{proof}

\subsection{$2n$-reflectivity of centrally-projective $n$-sided polygons}

Let $n\geq 3$ be an integer, $O,P_1,\ldots,P_n$ be points in $\RP^2$ such that no three of them are colinear. We consider the centrally-projective polygon based at $O,P_1,\ldots,P_n$.

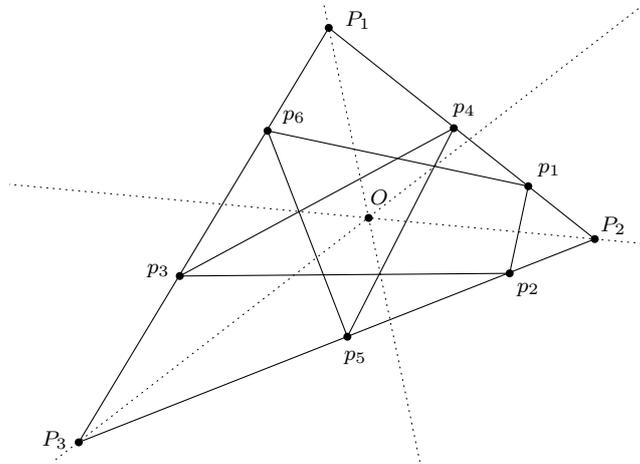
\begin{figure}[!h]
\centering
\begin{tikzpicture}[line cap=round,line join=round,>=triangle 45,x=3.5cm,y=3.5cm]
\clip(-1.2,-0.64) rectangle (1.2,1.09);
\draw (0,1)-- (-0.94,-0.57);
\draw (-0.94,-0.57)-- (1,0.2);
\draw (1,0.2)-- (0,1);
\draw [dotted,domain=-1.9:1.7] plot(\x,{(-0.15--0.72*\x)/-0.15});
\draw [dotted,domain=-1.9:1.7] plot(\x,{(--0.25-0.08*\x)/0.85});
\draw [dotted,domain=-1.9:1.7] plot(\x,{(-0.18-0.85*\x)/-1.09});
\draw (0.07,-0.17)-- (-0.23,0.61);
\draw (-0.23,0.61)-- (0.75,0.4);
\draw (0.75,0.4)-- (0.68,0.07);
\draw (0.68,0.07)-- (-0.56,0.06);
\draw (0.47,0.62)-- (0.07,-0.17);
\draw (-0.56,0.06)-- (0.47,0.62);
\begin{scriptsize}
\fill [color=black] (0,1) circle (1.5pt);
\draw[color=black] (0.11,1.03) node {$P_1$};
\fill [color=black] (-0.94,-0.57) circle (1.5pt);
\draw[color=black] (-1.03,-0.56) node {$P_3$};
\fill [color=black] (1,0.2) circle (1.5pt);
\draw[color=black] (1.07,0.25) node {$P_2$};
\fill [color=black] (0.07,-0.17) circle (1.5pt);
\draw[color=black] (0.1,-0.25) node {$p_5$};
\fill [color=black] (-0.23,0.61) circle (1.5pt);
\draw[color=black] (-0.13,0.66) node {$p_6$};
\fill [color=black] (0.75,0.4) circle (1.5pt);
\draw[color=black] (0.83,0.46) node {$p_1$};
\fill [color=black] (0.15,0.28) circle (1.5pt);
\draw[color=black] (0.19,0.36) node {$O$};
\fill [color=black] (0.68,0.07) circle (1.5pt);
\draw[color=black] (0.75,0.01) node {$p_2$};
\fill [color=black] (-0.56,0.06) circle (1.5pt);
\draw[color=black] (-0.64,0.08) node {$p_3$};
\fill [color=black] (0.47,0.62) circle (1.5pt);
\draw[color=black] (0.51,0.69) node {$p_4$};
\end{scriptsize}
\end{tikzpicture}
\caption{A $6$-periodic orbit $(p_k)_k$ on a centrally-projective triangle based at $0,P_1,P_2,P_3$. The dotted lines are representatives of the projective fields of lines on the sides of the triangle.}
\label{fig:centrally_projective_triangle}
\end{figure}

\begin{proposition}
\label{prop:centproj_odd_polygon_reflective}
Suppose that $n$ is odd. Then any $(p_1,p_2)\in P_1P_2\times P_2P_3$ with $p_1\neq p_2$ determines a $2n$-periodic orbit of the centrally-projective polygon based at $O,P_1,\ldots,P_n$. See Figure \ref{fig:centrally_projective_triangle}.
\end{proposition}

\begin{proof}
The idea of the proof is based on a construction which can be found in
\cite{taba_projectif}, to associate to a projective billiard a dual billiard. For a an introduction to dual billiards see for example \cite{taba_dual_billiards}. Consider a line $L_{\infty}\subset\RP^2$ which does not contain any of the points $O,P_1,\ldots,P_n$ (seen as the \textit{line at infinity}) and a polarity which sends $O$ to $L_{\infty}$ (that is a choice of a quadratic form for which $O$ is the pole of $L_{\infty}$; see Section \ref{section:general_properties_on_quadrics}). 

Consider the orbit $(p_j)_{j\in\ZZ}$ of $(p_1,p_2)\in P_1P_2\times P_2P_3$ with $p_1\neq p_2$. For each $j$, denote by $q_j$ the polar dual of the line $p_jp_{j+1}$ by $Q_j$ the polar dual of the line $P_jP_{j+1}$, and by $\omega_j$ the polar dual of the projective line at $p_j$ of the centrally-projective billiard, namely $Op_j$. Since the line $Op_j$ contains $O$, its polar dual $\omega_j$ belongs to $L_{\infty}$.

Given an integer $j$, the lines $p_{j-1}p_j$, $p_jp_{j+1}$, $P_jP_{j+1}$ and $Op_j$ form a harmonic quadruple of lines passing through the point $p_j$, hence the points $q_{j-1}$, $q_j$, $Q_j$, $\omega_j$ all belong to the same line (given by the polar dual of $p_j$) and they form a harmonic set of points. Since $\omega_j$ is on $L_{\infty}$, the harmonicity condition implies that $Q_j$ is at equal distance from $q_{j-1}$ and from $q_j$ in the open set $\RP^2\setminus L_{\infty}$ which is canonically diffeomorphic to $\RR^2$. This can be rewritten in terms of vectors of $\RP^2\setminus L_{\infty}\simeq\RR^2$ as
\begin{equation}
\label{eq:centproj_thales_relation}
\overrightarrow{q_{j-1}Q_j}=\overrightarrow{Q_jq_j}.
\end{equation}
Thus, we have transformed our problem into a simplified version of dual billiards in polygons, defined as follows:

\begin{definition}
Let $Q_1,\ldots,Q_n$ be distinct points in $\RR^2$ and $q_0\in\RR^2$. The \textit{virtual outer orbit} of $q_0$ associated to $(Q_1,\ldots,Q_n)$ is the sequence $(q_j)_{j\in\ZZ}$ of points of $\RR^2$ such that for each $j$, the point $Q_j$ is on the line $q_{j-1}q_j$ and is at equal distance from $q_{j-1}$ and from $q_j$.
\end{definition}

In our case, the problem is to show that any virtual outer orbit $q=(q_j)_{j\in\ZZ}$ of the above constructed points $Q_1,\ldots,Q_n$ is \textit{$2n-$periodic} in the sense that $q_{2n+j} = q_j$ for a certain $j$ (and thus for all $j$). By polar duality, we will recover that the corresponding virtual orbit of the original local projective billiard is $2n$-periodic.

\begin{lemma}
\label{lemma:centproj_thales_corollary}
For all $j\in\ZZ$ we have the relation
$$\overrightarrow{q_{j-1}q_{j+1}}=2\overrightarrow{Q_jQ_{j+1}}.$$
\end{lemma}

\begin{proof}
This relation comes from Relation \eqref{eq:centproj_thales_relation}, which defines a configuration as in the intercept theorem:
$$\overrightarrow{q_{j-1}q_{j+1}}=\overrightarrow{q_{j-1}Q_{j}}+\overrightarrow{Q_jq_{j}}+\overrightarrow{q_{j}Q_{j+1}}+\overrightarrow{Q_{j+1}q_{j+1}}=2\overrightarrow{Q_jq_{j}}+2\overrightarrow{q_{j}Q_{j+1}}=2\overrightarrow{Q_jQ_{j+1}}.$$
\end{proof}

We conclude the proof of of Proposition \ref{prop:centproj_odd_polygon_reflective} by showing that $q_{2n+1}=q_1$. Indeed, by Lemma \ref{lemma:centproj_thales_corollary} we have
$$\overrightarrow{q_{1}q_{2n+1}}=\sum_{j=1}^{n}\overrightarrow{q_{2j-1}q_{2j+1}}=2\sum_{j=1}^{n}\overrightarrow{Q_{2j}Q_{2j+1}}.$$
Since $n$ is odd, if we write $n=2m+1$ with an integer $m\geq 1$, the latter sum can be rewritten as
\begin{equation}
\label{eq:centproj_thales_relation_2}
\sum_{j=1}^{n}\overrightarrow{Q_{2j}Q_{2j+1}}=\sum_{j=1}^{m}\overrightarrow{Q_{2j}Q_{2j+1}}+\sum_{j=m+1}^{n}\overrightarrow{Q_{2j}Q_{2j+1}}=\sum_{j=1}^{m}\overrightarrow{Q_{2j}Q_{2j+1}}+\sum_{i=1}^{m+1}\overrightarrow{Q_{2i-1}Q_{2i}}
\end{equation}
where the last equality is obtained by the change of variables $j=i+k$ and the relation $Q_{i+n}=Q_i$. It is easy to see that the last quantity of \eqref{eq:centproj_thales_relation_2} equals $\overrightarrow{Q_{1}Q_{2m+2}}=\overrightarrow{Q_{1}Q_{1}}=0$. Hence $q_1=q_{2n+1}$ and therefore the lines $p_1p_2$ and $p_{2n+1}p_{2n+2}$ are the same which implies that the orbit $(p_j)_j$ is $2n$-periodic.
\end{proof}
	
\section{Billiards and Pfaffian systems}
	\label{section_pfaffian_systems}
	
In this section we present a strong link between $k$-reflective (eventually projective) billiards and integral surfaces of a certain distribution called \textit{Birkhoff's distribution}. The idea was developped in \cite{bary}, and interesting arguments are given in \cite{glutkud2} from which this section of the manuscript is inspired. Some of the results presented in this section are gathered in a preprint \cite{fierobe_triangular}.

\subsection{Classical Birkhoff's distribution}

Let $M$ be a smooth or analytic manifold and $k$ be a non-zero positive integer. We denote by $\grass{k}{TM}$ the fiber bundle over $M$ made by $k$-dimensional vector subspaces of $TM$, and by $\pi:\grass{k}{TM}\to M$ its natural projection.

\begin{definition}
\label{definition:distribution}
A \textit{$k$-dimensional distribution} on $M$ is a smooth (or analytic) section $\mathcal D:M\to\grass{k}{TM}$. An $\ell$-dimensional \textit{integral manifold (or surface) of $\mathcal D$} is a smooth (or analytic) submanifold $S\subset M$ of dimension $\ell$ such that for all $p\in S$
\begin{equation}
\label{eq:def_integral_surface}
T_pS\subset\mathcal{D}(p).
\end{equation}
An $\ell$-dimensional \textit{pseudo-integral manifold (or surface) of $\mathcal D$} is a smooth (or analytic) submanifold $S\subset M$ of dimension $\ell$ such that \eqref{eq:def_integral_surface} holds only for $p$ in a subset $V\subset S$ of non-zero Lebesgue measure, called \textit{integral set}.
\end{definition}

\begin{remark}
\label{remark:manifold_analytic_pseudoref_reflective}
A connected analytic pseudo-integral manifold of an analytic distribution is an integralble manifold. This result is implied by the Uniqueness Theorem for analytic extension.
\end{remark}



Now let us define the usual version of Birkhoff's distribution. We set $M=(\RR^d)^k$. Consider the open dense subset $U\subset M$ of $k$-tuples $p=(p_1,\ldots,p_k)$ such that for each $j=2,\ldots,k-1$ the points $p_{j-1},p_j,p_{j+1}$ do not lie on the same line of $\RR^d$. For each $j$ (modulo $k$), consider the interior bisector $L_j(p)\subset\RR^d$ of the oriented angle between the vectors $\overrightarrow{p_jp_{j-1}}$ and $\overrightarrow{p_jp_{j+1}}$ and denote by $H_j(p)\subset\RR^d$ its orthogonal hyperplane (with respect to the Euclidean metric of $\RR^d$). The hyperplanes $H_j(p)$ have the following simple property related to billiards which is simply due to the definition of the reflection law on the billiard.

\begin{lemma}
\label{lemma:classical_birkhoff_distribution_angles}
Let $\Omega$ be a (classical) billiard in $\RR^d$ and $p$ be a sequence of points $(p_1,\ldots,p_k)$ on the boundary $\partial\Omega$ such that $p\in U$. Then $p$ is a $k$-periodic orbit if and only if $T_{p_j}\partial\Omega=H_j(p)$.
\end{lemma}

Then we can identify $T_pU$ with $\oplus_{j=1}^k T_{p_j}\RR^d$, and consider the projections $\pi_j M\to\RR^d$ sending $p$ to $p_j$. We can consider the

\begin{definition}
\label{definition:classical_birkhoff_distribution}
The $k(d-1)$-dimensionnal analytic distribution on $U\subset(\RR^d)^k$ defined for all $p\in U$ by 
$$\mathcal{D}(p) = \oplus_{j=1}^k H_j(p)$$
is called \textit{Birkhoff's distribution}. We further say that an \textit{integral} (respectively \textit{pseudo-integral}) manifold $S$ of $\mathcal{D}$ is \textit{non-trivial} if the restriction of each $\pi_j$ to $S$ has rank $d-1$ for all $p\in S$ (respectively $p\in V$).
\end{definition}

Notice that if the restrictions of $\pi_j$ to $S$ have rank $d-1$ at $p$, then there is a small neighborhood $W$ of $p$ such that the $\pi_j(V)$ are submanifolds of $\RR^d$ of the same regularity than $S$. Birkhoff's distribution has then the following property which is related to the latter remark and to Lemma \ref{lemma:classical_birkhoff_distribution_angles}:

\begin{proposition}
\label{proposition:link_birkh_reflective}
 \textbf{1)} If $\mathcal{B}$ is a local $\mathcal{C}^{r}$-smooth (respectively analytic) $k$-reflective billiard, then there is a subset of $k$-periodic orbits of $\mathcal{B}$ which is a non-trivial $2(d-1)$-dimensional $\mathcal{C}^{r-1}$-smooth (respectively analytic) integral manifold of $\mathcal{D}$.\\
\textbf{2)} Conversely, if $S\subset U$ is a $\class^r$-smooth (respectively an analytic) non-trivial integral manifold of $\mathcal{D}$ of dimension $2(d-1)$, then for all $p\in S$ there is a neighborhood $W\subset S$ of $p$ for which $(\pi_1(W),\ldots,\pi_k(W))$ is a local $\class^r$-smooth (respectively analytic) $k$-reflective billiard.
\end{proposition}

\begin{proof}
1) Suppose that $\mathcal{B}$ is a local $\mathcal{C}^r$-smooth (respectively analytic) $k$-reflective billiard. We denote by $a_1,\ldots,a_k$ its classical boundaries. Consider $p=(p_1,\ldots,p_k)$ a $k$-periodic orbit of $\mathcal{B}$ such that any $(q_1,q_2)\in a_1\times a_2$ sufficiently close to $(p_1,p_2)$ can be completed into a $k$-periodic orbit of $\mathcal{B}$.\\
Given $(q_j,q_{j+1})\in a_j\times a_{j+1}$ close enough to $(p_j,p_{j+1})$, the line $q_{j}q_{j+1}$ is reflected at $q_{j+1}$ into a line intersecting $a_{j+1}$ at a certain point $q_{j+2}$ close to $p_{j+2}$. Therefore, one can define a $\class^{r-1}$ -smooth (respectively an analytic) map $B_j(q_j,q_{j+1})=(q_{j+1}, q_{j+2})$ locally on a neighborhood 
of $(p_j,p_{j+1})$, and which is a diffeomorphism onto its image (see Proposition \ref{proposition:billiard_regularity_rank}). Consider the set $S$ defined as the graph of the map $s:(q_1,q_2)\mapsto(q_3,\ldots,q_k)$ where each $q_{j+1}$ is defined as a map of $(q_1,q_2)$ by the relation $B_j\circ\cdots\circ B_1(q_1,q_2)=(q_{j+1},q_{j+2})$. By construction, $S$ is a $2(d-1)$-dimensional $\class^{r-1}$-smooth (respectively analytic) immersed submanifold of $U$, and the restriction of each $\pi_j$ to $S$ has rank $d-1$ since $B_j\circ\cdots\circ B_1$ is a local diffeomorphism. By assumptions,   one can suppose that $S$ contains only $k$-periodic orbits by shrinking the set of definition of $s$. By Lemma \ref{lemma:classical_birkhoff_distribution_angles}, for $q\in S$ and any $j$, $T_{q_j}a_j=H_j(q)$, hence $d\pi_j(T_q S)=H_j(q)$. Therefore $S$ is an integral manifold of $\mathcal{D}$.

2) Suppose that $S\subset U$ is a $\class^r$-smooth (respectively an analytic) non-trivial integral manifold of $\mathcal{D}$ and $p\in S$. Choose a neighborhood $W\subset S$ of $p$ for which $a_1:=\pi_1(W),\ldots,a_k:=\pi_k(W)$ are $\class^{r}$-smooth (respectively analytic) immersed submanifolds of $\RR^d$. Since $S$ is an integral manifold of $\mathcal{D}$, any $q=(q_1,\ldots,q_k)\in U$ satisfies
$$T_{q_j}a_j = d\pi_j(T_q S) = d\pi_j(\mathcal{D}(q))=H_j(q),$$
hence is a $k$-periodic orbit of $\mathcal{B}:=(a_1,\ldots,a_k)$ by Lemma \ref{lemma:classical_birkhoff_distribution_angles}. It remains to show that $\mathcal{B}$ is $k$-reflective. Consider the map $i:p\in S\mapsto (p_1,p_2)\in a_1\times a_2$. Let us show that $i$ is a local diffeomorphism in a neighborhood of $p$. The map $s$ of part \textit{1)} is such that $s\circ i(q) = q$ for all $q\in W$ since the latter are periodic orbits of $\mathcal{B}$. Therefore $di(p)$ is injective, and because $\dim S = \dim a_1\times a_2=2(d-1)$ the conclusion follows.
\end{proof}

Proposition \ref{proposition:link_birkh_reflective} has an analogue for pseudo-integral surfaces and $k$-pseudo-reflective billiards. In the following proposition, we say that a property is satisfied \textit{for almost all points} $p$ in a subset $V$ of a smooth manifold, if the set of points $p\in V$ for which it is not satisfied has zero Lebesgue measure.

\begin{proposition}
\label{proposition:link_birkh_pseudo_reflective}
 \textbf{1)} If $\mathcal{B}$ is a local $\mathcal{C}^{r}$-smooth $k$-pseudo-reflective billiard, then there is a subset of (not necessarily periodic) orbits of $\mathcal{B}$ which is a non-trivial $2(d-1)$-dimensional $\mathcal{C}^{r-1}$-smooth pseudo-integral manifold of $\mathcal{D}$.\\
\textbf{2)} Conversely, if $S\subset U$ is a $\class^r$-smooth non-trivial pseudo-integral manifold of $\mathcal{D}$ of dimension $2(d-1)$, then for almost all $p$ in the set $V$ of Definition \ref{definition:distribution} there is a neighborhood $W\subset S$ of $p$ for which $(\pi_1(W),\ldots,\pi_k(W))$ is a local $\class^r$-smooth $k$-pseudo-reflective billiard.
\end{proposition}

\begin{remark}
Notice that the analytic version of this result is given by Proposition \ref{proposition:link_birkh_reflective}, since $k$-pseudo-reflective analytic billiards are $k$-reflective, and connected analytic pseudo-integrable manifolds are integrable (see Remarks \ref{remark:billiard_analytic_pseudoref_reflective} and \ref{remark:manifold_analytic_pseudoref_reflective}).
\end{remark}

\begin{proof}
The proof is analogous to the proof of Proposition \ref{proposition:link_birkh_reflective}, except that we will work with so-called \textit{Lebesgue points}. Let us first recall defintions and results about them.

\begin{definition}
Let $V\subset \RR^d$ be a Lebesgue measurable set. A point $x\in V$ is said to be a \textit{Lebesgue point of $V$} if one has 
$$\lim_{r\to 0} \frac{\lambda(V\cap B(x,r))}{\lambda(B(x,r))} = 1$$
where $\lambda$ is the Lebesgue measure of $\RR^d$ and $B(x,r)$ is the Euclidean ball of radius $r>0$ centered at $x$. 
\end{definition}

This definition naturally extends to subset $V$ of smooth differentiable manifolds. We observe the following

\begin{theorem}[Lebesgue density theorem]
\label{theorem:leb_density_thm}
Let $V$ be a Lebesgue measurable set of a smooth differentiable manifold. Then almost all points of $V$ are Lebesgue poins of $V$.
\end{theorem}

\begin{lemma}[see \cite{glutkud2}]
\label{lemma:coincide_lebesgue_points}
Let $U\subset M$ be an open subset of a differentiable manifold $M$ and $f,g:U\to N$ be $\class^r$-smooth maps from $U$ to a differentiable manifold $N$. If $V$ is a subset of $U$ on which $f=g$ and $p$ is a Lebesgue point of $V$, then the $1$-jets of $f$ and $g$ at $p$ coincide.
\end{lemma}

\begin{proof}
Choosing a convenient set of coordinates, one can suppose that $U\subset M=\RR^d$, $N=\RR^k$, $p=0$ and also that $g=0$, by substituing $f-g$ to $f$. Consider the map $P:\RR^d\setminus\{0\}\to\partial B(0,1)$ defined by $P(x)=x/\|x\|$ where $\|\cdot\|$ is the Euclidean metric. We first prove that \textit{for $r>0$ $W_r:=P(V\cap B(0,r)\setminus\{0\})$ is dense in $\partial B(0,1)$}. Indeed, otherwise there would exist a non-empty open subset $U$ of $\partial B(0,1)$ included in $\partial B(0,1)\setminus W_r$. The cone $U':=P^{-1}(U)$ is open an satisfies for all $0<\rho\leq r$ that\\
\enum $\lambda(U'\cap B(0,\rho))=\rho^d\lambda(U'\cap B(0,1))$, since $U'\cap B(0,\rho)$ is obtained from $U'\cap B(0,1)$ by the dilatation $x\mapsto \rho x$;\\
\enum $U'\cap B(0,\rho)\subset B(0,\rho)\setminus V$, because $U$ and $W_r$ have empty intersection.\\
Hence for $0<\rho\leq r$
$$\frac{\lambda(V\cap B(0,\rho))}{\lambda(B(0,\rho))} = 1-\frac{\lambda(B(0,\rho)\setminus V)}{\lambda(B(0,\rho))}\leq 1-\frac{\lambda(U'\cap B(0,\rho))}{\lambda(B(0,\rho))}=1-\frac{\lambda(U'\cap B(0,1))}{\lambda(B(0,1))}<1$$
and the latter bound doesn't depend on $\rho$, which is impossible since $p=0$ is a Lebesgue point of $V$. Hence for $v\in \partial B(0,1)$, one can find a sequence of $v_n\in V$ such that $v_n\to 0$ and $P(v_n)\to v$. This implies that $df(0)\cdot v = \lim f(v_n)/\|v_n\| = 0$.
\end{proof}

Now we can prove Proposition \ref{proposition:link_birkh_pseudo_reflective}:

1) Suppose that $\mathcal{B}$ is a local $\mathcal{C}^r$-smooth $k$-pseudo-reflective billiard. We denote by $a_1,\ldots,a_k$ its classical boundaries and by $V'\subset a_1\times a_2$ a subset of non-zero measure of points $(q_1,q_2)$ which can be completed into a $k$-periodic orbit of $\mathcal{B}$.\\
Consider $p=(p_1,\ldots,p_k)$ a $k$-periodic orbit of $\mathcal{B}$ such that $(p_1,p_2)\in a_1\times a_2$ is a Lebesgue point of $V'$. The corresponding manifold $S$ defined in Proposition \ref{proposition:link_birkh_reflective} as the graph of a map $s$ does not only contain $k$-periodic orbits anymore. However by Lemma \ref{lemma:classical_birkhoff_distribution_angles}, all $q=(q_1,q_2)\in V'$ lying in an open subset $W'$ containing $p$ on which $s$ is defined is such that $s(q)$ is $k$-periodic. Then since $(p_1,p_2)\in W'\cap V'$ is a Lebesgue point of $V'$, $W'\cap V'$ has non-zero measure in $W'$. Hence the subset $s(W'\cap V')\subset S$ has non-zero measure in $S$ and contains only $k$-periodic orbits. 

2) Suppose that $S\subset U$ is a $\class^r$-smooth non-trivial pseudo-integral manifold of $\mathcal{D}$, $V$ the set of points $p\in S$ for which $T_pS\subset \mathcal{D}(p)$ and $p\in V$ a Lebesgue point. Choose a neighborhood $W\subset S$ of $p$ for which $a_1:=\pi_1(W),\ldots,a_k:=\pi_k(W)$ are $\class^{r}$-smooth immersed submanifolds of $\RR^d$. As in the proof of Proposition \ref{proposition:link_birkh_reflective}, any $q\in W\cap V$ is a $k$-periodic orbit of $\mathcal{B}:=(a_1,\ldots,a_k)$. Consider the map $i:p\in S\mapsto (p_1,p_2)\in a_1\times a_2$. The map $s$ of part \textit{1)} is such that $s\circ i(q) = q$ for all $q\in W\cap V$ since the latter are periodic orbits of $\mathcal{B}$. Therefore by Lemma \ref{lemma:coincide_lebesgue_points}, since $p$ is a Lebesgue point of $V$ we can write $ds(p_1,p_2)\circ di(p)=\id$ and the conclusion follows as before.
\end{proof}

\subsection{Prolongations of Pfaffian systems}
\label{subsection:prolongation_pfaffian_systems}

Let $M$ be an analytic manifold, $\mathcal{D}$ be an analytic distribution on $M$, and $k\in\{1,\ldots,\dim \mathcal{D}\}$. We denote by $\grass{k}{TM}$ the fiber bundle over $M$ made by $k$-dimensional vector subspaces of $TM$, with its natural projection $\pi:\grass{k}{TM}\to M$. 

One can define a natural analytic distribution $\mathcal{K}$ on $\grass{k}{TM}$, called \textit{contact distribution}, and defined for all $(x,E)\in\grass{k}{TM}$ by $\mathcal{K}(x,E) = d\pi^{-1}(E)$. In this subsection we introduce Pfaffian systems and their prolongations, as a way to link integral manifolds of $\mathcal{D}$ and integral manifolds of $\mathcal{K}$ contained in some submanifolds of $(x,E)\in\grass{k}{TM}$.

\begin{definition}[\cite{glutkud2}, definition 21]
Given a family of analytic distributions $(\mathcal{D}_i)_i$ on $M$, we call the data $\mathcal{P} = (M,\mathcal{D},k;(\mathcal{D}_i)_i)$ a \textit{Pfaffian system with transversality conditions}.\\
\enum A $k$-plane $E\in\grass{k}{TM}$ is said to be \textit{integral} if for any $1$-form $\omega$ vanishing on $\mathcal{D}$, $d\omega$ vanishes on $E$.\\
\enum An \textit{integral manifold} (or surface) of $\mathcal{P}$ is an integral manifold of $\mathcal{D}$ of dimension $k$ such that, for all $i$, its tangent subspaces either are transverse to $\mathcal{D}_i$, or intersect it by zero.\\
\enum An \textit{pseudo-integral manifold} (or surface) of $\mathcal{P}$ is a pseudo-integral manifold of $\mathcal{D}$ of dimension $k$ such that, for $x$ lying in its integral set $V$ (see Definition \ref{definition:distribution}) and for all $i$, $T_xS$ is either transverse to $\mathcal{D}_i$, or intersects it by zero.
\end{definition}

It follows immediately that the tangent planes to an integral manifold $S$ are integral. Notice also that if $S$ is a pseudo-integral manifold and $V$ is its integral set, then, due to Lemma \ref{lemma:coincide_lebesgue_points}, $T_xS$ is integral for any Lebesgue point $x$ of $V$.

In the following $\mathcal{P}=(M,\mathcal{D},k;(\mathcal{D}_i)_i)$ denote a Pfaffian system with transversality conditions. As described in \cite{glutkud2}, subsection 2.3, the set $\widetilde{M}_k$ of integral $k$-planes of $M$ is an analytic subset hence a stratified manifold: it is a locally finite and at most countable disjoint union of smooth analytically constructible subsets (see \cite{lojasiewicz}, section IV.8).

\begin{definition}[\cite{glutkud2}, definition 23]
Let $M'$ be a stratum of $\widetilde{M}_k$, $\mathcal{K}'$ the restriction of the contact distribution $\mathcal{K}$ to $M'$, and $\mathcal{D}'_i$ the pull-back of  $\mathcal{D}_i$ on $M'$ for each $i$. The Pfaffian system $\mathcal{P}' = (M',\mathcal{K}',k;(\mathcal{D}'_i)_i,\ker d\pi)$ is called \textit{a first Cartan prolongation} of $\mathcal{P}$.
\end{definition}

If $S\subset M$ is a $\class^{r}$-smooth submanifold of $M$ of dimension $k$, one can consider the subset $S^{(1)}\subset\grass{k}{TM}$ defined by
\begin{equation}
\label{equation:first_lift}
S^{(1)} = \ensemble{(x,T_xS)}{x\in S}.
\end{equation}
It is $\class^{r-1}$-smooth submanifold of $\grass{k}{TM}$, of dimension $\dim S$ and transverse to $\pi$. We call it the \textit{first lift} of $S$.

\begin{proposition}[\cite{glutkud2}, subsection 2.3, \cite{bryant_chern}, chapter VI]
\label{prop:integral_prolongation}
The lift $S^{(1)}$ of an integral surface $S$ of $\mathcal{P}$ contains an open dense subset such that each its connected component $S'$ lies in some stratum $M'$ of $\widetilde{M}_k$, and such that $S'$ is an integral surface of the first Cartan prolongation $\mathcal{P}' = (M',\mathcal{K}',k;(\mathcal{D}'_i)_i,\ker d\pi)$.

\end{proposition}

\begin{proof}
As explained, the tangent planes of an integral manifolds are integral, hence $S^{(1)}$ is contained in the set $\widetilde{M}_k$ of integral $k$-planes of $\mathcal{D}$. Now let $S'$ be a connected component of $S^{(1)}$ contained in a stratum $M'$ of $\widetilde{M}_k$. For all $p=(x,T_xS)\in S'$ we have $d\pi(T_pS^{(1)})=T_xS$, hence $T_pS^{(1)}\subset\mathcal{K}(p)$ and therefore $T_pS'\subset \mathcal{K}'(p)$. Moreover, the equality $d\pi(T_pS^{(1)})=T_xS$ implies that $d\pi$ is injective on $T_pS^{(1)}$ hence on $T_pS'$, and the transversality condition with $\ker d\pi$ is satisfied. One can easily check that he other transversality conditions are satisfied. 
\end{proof}

The converse result is also true, but only locally:

\begin{proposition}
\label{prop:projection_integral_manifolds}
Let $M'$ be a stratum of $\widetilde{M}_k$ and $S'$ be an integral manifold of the Pfaffian system $\mathcal{P}' = (M',\mathcal{K}',k;(\mathcal{D}'_i)_i,\ker d\pi)$ such that the intersection $TS'\cap\ker d\pi$ is $\{0\}$. Then for any $p\in S'$ there is an open subset $U\subset S'$ containing $p$ and   such that $S:=\pi(U)$ is an integral surface of $\mathcal{P}$ such that $S^{(1)}=S'$.
\end{proposition}

\begin{proof}
Since $TS'\cap\ker d\pi=\{0\}$ are transverse, there is a small neighborhood $U$ of $p$ such that $S:=\pi(U)$ is a $k$-dimensionnal manifold with $T_{\pi(q)}S=d\pi(T_qS')$ for any $q\in U$. Hence if $q=(x,E)\in U$, then $T_{\pi(q)}S=E$, because $S'$ is an integral manifold of the distribution $\mathcal{K}'$. Therefore $T_xS$ is an integral plane of $\mathcal{D}$, and thus $S$ is an integral manifold of $\mathcal{D}$. One can analogously check that $S$ satifies transversality conditions $\mathcal{D}_i$.
\end{proof}

The same constructions work also for pseudo-integral manifolds of $\mathcal{P}$:

\begin{proposition}[\cite{glutkud2}, subsection 2.3]
\label{prop:pseuo_integral_prolongation}
Let $S\subset M$ be a pseudo-intergral surface of $\mathcal{P}$, and $V$ be its integral set. Suppose $p=(x,T_xS)\in S^{(1)}$ is such that $x$ is a Lebesgue point of $V$. Then replacing $p$ by a Lebesgue point in $V$ arbitrarily close to $p$ (now denoted by $p$) one can achieve that there is a stratum $M'$ of $\widetilde{M}_k$ and a smooth submanifold $S'$ in an open subset of $M'$, such that:\\
\enum $S'$ contains $p$ and is tangent to $S^{(1)}$ at $p$;\\
\enum $S'$ is a pseudo-integral surface of the Cartan prolongation $\mathcal{P}' = (M',\mathcal{K}',k;(\mathcal{D}'_i)_i,\ker d\pi)$.
\end{proposition}

\begin{proof}
Denote by $\widetilde{V}\subset S^{(1)}$ the set of points $p=(x,T_xS)\in S^{(1)}$ such that $x\in V$: $\pi$ maps the Lebesgue points of $\widetilde{V}$ to the Lebesgue points of $V$. Denote by $\mathcal{L}(\widetilde{V})$ the set of Lebesgue points of $\widetilde{V}$. As explained, any Lebesgue point $p\in \mathcal{L}(\widetilde{V})$ of $\widetilde{V}$ belongs to $\widetilde{M}_k$ (since $T_xS$ is an integral plane, where $p=(x,T_xS)$).

Fix $p\in\mathcal{L}(\widetilde{V})$: $p$ is also a Lebesgue point of $\mathcal{L}(\widetilde{V})$ which is a simple consequence of Lebesgue density theorem (Theorem \ref{theorem:leb_density_thm}). In particular, any small neighborhood of $p$ in $\mathcal{L}(\widetilde{V})$ has non-zero measure in $S^{(1)}$. Hence one can choose a stratum $M'$ of $\widetilde{M}_k$ containing points arbitrarily close to $p$ such that $M'\cap \mathcal{L}(\widetilde{V})$ has non-zero Lebesgue measure in $S^{(1)}$ and the latter intersection (considered as a subset in $S^{(1)}$) has Lebesgue points arbitrarily close to $p$. From now on, $p$ will be one of the latter Lebesgue points.

On a small neighborhood $W$ of $p$ in $\grass{k}{TM}$, one can define a smooth map $s:W\to M'$ such that $s(q)=q$ for $q\in W\cap M'$ (take for example the orthogonal projection onto $M'$ in a set of coordinates). For $q\in \mathcal{L}(\widetilde{V})\cap M'$, we have $s(q)=q=i(q)$ where $i:S^{(1)}\to \grass{k}{TM}$ is the natural embedding of $S^{(1)}$. By Lemma \ref{lemma:coincide_lebesgue_points}, $ds(p) = di(p)$ hence $ds(p)$ is injective, and therefore one can suppose that $S':=s(S^{(1)}\cap W)$ is an $\ell$-dimensional submanifold of $M'$ (by shrinking $W$ if necessary). It is tangent to $S^{(1)}$ since $T_pS'=\im ds(p)=\im di(p) = T_pS^{(1)}$.

Let us show that $S'$ is a pseudo-integral manifold of the Cartan prolongation $\mathcal{P}'= $ $(M',\mathcal{K}',$ $k;(\mathcal{D}'_i)_i,\ker d\pi)$. Write $V'=\mathcal{L}(\widetilde{V})\cap W$, the Lebesgue points of $\widetilde{V}$ contained in $W$. Previous argument shows that $V'\subset S'$ and that for all $q\in V'$, $T_qS'=T_qS^{(1)}$, hence $d\pi(T_qS')=T_{\pi(q)} S$ and $S'$ is a pseudo-integral surface of $\mathcal{K}'$. The transversality conditions follow from the same argument.
\end{proof}

Propositions \ref{prop:integral_prolongation} and \ref{prop:pseuo_integral_prolongation} imply the existence of an integral (respectively a pseudo-integral) manifold $S'$  in the grassmanian as soon as there is an integral (respectively a pseudo-integral) manifold $S$ in $M$. Let us call \textit{$S'$ a first Cartan prolongation of $S$}. We deduce the following

\begin{corollary}
\label{cor:existence_of_prolongations}
Let $S\subset M$ be a $\mathcal{C}^r$-smooth integral (respectively pseudo-intergral) manifold of $\mathcal{P}$. Then there is a sequence $\left(\mathcal{P}^{(k)}\right)_{k=0\ldots r}$ of Pfaffian systems and a sequence $S_k$ of integral (respectively pseudo-integral) manifolds of $\mathcal{P}^{(k)}$, such that $\mathcal{P}^{(0)}=\mathcal{P}$ and such that for each $k<r$, $\mathcal{P}^{(k+1)}$ and $S_{k+1}$ are first Cartan prolongations of $\mathcal{P}^{(k)}$ and $S_k$.
\end{corollary}

We conclude this subsection by the folowing powerful result on prolongations of a Pfaffian system $\mathcal{P}$. It is cited in \cite{glutkud2}, theorem 24, and in \cite{bryant_chern}, chapter VI, paragraph 3. The original statement of this result can be found in \cite{rashevsky} which is in russian.

\begin{theorem}[E. Cartan \cite{cartan}, M. Kuranishi \cite{kuranishi}, P. K. Rachevsky \cite{rashevsky}]
\label{theorem:cartan_ku_ra}
Suppose that $\mathcal{P}$ has no analytic integral surfaces. Then for any sequence of Pfaffian systems $\mathcal{P}^{(k)}=(M^{(k)},\ldots)$ such that $\mathcal{P}^{(0)}=\mathcal{P}$ and $\mathcal{P}^{(k+1)}$ is a first Cartan prolongation of $\mathcal{P}^{(k)}$, one can find an integer $k_0>0$ for which $M^{(k_0)}=\emptyset$.
\end{theorem}

\subsection{$r$-jets approximation of integral manifolds}

Let $M$ be an analytic manifold, $\mathcal{D}$ be an analytic distribution on $M$, $k\in\{1,\ldots,\dim \mathcal{D}\}$ and $r>0$ an integer. 

Let $p\in M$. A \textit{germ of $\class^r$-smooth $k$-dimensional submanifold of $M$ at $p$} is the family of $\class^r$-smooth submanifolds of dimension $k$ of $M$ containing $p$ and satisfying: given to such submanifolds $S,S'$, there is an open subset $V$ of $M$ containing $p$ for which $S\cap V=S'\cap V$. Denote by $(S,p)$ the germs of submanifolds at $p$ containing $S$ and by $\mathcal{G}(M, k, r)$ the set of all germs of $\class^r$-smooth $k$-dimensional submanifolds of $M$ at $p$.

In the following we define a topology on $\mathcal{G}(M, k, r)$ (which is not Hausdorff). Given a $\class^r$-smooth submanifold $S\subset M$ containing $p$, there is an injective $\class^r$-smooth immersion $f$ defined on an open subset $U\subset\RR^k$ containing $0$ and such that $f(0)=p$ and $(f(U),p)=(S,p)$. Denote by $J^r_0(f)$ the $r$-jet at $0$ of a $\class^r$-smooth map $f:U\subset\RR^k\to M$ defined on an open subset $U$ of $\RR^k$ containing $0$, and by $J^r_k(M)$ the space of all such $r$-jets.

\begin{definition}
\label{definition:space_of_jet_topology}
Given an open subset $\Omega\subset J^r_k(M)$, we define $\mathcal{G}(\Omega)\subset\mathcal{G}(M, k, r)$ the set of germs $(S,p)$ for which one can find an injective immersion $f:U\subset\RR^k\to M$ satisfying $(f(U),p)=(S,p)$ with $f(0)=p$ and  $J^r_0(f)\in\Omega$.\\
The topology generated by all $\mathcal{G}(\Omega)$ will be called \textit{Whitney $\class^r$-topology on $\mathcal{G}(M, k, r)$}.
\end{definition}

The Whitney $\class^r$-topology on $\mathcal{G}(M, k, r)$ is not Hausdorff: if two germs $(S,p)$ and $(S',p)$ of $\class^r$-smooth submanifolds are parametrized by injective immersions having the same $r$-jets, then any neighborhood of $(S,p)$ contains $(S',p)$.

The following result shows that if prolongations of integral manifolds are close in the Whitney $\class^r$-topology, then the initial manifolds are also close in the Whitney $\class^{r+1}$-topology.

\begin{proposition}
\label{prop:surface_germs_approximation}
Let $(S,p)$ be a germ of $k$-dimensional $\class^{r+1}$-smooth submanifold. Then for any open subset $V_0\subset\mathcal{G}(M, k, r+1)$ containing $\left(S, p\right)$, there is an open subset $V_1\subset\mathcal{G}(\grass{k}{TM}, k, r)$ containing $\left(S^{(1)}, (p,T_pS)\right)$ such that 
$$\forall (S',q)\in\mathcal{G}(M, k, r+1) \qquad\left(S'^{(1)}, (q,T_{q}S')\right)\in V_1 \Rightarrow (S',q)\in V_0.$$
\end{proposition}

\begin{proof}
Given a $\class^{r+1}$-smooth injective immersion $f:U\to M$ defined on an open subset $U\subset\RR^k$, one can defined its \textit{first lift} to be the $\class^{r}$-smooth injective immersion
$$f^{(1)}:U\to\grass{k}{T\RR^d}$$
defined for all $x\in U$ by $f^{(1)}(x)=(x,\im df(x))$. Notice that if $S$ is a submanifold of $M$ parametrized by $f$, then the first lift $S^{(1)}$ of $S$ is parametrized by $f^{(1)}$. Therefore, by an appropriate choice of coordinates, we just have to show the following

\begin{lemma}
\label{lemma:convergence_jet_prol}
Let $f:U\subset\RR^k\to\RR^d$ be a $\class^{r+1}$-smooth injective immersion defined on an open subset $U\subset\RR^k$ containing $0$. Then for any neighborhood $V_0\subset J^{r+1}_k(\RR^d)$ containing $J^{r+1}_0(f)$, there is a neighborhood $V_1\subset J^r_k(\grass{k}{T\RR^d})$ containing $J^r_0(f^{(1)})$ at $0$, and verifying the following property:\\
for any $\class^{r+1}$-smooth injective immersion $g$ defined on an open subset $U'\subset\RR^k$ containing $0$, if $J^r_0(g^{(1)})\in V_1$, then there is a smooth diffeomorphism $\varphi:W\subset\RR^k \to W'\subset U$, sending $0$ to $0$ and for which $J^{r+1}_0(g\circ\varphi)\in V_0$.
\end{lemma}

\begin{proof}[Proof of Lemma \ref{lemma:convergence_jet_prol}]
Replacing $f$ and $g$ by $\psi\circ f\circ\psi'$ and by $\psi\circ g\circ\psi'$ for fixed $\class^{r+1}$-smooth diffeomorphisms $\psi:\RR^d\to\RR^d$ and $\psi':W\subset\RR^k \to W'\subset U$ with $\psi'(0)=0$, we can prove Lemma \ref{lemma:convergence_jet_prol} when $f$ is of the form $f:x\in U\mapsto (x,f_0(x))$, where $f_0:U\to\RR^{d-k}$ is a $\class^{r+1}$-smooth map.\\
Since $\pi\circ g^{(1)}=g$, we can choose a first neighborhood $V_1$ of $f^{(1)}$ such that if $J^r_0(g^{(1)})\in V_1$, then one can find a $\class^{\infty}$-smooth diffeomorphism $\varphi:W\subset\RR^k \to W'\subset U$, sending $0$ on $0$, such that $g\circ\varphi$ is of the form $x\in W\mapsto (x,g_0(x))$, where $g_0:W\to\RR^{d-k}$ is a $\class^{r+1}$-smooth map. Notice that $\im d(g\circ\varphi) = \im dg$, hence $\im dg$ is generated by the vectors $(e_i,\partial_i g_0)$, $i=1\ldots k$, where $B=(e_1,\ldots,e_d)$ is the canonical basis of $\RR^d$ and $\partial_i$ is the $i$-th partial derivative. Similarly $\im df$ is generated by the vectors $(e_i,\partial_i f_0)$, $i=1\ldots k$.\\
The canonical basis $B$ defines a set of coordinates in $\grass{k}{T\RR^d}$ in which the coordinates of $\im df$ are the coordinates of $(\partial_1 f_0,\ldots,\partial_k f_0)$ in $(e_{d-k},\ldots,e_d)$, and the coordinates of $\im dg$ are the coordinates of $(\partial_1 g_0,\ldots,\partial_k g_0)$ in $(e_{d-k},\ldots,e_d)$. Therefore, saying that the $r$-jet of $\im dg$ at $0$ is close to the $r$-jet of $\im df$ at $0$ means that the same holds for the $r$-jets at $0$ of the partial derivatives $(\partial_1 g_0,\ldots,\partial_k g_0)$ and $(\partial_1 f_0,\ldots,\partial_k f_0)$. And with the additionnal assumption that $g(0)$ is close to $f(0)$, this means that the $(r+1)$-jets of $f$ and $g\circ\varphi$ at $0$ are close.
\end{proof}
This concludes the proof of Proposition \ref{prop:surface_germs_approximation}.
\end{proof}

Given a Pfaffian system $\mathcal{P}=(M,\mathcal{D},k;(\mathcal{D}_i)_i)$, a $\class^{\infty}$-smooth integral manifold $S$ of $\mathcal{P}$ and $p\in S$, the following result establishes the existence of germs of analytic integral manifolds of $\mathcal{P}$ arbitrarily close to $(S,p)$ in the Whitney $\class^r$-topology.

\begin{proposition}
\label{prop:analytic_approximation}
Let $S\subset M$ be a $\class^{\infty}$-smooth integral manifold of a Pfaffian system $\mathcal{P}=(M,\mathcal{D},k;(\mathcal{D}_i)_i)$, $p\in S$ and a positive integer $r$. Then for any open subset $V\subset\mathcal{G}(M, k, r)$ containing $(S,p)$ one can find an analytic integral manifold $S_a\subset M$ of $\mathcal{P}$ and $p_a\in S_a$ such that $(S_a,p_a)\in V$.
\end{proposition}

\begin{proof}
We first prove Proposition \ref{prop:analytic_approximation} for $r=0$. We need to show that for any open subset $V\subset M$ containing $p$, one can find an analytic integral manifold $S'$ of $\mathcal{P}$ intersecting $V$. The set $V$ contains the $\class^{\infty}$-smooth integral manifold $S\cap V$ of $\mathcal{P}$. Hence, one can find a sequence of prolongations $\mathcal{P}^{(r)}=(M^{(r)},\ldots)$ of the Pfaffian system $\mathcal{P}=(M,\mathcal{D},k;(\mathcal{D}_i)_i)$ such that $M^{(r)}\neq \emptyset$ for all $r$ (Corollary \ref{cor:existence_of_prolongations}). Therefore, Theorem \ref{theorem:cartan_ku_ra} implies the existence of an analytic integral surface of $\mathcal{P}$ in $V$, which concludes the proof in the case when $r=0$.\\
We conclude the proof by induction. Let $r>0$ and $V\subset\mathcal{G}(M, k, r)$ be an open subset containing $(S,p)$. By Proposition \ref{prop:surface_germs_approximation}, there is an open subset $V_1\subset\mathcal{G}(\grass{k}{TM}, k, r-1)$ containing $\left(S^{(1)}, (p,T_pS)\right)$ such that for all $(\tilde{S},q)\in\mathcal{G}(M, k, r)$ satisfying $\left(\tilde{S}^{(1)}, (q,T_{q}\tilde{S})\right)\in V_1$ then $(\tilde{S},q)\in V_0$. Hence given a prolongation $\mathcal{P}' = (M',\mathcal{K}',k;(\mathcal{D}'_i)_i,\ker d\pi)$ of $\mathcal{P}$ on a stratum $M'$ containing a prolongation $S'$ of $S$, the germs of $V_1$ contained in $\mathcal{G}(M', k, r-1)$ define an open set $V_1'$ of $\mathcal{G}(M', k, r-1)$ containing $(S',(p,T_pS))$. By shrinking $V_1'$ if necessary, one can further suppose that if $(\tilde{S},q)\in V_1'$, then $T\tilde{S}\cap \ker d\pi$ is $\{0\}$. By induction, one can find an analytic integral manifold $S'_a$ of $\mathcal{P}'$ and a point $q_a\in S_a$ such that $(S'_a,q_a)\in V_1'$. The conclusion follows immediately from Proposition \ref{prop:projection_integral_manifolds}.
\end{proof}

\subsection{From smooth to analytic $k$-reflective classical billiards}

In this subsection, we show that the existence of a $k$-pseudo-reflective $\class^{\infty}$-smooth classical billiard implies the existence of a $k$-reflective analytic classical billiard. We also prove that $k$-reflective $\class^{\infty}$-smooth classical billiards can be approximated by $r$-jets of $k$-reflective analytic billiards (a more precise meaning will be given).

To establish these results, we translate Propositions \ref{proposition:link_birkh_reflective} and \ref{proposition:link_birkh_pseudo_reflective} in terms of Pfaffian systems with transversality conditions. Let $\mathcal{D}$ be the classical Birkhoff's distribution defined on the subset $U$ of $M=(\RR^d)^k$ constituted by all $p=(p_1,\ldots,p_k)$ such that $p_{j-1}$, $p_j$, $p_{j+1}$ do not lie on the same line for each $j$ modulo $k$ (see Definition \ref{definition:classical_birkhoff_distribution}). Denote by $\pi_1,\ldots,\pi_k:M\to\RR^d$ the maps given by $\pi_j(p)=p_j$ for each $j$. Given a $2(d-1)$-dimensional vector space $E\subset\mathcal{D}(p)$, we have $\rk d\pi_{|E}\leq d-1$ by construction of Birkhoff's distribution. We consider the following

\begin{lemma}
\label{lemma:transversality_planes}
A $2(d-1)$-dimensional vector space $E\subset\mathcal{D}(p)$ is transverse to $\ker d\pi_j$ if and only if $\rk d\pi_{j|E}= d-1$.
\end{lemma}

\begin{proof}
By Grassmann's formula, we get $\dim (\ker d\pi_{j|\mathcal{D}}+ E)=\rk d\pi_{j|E} +\dim\mathcal D-\rk d\pi_{j|\mathcal{D}}$. Hence $\ker d\pi_{j|\mathcal{D}}+ E=\mathcal{D}$ if and only if $\rk d\pi_{j|E} +\dim\mathcal D=\rk d\pi_{j|\mathcal{D}}=d-1$.
\end{proof}

Lemma \ref{lemma:transversality_planes} implies that in Propositions \ref{proposition:link_birkh_reflective} and \ref{proposition:link_birkh_pseudo_reflective} we can replace the terms non-trivial $2(d-1)$-dimensional integral (respectively pseudo-integral) manifolds of Birkhoff's distribution $\mathcal{D}$ by integral (respectively pseudo-integral) manifolds of the Pfaffian system with transversality conditions given by $\mathcal{P}=(U,\mathcal{D},2(d-1);\ker d\pi_1,\ldots,\ker d\pi_k)$.

\begin{theorem}
Suppose that one can find a local $\mathcal{C}^{\infty}$-smooth $k$-pseudo-reflective classical billiard $(a_1,\ldots,a_k)$. Then given open subsets $V_1,\ldots,V_k$ of $\RR^d$ containing respectively $a_1,\ldots,a_k$, one can find a local analytic $k$-reflective classical billiard $(b_1,\ldots,b_k)$ such that $b_j\subset V_j$ for all $j$.
\end{theorem}

\begin{proof}
By proposition \ref{proposition:link_birkh_pseudo_reflective} and Lemma \ref{lemma:transversality_planes}, one can find a $\class^{\infty}$-smooth pseudo-integral manifold $S$ of the above defined Pfaffian system $\mathcal{P}$, which is contained in $V:=V_1\times\ldots\times V_k$. In particular, $S$ is a $\class^{\infty}$-smooth pseudo-integral manifold of the Pfaffian system $\restreint{\mathcal{P}}{V}:=(V,\mathcal{D},2(d-1);\ker d\pi_1,\ldots,\ker d\pi_k)$. Now by Corollary \ref{cor:existence_of_prolongations} and the theorem of Cartan-Kuranishi-Rachevsky (Theorem \ref{theorem:cartan_ku_ra}), $\restreint{\mathcal{P}}{V}$ posseses an analytic integral manifold. The conclusion follows from Proposition \ref{proposition:link_birkh_reflective} and Lemma \ref{lemma:transversality_planes}.
\end{proof}

In the following, we name by \textit{$k$-reflective set} of a local $k$-reflective billiard $(a_1,\ldots,a_k)$ the set of its $k$-periodic orbits $p=(p_1,p_2,\ldots,p_k)$ for which there is an open subset $U\subset a_1\times a_2$ containing $(p_1,p_2)$ and whose elements $(q_1,q_2)\in U$ can be completed into a $k$-periodic orbit close to $p$.

\begin{theorem}
\label{prop:billiard_analytic_approx}
Suppose that one can find a local $\mathcal{C}^{\infty}$-smooth $k$-reflective classical billiard $(a_1,\ldots,a_k)$. Then given any integer $r>0$, any orbit $p=(p_1,\ldots,p_k)\in a_1\times\ldots\times a_k$ in its $k$-reflective set, and any neighborhoods $V_1,\ldots,V_k\subset\mathcal{G}(\RR^d, d-1, r)$ containing respectively the germs of hypersurfaces $(a_1,p_1),\ldots,(a_k,p_k)$, one can find a local analytic $k$-reflective classical billiard $(b_1,\ldots,b_k)$ and an orbit $q=(q_1,\ldots,q_k)$ in its $k$-reflective set such that $(b_j,q_j)\in V_j$ for all $j$.
\end{theorem}

\begin{remark}
It might be possible that Theorem \ref{prop:billiard_analytic_approx} remains valid for a $k$-pseudo-reflective billiard $(a_1,\dots,a_k)$. The answer to this problem can possibly be found using the results of this manuscript.
\end{remark}

\begin{proof}
Choose $p$ and $V_1,\ldots,V_d$ as in the statement of Theorem \ref{prop:billiard_analytic_approx}. By Proposition \ref{proposition:link_birkh_reflective} and Lemma \ref{lemma:transversality_planes}, one can find a $\class^{\infty}$-smooth integral manifold $S$ of the above defined Pfaffian system $\mathcal{P}$ such that $(\pi_j(S),\pi_j(p))\in V_j$ for each $j$. Let an open set $W\subset\mathcal{G}((\RR^d)^k,2(d-1),r)$ containing the germ $(S,p)$ be such that any germ $(S',q)\in W$ transverse to all $\ker d\pi_j$ satisfies $(\pi_j(S'),\pi_j(q))\in V_j$. By Proposition \ref{prop:analytic_approximation}, the Pfaffian system $\mathcal{P}:=((\RR^d)^k,\mathcal{D},2(d-1);\mathcal{D}_1,\ldots,\mathcal{D}_k)$ posseses an analytic integral manifold $S_a$ for which one can find $p_a\in S_a$ verifying $(S_a,p_a)\in W$. The conclusion follows from Proposition \ref{proposition:link_birkh_reflective}.
\end{proof}

\subsection{From smooth to analytic $k$-reflective projective billiards}

In this subsection, we extend the classical Birkhoff's distribution and its links with the $k$-reflective classical billiards to projective billiards. In the case of projective billiards, the natural space on which the distribution is defined has to be replaced to take into account each field of projective transverse lines.
	
\subsubsection{Projective Birkhoff distribution}

Consider the space $\mathbb{P}(T\RR^d)$, which can be identified as the set of $(p,L)$ such that $p\in\RR^d$ is a point and $L\subset\RR^d$ is a line containing $p$, together with the natural projection $\Pi:\mathbb{P}(T\RR^d)\to\RR^d$. Consider the manifold $M=\mathbb{P}(T\RR^d)^k$ and the constructible subset $U\subset M$ of elements $(p_1,L_1,\ldots p_k, L_k)$, such that for each $j$ (modulo $k$), $(p_j,L_j)\in\mathbb{P}(T\RR^d)$, $p_{j-1}$, $p_j$, $p_{j+1}$ do not lie on the same line, $L_j$ belongs to the plane $p_{j-1}p_jp_{j+1}$ and doesn't coincide with the lines $p_{j-1}p_j$ and $p_{j+1}p_j$. Note that when $d=2$, $U$ is an open dense subset of $M$.

We define for each $j$ the maps $\proj_j:M\to\mathbb{P}(T\RR^d)$ and $\pi_j:M\to\RR^d$ by
$$\proj_j(p_1,L_1,\ldots p_k, L_k)=(p_j,L_j)\qquad\text{ and }\qquad\pi_j(p_1,L_1,\ldots p_k, L_k)=p_j.$$
In what follows, we introduce the analogue of Birkhoff's exterior bissectors for elements in $U$. If $P=(p_1,L_1,\ldots p_k, L_k)\in U$, one can define for each $j$ (modulo $k$) the line $T_j(P)\subset\RR^d$ containing $p_j$ and such that the four lines $p_{j-1}p_j$, $p_{j+1}p_j$, $L_j$ $T_j(P)$ are in the same plane and form a harmonic set of lines. This induces an analytic map $U\to\mathbb{P}(T\RR^d)^k$ which associates to a $P\in U$ the element $(p_1,T_1(P),\ldots,p_k,T_k(P))$.

The analogue of Lemma \ref{lemma:classical_birkhoff_distribution_angles} is given by the following

\begin{lemma}
\label{lemma:projective_birkhoff_distribution_angles}
Let $\mathcal{B}=(\alpha_1,\ldots,\alpha_k)$ be a local projective billiard in $\PP{}{R\RR^d}$ with classical boundaries $a_1,\ldots,a_k$ and $P=(p_1,L_1,\ldots,p_k,L_k)\in U$ such that $(p_j,L_j)\in\alpha_j$ for all $j$. Then $p=(p_1,\ldots,p_k)$ is a $k$-periodic orbit of $\mathcal{B}$ if and only if $T_j(P)\subset T_{p_j}a_j$ for all $j$.
\end{lemma}

\begin{proof}
Fix $j$ and denote by $H$ the plane containing $p_{j-1}$, $p_j$ and $p_{j+1}$. The line $p_{j-1}p_j$ is reflected into the line $p_jp_{j+1}$ by the projective law of reflection at $p_j\in a_j$ if and only if the lines $p_{j-1}p_j$, $p_jp_{j+1}$, $L_j$, $T_{p_j}a_j\cap H$ form a harmonic set of lines in $H$. This is the same as saying that $T_j=T_{p_j}a_j\cap H\subset T_{p_j}a_j$.
\end{proof}

Notice that when $d=2$, the inclusion $T_j\subset T_{p_j}a_j$ is in fact an equality.

\begin{definition}
\label{definition:projective_birk_distrib}
The \textit{projective Birkhoff's distribution} is the analytic map $\mathcal{D}_{proj}:U\to \grass{k}{TM}$ defined for all $P\in U$ by 
$$\mathcal{D}_{proj}(P)=T_1(P)\oplus\ldots\oplus T_k(P).$$
\end{definition}

\subsubsection{Local projective billiards and integral manifolds}

The proofs of Propositions \ref{proposition:link_birkh_reflective} and \ref{proposition:link_birkh_pseudo_reflective} cannot be immediately applied for projective billiards since a $3$-reflective local projective billiard does not correspond to a integral surface of $\mathcal{D}_{proj}$ anymore by Lemma \ref{lemma:projective_birkhoff_distribution_angles}, except in the case when $d=2$. To solve this problem, we consider a version of the $k$-reflective billiard problem in the grassmannian bundle.

Denote by $\restreint{TM}{U}$ the set of $(P,E)\in U\times\grass{2(d-1)}{T_PM}$. We consider the set $M'\subset\grass{2(d-1)}{\restreint{TM}{U}}$ of $2(d-1)$
-dimensional vector spaces $(P,E)\in U\times\grass{2(d-1)}{T_PM}$ satisfying for all $j$ the following conditions:\\
\enum $T_j(P)\subset d\pi_j(E)$;\\
\enum $\rk \restreint{d\pi_j}{E}=d-1$.\\
\enum $\rk \restreint{d\proj_j}{E}=d-1$.\\
In the space of coordinates, the first condition can be expressed as a closed algebraic condition and the second and third ones define a Zariski open subset in an algebraic set. Hence $M'$ is a constructible subset of $\grass{2(d-1)}{\restreint{TM}{U}}$.

We can endow $M'$ with the restriction of the contact distribution defined by $\mathcal{K}'(p,E)=\inverse{d\pi}(E)\cap T_{(p,E)}M'$ for all $(p,E)\in M'$, where $\pi:\grass{2(d-1)}{TM}\to M$ is the natural projection. Consider the Pfaffian system $\mathcal{P}':=(M',\mathcal{K}',2(d-1);\ker d\pi)$.

\begin{proposition}[Analogue of Proposition \ref{proposition:link_birkh_reflective} for projective billiards]
\label{proposition:link_projbill_integral_surfaces}
\textbf{1)} Let $\mathcal{B}$ be a local $\class^{r}$-smooth (resp. analytic) $k$-reflective projective billiard. Then the lifting to $U$ of the set of $k$-periodic orbits of $\mathcal{B}$ contains a $2(d-1)$-dimensional $\class^{r-1}$-smooth (resp. analytic) submanifold $S$ of $U$. The first lift $S^{(1)}$ of $S$ to $\grass{2(d-1)}{TM}$ contains an open dense subset $S'$ which is a $\class^{r-2}$-smooth (resp. analytic) integral manifold of the Pfaffian system $\mathcal{P}'=(M',\mathcal{K}',2(d-1);\ker d\pi)$.

\textbf{2)} Suppose that one can find a $\class^{r}$-smooth (resp. an analytic) integral manifold $S'$ of the Pfaffian system $\mathcal{P}'$ such that the intersection $\ker d\pi$ with $TS'$ is $\{0\}$. Then for $q\in S'$, there is an open subset $W\subset S'$ containing $q$ and such that $(\proj_1\circ\pi(W), \ldots, \proj_k\circ\pi(W))$ is a local $\class^{r}$-smooth (resp. analytic) $k$-reflective projective billiard. 
\end{proposition}

\begin{proof}
\textbf{1)} Write $\mathcal{B}=(\alpha_1,\ldots,\alpha_k)$, and let $a_1,\ldots,a_k$ be its classical boundaries. For each $j$, denote by $L_j(p)$ the projective line of $\alpha_j$ at a point $p\in a_j$. As in Proposition \ref{proposition:link_birkh_reflective}, we can consider the $\class^{r-1}$-smooth (resp. analytic) map $s:(p_1,p_2)\in a_1\times a_2\mapsto (p_3,L_3(p_3),\ldots,p_k,L_k(p_k))\in \alpha_3\times\ldots\times \alpha_k$ such that for each $1<j<k$, $\mathcal{B}_j(p_{j-1},p_j)=(p_j,p_{j+1})$, where $\mathcal{B}_j:a_{j-1}\times a_j\to a_j\times a_{j+1}$ is the projective billiard map. Let $W\subset a_1\times a_2$ be an open subset such that for any $(p_1,p_2)\in W$, the set $(p_1,p_2,p_3,\ldots,p_k)$ is a $k$-periodic orbit of $\mathcal{B}$. Then $S=\ensemble{(p_1,L_1(p_1),p_2,L_2(p_2),s(p_1,p_2))}{(p_1,p_2)\in W}$ is a $2(d-1)$-dimensional $\class^{r-1}$-smooth (resp. analytic) injectively immersed submanifold of $U$. Let $S^{(1)}$ be its first lift to $\grass{2(d-1)}{TM}$. If $(P,E)\in S^{(1)}$, then $P\in S$ and $E=T_PS$. Since the billiard map is a local diffeomorphism (see Proposition \ref{proposition:billiard_regularity_rank}), for each $j$ the maps $\restreint{d\proj_j}{T_PS}$, $\restreint{d\pi_j}{T_PS}$ have rank $d-1$ and $d\pi_j(T_PS)=T_{p_j}a_j$. By Lemma \ref{lemma:projective_birkhoff_distribution_angles}, we can write $T_j(P) \subset T_{p_j}a_j=d\pi_j(T_PS)=d\pi_j(E)$ since $P$ corresponds to a $k$-periodic orbit. Hence $S^{(1)}\subset M'$, and the rest of the proof follows easily.

\textbf{2)} Let $q\in S'$. The transversality condition with $d\pi$ implies that there is an open subset $W'\subset S'$ such that $S:=\pi(W')$ is a $2(d-1)$-dimensional $\class^r$-smooth (resp. analytic) submanifold of $U$. Since $S'$ is an integral manifold of the contact distribution $\mathcal{K}'$, we can write $S^{(1)}=W'$ and for $P\in S$ we have $(P,T_PS)\in M'$. We conclude that $\restreint{d\proj_j}{T_PS}$ and $\restreint{d\pi_j}{T_PS}$ have rank $d-1$ for each $j$, and that $T_j(P)\subset d\pi_j(T_PS)$. The rank condition implies the existence of an open subset $W\subset W'$ containing $q$ such that for each $j$, $\alpha_j=\proj_j\circ\pi(W)$ is a line-framed hypersurface over a hypersurface $a_j=\pi_j\circ\pi(W)$. If $P=(p_1,L_1,\ldots,p_k,L_k)\in \pi(W)$, then $(p_1,\ldots,p_k)$ is a $k$-periodic orbit of $(\alpha_1,\ldots,\alpha_k)$ since for all $j$ we have $(p_j,L_j)\in\alpha_j$, $T_j(P)\subset T_{p_j}a_j=d\pi_j(T_P S)$. Finally the same argument as in the proof of Proposition \ref{proposition:link_birkh_reflective} shows that the map $P\in \pi(W)\mapsto(p_1,p_2)\in a_1\times a_2$ is a local diffeomorphism, hence that $(\alpha_1,\ldots,\alpha_k)$ is $k$-reflective.
\end{proof}

\begin{proposition}[Analogue of Proposition \ref{proposition:link_birkh_pseudo_reflective} for projective billiards]
\label{proposition:link_projbill_pseudo_integral_surfaces}
\textbf{1)} Let $\mathcal{B}$ be a local $\class^{r}$-smooth $k$-pseudo-reflective billiard. Then the 
Pfaffian system $\mathcal{P}'=(M',\mathcal{K}',2(d-1);\ker d\pi)$ has a $\class^{r-2}$-smooth pseudo-integral manifold.

\textbf{2)} Suppose that one can find an analytic pseudo-integral manifold $S'$ of the Pfaffian system $\mathcal{P}'$ such that the intersection $\ker d\pi(p)\cap T_pS'$ is $\{0\}$ for every $p\in S'$. Then for almost all $q$ in the set $V$ of Definition \ref{definition:distribution}, there is an open subset $W\subset S'$ containing $q$ and such that $(\proj_1\circ\pi(W), \ldots, \proj_k\circ\pi(W))$ is a local $\class^{r}$-smooth $k$-pseudo-reflective projective billiard. 
\end{proposition}

\begin{remark}
Notice that the analytic version of this result is given by Proposition \ref{proposition:link_projbill_integral_surfaces}, since $k$-pseudo-reflective analytic projective billiards are $k$-reflective, and connected analytic pseudo-integrable manifolds are integrable (see Remarks \ref{remark:billiard_analytic_pseudoref_reflective} and \ref{remark:manifold_analytic_pseudoref_reflective}).
\end{remark}

\begin{proof}
\textbf{1)} Write $\mathcal{B}=(\alpha_1,\ldots,\alpha_k)$ and denote by $a_1,\ldots,a_k$ its classical boundaries. Let $V_0\subset a_1\times a_2$ be a set of non-zero measure be such that all $(q_1,q_2)\in V_0$ can be completed in a $k$-periodic orbit of $\mathcal{B}$. Let $p=(p_1,\ldots,p_k)$ be a $k$-periodic orbit of $\mathcal{B}$ such that $(p_1,p_2)\in V_0$ is a Lebesgue point of $V_0$. Similarly to the proofs of Propositions \ref{proposition:link_birkh_pseudo_reflective} and \ref{proposition:link_projbill_integral_surfaces}, there is an open subset $U_{(p_1,p_2)}\subset a_1\times a_2$ containing $(p_1,p_2)$ such that the set $S$ of elements $(q_1,L_1(q_1),\ldots,q_k,L_k(q_k))\in\alpha_1\times\ldots\times\alpha_k$ for which $(q_1,\ldots,q_k)$ is a (non-necessarily periodic) orbit of $\mathcal{B}$ with $(q_1,q_2)\in U$ is a $2(d-1)$-dimensional submanifold of $M=\PP{}{T\RR^d}^k$ diffeomorphic to $U_{(p_1,p_2)}$. Let $V\subset S$ be the set of non-zero measure corresponding to $V_0$ in $S$. For $Q=(q_1,L_1(q_1),\ldots,q_k,L_k(q_k))\in S$, the maps $\restreint{d\proj_j}{T_QS}$ and $\restreint{d\pi_j}{T_QS}$ have rank $d-1$ by Lemma \ref{proposition:billiard_regularity_rank}, and if $Q\in V$ we have $Q\in U$ and $T_j(Q)\subset T_{q_j}a_j$ for all $j$ by $k$-periodicity. Therefore any $Q\in V$ is such that $(P,T_PS)\in M'$, hence the first lift $S^{(1)}\grass{2(d-1)}{TM}$ of $S$ contains a subset $V':=\inverse{\pi}(V)\cap S^{(1)}$ of non-zero measure included in $M'$. Now as in the proof of Proposition \ref{prop:pseuo_integral_prolongation}, we can project a neighborhood of $S^{(1)}$ containing a Lebesgue point of $V'$ on an pseudo-integral manifold $S'\subset M'$ of the desired Pfaffian system $\mathcal{P}'$.

 \textbf{2)} As in the proof of Proposition \ref{proposition:link_projbill_integral_surfaces}, if $q$ is a Lebesgue point of $V$, the transversality conditions implies that there is an open subset $W\subset S'$ containing $q$ such that $S:=\pi(W)$ is a $2(d-1)$-dimensional $\class^r$-smooth submanifold of $U$. Let $V_1\subset S$ be the image by $\pi$ of the set $V\cap W$. For $P\in V$, if $(P,E)\in S'$, we have $T_{(P,E)}S'\subset\mathcal{K}'(P,E)$ hence $T_PS=d\pi(T_{(P,E)}S')\subset E$. This implies that for all $j$ and all $P\in V_1$, $\rk\restreint{d\pi_j}{T_PS}=d-1$, $\rk\restreint{d\proj_j}{T_PS}=d-1$ and $T_j(P)\subset d\pi_j(T_PS)$. The first two rank conditions are analytically open conditions satisfied on the subset $V_1\subset S$ of non-zero measure, hence are satisfied on an open dense subset of $S$. Hence by shrinking $W$, one can suppose that for all $j$, the set $\alpha_j:=\proj_j(S)$ is a line-framed hypersurface over the hypersurface $a_j:=\pi_j(S)$. If $Q=(q_1,L_1(q_1),\ldots,q_k,L_k(q_k))\in V_1$, then $(q_1,\ldots,q_k)$ is a $k$-periodic orbit of $(\alpha_1,\ldots,\alpha_k)$. Since $V_1$ has non-zero measure in $S$, the conclusion follows from the same argument as in the proof of Proposition \ref{proposition:link_birkh_pseudo_reflective}.
\end{proof}

\subsubsection{From smooth to analytic $k$-reflective projective billiards}

\begin{theorem}
\label{prop:proj_billiard_pseudo_to_analytic}
Suppose that one can find a local $\mathcal{C}^{\infty}$-smooth $k$-pseudo-reflective projective billiard $(\alpha_1,\ldots,\alpha_k)$. Then given open subsets $V_1,\ldots,V_k$ of $\PP{}{T\RR^d}$ containing respectively $\alpha_1,\ldots,\alpha_k$, one can find a local analytic $k$-reflective projective billiard $(\beta_1,\ldots,\beta_k)$ such that $\beta_j\subset V_j$ for all $j$.
\end{theorem}

\begin{proof}
As in Proposition \ref{proposition:link_projbill_pseudo_integral_surfaces} 1), one can find a $\class^{\infty}$-smooth pseudo-integral manifold $S'$ of the Pfaffian system $\mathcal{P}'=(M',\mathcal{K}',2(d-1);\ker d\pi)$ defined in the proposition, which is contained in the fiber over $V:=V_1\times\ldots\times V_k$. Now by Corollary \ref{cor:existence_of_prolongations} and the theorem of Cartan-Kuranishi-Rachevsky (Theorem \ref{theorem:cartan_ku_ra}), $\restreint{\mathcal{P}'}{V}$ posseses an analytic integral manifold. The conclusion follows from Proposition \ref{proposition:link_projbill_pseudo_integral_surfaces} 2).
\end{proof}

In the following, we name by \textit{$k$-reflective set} of a local $k$-reflective billiard $(\alpha_1,\ldots,\alpha_k)$ the set of its $k$-periodic orbits $p=(p_1,p_2,\ldots,p_k)$ for which there is an open subset $U\subset a_1\times a_2$ containing $(p_1,p_2)$ and whose elements $(q_1,q_2)\in U$ can be completed into a $k$-periodic orbit close to $p$.

\begin{theorem}
\label{prop:proj_billiard_analytic_approx}
Suppose that one can find a local $\mathcal{C}^{\infty}$-smooth $k$-reflective projective billiard $(\alpha_1,\ldots,\alpha_k)$. Then given any integer $r>0$, any element $P=(p_1,L_1(p_1),\ldots,p_k, L_k(p_k))\in \alpha_1\times\ldots\times \alpha_k$ such that $(p_1,\ldots,p_k)$ lies in its $k$-reflective set, and any neighborhoods $V_1,\ldots,V_k\subset\mathcal{G}(\PP{}{T\RR^d}, d-1, r)$ containing respectively the germs of hypersurfaces $(\alpha_1,(p_1,L_1(p_1)))$, $\ldots$, $(a_k,(p_k,L_k(p_k)))$, one can find a local analytic $k$-reflective projective billiard $(\beta_1,\ldots,\beta_k)$ and an orbit $q=(q_1,\ldots,q_k)$ in its $k$-reflective set such that $(\beta_j,(q_j,L_j(q_j)))\in V_j$ for all $j$.
\end{theorem}

\begin{remark}
It might be possible that Theorem \ref{prop:proj_billiard_analytic_approx} remains valid for a $k$-pseudo-reflective projective billiard $(\alpha_1,\dots,\alpha_k)$. The answer to this problem can possibly be found using the results of this manuscript.
\end{remark}

\begin{proof}
Choose $P$ and $V_1,\ldots,V_d$ as in Theorem \ref{prop:proj_billiard_analytic_approx}. As in Proposition \ref{proposition:link_projbill_integral_surfaces} 1), one can find a $\class^{\infty}$-smooth pseudo-integral manifold $S'$ of the Pfaffian system $\mathcal{P}'=(M',\mathcal{K}',2(d-1);\ker d\pi)$ defined in the proposition, such that the $2(d-1)$-dimensional manifold $S:=\pi(S')$ is a set containing $k$-periodic orbits of the projective billiard $(\alpha_1,\ldots,\alpha_k)$. Hence $(\proj_j\circ\pi(S'),\proj_j(P))\in V_j$ for each $j$. Now for $P\in S$, we can choose an open set $W\subset\mathcal{G}(M',2(d-1),r)$ containing the germ of $S'$ at $(P,T_PS)$ such that any germ $(S'_1,P_1)\in W$ is transverse to $\ker d\pi$ and satisfies $(\proj_j\circ\pi(S'),\proj_j(P_1))\in V_j$. By Proposition \ref{prop:analytic_approximation}, the Pfaffian system $\mathcal{P}'$ posseses an analytic integral manifold $S'_a$ for which one can find $(P_a,E)\in S'_a$ verifying $(S'_a,(P_a,E))\in W$. The conclusion follows from Proposition \ref{proposition:link_projbill_integral_surfaces}.
\end{proof}

\section{Triangular orbits of projective billiards}
	\label{section_3_reflective_proj_billiards}
	
In this section, we study the particular case of triangular orbits of projective billiards. More precisely, we investigate the question of classifying the $3$-reflective and $3$-pseudo-reflective local projective billiards. 

As shown in Section \ref{section:examples_reflective_billiards}, given three non-colinear points $P_1,P_2,P_3$ in the Euclidean plane, the right-spherical billiard based at $P_1,P_2,P_3$ is $3$-reflective (see Proposition \ref{prop:right_spherical_3_reflective}). In this section, we present a result classifying the other $3$-reflective projective billiards, see Theorem \ref{theorem:classification_3_projective_billiards}. Let us say that a local projective billiard $\mathcal{B}=(\alpha_1,\alpha_2,\alpha_3)$ of $\PP{}{T\RR^2}$ is \textit{right-spherica}l if it has a $3$-periodic orbit $p=(p_1,p_2,p_3)$ such that each $\alpha_j$ contains a neighborhood of $p_j$ which coincides with the boundary of a right-spherical billiard. Notice that in the analytic case, this is the same as saying that each $\alpha_j$ should be \textit{contained} in the boundary of a right-spherical billiard. This section is devoted to prove the following

\begin{theorem}
\label{theorem:classification_3_projective_billiards}
\textbf{1)} (Planar billiards) If a $\class^{\infty}$-smooth local projective billiard in the Euclidean plane is $3$-reflective, then it is right-spherical.\\
\textbf{2)} (Multidimensional billiards) There are no $\class^{\infty}$-smooth $3$-pseudo-reflective local projective billiards in $\RR^d$ with $d\geq 3$.
\end{theorem}

\subsection{Complex projective billiards}

We can define a complex version for local projective billiards. It consists of definitions analogous to the real case and taking place in the space $\PP{}{T\CC^d}$, considered as the space of pairs $(p,L)$ where $p\in\CC^d$ and $L$ is a complex line of $\CC^d$ containing $p$. Denote by $\pi$ the map $\PP{}{T\CC^d}\to\CC^d$ which associates to a pair $(p,L)\in\PP{}{T\CC^d}$ the point $p$.

\begin{definition}
\label{definition:complex_lf_hypersurfaces}
A \textit{complex line-framed hypersurface} is a $(d-1)$-dimensional connected complex submanifold $\Sigma$ of $\PP{}{T\CC^d}$ such that:\\
\enum $\pi$ is a biholomorphism between $\Sigma$ and a complex hypersurface $S\subset\CC^d$;\\
\enum any pair $(p,L)\in\Sigma$ is such that $L$ is transverse to $S$ at $p$.\\
We say that $\Sigma$ is a line-framed hypersurface \textit{over} $S$. In the case when $d=2$, we say that $\Sigma$ is \textit{complex line-framed curve}.
\end{definition}

\begin{remark}
We can also consider the analogous definition of a complex line-framed hypersurface of $\PP{}{T\CP^d}$ over a complex hypersurface of $\CP^d$.
\end{remark}

The projective law of reflection (see Definition \ref{definition:projective_reflection_law}) can be analogously defined in $\CC^d$ using the same harmonicity conditions on complex lines (see Section \ref{subsection:harmonicity}):

\begin{definition}
\label{definition:projective_reflection_law}
Let $\Sigma$ be a complex line-framed hypersurface over $S$. Let $p\in S$ and $\ell,\ell'$ be complex lines intersecting $S$ at $p$. We say that $\ell'$ is \textit{obtained from $\ell$ by the projective reflection law on $\Sigma$ at $p$} if\\
\enum the lines $\ell$, $\ell'$, $L(p)$ are contained in a complex plane $\mathcal{P}$;\\
\enum the quadruple of lines $\ell,\ell',L(p),T_pS\cap \mathcal{P}$ is harmonic in $\mathcal{P}$.
\end{definition}

\begin{definition}
A complex local projective billiard $\mathcal{B}$ is a collection of complex line-framed hypersurfaces $(\alpha_1,\ldots,\alpha_k)$ over complex hypersurfaces $a_1,\ldots,a_k$ of $\CC^d$ (or $\CP^d$) called \textit{classical boundaries} of $\mathcal{B}$.
\end{definition}

We can define analogously \textit{complex orbits} and complex periodic orbits of $\mathcal{B}$ as in Definition \ref{definition:local_projective_billiard} without the statement on orientation of lines. The notions of \textit{$k$-reflective} complex local projective billiard and $k$-reflective set of such billiard admit also a similar definition. Finally, right-spherical billiards in $\CC^2$ can be defined exactly as in the real case by considering complex lines instead of real lines.

The reason why we introduce complex versions of line-framed hypersurfaces and of local projective billiards is the following: given an analytic line-framed hypersurface $\Sigma$ of $\PP{}{T\RR^d}$, we can consider its complexification $\widehat{\Sigma}$ which is a complex line-framed hypersurface of $\PP{}{T\CC^d}$. Hence given an analytic local projective billiard $\mathcal{B}=(\alpha_1,\ldots,\alpha_k)$ of $\RR^d$ and an orbit $p=(p_1,\ldots,p_k)$ of $\mathcal{B}$, the complexification $\widehat{\alpha_j}$ of each $\alpha_j$ defines a complex line-framed hypersurface in a neighborhood of $\inverse{\pi}(p_j)\cap\alpha_j$. Now if $\mathcal{B}$ is $k$-reflective and $p$ is a periodic orbit and in the $k$-reflective set of $\mathcal{B}$, then by analyticity the complex local projective billiard $\widehat{\mathcal{B}}:=(\widehat{\alpha_1},\ldots,\widehat{\alpha_k})$ is also $k$-reflective.

\subsection{$3$-reflective projective billiards supported by lines}

In this section we show that if a local projective billiard in $\RR^2$ or $\CC^2$ has its classical boundary supported by lines and is $3$-reflective, then it is a right-spherical billiard. We first prove the complex version and then we deduce the real case.

\begin{proposition}
\label{proposition:3_reflective_billiards_on_lines}
Let $\mathcal{B}=(\alpha_1,\alpha_2,\alpha_3)$ be a complex local projective billiard of $\CC^2$ such that its classical boundaries are included in complex lines. If $\mathcal{B}$ is $3$-reflective then it is right-spherical.
\end{proposition}

\begin{proof}
For each $j=1,2,3$, let $\ell_j$ be the line of $\CC^2$ such that $a_j=\pi(\alpha_j)$ is included in $\ell_j$. 

We first show that each $\alpha_j$ can be extended into a complex line-framed hypersurface $\alpha_j'$ over the whole line $\ell_j$. Notice that if such an extension exists it is unique by analyticity. Let $(p_1,p_2,p_3)$ be a $3$-periodic orbit of $\mathcal{B}$ such that $p_1$ is not contained on $\ell_2$ nor $\ell_3$. Given a point $q_2\in\ell_2$, we construct a point $q_3\in\ell_3$ as follows: the line $p_1q_2$ is reflected into a line intersecting $\ell_3$ at a point $q_3$ by the projective law of reflection at $p_1$ with respect to $\alpha_1$. Here maybe $q_3$ lies at infinity with respect to an embedding $\CC^2\subset\CP^2$. Yet the line $q_2q_3$ is well-defined and depends analytically on $q_2$ since $q_3$ depends analytically on $q_2$ by the implicit function theorem. If $q_2$ is such that $p_1,q_2,q_3$ are not on the same line, we can define a unique line $L_2(q_2)$ containing $q_2$ such that the four lines $p_1q_2$, $q_2q_3$, $\ell_2$, $L_2(q_2)$ form a harmonic set of lines. This defines a meromorphic map $s_2:q_2\in\ell_2\mapsto (q_2,L_2(q_2))\in \PP{}{T\CC^2}$. Identifying $\PP{}{T\CC^2}$ with $\CC^2\times\CP^1$, the map $L_2$ can be seen as a holomorphic map $\ell_2\to\CP^1$ hence is defined everywhere. Define $\alpha_2'$ to be the image of $s_2$. Since $\alpha_2'$ coincide with $\alpha_2$ on an open subset by $3$-reflectivity, it contains $\alpha_2$. We can do the same with $\alpha_1$ and $\alpha_3$ defining $\alpha_1'$ and $\alpha_3'$.

The projective billiard maps of $\mathcal{B}'=(\alpha_1',\alpha_2', \alpha_3')$ denoted by $B_1:\ell_1\times \ell_2\to \ell_2\times \ell_3$, $B_2:\ell_2\times \ell_3\to \ell_3\times \ell_1$, $B_3:\ell_3\times \ell_1\to \ell_1\times \ell_2$ are analytic and satisfy $B_3\circ  B_2\circ B_1=Id$ on an open subset of $\ell_1\times \ell_2$ hence on an open dense subset of $\ell_1\times \ell_2$ where the relation is well-defined. This means that all orbits $(p_1,p_2,p_3)$ of $\mathcal{B}'$ are $3$-periodic.

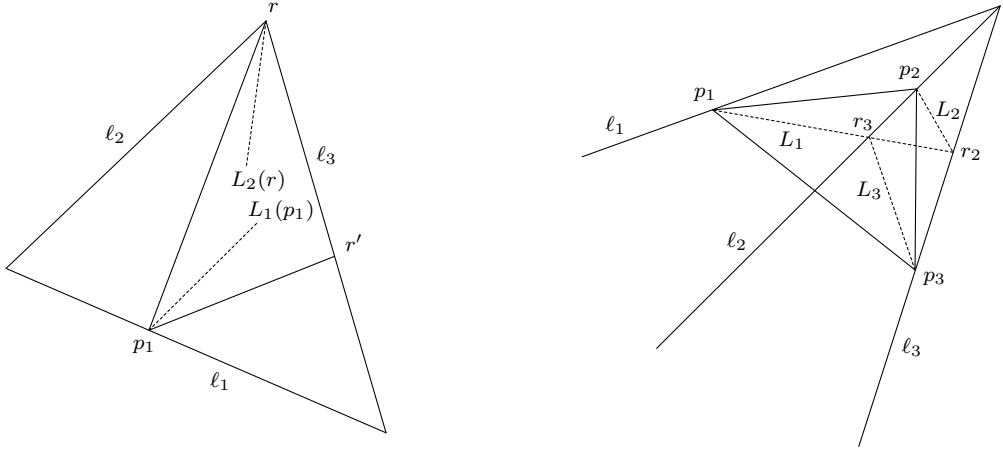
\begin{figure}[!h]
\centering
\begin{tikzpicture}[line cap=round,line join=round,>=triangle 45,x=2cm,y=2cm]
\clip(-1.5,-1.5) rectangle (1.5,2);
\draw (-1.31,0.12)-- (0.4,1.76);
\draw (0.4,1.76)-- (1.19,-0.97);
\draw (1.19,-0.97)-- (-1.31,0.12);
\draw (-0.37,-0.29)-- (0.4,1.76);
\draw (-0.37,-0.29)-- (0.85,0.2);
\draw[dash pattern=on 1pt off 1pt] (-0.37,-0.29)-- (0.35,0.43);
\draw[dash pattern=on 1pt off 1pt] (0.4,1.76)-- (0.27,0.79);
\begin{scriptsize}
\draw[color=black] (0.45,1.85) node {$r$};
\draw[color=black] (-0.4,-0.4) node {$p_1$};
\draw[color=black] (0.98,0.29) node {$r'$};
\draw[color=black] (0.5,0.5) node {$L_1(p_1)$};
\draw[color=black] (0.35,0.7) node {$L_2(r)$};
\draw[color=black] (-0.59,1.01) node {$\ell_2$};
\draw[color=black] (0.8,0.87) node {$\ell_3$};
\draw[color=black] (0.1,-0.63) node {$\ell_1$};
\end{scriptsize}
\end{tikzpicture}\hspace*{1cm}\begin{tikzpicture}[line cap=round,line join=round,>=triangle 45,x=2.0cm,y=2.0cm]
\clip(-2,-1.5) rectangle (1.3,2);
\draw (-1.52,0.86)-- (1.23,1.86);
\draw (1.23,1.86)-- (-1.03,-0.41);
\draw (1.23,1.86)-- (0.3,-1.06);
\draw (-0.66,1.17)-- (0.68,1.31);
\draw (0.68,1.31)-- (0.67,0.11);
\draw (0.67,0.11)-- (-0.66,1.17);
\draw [dash pattern=on 1pt off 1pt] (0.68,1.31)-- (0.92,0.89);
\draw [dash pattern=on 1pt off 1pt] (0.92,0.89)-- (-0.66,1.17);
\draw [dash pattern=on 1pt off 1pt] (0.67,0.11)-- (0.37,0.99);
\begin{scriptsize}
\draw[color=black] (-1.3,1.07) node {$\ell_1$};
\draw[color=black] (-0.5,0.3) node {$\ell_2$};
\draw[color=black] (0.65,-0.4) node {$\ell_3$};
\draw[color=black] (-0.71,1.27) node {$p_1$};
\draw[color=black] (0.64,1.42) node {$p_2$};
\draw[color=black] (0.8,0.05) node {$p_3$};
\draw[color=black] (1.04,0.88) node {$r_2$};
\draw[color=black] (0.89,1.17) node {$L_2$};
\draw[color=black] (-0.14,0.97) node {$L_1$};
\draw[color=black] (0.33,1.08) node {$r_3$};
\draw[color=black] (0.37,0.64) node {$L_3$};
\end{scriptsize}
\end{tikzpicture}
\caption{The two cases of Proposition \ref{proposition:3_reflective_billiards_on_lines}: on the left, the lines $\ell_1$, $\ell_2$, $\ell_3$ do not intersect at the same point; on the right, they do intersect at the same point. Each transverse line is represented as a dotted line (on the right $L_j$ stands for $L_j(p_j)$).}
\label{figure:proof_prop_3reflect_lines}
\end{figure}

We can consider two cases depending on the position of the lines $\ell_1$, $\ell_2$, $\ell_3$ (see Figure \ref{figure:proof_prop_3reflect_lines}):

\textbf{First case.} \textit{Suppose that $\ell_1$, $\ell_2$, $\ell_3$ intersect at the same point.}\\
Fix $p=(p_1,p_2,p_3)$ a periodic orbit in the $3$-reflective set of $\mathcal{B}$. Let $L_1(p_1)$, $L_2(p_2)$, $L_3(p_3)$ be the respective projective lines of $\alpha_1$, $\alpha_2$, $\alpha_3$ over $p_1$, $p_2$, $p_3$. The lines $L_1(p_1)$ and $L_2(p_2)$ intersects $\ell_3$ at the same point $r_2$ since both quadruples of lines $(p_1p_3,p_1p_2,\ell_1,L_1(p_1))$ and $(p_2p_3,p_2p_1,\ell_2,L_2(p_2))$ are harmonic and the three first lines of one quadruple intersect $\ell_3$ at the same points as the three first lines of the other quadruple. The same argument shows that the lines $L_1(p_1)$ and $L_3(p_3)$ intersects $\ell_2$ at the same point $r_3$. Since $p$ is in the $3$-reflective set, one can find $3$-periodic orbits of the form $(q_1,p_2,q_3)$ and $(q_1,q_2,p_3)$, with for all $j$, $q_j\in\ell_j$ close to $p_j$. In the first case $r_2$ is constant since $p_2$ is fixed, and $r_2$ is contained in $L_1(q_1)$. In the second case, $r_3$ is also constant and is contained in $L_1(q_1)$. Hence, for $q_1$ close to $p_1$, the line $L_1(q_1)$ is constant which is impossible since it should contain $q_1$.

\textbf{Second case.} \textit{Suppose that $\ell_1$, $\ell_2$, $\ell_3$ do not intersect at a the same point.}\\
Let $r$ the point of intersection of $\ell_2$ with $\ell_3$. We show that given $p_1\in\ell_1$, the projective line $L_1(p_1)$ contains $r$. Suppose the contrary: $r\notin L_1(p_1)$. Further suppose that $p_1$ is not contained in $\ell_2$ nor $\ell_3$, and that the quadruple of lines $(p_1r,\ell_3,L_2(r),\ell_2)$ is not harmonic. Consider the intersection point $r'\in\ell_3$ of the line containing $p_1$ and reflected from the line $p_1r$ by the projective law of reflection at $p_1$. Since $r$ is not in the line $L_1(p_1)$, the points $r'$ and $r$ are distinct. Approching $(p_1,r,r')$ by a $3$-periodic orbit, we see that the quadruple of lines $(p_1r,r'r,L_2(r),\ell_2)$ is harmonic. Two cases can happen: either $L_2(r)=\ell_2$, or $L_2(r)\neq\ell_2$. In the first case, three lines of the harmonic quadruple should be the same (see Remark \ref{remark:degenerate_harmonicity}). In the second case, the four lines are pairwise distinct, and the line $p_1r$ is completely determined by the triple of lines $(\ell_2,\ell_3,L_2(r))$, hence does not depend on $p_1$. We get a contradiction in both cases and we conclude that $r$ is contained in $L_1(p_1)$. Using Lemma \ref{proposition:billiard_regularity_rank}, the same argument applied to $p_2$ and $p_3$ shows that $\mathcal{B}$ is right-spherical.
\end{proof}

\begin{corollary}
\label{corollary:regular_proj_billiard_on_lines}
Let $\mathcal{B}=(\alpha_1,\alpha_2,\alpha_3)$ be a $\class^1$-smooth local projective billiard of $\RR^2$ such that its classical boundaries are included in lines. If $\mathcal{B}$ is $3$-reflective and $p=(p_1,p_2,p_3)$ is a $3$-periodic orbit in its $3$-reflective set, then each $\alpha_j$ coincides with the boundary of a right-spherical billiard in a neighborhood of $\inverse{\pi}(p_j)\cap \alpha_j$. 
\end{corollary}

\begin{proof}
If the classical boundaries of $\alpha_1,\alpha_2,\alpha_3$ are contained in lines $\ell_1$, $\ell_2$, $\ell_3$, then as in the proof of Proposition \ref{proposition:3_reflective_billiards_on_lines} one can define for each $j$ an analytic map $s_j:\ell_j\to\PP{}{T\RR^2}$ such that $\im s_j$ and $\alpha_j$ coincide in a neighborhood of $p_j$. Hence $\alpha_j$ is analytic in a neighborhood of $p_j$ and we can consider its complexification $\widehat{\alpha_j}$. The complex local projective biliard $\widehat{\mathcal{B}}:=(\widehat{\alpha_1},\widehat{\alpha_2},\widehat{\alpha_3})$ is also $3$-reflective and its classical boundaries are contained in lines. Hence it is a right-spherical billiard by Proposition \ref{proposition:3_reflective_billiards_on_lines}. This concludes the proof.
\end{proof}

\subsection{Space of $3$-periodic orbits attached to a curve}
	\label{subsection:space_3_periodic_orbits}
	
In this section, we study an anologue idea to \cite{glut} which can be described as follows. Given a complex local projective billiard $\mathcal{B}=(\alpha_1,\alpha_2,\alpha_3)$ in $\CC^2$, we can consider other complex local projective billiards of the form $\mathcal{B}'=(\alpha_1,\alpha_2',\alpha_3')$, that is just with the same first boundary $\alpha_1$. Let us say that such a billiard $\mathcal{B}'$ is a local projective billiard \textit{attached to $\alpha_1$}.

Suppose now that $\mathcal{B}$ is $3$-reflective: we can ask if there is a billiard $\mathcal{B}'$ attached to $\alpha_1$ which is $3$-reflective but different from $\mathcal{B}$. We show in fact that this is the case, and that particularly interesting such billiards appear. The main arguments of this subsection are taken from the theory of analytic distribution and of an analogue of Birkhoff's distribution in the complex projective case.

\subsubsection{Singular analytic distributions}

\newcommand{\sing}{\text{Sing}}

We recall some definitions and properties of singular analytic distributions, which can be found in \cite{glut}.

\begin{definition}[\cite{glut}, Lemma $2.27$]
Let $W$ be a complex manifold, $\Sigma\subset W$ a nowhere dense closed subset, $k\in\{0,\ldots,n\}$ and $\mathcal{D}$ an analytic field of $k$-dimensional planes defined on $W\setminus\Sigma$. We say that $\mathcal{D}$ is a \textit{singular analytic distribution of dimension $k$ and singular set $\sing(\mathcal{D})=\Sigma$} if $\mathcal{D}$ extends analytically to no points in $\Sigma$ and if for all $x\in W$, one can find holomorphic $1$-forms $\alpha_1$,..., $\alpha_p$ defined on a neighborhood $U$ of $x$ and such that for all $y\in U\setminus\Sigma$,
$$\mathcal{D}(y) = \bigcap_{i=1}^p \ker\alpha_i(y).$$
\end{definition}

\noindent Singular analytic distributions can be restricted to analytic subsets:

\begin{proposition}[\cite{glut}, Definition $2.32$]
Let $W$ be a complex manifold, $M$ an irreducible analytic subset of $W$ and $\mathcal{D}$ a singular analytic distibution on $W$ with $M\nsubseteq\sing(\mathcal{D})$. Then there exists an open dense subset $M^o_{reg}$ of point $x\in M_{reg}$ for which
$$\mathcal{D}_{|M}(x) := \mathcal{D}(x)\cap T_xM$$
has minimal dimension. We say that $\mathcal{D}_{|M}$ is a singular analytic distribution on $M$ of singular set $\sing(\mathcal{D}):=M\setminus M^o_{reg}$.
\end{proposition}

\begin{remark}
When $M$ is not irreducible anymore, we still can restrict $\mathcal{D}$ to $M$ by looking at its restriction to each of the irreducible components of $M$.
\end{remark}

\noindent As in the smooth case, we can look for integral surfaces defined by the following

\begin{definition}[\cite{glut}, Definition $2.34$]
Let $\mathcal{D}$ be a $k$-dimensional analytic distribution on an irreducible analytic subset $M$ and $\ell\in\{0,\ldots,k\}$. An \textit{integral $\ell$-surface of $\mathcal{D}$} is a submanifold $S\subset M\setminus\sing(\mathcal{D})$ of dimension $\ell$ such that for all $x\in S$, we have the inclusion $T_xS\subset \mathcal{D}(x)$. The analytic distribution $\mathcal{D}$ is said to be \textit{integrable} if each $x\in M\setminus\sing(\mathcal{D})$ is contained in an integral $k$-surface. (In this case the $k$-dimensional integral surfaces form a holomorphic foliation of the manifold $M\setminus\sing(\mathcal D)$, by Frobenius theorem.)
\end{definition}

\noindent We can finally introduce the following lemma, which will be used in a key result (Corollary \ref{cor:birkhoff_integrable}). We recall here that the analytic closure of a subset $A$ of a complex manifold $W$, is the smallest analytic subset of $W$ containing $A$. We denote it by $\overline{A}^{an}$.

\begin{lemma}[\cite{glut}, Lemma $2.38$]
\label{lemma:sad_frobenius_integrable}
Let $\mathcal{D}$ be a $k$-dimensional singular analytic distribution on an analytic subset $N$ and $S$ be a $k$-dimensional integral surface of $\mathcal{D}$. Then the restriction of $\mathcal{D}$ to $\overline{S}^{an}$ is an integrable analytic distribution of dimension $k$.
\end{lemma}

\noindent The proof is the same as in \cite{glut}:

\begin{proof}
Write $M=\overline{S}^{an}$. First, let us prove that $\mathcal{D}_{|M}$ is $k$-dimensional. Consider the subset 
$$A := \ensemble{x\in M\setminus\sing(\mathcal{D}_{|M})}{ \mathcal{D}(x)\subset T_xM}.$$
It contains $S\setminus\sing(\mathcal{D}_{|M})$, hence its closure, which is an analytic subset of $M$, contains $S$. By definition, $\overline{A}^{an}=M$ which implies that $\mathcal{D}_{|M}$ is $k$-dimensional.

Now let us show that $\mathcal{D}_{|M}$ is integrable. The argument is similar: define the subset $B$ of those $x\in M\setminus\sing(\mathcal{D}_{|M})$ such that the Frobenius integrability condition is satisfied. $B$ contains $S\setminus\sing(\mathcal{D}_{|M})$ and its closure is an analytic subset of $M$ containing $S$, hence it is the whole $M$. Thus Frobenius theorem can be applied on the manifold $M\setminus\sing(\mathcal{D}_{|M})$, which implies the result.
\end{proof}

\subsubsection{Birkhoff's distribution and the $3$-reflective billiard problem}
\label{subsec:birk_distrib_1}

In this section we define an analogue of Birkhoff's distribution in the case of complex local projective billiards attached to a fixed line-framed curve $\alpha$. We give an analogue of Proposition \ref{proposition:link_birkh_reflective} for such billiards at Proposition \ref{prop:_birkh_billiard_attached}.

We first define the \textit{space of the distribution}. Let $\mathcal{L}$ be the fiber bundle 
$$\mathcal{L}= \PP{}{T\CP^2}\underset{\CP^2}{\times}\PP{}{T\CP^2}$$ 
that is the set of triples $(p,L,T)$ where $p\in\CP^2$ and $L,T$ are lines in $T_p\CP^2$. Consider the space $\alpha\times\mathcal{L}\times\mathcal{L}$ of triples $P=(P_1,P_2,P_3)$ where $P_1=(p_1,L_1)\in\alpha$, $P_2=(p_2,L_2,T_2)\in\mathcal L$, $P_3=(p_3,L_3,T_3)\in\mathcal{L}$. 

Consider the subspace $M_{\alpha}^{0}$ of $3$-periodic billiard orbits having one reflection in $\alpha$, that is the set of elements $P\in\alpha\times\mathcal{L}\times\mathcal{L}$ such that the points $p_1,p_2,p_3$ do not lie on the same line and the quadruples of lines $(p_1p_2,p_1p_3,L_1,T_{p_1}\alpha)$, $(p_2p_3,p_2p_1,L_2,T_2)$, $(p_3p_1,p_3p_2,L_3,T_3)$ form harmonic sets of distinct lines. Denote by $M_{\alpha}$ the analytic closure of $M_{\alpha}^{0}$. 

We use the same notations for the different projections as in Section \ref{section_pfaffian_systems}:\\
\enum $\proj_1:\alpha\times\mathcal{L}^2\to\alpha$, $\proj_2,\proj_3:\alpha\times\mathcal{L}^2\to\PP{}{T\CP^2}$ defined for all $P\in \alpha\times\mathcal{L}^2$ and all integers $j=1,2,3$ by $\proj_j(P)=(p_j,L_j)$;\\
\enum $\pi_1,\pi_2,\pi_3:\alpha\times\mathcal{L}^2\to\CP^2$ defined for all $P\in \alpha\times\mathcal{L}^2$ and all integer $j=1,2,3$ by $\pi_j(P)=p_j$.
\begin{definition}
We call \textit{Birkhoff's distribution attached to $\alpha$} the restriction $\mathcal{D}_{\alpha}$ to $M_{\alpha}$ of the analytic distribution $\mathcal{D}$ defined for all $P\in\alpha\times\mathcal{L}^2$ by
$$\mathcal{D}(P) = d\inverse{\pi_2}(T_2)\cap d\inverse{\pi_3}(T_3).$$
\end{definition}

\begin{proposition}[Analogue of Proposition \ref{proposition:link_birkh_reflective} for $\mathcal{D}_{\alpha}$]
\label{prop:_birkh_billiard_attached}
Let $P\in M_{\alpha}^0$ such that one can find a $2$-dimensional integral analytic surface $S$ of $D_{\alpha}$ containing $P$. Suppose that for each $j=1,2,3$ the restrictions of $\proj_j$ and $\pi_j$ to $S$ have rank $1$. Then there exists a neighborhood $U$ of $P$ in $S$ such that the complex local projective billiard $(\alpha,\proj_2(U),\proj_3(U))$ is $3$-reflective.
\end{proposition}

\begin{proof}
Let $U\subset S$ be an open subset such that $\alpha_2:=\proj_2(U)$ and $\alpha_3:=\proj_3(U)$ are complex line-framed curves over the complex curves $a_2:=\pi_2(U)$ and $a_3:=\pi_3(U)$. Since $S$ is an integral surface of $\mathcal{D}$, for $Q=(q_1,L_1,q_2,L_2,T_2,q_3,L_3,T_3)\in U$ we can write $T_{q_2}a_2=d\pi_2(T_QS)\subset T_2$, hence $T_{q_2}a_2=T_2$ and it follows that the quadruple of lines $(q_2q_1,q_2q_3,L_2,T_{q_2}a_2)$ is harmonic, and a similar argument can be applied to the lines through $q_3$. Hence $(q_1,q_2,q_3)$ is a $3$-periodic orbit of $\mathcal{B}:=(\alpha,\alpha_2,\alpha_3)$.

It remains to show that $\mathcal{B}$ is $3$-reflective. Indeed, write $a=\pi_1(\alpha)$ and let us show that the projection $j:U\to a\times a_2$ onto $(q_1,q_2)$ has rank $2$ in a neighborhood of $P$. Denote by $s:a\times a_2\to \alpha\times\mathcal{L}^2$ the map defined by $s(q_1,q_2)=(q_1,L_1(q_1),q_2,L_2(q_2),T_{p_2}a_2,p_3,L_3(q_3),T_{q_3}a_3)$ where $q_3$ is the point of intersection with $a_3$ of the line reflected from $q_1q_2$ by the projective reflection on $\alpha_2$, $L_1(q_1)$ is the projective line of $\alpha$ at $q_1$, $L_2(q_2)$ is the projective line of $\alpha_2$ at $q_2$ and $L_3(q_3)$ is the projective line of $\alpha_3$ at $q_3$. The map $s$ is defined in a neighborhood of $(p_1,p_2)$ and satisfies $s\circ j(Q)=Q$ for all $Q\in U$ close to $P$. Hence $j$ has rank $2$ in a neighborhood of $P$ and therefore $\mathcal{B}$ is $3$-reflective.  
\end{proof}

\subsubsection{Reduction of the space of orbits}
\label{subsec:reduction_space_orbits}

In this subsection, we suppose that we are given a complex local projective billiard $\mathcal{B}=(\alpha_1,\alpha_2,\alpha_3)$ with classical boundaries $a_1$, $a_2$, $a_3$, which is $3$-reflective, and we investigate the structure of complex local projective billiards attached to $\alpha:=\alpha_1$.

Since $\mathcal{B}$ is $3$-reflective, there is a $2$-dimensional integral surface $S$ of $\mathcal{D}_{\alpha}$ in $M_{\alpha}^0$ such that for each $P\in S$, $(p_1,p_2,p_3)=(\pi_1(P),\pi_2(P),\pi_3(P))$ is a $3$-periodic orbit of $\mathcal{B}$ (this is an easy consequence of the arguments detailed in the proof of Proposition \ref{prop:_birkh_billiard_attached}). Denote by $\hat{S}$ the analytic closure of $S$ in $M_{\alpha}$. In this subsection we want to prove that $\dim \hat{S}\leq 4$.

To achive this result on dimension, we first construct two analytic subsets $M_{\alpha,2}$ and $M_{\alpha,3}$ of $\alpha\times\mathcal{L}^2$ containing $S$. Consider an element $P=(p_1,L_1,p_2,L_2,T_2,p_3,L_3,T_3)\in M_{\alpha}^{0}$. By the implicit function theorem, we can define an analytic map $j_{(p_2,L_2,T_2)}$ on a neighborhood of $p_3$ in $a:=\pi(\alpha)$, with values in $\CP^2$, as follows: if $q_1\in a$, the line $\ell_1$ obtained from $q_1p_2$ by the projective law of reflection at $p_1$, and the line $\ell_2$ such that the quadruple of lines $(q_1p_2,\ell_2,L_2,T_2)$ is harmonic, intersect at a point $q_3$, and we set $j_{(p_2,L_2,T_2)}(q_1)=q_3$. The map $j_{(p_2,L_2,T_2)}$ is analytic and obviously non-constant with $j_{(p_2,L_2,T_2)}(p_1)=p_3$. Its image is an irreducible germ of analytic curve at $p_3$, and we can consider the latter's tangent line at $p_3$ denoted by $T_{p_1}j_{(p_2,L_2,T_2)}$. 

Let $M_{\alpha,3}\subset M_{\alpha}$ be the analytic closure of the set $\ensemble{P\in M_{\alpha}^0}{T_3=T_{p_1}j_{(p_2,L_2,T_2)}}$. We can analogously define $M_{\alpha,2}\subset M_{\alpha}$ by exchanging the roles of $p_3$ and $p_2$.

\begin{proposition}
The analytic closure $\hat{S}$ of $S$ is contained in $M_{\alpha,2}\cap M_{\alpha,3}$.
\end{proposition}

\begin{proof}
We only have to show that $S$ is contained in $M_{\alpha,2}\cap M_{\alpha,3}$. If $P\in S$, the image of a neighborhood of $p_1$ by $j_{(p_2,L_2,T_2)}$ is contained in $a_3=\pi(\alpha_3)$ by $3$-reflectivity of the local projective billiard $\mathcal{B}$. Hence $T_{p_1}j_{(p_2,L_2,T_2)}$ coincide with the tangent line $T_{p_3}a_3=T_3$, and therefore $P\in M_{\alpha,3}$. The same argument applied to $M_{\alpha,2}$ implies the result.
\end{proof}

We can now prove that $\dim\hat{S}\leq 4$. Consider the set $F$ of triples $(p_1,p_2,p_3)\in a\times\CP^2\times\CP^2$ such that the points $p_1$ and $p_2$ are contained in a line which is reflected into a line containing $p_3$ by the projective reflection law at $p_1$ on $\alpha$. $F$ is an analytic set of dimension $4$, as one can easily see. Now the map $s:\alpha\times\mathcal{L}^2\to a\times\CP^2\times\CP^2$ which associates to $P$ the triple $p=(p_1,p_2,p_3)$ is such that $s(S)\subset F$, hence $s(\hat{S})$ is an analytic subset of~$F$. 

\begin{proposition}
\label{prop:finite_fibers}
The map $s:\hat{S}\to F$ has generically finite fibers, in the following sense: there exists an open dense subset $U\subset s(\hat{S})$ (a complement to a proper analytic subset) such that $\inverse{s}(p)$ is finite for every $p\in U$. In particular $\dim \hat{S}\leq 4$.
\end{proposition}

\begin{proof}
Consider the open dense subset $U\subset F$ of triples $(p_1,p_2,p_3)$ in $F$ for which $p_1,p_2,p_3$ are not on the same line and the points $p_2,p_3$ are not contained in the line $T_{p_1} a$ nor the projective line $L_1(p_1)$ of $\alpha$ at $p_1$.

Consider $p=(p_1,p_2,p_3)\in U$ and suppose that the fiber $\inverse{s}(p)$ is not finite. By construction, $\inverse{s}(p)$ is an analytic subset contained in $\{(p_1,L_1(p_1))\}\times\{p_2\}\times\PP{}{T_{p_2}\CP^2}^2\times\{p_3\}\times\PP{}{T_{p_3}\CP^2}^2$, hence it is algebraic by Chow's theorem (see \cite{GH}). Since $\inverse{s}$ is not finite, at least one of the projection from $\inverse{s}(p)$ to $L_2$, $T_2$, $L_3$ or $T_3$ is infinite. Without loss of generality we suppose that the projection to $T_2$ is infinite: the image of such projection is an infinite analytic subset in $\PP{}{T_{p_2}\CP^2}\simeq\CP^1$, hence is the whole $\CP^1$ by Chow's theorem.

Therefore we can consider an element $P=(p_1,L_1,p_2,L_2,T_2,p_3,L_3,T_3)\in\inverse{s}(p)$ for which $T_2$ is different from the lines $p_1p_2$ and $p_2p_3$ with usual identification. By the same argument, we can consider another element $P'\in\inverse{s}(p)$ of the form $P'=(p_1,L_1,p_2,L_2',L_2,p_3,L_3',T_3')$, that is the projection on $T_2$ of $P'$ gives previous $L_2$. Since $P$ and $P'$ are in $M_{\alpha}$, the quadruple of lines $(p_1p_2,p_2p_3, L_2,T_2)$ and $(p_1p_2,p_2p_3, L_2',L_2)$ are harmonic, hence $L_2'=T_2$.

Since $P\in M_{\alpha,3}$, $T_3$ is defined by the relation $T_3=T_{p_1}j_{(p_2,L_2,T_2)}$. Applying the same argument to $P'$ we get $T_3'=T_{p_1}j_{(p_2,T_2,L_2)}$. Now, permuting $L_2$ and $T_2$ doesn't change the projective reflection law at $p_2$, and $j_{(p_2,L_2,T_2)}=j_{(p_2,T_2,L_2)}$, therefore $T_3=T_3'$. The harmonicity conditions at $p_3$ implies that $L_3=L_3'$.

Thus we just proved that $P=(p_1,L_1,p_2,L_2,T_2,p_3,L_3,T_3)$ and $P'=(p_1,L_1,p_2,T_2,L_2,p_3,L_3,T_3)$. But if we consider now that $P,P'\in M_{\alpha,2}$, by the same arguments we get that $T_2=T_{p_1}j_{(p_1,L_3,T_3)}=L_2$. This contradicts the harmonicity condition of the quadruple of lines $(p_1p_2,p_1p_3,L_2,T_2)$. Hence $\inverse{s}(p)$ is finite.
\end{proof}

Given a point $p_1\in a=\pi(\alpha_1)$, we denote by $\hat{S}_{p_1}$ the set $\inverse{\pi_1}(p_1)\cap\hat{S}$. It is algebraic by Chow's theorem.

\begin{lemma}
\label{lemma:epimorphic_projections}
Suppose $\dim \hat{S}\geq 3$. Then for all $p_1$ lying outside a discrete subset of $a$ we have either $\pi_2(\hat{S}_{p_1})=\CP^2$ or $\pi_3(\hat{S}_{p_1})=\CP^2$.
\end{lemma}

\begin{proof}
For $j=2,3$ and any $p_1\in a$, the set $\pi_j(\hat{S}_{p_1})\subset\CP^2$ is algebraic by Chow's theorem and contains the classical boundary $a_j=\pi(\alpha_j)$, hence it has dimension at least $1$. Now since the map $\pi_1:\hat{S}\to a$ is surjective, for $p_1$ lying outside a discrete subset $a^{\ast}$ of $a$ the algebraic set $\hat{S}_{p_1}$ has dimension at least $2$. And exactly as in the proof of Proposition \ref{prop:finite_fibers}, the restriction of $s$ to $\hat{S}_{p_1}$ has generically finite fibers. Hence if $p_1\notin a^{\ast}$ we have $\dim s(\hat{S}_{p_1})\geq 2$.

Suppose that $\pi_2(\hat{S}_{p_1})$ has dimension $1$. Hence for all points $p_2$ in an open and dense subset of $\pi_2(\hat{S}_{p_1})$, the fiber over $p_2$ of the projection $s(\hat{S}_{p_1})\to\pi_2(\hat{S}_{p_1})$ sending $(p_1,p_2,p_3)$ to $p_2$ has dimension $1$. By definition of $F$, this fiber over a fixed point $p_2$ is contained in the set of triples $(p_1,p_2,p_3)$ for which $p_3$ is in a line $\ell$ determined by $p_1$ and $p_2$: $\ell$ is obtained from the line $p_1p_2$ by the projective law of reflection on $\alpha$ at $p_1$. Hence the fiber over $p_2$ contains all triples $(p_1,p_2,p_3)$ where $p_3\in\ell$. Therefore $\pi_3(\hat{S}_{p_1})$ contains all lines $\ell$ obtained by the projective law of reflection from a line $p_1p_2$ where $p_2\in\pi_2(\hat{S}_{p_1})$.

If $a_2$ is not contained in a line, there are infinitely many such lines $\ell$ and we get $\pi_3(\hat{S}_{p_1})=\CP^2$. If $a_2$ is contained in a line, we get the same result by choosing $p_1$ outside this line.
\end{proof}

\subsubsection{Integrability of Birkhoff's distribution on $\hat{S}$}

In this subsection, we suppose that we are given a complex local projective billiard $\mathcal{B}=(\alpha_1,\alpha_2,\alpha_3)$ with classical boundaries $a_1$, $a_2$, $a_3$, which is $3$-reflective. We consider the restriction to $\hat{S}$ of Birkhoff's distribution attached to $\alpha_1$, denoted by $\mathcal{D}_{\hat{S}}$. We first compute the dimension of $\mathcal{D}_{\hat{S}}$ and then we show that it is integrable.

\begin{proposition}
\label{prop:dimension_birkhoff_ditrib}
The singular analytic distribution $\mathcal{D}_{\hat{S}}$ is $2$-dimensional.
\end{proposition}

\begin{proof}
We first have $\dim \mathcal{D}_{\hat{S}}\geq\dim S = 2$ since $T_PS\subset \mathcal{D}_{\hat{S}}(P)$ for $P\in S$. By Proposition \ref{prop:finite_fibers}, $2\leq\dim\hat{S}\leq 4$ and so is $\dim \mathcal{D}_{\hat{S}}$. We consider two cases: $\dim \hat{S}=3$ and $\dim \hat{S}=4$. In both cases, we consider a regular point $P=(p_1,L_1,p_2,L_2,T_2,p_3,L_3,T_3)$ of $\hat S$ such that $\dim \mathcal{D}_{\hat{S}}(P)$ is minimal. We can further suppose that $P\in M_{\alpha}^0$ since $\hat S\cap M_{\alpha}^0$ is an open dense subset of $\hat S$.

\textbf{Case when $\dim \hat{S}=3$.} We have to find one vector $u\in T_P \hat S$ which is not in $\mathcal{D}_{\hat{S}}(P)$. By Lemma \ref{lemma:epimorphic_projections} we can suppose that $p_1$ satisfies without loss generality $\pi_{2}(\hat{S}_{p_1}) = \CP^2$. Hence there is a path $u(t)\in \hat{S}$ such that $u(0)=P$ and $\pi_2\circ u(t)$ is contained in the the line $p_1p_2$ with non-zero derivative at $0$. The vector $u'(0)$ of $T_P\hat S$ is such that $d\pi_2\cdot u'(0)$ is a non-zero vector directed along the line $p_1p_2$. Hence
$d\pi_2\cdot u'(0)\notin T_2$ since $T_2\neq p_1p_2$ because $P\in M_{\alpha}^0$. We conclude that $u'(0)\notin \mathcal{D}_{\hat{S}}(P)$.

\textbf{Case when $\dim \hat{S}=4$.} Let us find two linearly independent vectors $u,v\in T_P\hat{S}$ such that $\mathcal{D}_{\hat{S}}(P)$ and the plane spanned by $(u,v)$ intersect by $\{0\}$. We can suppose that $p_1$ is such that $\hat{S}_{p_1} = 3$, and as in the proof of Proposition \ref{prop:finite_fibers} that $s(\hat{S}_{p_1})=F\cap\left(\{p_1\}\times (\CP^2)^2\right)$. Let $u(t)\in\hat S$ be a path such that $u(0)=P$, $\pi_2\circ u(t)$ belongs to the line $p_1p_2$  with non-zero derivative at $0$ and $\pi_3\circ u(t)$ is contant equal to $p_3$. Now exchange the role of $p_2$ and $p_3$ and define similarly a path $v(t)$ such that $v(0)=P$, $\pi_3\circ v(t)\in p_1p_3$ with non-zero derivative at $0$ and $\pi_2\circ v(t)=p_2$.

Then $u'(0)$ and $v'(0)$ are linearly independant since $d\pi_2\cdot u'(0)\neq 0$ and $d\pi_3\cdot v'(0)\neq 0$ while $d\pi_3\cdot u'(0)=0$ and $d\pi_2\cdot v'(0)=0$. Moreover, if there are $\lambda,\mu\in\CC$ such that $\lambda u'(0)+\mu v'(0)\in \mathcal{D}_{\hat{S}}(P)$, then $\lambda d\pi_2\cdot u'(0)= d\pi_2(\lambda u'(0)+\mu v'(0))\in T_2$ by the definition of $\mathcal{D}_{\hat{S}}$. Thus $\lambda=0$ since $p_1p_2\neq T_2$. Similarly $\mu=0$ and this concludes the proof.
\end{proof}

\noindent As a consequence of Proposition \ref{prop:dimension_birkhoff_ditrib} and Lemma \ref{lemma:sad_frobenius_integrable}, we can state the following

\begin{corollary}
\label{cor:birkhoff_integrable}
The singular analytic distribution $\mathcal{D}_{\hat{S}}$ is integrable.
\end{corollary}

\subsection{Proof of Theorem \ref{theorem:classification_3_projective_billiards} for analytic planar billiards}
\label{sec:border_lines}

In this section, we prove the following result:

\begin{proposition}
\label{proposition:reflective_cmplx_proj_billiards}
Let $\mathcal{B}$ be a complex local projective billiard of $\CC^2$. If $\mathcal{B}$ is $3$-reflective, then its classical boundaries are contained in lines.
\end{proposition}

We deduce a proof of Theorem \ref{theorem:classification_3_projective_billiards} in the case of complex local projective billiards from this result and from Proposition \ref{proposition:3_reflective_billiards_on_lines}:

\begin{theorem}[Complex version of Theorem \ref{theorem:classification_3_projective_billiards} case 1.]
\label{theorem:classification_3_projective_analytic_billiards}
Let $\mathcal{B}$ be a complex local projective billiard of $\CC^2$ or an analytic local projective billiard of $\RR^2$. If $\mathcal{B}$ is $3$-reflective, then it is right-spherical.
\end{theorem}

\begin{proof}
By complexification, we can suppose that $\mathcal{B}$ is a complex local projective billiard of $\CC^2$ which is $3$-reflective. By Proposition \ref{proposition:reflective_cmplx_proj_billiards} the classical boundaries of $\mathcal{B}$ are contained in lines. This implies that $\mathcal{B}$ is right-spherical by Proposition \ref{proposition:3_reflective_billiards_on_lines}.
\end{proof}

The idea of the proof is as follows: let $\mathcal{B}=(\alpha_1,\alpha_2,\alpha_3)$ be a complex local projective billiard with classical boundaries denoted by $a_1,a_2,a_3$ and suppose that $\mathcal{B}$ is $3$-reflective. If one of the classical boundaries, say $a_1$, is not contained in a line, then we consider the complex local projective billiards \textit{attached} to $\alpha_1$ (as defined at the beginning of Section \ref{subsection:space_3_periodic_orbits}). We show that the existence of $\mathcal{B}$ implies the existence of a complex local projective billiard attached to $\alpha_1$ having what we call \textit{a one-parameter family of flat orbits} defined below (Definition \ref{definition:3_reflective_flat_orbit}). We finally show that if such a billiard has this property, then $a_1$ is contained in a line.

\begin{definition}
\label{definition:3_reflective_flat_orbit}
Let $\mathcal{B} = (\alpha_1,\alpha_2,\alpha_3)$ be a complex local projective billiard with classical boundaries $a_1,a_2,a_3$.\\
\enum We say that $\mathcal{B}$ has \textit{a one-parameter family of flat orbits} if there is a an open subset $V\subset a_1$ such for all points $p_1\in V$, the tangent line $T_{p_1}a_1$ intersects $a_2$ at a point $p_2$ and $a_3$ at a point $p_3$ depending continuously on $p_1$, and verifying the following property: there is a sequence of $3$-periodic orbits of $\mathcal{B}$ of the form $(p_1,q_2^n,q_3^n)$ converging to $(p_1,p_2,p_3)$ and belonging to the $3$-reflective set of $\mathcal{B}$.\\
\enum Any triple $(p_1,p_2,p_3)$  as above is called an \textit{$\alpha_1$-flat orbit of $\mathcal{B}$.}
\end{definition}

Note that the property of having a one-parameter family of flat orbits can be found on right-spherical billiards. They are the only such analytic billiards as Theorem \ref{theorem:classification_3_projective_billiards} shows.

\subsubsection{Existence of a particular $3$-reflective local projective billiard}
\label{subsec:existence_integral}

Let $\mathcal{B}=(\alpha_1,\alpha_2,\alpha_3)$ be a complex local projective billiard with classical boundaries $a_1,a_2,a_3$. Suppose that $\mathcal{B}$ is $3$-reflective and that $a_1$ is not contained in a line. 

In what follows, we use the following definition: given two curves $\gamma,\gamma'\subset\CP^2$ and point $p\in\gamma$, $p,p'\in\gamma'$, we say that \textit{the germs $(\gamma,p)$ and $(\gamma',p')$ coincide} if $p=p'$ and there is an open subset $U\subset\CP^2$ containing $p=p'$ such that $\gamma\cap U=\gamma'\cap U$.

\begin{figure}[h]
\centering
\input{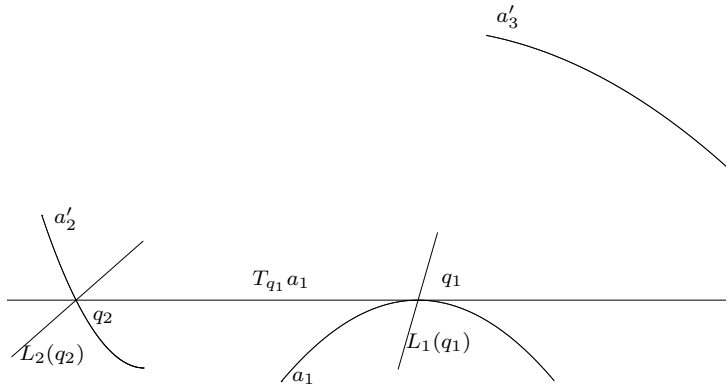}
\caption{The local projective billiard in the second case of Proposition \ref{prop:first_step_flat_orbit}: $T_{q_1}a_1$ intersects $a_2'$ transversally at $q_2$.}
\label{figure:classification_1}
\end{figure}

\begin{proposition}
\label{prop:first_step_flat_orbit}
There is a $3$-reflective complex local projective billiard $\mathcal{B}'=(\alpha_1,\alpha_2',\alpha_3')$ with classical boundaries $a_1$, $a_2'$, $a_3'$, and points $q_1\in a_1$, $q_2\in a_2'$, such that one of the following cases holds:\\
1) The germs of curves $(a_1,q_1)$ and $(a_2',q_2)$ coincide;\\
2) $q_1\neq q_2$ and $T_{q_1}a_1$ intersects $a_2'$ transversally at $q_2$ (see Figure \ref{figure:classification_1}).
\end{proposition}

\begin{proof}
Let $S$ be the $2$-dimensional integral surface $S$ of $\mathcal{D}_{\alpha_1}$ in $M_{\alpha_1}^0$ such that for each $P\in S$, $(p_1,p_2,p_3)=(\pi_1(P),\pi_2(P),\pi_3(P))$ is a $3$-periodic orbit of $\mathcal{B}$, and denote by $\hat{S}$ the analytic closure of $S$ in $M_{\alpha}$. By Corollary \ref{cor:birkhoff_integrable}, the restriction $\mathcal{D}_{\hat S}$ of $\mathcal{D}_{\alpha_1}$ to $\hat S$ is integrable.

Consider the subset $\hat{S}^0\subset \hat S$ consisting of points $P$ of $\hat S\cap M_{\alpha}^0$ outside the singular set of $\mathcal{D}_{\hat S}$ for which the restrictions of $d\proj_j(P)$ and $d\pi_j(P)$ is of rank $1$: $\hat{S}^0$ is an Zariski-open dense subset of $\hat S$ since these conditions are given by anlytically open relations which are satisfied on $S$. 

If $p_1\in a_1$, the set $\hat{S}_{p_1}^0:=\hat{S}_{p_1}\cap\hat{S}^0$ is also a Zariski-open dense subset of $\hat{S}_{p_1}$ satisfying the following property resulting from Propositions \ref{cor:birkhoff_integrable} and \ref{prop:_birkh_billiard_attached}: if $P\in\hat{S}_{p_1}^0$, there is a local projective billiard of the form $(\alpha_1,\alpha_2',\alpha_3')$ which is $3$-reflective and for which $(\pi_1(P)=p_1,\pi_2(P),\pi_3(P))$ is a $3$-periodic orbit.

\begin{lemma}
Let $p_1\in a_1$. There is a point $q_1\in a_1$ which can be chosen arbitrary close to $p_1$ such that $T_{q_1}a_1$ intersects $\pi_2(\hat{S}_{p_1}^0)$ at a point $q_2$ distinct from $q_1$.
\end{lemma}

\begin{proof}
By Chow's theorem, $\pi_2(\hat{S}_{p_1})$ is an algebraic subset of $\CP^2$ which contains the classical boundary $a_2$ of $\alpha_2$. By Chevalley's theorem, $\pi_2(\hat{S}_{p_1}^0)$ is a constructible dense subset of $\pi_2(\hat{S}_{p_1})$. Now since $a_1$ is not a line, the map $q_1\mapsto T_{q_1}a_1$ is not a constant map. Hence if $\pi_2(\hat{S}_{p_1})=\CP^2$, there are points $q_1$ arbitrary close to $p_1$ such that $T_{q_1} a_1$ contains an open dense subset of points $q_2\in\pi_2(\hat{S}_{p_1}^0)$ which are not in $a_1$. If $\dim\pi_2(\hat{S}_{p_1})=1$, the algebraic set $\pi_2(\hat{S}_{p_1})\smallsetminus \pi_2(\hat{S}_{p_1}^0)$ is finite. Since $a_1$ is not a line, we can choose $q_1$ close to $p_1$ such that $T_{q_1}a_1$ doesn't intersect $\pi_2(\hat{S}_{p_1})\smallsetminus \pi_2(\hat{S}_{p_1}^0)$. By Bezout's theorem, $T_{q_1}a_1$ intersects $\pi_2(\hat{S}_{p_1})$ hence $\pi_2(\hat{S}_{p_1}^0)$.
\end{proof}

Choose a point $q_1\in a_1$ close to $p_1$ and a point $P\in\hat{S}_{p_1}^0$ such that $p_2:=\pi_2(P)$ is contained in $T_{q_1} a_1$. By Propositions \ref{cor:birkhoff_integrable} and \ref{prop:_birkh_billiard_attached} there is a local projective billiard of the form $(\alpha_1,\alpha_2',\alpha_3')$ with projective boundaries $(a_1,a_2',a_3')$ which is $3$-reflective and for which $(p_1,p_2,p_3):=(\pi_1(P),\pi_2(P),\pi_3(P))$ is a $3$-periodic orbit. By construction, $p_2$ is contained in $T_{q_1}a_1$ and $a_2'$ at the same time. If the germs of curves $(a_1,q_1)$ and $(a_2',p_2)$ coincide, set $q_2:=p_2$ and there is nothing more to do. Otherwise, we can change $q_1$ for a point arbitrary close to $q_1$ such that $q_1\neq p_2$ and $T_{q_1}a_1$ intersects $a_2'$ transversally at a point $q_2$ close to $p_2$ and different from $q_1$. This concludes the proof.
\end{proof}

\begin{figure}[h]
\centering
\input{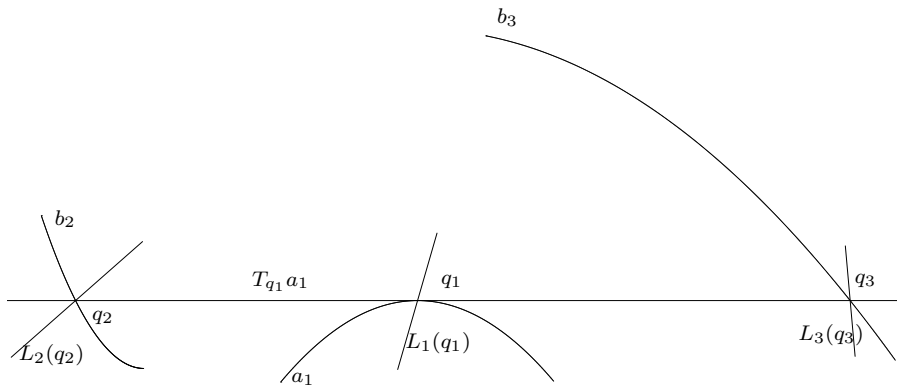}
\caption{The local projective billiard of Proposition \ref{prop:existence_flat_orbit}. Here the three points $q_1$, $q_2$, $q_3$ are pairwise distinct.}
\label{figure:classification_2}
\end{figure}

\begin{proposition}
\label{prop:existence_flat_orbit}
There is a $3$-reflective complex local projective billiard $(\alpha_1,\beta_2,\beta_3)$ attached to $\alpha_1$, with classical boundaries $a_1$, $b_2$, $b_3$, having a one-parameter family of flat orbits. We can find a $\alpha_1$-flat orbit $(q_1,q_2,q_3)$ with the following properties (see Figure \ref{figure:classification_2}):

1) The points $q_1$, $q_2$, $q_3$ lies on $T_{p_1}a_1$.

2) If two points among $\{q_1,q_2,q_3\}$ coincide, then the corresponding classical borders coincide. 

3) $T_{p_1}a_1$ intersect $b_2$ transversally at $q_2$ if $q_1\neq q_2$, and $b_3$ transversally at $q_3$ if $q_1\neq q_3$.
\end{proposition}

\begin{proof}
Let $\mathcal{B}'=(\alpha_1,\alpha_2',\alpha_3')$ be the local projective billiard from Proposition \ref{prop:first_step_flat_orbit}, with classical boundaries $a_1$, $a_2'$, $a_3'$. Let $q_1,q_2$ be the points from Proposition \ref{prop:first_step_flat_orbit}.

We first define a meromorphic map from $a_1\times a_2'$ to $\CP^2$ as follows: let $(p_1,p_2)\in a_1\times a_2'$ such that $p_1\neq p_2$ and consider the point $p_3\in\CP^2$ of intersection of the lines $\ell_1$ and $\ell_2$, where $\ell_1$ is the line reflected from $p_1p_2$ by the projective law of reflection on $\alpha_1$ at $p_1$, and $\ell_2$ is the line reflected from $p_1p_2$ by the projective law of reflection on $\alpha_2'$ at $p_2$. The map $p_3:(p_1,p_2)\mapsto p_3(p_1,p_2)\in\CP^2$ is a meromorphic map, hence is well-defined and analytic outsie a discret subset of $a_1\times a_2'$ (an analytic subset of codimension $2$). Hence by eventually moving $q_1$ a little, one can suppose that the map $p_3$ is analytic at $(q_1,q_2)$ and we write $q_3=p_3(q_1,q_2)$.

Moreover, the map $p_3(p_1,p_2)$ has rank one on an open subset of $a_1\times a_2'$ since $\mathcal{B'}$ is $3$-reflective, hence it is of rank one on an open dense subset of $a_1\times a_2'$ and sends a small neighborhood of $(q_1,q_2)$ into an analytic curve $b_3$ of $\CP^2$ intersecting $T_{q_1} a_1$ at $q_3$. Hence for $(p_1,p_2)$ in a neighborhood of $(q_1,q_2)$, we can define the line $L_3(p_3)$ containing $p_3$ and such that the quadruple of lines $(p_1p_3,p_2p_3,L_3(p_3),T_{p_3}b_3)$ is harmonic. By the same argument, on an open dense subset the map has rank one and $L_3(p_3)\neq T_{p_3}b_3$. Again by moving $q_1$ a little we can suppose that $L_3(q_3)\neq T_{q_3}b_3$ and that if the germs $(a_1,p_1)$ and $(b_3,q_3)$ do not coincide, then $p_1\neq q_3$, and the same with $(a_2',q_2)$ instead of $(a_1,p_1)$. Hence the image $\beta_3$ of the map $(p_1,p_2)\mapsto (p_3,L_3)$ is a complex line-framed curve with classical boundary $b_3$.

By construction, if we denote by $\beta_2=\alpha_2'$ the line-framed curve over $b_2:=a_2'$, then $(\alpha_1,\beta_2,\beta_3)$ is the desired $3$-reflective complex local projective billiard.
\end{proof}

\subsubsection{The $3$-reflective local projective billiard of Proposition \ref{prop:existence_flat_orbit} cannot exist}
\label{subsec:integral_not_reflective}

Let $\mathcal{B}=(\alpha_1,\alpha_2,\alpha_3)$ be a complex local projective billiard with classical boundaries $a_1,a_2,a_3$. Suppose that $\mathcal{B}$ is $3$-reflective and that $a_1$ is not contained in a line. Let $\mathcal{B}_0=(\alpha_1,\beta_2,\beta_3)$ be the $3$-reflective local projective billiard from Proposition \ref{prop:existence_flat_orbit} with classical boundaries $a_1$, $b_2$, $b_3$. In this subsection we show that the existence of $\mathcal{B}_0$ is impossible (under the already made assumption that $a_1$ is not a line).

\begin{figure}[h]
\centering
\input{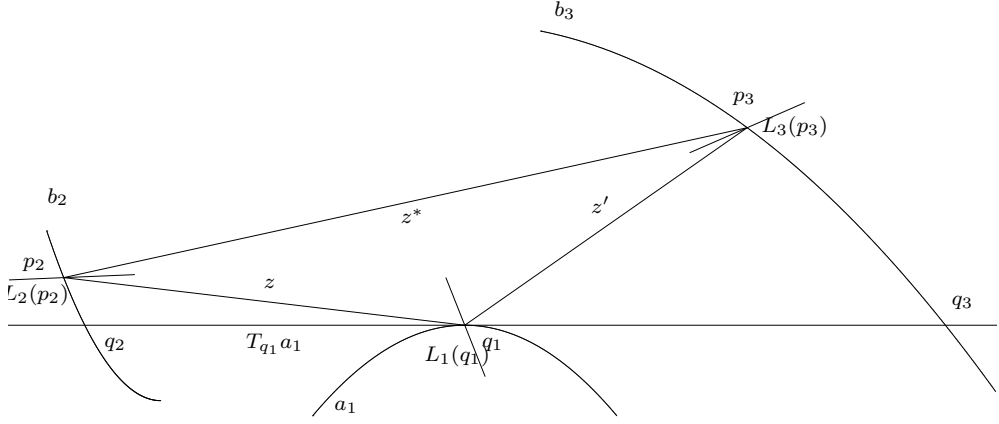}
\caption{The local projective billiard $\mathcal{B}_0$ with an orbit $(q_1,p_2,p_3)$.}
\label{figure:classification_3}
\end{figure}

Let $(q_1,q_2,q_3)$ be the $\alpha_1$-flat orbit of Proposition \ref{prop:existence_flat_orbit}. Denote by $L_1(p_1)$, $L_2(p_2)$, $L_3(p_3)$ the fields of projective lines respectively on $a_1$, $b_2$, $b_3$. Choose an affine chart $\mathbb C^2\subset\mathbb{CP}^2$ containing the points $q_1$, $q_2$, $q_3$ and a coordinate on the line $L_{\infty}=\CP^2\smallsetminus\CC^2$ such that 
$$\az(T_{q_1}a_1)=0\qquad\text{ and }\qquad\infty\notin\{\az(T_{q_2}b_2),\az(T_{q_3}b_3),\az(L_1(q_1))\}$$
where $\az(\ell)$ is the coordinate of the intersection point $L\cap L_{\infty}$ of a line $\ell$ with $L_{\infty}$ (see Subsection \ref{subsection:harmonicity}). When considering a $3$-periodic orbit of the form $(q_1,p_2,p_3)$, we will write
$$z=\az(q_1p_2), \quad z^{\ast}=\az(p_2p_3),\quad z'=\az(q_1p_3)$$
and find asymptotic relations on $z,z^{\ast},z'$ when $(q_1,p_2,p_3)$ is close to $(q_1,q_2,q_3)$ (see Figure \ref{figure:classification_3}, and section \ref{subsection:harmonicity} for further details on azimuths).

\begin{proposition}
\label{prop:three_equivalences}
When $(q_1,p_2,p_3)$ is close to $(q_1,q_2,q_3)$, the following asymptotic equivalence relations are satisfied: 
$$z' \sim (-z),\qquad z^{\ast} \sim (2I_2-1)z,\qquad z^{\ast}\sim (2I_3-1)z'$$
where $I_2$ (respectively, $I_3$) is the intersection index of $b_2$ (respectively $b_3$) with the tangent line $T_{q_1}a_1$ at $q_2$ (respectively $q_3$).
\end{proposition}

\noindent From Proposition \ref{prop:three_equivalences}, we deduce that $2I_3-1=-(2I_2-1)$ which is impossible since $2I_2-1$ and $2I_3-1$ are strictly positive integers. Hence $\mathcal{B}_0$ cannot exist. 

We will prove the three asymptotic relations of Proposition \ref{prop:three_equivalences} in what follows, separated in three propositions (Propositions \ref{prop:equivalence_A}, \ref{prop:intersection_index} and \ref{prop:equivalence_azimuths}).

\begin{proposition}[$z'\sim (-z)$]
\label{prop:equivalence_A}
When $p_2\in b_2$ goes to $q_2$, we have
$$z'\sim (-z).$$
\end{proposition}

\begin{proof}
Equation \eqref{equation:projective_transformation_harmonic} in Section \ref{subsection:harmonicity} implies that 
$$z'=\frac{(\ell+t)z-2\ell t}{2z-(\ell+t)}$$
where $t = \az(T_{q_1}a_1)$, $\ell=\az(L_{1}(q_1))$. In the chosen set of coordinates we have, when $p_2\to q_2$,
$$z'=\frac{\ell z}{2z-\ell}\sim\frac{\ell z}{-\ell}=-z.$$
\end{proof}

\begin{proposition}
\label{prop:intersection_index}
If $q_1=q_2$ then when $p_2\in b_2$ goes to $q_2$, we have
$$z^{\ast}\sim (2I_2-1)z$$
where $I_2\geq2$ is the index of intersection of $a$ with the tangent line $T_{q_1}a_1$ at $A_1$.
\end{proposition}

\begin{proof}
In the case when $q_2=q_1$, the germs $(b_2,q_2)$ and $(a_1,q_1)$ coincide as prescribed in Proposition \ref{prop:existence_flat_orbit}. Take a $3$-periodic orbit of the form $(q_1,p_2,p_3)$ close to $(q_1,q_2,q_3)$. Write $t = \az(T_{p_2}b_2)$, $\ell=\az(L_2(p_2))$. Equation \eqref{equation:projective_transformation_harmonic} in Section \ref{subsection:harmonicity} implies that 
$$\frac{z^{\ast}}{z}=\frac{(\ell+t)z-2\ell t}{z(2z-(\ell+t))}.$$
Now, when $p_2\to q_2$, since $a_1$ and $b_2$ coincide in a neighborhood of $q_2$, we can compute that $t\sim Iz$. Thus
$$\frac{z^{\ast}}{z}\sim \frac{(1-2I)\ell z}{-\ell z}=2I_2-1.$$
\end{proof}

\begin{lemma}
\label{lemma:distinct_points1}
If $q_2=q_3$, then the germs $(a_1,q_1)$, $(b_2,q_2)$ and $(b_3,q_3)$ coincide.
\end{lemma}

\begin{proof}
Suppose that the three germs do not coincide and that $q_2=q_3$: as prescribed in Proposition \ref{prop:existence_flat_orbit} we should have $(b_3,q_3)=(b_2,q_2)$ but $q_1\neq q_2$ with $T_{p_1}a_1$ intersecting $b_2$ transversally at $q_2$ by Proposition \ref{prop:existence_flat_orbit}. Consider a $3$-periodic orbit of the form $(q_1,p_2,p_3)$ close to $(q_1,q_2,q_3)$. Then write $t = \az(T_{p_2}b_2)$ and $\ell=\az(L_2(p_2))$. Remark \ref{remark:order_harmonicity} and Equation \eqref{equation:projective_transformation_harmonic} in Section \ref{subsection:harmonicity} imply that 
$$\ell=\frac{(z+z^{\ast})t-2zz^{\ast}}{2t-(z+z^{\ast})}.$$
Now, when $p_2\to q_2$, we have $z\to 0$ and $t\to t_0$ where $t_0:=\az(T_{q_2}b_2)\notin\{0,\infty\}$ by transversality of $b_2$ with $T_{q_1}a_1$ at $q_2$ and by choice of coordinates. But we also have $z^{\ast}\to t_0$
because $p_2p_3\to T_{q_2}b_2$ since $p_2$,$p_3$ are distinct points of the same irreducible germ of curve $b_2=b_3$ converging to the same point $q_2=q_3$. Hence, when $p_2\to q_2$,
$$\ell\to \frac{t_0^2}{t_0}=t_0$$
which means that $L_{2}(q_2)=T_{q_2}b_2$. But this is not the case by Proposition \ref{prop:existence_flat_orbit}, contradiction.
\end{proof}

\begin{proposition}
\label{prop:equivalence_azimuths}
Suppose that $q_2\neq q_1$. Then when $p_2\in b_2$ goes to $q_2$, we have
$$z^{\ast}\sim z$$
which allows to extend the formula of Proposition \ref{prop:intersection_index} by setting $I_2=1$ in this case (transverse intersection).
\end{proposition}

\begin{proof}
First, let us prove the following lemma, which gives the form of the projective field of lines locally around $q_2$:

\begin{lemma}
\label{lemma:tangent_field_of_lines}
Suppose that $q_2\neq q_1$. Then when $p_2\in b_2$ is close to $q_2$, there is a $p_1\in a_1$ close to $q_1$ for which $L_2(p_2)$ is tangent to $a_1$ at $p_1$.
\end{lemma}

\begin{proof}
Proposition \ref{prop:existence_flat_orbit} implies that $T_{q_1}a_1$ intersects $b_2$ transversally at $q_2$. By the implicit function theorem, there is an analytic map which associates to any $p_1$ close to $q_1$ a point $p_2\in b_2$ close to $q_2$ which is contained in $T_{p_1}a_1$. Since $a_1$ is not a line, this map is not constant, hence is open and thus parametrizes (maybe non-bijectively) the germ of $b_2$ at $q_2$. We choose $p_1$ in the neighborhood of $q_1$, and denote by $p_2$ the corresponding point on $b_2$ obtained via the above parametrization. 

We can suppose that that $T_{p_1}a_1$ is transverse to $b_2$ at $p_2$ and that $p_2\notin b_3$ (possible by Lemma \ref{lemma:distinct_points1}). Consider a $3$-periodic orbit of the form $(p_1,p_2',p_3')$ such that $p_2'$ converges to $p_2$. The line $p_1p_2'$ converges to $T_{p_1}a_1$ and by the projective reflection law at $p_1$ we get that the line $p_1p_3'$ also converges to $T_{p_1}a_1$, hence the limit $p_3$ of $p_3'$ lies on $T_{p_1}a_1$. 
We also have that $p_3\neq p_1$ (Lemma \ref{lemma:distinct_points1}). 
Hence $p_2p_3=T_{p_1}a=p_1p_2$: $T_{p_1}a_1$ is invariant by the projective reflection law of $\beta_2$ at $p_2$. Since the latter tangent line is transverse to $b_2$, we have $T_{p_1}a_1=L_{2}(p_2)$, and this concludes the proof.
\end{proof}

We now conclude the proof of Proposition \ref{prop:equivalence_azimuths}. As in Lemma \ref{lemma:tangent_field_of_lines}, when $p_2\in b_2$ is close to $q_1$, $L_2(q_2)$ is tangent to $a_1$ at a point $p_1$ close to $q_1$. Write $t = \az(T_{p_2}b_2)$, $\ell=\az(L_2(p_2))$. We have, by Equation \eqref{equation:projective_transformation_harmonic} in Section \ref{subsection:harmonicity}, 
\begin{equation}
\label{eq:limit_symetry}
z^{\ast} = \frac{(\ell+t)z-2\ell t}{2z-(\ell+t)}.
\end{equation}
Now in this configuration, we easily compute using Lemma \ref{lemma:tangent_field_of_lines} that, when $p_2\to q_2$, 
$$\ell\sim z.$$ 
Here besides Lemma \ref{lemma:distinct_points1}, we essentially use the inequality $q_2\neq q_1$. This allows to use the following argument. As $p_2$ tends to $q_2$, the lines $L_2(p_2)$ and $L_{2}(q_2)=T_{q_1}a_1$ intersect at a point converging to $q_1$, while $p_2$ remains distant from $q_1$. This implies the required asymptotic equivalence of azimuths.

But we have also $t\to t_0$ where $t_0=\az(T_{q_2}b_2)\notin\{0,\infty\}$ (by the transversality condition of the intersection with $T_{q_1}a_1$). Hence, Equation \eqref{eq:limit_symetry} implies, when $p_2\to q_2$, that
$$\frac{z^{\ast}}{z} = \frac{(\ell+t)z-2t\ell}{z(2z-(\ell+t))}\sim\frac{-t_0z}{-t_0z}=1.$$
\end{proof}

\begin{proof}[Proof of Proposition \ref{prop:three_equivalences}]
The first asymptotic equivalence, $z\sim-z'$, comes from Proposition \ref{prop:equivalence_A}. The second one comes from Proposition \ref{prop:intersection_index}, when $I_2\geq 2$, and from Proposition \ref{prop:equivalence_azimuths}, when $I_2=1$. Finally the third one can be deduced from the second one by interchanging the germ of curves $(b_2,q_2)$ and $(b_3,q_3)$.
\end{proof}

\subsubsection{Proof of Proposition \ref{proposition:reflective_cmplx_proj_billiards}}

We finally prove Proposition \ref{proposition:reflective_cmplx_proj_billiards} which will complete the proof of Theorem \ref{theorem:classification_3_projective_analytic_billiards}.

Let $\mathcal{B}=(\alpha_1,\alpha_2,\alpha_3)$ be a complex local projective billiard with classical boundaries $a_1,a_2,a_3$. Suppose that $\mathcal{B}$ is $3$-reflective and that one curve among $\{a_1,a_2,a_3\}$, say for example $a_1$, is not contained in a line. Let $\mathcal{B}_0=(\alpha_1,\beta_2,\beta_3)$ be the $3$-reflective local projective billiard from Proposition \ref{prop:existence_flat_orbit}.

Using the same notations as in Proposition \ref{prop:three_equivalences}, we deduce that $2I_3-1=-(2I_2-1)$ which is impossible since $2I_2-1$ and $2I_3-1$ are strictly positive integers. Hence $\mathcal{B}_0$ cannot exist, contradiction: $a_1$ is contained in a line. By symmetry of previous argument, $a_1$, $a_2$, an $a_3$ are contained in lines, which proves Proposition \ref{proposition:reflective_cmplx_proj_billiards}.

\subsection{Proof of Theorem \ref{theorem:classification_3_projective_billiards}: planar case}

In this section we give a proof of Theorem \ref{theorem:classification_3_projective_billiards} in the case of planar projective billiards, that is case 1). Let $(\alpha_1,\alpha_2,\alpha_3)$ be a $\class^{\infty}$-smooth local projective billiard of $\RR^2$ which is $3$-reflective. Let $a_1,a_2,a_3$ be its classical boundaries. 

Let $p=(p_1,p_2,p_3)$ be in its $3$-reflective set. By Theorem \ref{prop:proj_billiard_analytic_approx}, one can find an analytic $3$-reflective local projective billiard $\mathcal{B}_a$ with classical boundaries $b_1,b_2,b_3$ and a $3$-periodic orbit $q$ of $\mathcal{B}_a$ such that for each $j$, the germs of curves $(b_j,q_j)$ and $(a_j,p_j)$ are arbitrary close in the Withney $\class^r$-topology (see Definition \ref{definition:3_reflective_flat_orbit}), for a fixed integer $r>0$. By Theorem \ref{theorem:classification_3_projective_analytic_billiards}, the germs $(b_j,q_j)$ are germs of lines. 

Since this is also true for points close to $p_1$, $p_2$ and $p_3$ in $a_1$, $a_2$ and $a_3$ respectively (by definition of $3$-reflectivity and by Proposition \ref{proposition:billiard_regularity_rank}), the latter curves coincide with lines on open subsets. The conclusion follows from Corollary \ref{corollary:regular_proj_billiard_on_lines}.

\subsection{Proof of Theorem \ref{theorem:classification_3_projective_billiards}: multidimensional case}
\label{sec:multidim}

In this section we give a proof of Theorem \ref{theorem:classification_3_projective_billiards} in the case of local projective billiards in dimension $d\geq 3$, that is case 2). As in the planar case, we first prove a complex analytic version of the theorem, and then we use Pfaffian systems to extend the result to $\class^{\infty}$-smooth local projective billiards.

\begin{theorem}[Complex version of Theorem \ref{theorem:classification_3_projective_billiards} case 2.]
\label{theorem:classification_pseudo_reflective_multi}
Let $d\geq 3$. There are no complex local projective billiards in $\CC^d$ which are $3$-reflective.
\end{theorem}

If we admit this theorem, we can prove the multidimensional case of Theorem \ref{theorem:classification_3_projective_billiards}:

\begin{proof}[Proof of Theorem \ref{theorem:classification_3_projective_billiards}]
Suppose that one can find a $\class^{\infty}$-smooth local projective billiard $(\alpha_1,\alpha_2,\alpha_3)$ in $\RR^d$ which is $3$-pseudo-reflective. Theorem \ref{prop:proj_billiard_pseudo_to_analytic} implies the existence of an analytic local projective billiard which is $3$-reflective. By complexification, this contradicts Theorem \ref{theorem:classification_pseudo_reflective_multi}.
\end{proof}

We first prove this auxiliary lemma:

\begin{lemma}
\label{lemma:lines_and_tangent}
Let $W\subset\CC^d$ be a complex hypersurface, $p\in W$ and $U$ a non-empty open subset of $\PP{}{T_pW}$. Suppose that for any $v\in T_pW$ with $\pi(v)\in U$, the hypersurface $W$ contains the points $p+tv$ for all $t\in \CC$ in a neighborhood of $0$ depending on $v$. Then $W$ is a hyperplane.
\end{lemma}

\begin{proof}
We can suppose that $p=0$, $T_pW={z_d=0}$ and $W$ is locally the graph of an analytic map $f:V\to\CC$ where $V\subset\CC^{d-1}$ is an open subset containing $0$. Let $v\in\CC^{d-1}$ be a non-zero vector such that $\PP{}{v}\in U$. By assumption, for $t$ close to $0$ we have $g_v(t):=f(tv)=0$. Since $g_v$ is analytic, it is $0$ everywhere where it is defined. Yet the set $\{tv| t\in\RR, \PP{}{v}\in U\}$ contains a non-empty open subset of $V$, on which $f$ should vanish. By analyticity $f=0$ and $W$ is the hyperplane defined by the equation $z_d=0$.
\end{proof}

\begin{proof}[Proof of Theorem \ref{theorem:classification_pseudo_reflective_multi}]
Suppose that we can find a $3$-reflective complex local projective billiard $\mathcal{B}=(\alpha_1,\alpha_2,\alpha_3)$ in $\PP{}{T\CC^d}$ with classical boundaries $a_1,a_2,a_3$. Denote by $L_1(p_1)$, $L_2(p_2)$ and $L_3(p_3)$ the field of projective lines respectively of $\alpha_1$, $\alpha_2$ and $\alpha_3$ at $p_1\in a_1$, $p_2\in a_2$, $p_3\in a_3$. Let $U\times V\subset a_1\times a_2$ be an open subset such that all $(p_1,p_2)\in U\times V$ can be completed in $3$-periodic orbits of $\mathcal{B}$. Let us state the following obvious result:

\begin{lemma}
\label{lemma:planar_orbits}
Let $(p_1,p_2,p_3)$ be a $3$-periodic orbit of $\mathcal{B}$. Then all lines $p_1p_2$, $p_2p_3$, $p_3p_1$, $L_1(p_1)$, $L_2(p_2)$, $L_3(p_3)$ belong to the plane $p_1p_2p_3$, which is transverse to $a_1,a_2,a_3$ at $p_1,p_2,p_3$ respectively.
\end{lemma}

First let us show the

\begin{lemma}
\label{lemma:hypersurf_planes}
The hypersurfaces $a_1$ and $a_2$ are contained in hyperplanes.
\end{lemma}

\begin{proof}
By symmetry, let us just show that $a_1$ is supported by a hyperplane. Fix $p_1\in U$. For $p_2\in V$, consider the plane $\mathcal{P}(p_2)$ containing the triangular orbit starting by $(p_1,p_2)$, as in Lemma \ref{lemma:planar_orbits}. Consider $a_1(p_2)$, $a_2(p_2)$, $a_3(p_2)$ to be the intersections of $\mathcal{P}(p_2)$ respectively with $a_1$, $a_2$, $a_3$: by transversality, and shrinking them if needed, we can suppose that they are immersed curves of $\mathcal{P}(p_2)$.

Now consider for each $j=1,2,3$ the curve $\alpha_j(p_2)=\pi^{-1}(a_j(p_2))\cap\alpha_j$. Let us show that $\mathcal{B}(p_2)=(\alpha_1(p_2),\alpha_2(p_2),\alpha_3(p_2))$ is a planar $3$-reflective projective billiard. Consider the open subsets $U'=U\cap a_1(p_2)$ of $a_1(p_2)$ and $V'=V\cap a_2(p_2)$ of $a_2(p_2)$. Any $q_2\in V'$ is such that $(p_1,q_2)$ can be completed in a $3$-periodic orbit $(p_1,q_2,q_3)$ of $\mathcal{B}$ and by Lemma \ref{lemma:planar_orbits}, $p_1q_2q_3$ is a plane containing $L_1(p_1)$ and $p_1q_2$, which are intersecting lines inside $\mathcal{P}(p_2)$. Hence $p_1q_2q_3=\mathcal{P}(p_2)$ and thus $\alpha_2(p_2)$ is an analytic curve such that for all $q_2\in V'$, the point $q_2$ and $L_2(q_2)$ are in $\mathcal{P}(p_2)$. The same argument work for $\alpha_1(p_2)$, and also for $\alpha_3(p_2)$ by \ref{proposition:billiard_regularity_rank}. This implies that $\mathcal{B}(p_2)$ is a $3$-reflective local projective billiard inside $\mathcal{P}(p_2)$.

In particular, by Theorem \ref{theorem:classification_3_projective_analytic_billiards}, $a_1(p_2)$ is contained in a line denoted by $\ell(p_2)$ which is itself included in $T_{p_1}a_1$ (since the tangent space of $a_1(p_2)$ is included in the tangent space of $a_1$) and in $\mathcal{P}(p_2)$. Hence $a_1$ intersect $\ell(p_2)$ in an open subset of $\ell(p_2)$ containing $p_1$. This result is true for any $p_2\in V$, implying the same result for lines in a neighborhood of $\ell(p_2)$ in $T_{p_1}a_1$: hence by Lemma \ref{lemma:lines_and_tangent}, $a_1$ is supported by an hyperplane, which concludes the proof.
\end{proof}

Let $H_1$ be the hyperplane containing $a_1$ and $H_2$ be the hyperplane containing $a_2$.

\begin{lemma}
\label{lemma:3reflect_analytic_multi_projlines}
There is a point $q_2\in H_2$ such that for all $p_1\in a_1$ the line $L_1(p_1)$ goes through $q_2$.
Similarly, there is a point $q_1\in H_1$ such that for all $p_2\in a_2$ the line $L_2(p_2)$ goes through $q_1$.
\end{lemma}

\begin{proof}
Let us show the existence of $q_2$, the existence of $q_1$ being analogous. Fix $p_1\in U$, and consider the point $q_2\in H_2$ of intersection of $L_1(p_1)$ with $H_2$. For $p_2\in V$, consider the plane $\mathcal{P}(p_2)$ containing the triangular orbit starting by $(p_1,p_2)$, as in Lemma \ref{lemma:planar_orbits}: define $a_1(p_2)$, $a_2(p_2)$, $a_3(p_2)$, $\alpha_1(p_2)$, $\alpha_2(p_2)$, $\alpha_3(p_2)$, $U'$, $V'$  as in the proof of Lemma \ref{lemma:hypersurf_planes}. One has $q_2\in a_2(p_2)\subset \mathcal{P}(p_2)$, by Lemma \ref{lemma:planar_orbits}. We recall that $(\alpha_1(p_2),\alpha_2(p_2),\alpha_3(p_2))$ is a planar $3$-reflective complex local projective billiard. 

By Theorem \ref{theorem:classification_3_projective_analytic_billiards} it is a right-spherical billiard, hence each $p_1'\in U'$ is such that $L_1(p_1')$ and $L_1(p_1)$ intersect $a_2(p_2)$ at the same point which is $q_2=L_1(p_1)\cap a_2(p_2)$ by construction. Therefore, any $p_1'\in U'$ is such that $L_1(p_1')$ passes through $q_2$. Hence by analyticity, if $\ell(p_2)$ is the line of intersection of $\mathcal{P}(p_2)$ with $a_1$, every $p_1'\in a_1\cap\ell(p_2)$ is such that $L_1(p_1')$ passes through $q_2$.

Now the union of all $\ell(p_2)$ for $p_2\in V$ contains a non-empty open subset $\Omega$ of $a_1$, which by construction has the following property: all $p_1'\in \Omega$ is such that $L_1(p_1')$ passes through $q_2$. By analyticity, this is also true for all $p_1'\in a_1$, and the proof is complete.
\end{proof}

Now we can finish the proof of Theorem \ref{theorem:classification_pseudo_reflective_multi}. Indeed, any $p=(p_1,p_2)\in U\times V$, can be completed in a $3$-periodic orbit which lies in a plane $\mathcal{P}(p)$. This plane $\mathcal{P}(p)$ contains $L_1(p_1)$ and $L_2(p_2)$ (Lemma \ref{lemma:planar_orbits}), hence goes through $q_1$ and $q_2$ as in Lemma \ref{lemma:3reflect_analytic_multi_projlines}. 
If $q_1\neq q_2$, $P(p) = p_1q_1q_2$, but this is impossible since in this case $P(p_2)$ doesn't depend on $p_2)$, which therefore can be chosen outside the plane $p_1q_1q_2$. Hence $q_1=q_2\in H_1$, implying that all $p_1\in U$ are such that $L_1(p_1)\subset T_{p_1}a_1$. This contradicts the definition of $\alpha_1$, and the result is proved.
\end{proof}


\begin{thebibliography}{20}
\bibitem{bary}
Y.M. Baryshnikov, V. Zharnitsky, Billiards and nonholonomic distributions, \emph{J. Math. Sci.} 128 (2005), 2706–2710.

\bibitem{bary2}
Y.M. Baryshnikov, V. Zharnitsky, Sub-Riemannian geometry and periodic orbits in classical billiards, \emph{Math. Res. Lett.} 13 (2006), no. 4, 587–598.

\bibitem{bary_introuvable}
Y.M. Baryshnikov, Spherical billiards with periodic orbits, \textit{preprint}.

\bibitem{beltrami}
Beltrami, E. Risoluzione del problema: Riportare i punti di una superficie sopra un piano in modo che le linee geodetiche vengano rappresentate da linee rette. \textit{Annali di Matematica pura ed applicata} 7, 185–204 (1865).

\bibitem{berger_geometry}
M. Berger, \textit{Geometry}, Volumes I and II, Springer-Verlag, 1987.

\bibitem{berger_caustics}
M. Berger, Seules les quadriques admettent des caustiques, \emph{Bull. Soc. math. France} 123 (1995), 107–116.

\bibitem{bialy}
M. Bialy, Convex billiards and a theorem by E. Hopf, \textit{Math. Z.} 214, No. 1 (1993), 147–15.

\bibitem{bialymironov_poly}
M. Bialy, A. E.  Mironov,  Angular  billiard  and  algebraic  Birkhoff  conjecture, \textit{Adv. Math.} 313 (2017), 102–126.

\bibitem{bialymironov_poly_bis}
M. Bialy, A. E.  Mironov,  Algebraic Birkhoff conjecture for billiards on Sphere and Hyperbolic plane, \textit{J. Geom. Phys.} 115 (2017), 150–156.

\bibitem{VKNZ}
V. Blumen, K. Kim, J. Nance, V. Zharnitsky, Three-Period Orbits in Billiards on the Surfaces of Constant Curvature, \emph{International Mathematics Research Notices} (2012), 10.1093/imrn/rnr228. 

\bibitem{bryant_chern}
R. L. Bryant, S. S. Chern, R. B. Gardner, H. L. Goldschmidt, P. A. Griffiths, \textit{Exterior Differential Systems}, Springer Science and Business Media, 2013.

\bibitem{cartan}
E. Cartan, \textit{Les systèmes différentiels extérieurs et leurs applications géométriques}, Paris, 1945.

\bibitem{CKS}
S.J. Chang, B. Crespi, K.J. Shi, Elliptical billiard systems and the full Poncelet’s theorem in $n$ dimensions, Journal of Mathematical Physics 34, 2242 (1993).

\bibitem{DragRad_bicent}
V. Dragović, M. Radnović, Bicentennial of the Great Poncelet Theorem (1813-2013): Current Advances, \textit{Bulletin of the American Mathematical Society} 51 (2013), No. 3, 373–445.

\bibitem{DragRad_minkowski2}
V. Dragović, M. Radnović, Ellipsoidal billiards in pseudo-Euclidean spaces and relativistic quadrics, Advances in Mathematics 231 (2012), 1173–1201.

\bibitem{DragRad_minkowski}
A.K. Adabrah, V. Dragović, M. Radnović, Periodic billiards within conics in the Minkowski plane and Akhiezer polynomials, \textit{Regul. Chaotic Dyn.} 24 (2019), 464–501.

\bibitem{DragRad}
V. Dragović, M. Radnović, \emph{Poncelet Porisms and Beyond}, Springer Basel, 2011.

\bibitem{fierobe_circumcenters}
C. Fierobe, On the Circumcenters of Triangular Orbits in Elliptic Billiard, \textit{J. Dyn. Control Syst.} (2021). \url{https://doi.org/10.1007/s10883-021-09537-2}

\bibitem{fierobe_caustics}
C. Fierobe, Complex Caustics of the Elliptic Billiard, \textit{Arnold Math. J.} 7 (2021), 1–30. \url{https://doi.org/10.1007/s40598-020-00152-w}

\bibitem{fierobe_triangular}
C. Fierobe, On projective billiards with open subsets of triangular orbits, \url{https://arxiv.org/abs/2005.02012}

\bibitem{fierobe1}
C. Fierobe, Examples of reflective projective billiards, \url{https://arxiv.org/pdf/2002.09845.pdf}

\bibitem{flatto}
L. Flatto, \textit{Poncelet's Theorem}, AMS, 2009.

\bibitem{glut}
A. A. Glutsyuk, On 4-reflective complex analytic billiards, \emph{Journal of Geometric Analysis} 27 (2017), 183–238.

\bibitem{glut1}
A. A. Glutsyuk, On quadrilateral orbits in complex algebraic planar billiards, \emph{Moscow Mathematical Journal} 14 (2014), 239–289

\bibitem{glut2}
A. A. Glutsyuk, On Odd-periodic Orbits in complex planar billiards, \emph{Journal of Dynamical and Control Systems} 20 (2014), 293–306.

\bibitem{glutkud1}
A.A. Glutsyuk, Yu.G. Kudryashov, On quadrilateral orbits in planar
billiards, \emph{Doklady Mathematics} 83 (2011), No. 3, 371–373.

\bibitem{glutkud2}
A. A. Glutsyuk, Yu. G. Kudryashov, No planar billiard possesses an
open set of quadrilateral trajectories, \emph{J. Modern Dynamics} 6 (2012), No. 3, 287–326.

\bibitem{glut_caustics}
A. A. Glutsyuk, On commuting billiards in higher-dimensional spaces of constant curvature, \emph{Pacific Journal of Mathematics} 305 (2020), No. 2, 577–595.

\bibitem{glut_integrability}
A. A. Glutsyuk, On Two-Dimensional Polynomially Integrable Billiards on Surfaces of Constant Curvature, \textit{Dokl. Math.} 98 (2018), 382–385.

\bibitem{glut_integrability_bis}
A. A. Glutsyuk, On polynomially integrable Birkhoff billiards on surfaces of constant curvature, \textit{J. Eur. Math. Soc.} \url{https://www.ems-ph.org/journals/of_article.php?jrn=jems&doi=1027}

\bibitem{GH}
Ph. Griffiths and J. Harris, \emph{Principles of algebraic geometry}, John Wiley and Sons, 1978.

\bibitem{GHponc}
Ph. Griffiths and J. Harris, Cayley's explicit solution to Poncelet's porism, \emph{L'Enseignement Mathématiques} 24 (1978), 31–40.

\bibitem{ivrii}
V. Y. Ivrii, The second term of the spectral asymptotics for a Laplace-Beltrami operator on manifolds with boundary, \textit{Func. Anal. Appl.}, 14 (1980), 98–106.

\bibitem{izmestiev}
I. Izmestiev, Spherical and hyperbolic conics, \emph{Eighteen Essays in Non-Euclidean Geometry}, EMS, 2019.

\bibitem{kac}
M. Kac, Can One Hear the Shape of a Drum?, The American Mathematical Monthly 73, No. 4 (1966), 1–23.

\bibitem{kaloshinsorrentino}
V. Kaloshin, A. Sorrentino, On the local Birkhoff conjecture for convex billiards, \textit{Annals of Mathematics} 188, No. 1 (2018), 315–380.

\bibitem{treshchev}
V. V. Kozlov, D. V. Treshchëv, \textit{Billiards, A Genetic Introduction to the Dynamics of Systems with Impacts}, Translations of Mathematical Monographs vol. 89, Amercian Mathematical Society, 1991. 

\bibitem{krattenthaler}
C. Krattenthaler, Advanced determinant calculus: a complement, \textit{Linear Algebra Appl.} 411 (2005), 68–166.

\bibitem{krattenthaler_cat}
C. Krattenthaler, Determinants of (generalised) Catalan numbers, \textit{Journal of Statistical Planning and Inference} 140 8 (2010), 2260–2270.

\bibitem{khesin_taba}
B. Khesin, S. Tabachnikov, Pseudo-Riemannian geodesics and billiards, \textit{Advances in Mathematics} 221 (2009), 1364–1396.

\bibitem{klein26}
F. Klein, \emph{\"Uber höhere Geometrie}, Springer, 1926

\bibitem{kuranishi}
M. Kuranishi, On E. Cartan's Prolongation Theorem of Exterior Differential Systems, \emph{The Johns Hopkins University Press} 79 (1957), 1–47. 

\bibitem{lojasiewicz}
S. Lojasiewicz, \emph{Introduction to Complex Analytic Geometry}, Springer, 1991.


\bibitem{matveev}
V.S. Matveev, Geometric explanation of the Beltrami theorem, \textit{International Journal of Geometric Methods in Modern Physics} 3, No. 03 (2006), 623–629.

\bibitem{petkovstojanov}
V. Petkov, L. Stojanov, On the number of periodic reflecting rays in generic domains, \textit{Ergodic Theory and Dynamical Systems}, 8 (1988), 81–91.

\bibitem{poncelet}
J. V. Poncelet, \textit{Traite des propriétés projectives des figures: vol. 2}, Gauthier Villars, 1866.

\bibitem{rashevsky}
P. K. Rachevsky, \textit{Geometrical theory of partial differential equations}, OGIZ Gostehizdat, 1947.

\bibitem{reznik_youtube}
D. Reznik, \emph{http://www.youtube.com/watch?v=BBsyM7RnswA}

\bibitem{reznik_github}
D. Reznik, R. Garcia and J. Koiller, New Properties of Triangular Orbits in Elliptic Billiards, \emph{https://dan-reznik.github.io/Elliptical-Billiards-Triangular-Orbits/} (april 2019)

\bibitem{romaskevich2014}
O. Romaskevich, On the incenters of triangular orbits in elliptic billiard, \emph{L'Enseignement Mathématiques} 60 (2014), 247-255

\bibitem{ronaldo}
R. Garcia, Elliptic Billiards and Ellipses Associated to the $3$-Periodic Orbits, \textit{The American Mathematical Monthly} 126 (2019), 491–504.


\bibitem{rychlik}
M. R. Rychlik, Periodic points of the billiard ball map in a convex domain, \emph{Journal of Differential Geometry} 30 (1989), 191–205.

\bibitem{stojanov}
L. Stojanov, Note on the periodic points of the billiard, \emph{Journal of Differential Geometry} 34 (1991), 835–837.

\bibitem{taba_centers}
R. Schwartz and S. Tabachnikov, Centers of mass of Poncelet polygons, 200 years after, \emph{https://math.psu.edu/tabachni/prints/Poncelet5.pdf}

\bibitem{shaidenko_vivaldi}
A. V. Shaidenko, F. Vivaldi, Global stability of a class of discontinuous dual billiards, \textit{Comm. Math. Phys.} 110 (1987), 625–640. 

\bibitem{taba_commut}
S. Tabachnikov, Commuting dual billiard maps, \textit{Geometriae Dedicata} 53 (1994), 57–68.

\bibitem{taba_proj_metrics}
S. Tabachnikov, Ellipsoids, complete integrability and hyperbolic geometry, \emph{Moscow Mathematical Journal} 2 (2002), 185–198.

\bibitem{taba_projectif_ball}
S. Tabachnikov, Exact transverse line fields and projective billiards in a ball, \emph{GAFA Geom. funct. anal.} 7 (1997), 594–608.

\bibitem{taba_book}
S. Tabachnikov, \textit{Geometry and Billiards}, American Mathematical Sociecty, 2005. 

\bibitem{taba_projectif}
S. Tabachnikov, Introducing projective billiards, \emph{Ergodic Theory and Dynamical Systems} 17 (1997), 957–976.

\bibitem{taba_dual_billiards}
S. Tabachnikov, On the dual billiard problem, \emph{Adv. in Math.} 115 (1995), 221–249.

\bibitem{vasiliev}
D. Vasiliev, Two-term asymptotics of the spectrum of a boundary value problem in interior reflection of general form, \emph{Funct Anal Appl.} 18 (1984), 267–277.

\bibitem{vorobets}
Ya. B. Vorobets, On the measure of the set of periodic points of the billiard, \emph{Math. Notes} 55 (1994), 455–460.

\bibitem{weyl}
H. Weyl, Das asymptotische Verteilungsgesetz der Eigenwerte linearer partieller Differentialgleichungen (mit einer Anwendung auf die Theorie der Hohlraumstrahlung), \emph{Math. Ann.} 71 (1912), 441–479.

\bibitem{whitney_diffmani}
H. Whitney, Differentiable manifolds, \emph{Annals of Mathematics} 37 (1936), 645–680.

\bibitem{wojtkowski}
M.P. Wojtkowski, Two applications of Jacobi fields to the billiard ball problem, \emph{Journal of Differential Geometry} 40 (1994), 155–164.

\bibitem{zaslawski}
A. Zaslavsky, D. Kosov, and M. Muzafarov, Trajectories of remarkable points of the Poncelet triangle (in Russian), \emph{Kvanto} 2 (2003), 22–25.

%
%
%
%
%
%
%
%
%
%
%
%
%
%
%
%
%
%
%
%
%
%
\end{thebibliography}
\end{document}